\documentstyle[epsfig]{article}

\input xypic

\input{amssymb.sty}

\makeindex

\title{Equivariant Seiberg-Witten Floer Homology}
\author{Matilde Marcolli, Bai-Ling Wang}
\date{}

\newtheorem{thm}{Theorem}[section]
\newtheorem{lem}[thm]{Lemma}
\newtheorem{corol}[thm]{Corollary}
\newtheorem{defin}[thm]{Definition}

\newtheorem{prop}[thm]{Proposition}
\newtheorem{rem}[thm]{Remark}

\newcommand{\beq}{\begin{equation}}
\newcommand{\eeq}{\end{equation}}

\newcommand{\M}{{\cal M}}

 \def\AA{{\mathbb A}}
 \def\RR{{\mathbb R}}
 \def\CC{{\mathbb C}}
 \def\ZZ{{\mathbb Z}}
 
 \def\HH{{\mathbb H}}

\begin{document}
\maketitle

\tableofcontents

\section{Introduction}

\begin{verse} {\em
und das hat mit ihrem Singen \\
die Loreley getan. \\
\vskip .1in
(H. Heine) }
\end{verse}

Floer homology, as a gauge theoretic invariant of three-manifolds,
made its first appearance in the seminal work of Floer
\cite{Fl} on what is now referred to as the Yang--Mills, or instanton, 
Floer homology. Instanton Floer homology came to play an essential
conceptual role in Atiyah's Topological Quantum Field Theory
formulation of Donaldson theory \cite{A}, and in the definition of
relative 
Donaldson invariants of 4-manifolds with boundary and their gluing 
formulae, see for instance the work of Taubes \cite{Taubes}.
Soon after the introduction of the new Seiberg--Witten gauge theory,
it became clear that an analogue of Floer homology existed in the
Seiberg--Witten context. In the last four years, there has been a
number of significant papers using Seiberg--Witten Floer homology and
gluing theorems for Seiberg--Witten invariants. It is now known
(see e.g. \cite{Kronheimer}) that a good understanding of the relation
between the instanton and the Seiberg--Witten Floer homology can lead
to very striking topological consequences. However, in the
literature (both in the instanton and in the Seiberg--Witten context)
there seems to be a fundamental need for a satisfactory foundational
work which presents a detailed and careful construction of the Floer
homology. 
It is important to mention that, perhaps surprisingly, the
technical aspects involved in the Seiberg--Witten gauge theory are
substantially different from the Yang--Mills case and require
different techniques. An important issue, which has no analogue in the
Yang--Mills case and to which much of this paper
is dedicated, is the metric dependence of the Seiberg--Witten Floer
homology in the case of homology 3-spheres. Providing a metric
independent version of Floer homology is essential in all the
important applications which involve surgery formulae and gluing of
relative invariants. We develop an equivariant version of
Seiberg--Witten Floer homology precisely to the purpose of avoiding
the metric dependence problem. The main technical issues involved
center around the fine structure of the compactification of the moduli 
spaces of flow lines of the Chern--Simons--Dirac functional, and
around the use of Taubes' obstruction bundles in the course of the
proof of the topological invariance of the equivariant Seiberg--Witten 
Floer homology. 

It is not easy to give the skeptical reader an introduction to this
paper that justifies the more than hundred and fifty pages to come. To be
perfectly honest, it took a long time to convince ourselves as well of 
the real need to unravel all the detailed technical issues that appear
in the corpus of this paper. Nonetheless, we hope the introductory
paragraph above gives sufficiently clear an idea of why it is
genuinely useful to write a paper that covers the material presented
here. This paper has been under preparation since 1995. It circulated
in different draft versions in 1996 and 1997. Regrettably,
more than two years passed between 
the time when the paper was submitted and the time when a referee
report was finally made available to us. We should signal to the
reader the appearance, in the meanwhile, of other papers 
that have a more or less extended overlap with the present one, most
notably \cite{Iga} and \cite{Ye}. In that respect, we are certainly
guilty of not having circulated, through the customary electronic
distribution, the last revision of this paper dating back to November
1997. We trust that this paper, in the present and hopefully last of
its long series of {\em avatars}, will soon see the light as a printed
form incarnation.  

We summarize, in the rest of the Introduction, the main
results of the various sections of this paper, pointing to precise
references to Theorems, Sections, etc. where the various statements
are presented and proved. 

The dimensional reduction of the Seiberg-Witten equations on
4-manifolds leads to equations on a compact oriented
3-manifold $Y$, obtained by considering translation invariant
solutions of the original equations in a temporal gauge on an infinite
cylinder $Y\times \RR$ with translation invariant metric. The
solutions modulo gauge transformations of the reduced equations on $Y$
can be regarded  
as the critical points of a Chern--Simons--Dirac functional, defined on
the configuration 
space of $U(1)$-connections and spinors. There is an associated
Floer-type homology, which depends on the choice of the
$Spin_c$-structure.  Some of the properties of this
Seiberg-Witten Floer homology have been discussed in \cite{CMW},
\cite{Wa}, \cite{Iga}, \cite{Ma}, \cite{MW2}, \cite{Ni4}  
\cite{Wa2}, \cite{GWa}, \cite{Ye}. 
The properties are different according to whether the
manifold $Y$ is a rational homology sphere or a manifold with
$b^1(Y)>0$. In particular, if $Y$ is a rational homology sphere, there
is always a reducible point in the moduli space, namely solutions with 
vanishing spinor and a non-trivial stabilizer of the gauge group
action. Unlike
the case with non-zero Betti number, the reducible cannot be perturbed
away just by adding a co-closed 1-form to the curvature equation.
This gives rise to an interesting metric dependence phenomenon.
The Seiberg-Witten Floer homology for manifolds with
$b^1(Y)>0$ presents a periodicity of the grading and is graded over
the integers only after passing to a cover. In this case, however,
with the formulation used in this paper, it
is no longer finitely generated. There are different ways of defining
the Floer homology when $b^1(Y)>1$, some of which still give rise to a
finitely generated complex. Some of these different
constructions are summarized briefly in Section 4, in Remark
\ref{correction}, an a more detailed
account will be given in \cite{MW4}. For the purpose of this paper,
the important issue is that, for $b^1(Y)>1$, the Floer homology
is independent of the metric and perturbation. Since we consider
perturbations that are cohomologically trivial, the metric and
perturbation independence is achieved in the case $b^1(Y)=1$ as well.
Thus, we have Seiberg-Witten Floer homologies with substantially
different properties depending on the underlying three-manifold $Y$.
The Floer homology always depends on a choice of the
$Spin_c$-structure, as does the construction of Seiberg--Witten
invariants. Thus, we really obtain, over a 3-manifold $Y$, a family of 
Floer homologies $SWF_*(Y,{\bf s})$, parameterized by the set of
$Spin_c$-structures ${\bf s}\in {\cal S}(Y)$. We shall always consider 
a fixed $Spin_c$-structure, hence, for simplicity of notation, we
shall not explicitly mention the ${\bf s}$ dependence in the
following. 

We give a unifying approach to Seiberg-Witten Floer
theory, by introducing a version of Seiberg-Witten Floer
homology that is defined for all three-manifolds and is always metric
independent. In the case of manifolds with $b^1(Y)>0$ we shall only
consider $Spin_c$-structures with $c_1(L)\neq 0$. The case $b^1(Y)>0$
and $c_1(L)=0$ has a more subtle behavior with respect to the choice of the
perturbation and the compactification of 
the moduli spaces of flow lines (cf. Remark \ref{Ybadcase}). Since
this case arises in important problems connected to surgery formulae,
we are going to deal with it separately in \cite{MW4}. 
Along the lines of the theory Austin and Braam
\cite{AB} developed for instanton homology, we 
construct an equivariant version of the Seiberg-Witten Floer
homology. In the case of rational homology spheres, this approach
counts all the contributions, both from the
reducible point and from the irreducibles. The equivariant Floer
homology also depends on the choice of the $Spin_c$-structure. 
We prove the invariance of the equivariant Floer homology with respect to
the metric and perturbation. The argument involves Taubes' obstruction
bundle technique and an appropriate gluing theorem.
In the two distinct cases of rational homology spheres and of manifolds with
non-trivial Betti number, we compare the equivariant theory with the
constructions of non-equivariant Seiberg-Witten Floer theories,
as presented in the work of K. Iga, M. Marcolli, B.L. Wang, and
R.G. Wang. In the 
case of manifolds with non-trivial Betti number, we obtain an
isomorphism of the equivariant and the non-equivariant theories (both
considered for a fixed $Spin_c$-structure) and in the case of rational
homology spheres we construct some exact sequences that relate the
equivariant and the non-equivariant theory. The results are obtained
by means of the spectral sequences associated to a filtration
of the complexes.
The metric independence of the equivariant Floer homology 
together with the exact sequences lead to an algebraic
proof of the wall crossing formula for the associated Casson-type
invariant. We also provide a geometric proof of the same
formula, through the analysis of the local structure of the moduli
spaces. 

The paper is organized as follows. Section 1 and 2 present preliminary
material, starting, in Propositions \ref{compactM},
\ref{b1>0}, and \ref{nu},
with a brief account of Seiberg-Witten gauge theory
on three-manifolds and a brief overview of the properties and the
local structure of the moduli space ${\cal M}$ of critical points of
the Chern--Simons--Dirac functional, and then continuing in Sections 2.2 and 
2.3 with an account of the properties
of the moduli spaces ${\cal
M}(a,b)$ of flow lines connecting critical points. These are the 
moduli spaces that appear in the construction of the non-equivariant
versions of Seiberg-Witten Floer theory, as  used for instance in \cite{Wa},
\cite{Ma}, \cite{Wa2}, 
\cite{GWa}. The main result in Sections 2.2 and 2.3 is the construction of a
suitable class of perturbations of the flow lines equations, in
Definition \ref{calP}, and the proof of transversality in Proposition
\ref{transverse}. 
We also discuss the notion of relative Morse index of critical points
of the Chern--Simons--Dirac functional in Proposition \ref{relmorseind}
and the orientation of the moduli spaces ${\cal M}(a,b)$ in
Proposition \ref{orientation}. 
The construction and properties of the non-equivariant
Seiberg-Witten Floer theory are briefly summarized in Section 2.4, the 
main result being the gluing formula of Lemma \ref{glue},
which is essential in establishing the property $\partial\circ
\partial=0$ for the boundary of the Floer complex. We do not present a 
complete treatment of the non-equivariant Floer homology, since all the
necessary results follow from the equivariant case of Sections 4, 5,
and 6, and the comparison results of Section 7.
However, we discuss some essential aspects of the non-equivariant
theory, for instance the phenomenon of metric dependence in the case
of homology spheres. In the proof of Theorem \ref{Morse} it is shown
that, in the case of homology spheres, the unique
reducible point is isolated whenever the metric and perturbation are
chosen so that the twisted Dirac operator has trivial kernel. This
condition on the choice of metric and perturbation determines a
chamber structure with codimension one wall which is analyzed in
Section 2.1. We prove in Theorem \ref{codim:1} that the condition on
the Dirac operator is 
generic and that the space of metrics and perturbations breaks into
chambers with codimension one walls. We analyze the structure of the
walls in Theorem \ref{wallstrata}.

In Section 3 we introduce the framed moduli space ${\cal M}^0$ in the
framed configuration space ${\cal B}^0$. The results of Lemma
\ref{fixedpoints}, Lemma \ref{free:orbits}, and Theorem \ref{Morse}
can be summarized as follows: the framed moduli space ${\cal M}^0$ is
a smooth 
manifold with a $U(1)$-action: the action is free in the case of
manifolds with $b^1(Y)>0$. When $b^1(Y)=0$, there is a unique fixed
point in ${\cal M}^0$, which corresponds to the reducible point in
${\cal M}$. In the case of rational homology spheres, the finitely
many inequivalent flat $U(1)$ connections correspond to the unique
fixed point in ${\cal M}^0$ for all the different possible choices of
the $Spin_c$-structure. 
In Theorem \ref{Morse} we also show that the
critical orbits are all isolated and the Hessian is non-degenerate in
the directions orthogonal to the orbits. That is, the 
Chern--Simons--Dirac functional satisfies the Morse-Bott condition.

In the remaining of Section 3 we analyze the moduli spaces ${\cal
M}(O_a,O_b)$ of gradient flow lines connecting critical orbits $O_a$
and $O_b$. In Section 3.3, Theorem \ref{relmorseind2}, we discuss the
relative Morse index of critical orbits, which requires
introducing the relevant Fredholm theory. The analysis of
Section 3.4 is aimed at proving that all finite energy solutions of
the Seiberg--Witten equations on the manifold $Y\times \RR$ decay
asymptotically to critical orbits, with an exponential weight which is 
determined by the smallest absolute value of the non-trivial
eigenvalues of the Hessian. Thus, every finite energy solution lies in 
some moduli space ${\cal M}(O_a,O_b)$. The exponential decay to the
endpoints of solutions in a moduli space ${\cal M}(O_a,O_b)$ is proved 
in Theorem \ref{decay}, using the result of Lemma \ref{PalaisSmale} (a
Palais--Smale condition), and the estimate of Lemma
\ref{near:critical}. The argument of Theorem \ref{decay} is basically
a ``finite energy implies finite length'' type result
(cf. \cite{Si}). Finally, we 
give the transversality result for the moduli spaces ${\cal
M}(O_a,O_b)$ in Proposition \ref{flowlines}.

Section 4 is dedicated to the existence and properties of the
compactification of the moduli spaces  
$\hat {\cal M}(O_a,O_b)$ of unparameterized flow lines. We prove that
a compactification can be obtained by adding boundary strata of broken
trajectories. The codimension k strata in the boundary are of the form
\beq
\bigcup_{c_1, \cdots c_k} \hat{\cal M}(O_a, O_{c_1}) \times _{O_{c_1}}
\hat{\cal M}(O_{c_1}, O_{c_2})\times \cdots \times_{O_{c_k}}
\hat{\cal M}(O_{c_k}, O_b),
\label{ideal:bound} \eeq
where the union is over all possible sequences of the critical points
$c_1, \cdots c_k$ with decreasing indices. We also show that the
compactification has a fine structure of a smooth manifold with
corners. 
The analysis in Section 4.1 shows that a compactification exists, and
that the points in the ideal boundary consist of a certain subset of
the set (\ref{ideal:bound}) of
broken trajectories through intermediate critical points. This is the
main result of Theorem \ref{seqcompact}. Unlike other problems of
compactification in gauge theory, the construction of Floer homology
requires a more detailed analysis of the compactification. In fact,
the property the the boundary $D$ of the (equivariant) Floer complex
satisfies $D\circ D =0$ requires non only to know that a
compactification of the moduli spaces $\hat {\cal M}(O_a,O_b)$ exists, 
but also that {\em all} the broken trajectories (\ref{ideal:bound}) in the
ideal boundary occur in the actual boundary of the
compactification. Moreover, since the argument in the proof of $D\circ
D =0$ is based on a version of Stokes' theorem for manifolds with
corners, it is essential to analyze the fine structure of the
compactification at the boundary strata (\ref{ideal:bound}).

In Section 4.2 we prove that in fact all broken trajectories through
intermediate critical points listed in (\ref{ideal:bound}) actually
occur in the compactification, the main result being the gluing theorem
\ref{equivgluing}. This theorem only deals with the codimension one
boundary, namely with the trajectories that break through one intermediate
critical orbit. The gluing construction which proves Theorem
\ref{equivgluing} involves several technical aspects. The main idea 
is to introduce a pre-gluing procedure, in which an approximate
solution to the Seiberg--Witten equations is obtained by splicing
together with cutoff functions a pair of solutions in the fibered product
$$ \hat{\cal M}(O_a,O_b)\times_{O_b} \hat{\cal M}(O_b,O_c). $$ 
The gluing construction then takes place by proving, via a fixed point 
argument (Remark \ref{fixedpoint}), that close enough to the
approximate solution there is a 
unique actual solution in $\hat{\cal M}(O_a,O_c)$.
The central technical issues connected to this argument revolve around
the analysis of the eigenspace of small 
eigenvalues of the linearization of the Seiberg--Witten equations at
the approximate solutions. Lemma \ref{uniformslice} and Lemma
\ref{slicelemma} deal with the slices of the gauge action. These are
necessary in order to introduce the pre-gluing construction in Lemma
\ref{inducedgluing}. The analysis of the small eigenvalues is
developed essentially in Lemma \ref{split:convergence}, Lemma
\ref{surj2}, and Lemma \ref{apprkercoker}, and then recalled, in
Section 6, in Lemma \ref{spectral:lemma}. In Lemma 
\ref{apprkercoker} we also provide an estimate of the rate of decay of 
the small eigenvalues. The eigenspaces of small
eigenvalues, which give the normal bundle for the gluing construction, 
are introduced in Definition \ref{defapprker} and in Proposition
\ref{apprkercoker2}, cf. also Corollary \ref{indexsplit} on the
splitting of the index, and Remark \ref{rem:isom}.

In Section 4.3 we extend the results of Section 4.2, in order to deal
with the strata of higher codimension in the compactification, and
to show the fine structure of the compactified moduli spaces, namely
the fact 
that the moduli spaces $\hat{\cal M}(O_a,O_b)$ compactify to smooth
manifolds with corners, in the sense of Melrose \cite{Melrose}. The
main result on the corner structure, Theorem \ref{corners}, is based
on Proposition \ref{equivgluing:k}, which is an inductive
generalization of Theorem \ref{equivgluing}, and shows the existence
of a smooth atlas of charts with corners.

Throughout all the gluing construction in Section 4, we make essential 
use of the transversality result, namely of the fact that the
linearizations at the solutions in $\hat{\cal M}(O_a,O_b)$ and
$\hat{\cal M}(O_b,O_c)$ have trivial cokernel. We say in such cases
that the gluing is ``unobstructed''. The spectral analysis in Lemma
\ref{split:convergence}, Lemma \ref{surj2}, and Lemma
\ref{apprkercoker}, however, is formulated in more general terms that
adapt to the case where cokernels are present. This gluing theory
``with obstructions'' is elaborated in Section 6, with the purpose of
proving the topological invariance of the equivariant Floer homology. 

As we discuss in the beginning of Section 4, the fine structure of the 
compactification of the moduli spaces of flow lines is 
necessary in order to establish the existence of the Floer complex,
namely the fact that the boundary operator $D$ of the equivariant
Floer complex satisfies $D^2=0$. Section 5 introduces the equivariant
complex, the boundary operator, and the equivariant Floer
homology. With the essential use of the results of Section 4, we prove 
the property $D^2=0$ in Theorem \ref{DD=0}. We then give an explicit
description of the boundary operator of the equivariant Floer complex
in Proposition \ref{boundaryterms}. 

Section 6 contains the proof of the invariance of the
equivariant Seiberg--Witten Floer homology with respect to the choice
of the metric and perturbation. The proof
of the invariance is obtained by defining a chain map $I$ connecting the
equivariant Floer complexes associated to choices $(g_0,\nu_0)$ and
$(g_1,\nu_1)$, a similar chain map $J$ in the opposite direction, and
then showing that there is a chain homotopy $H$, satisfying
$id-JI=DH+HD$, that induces an  
isomorphism on the level of cohomology. The chain map is constructed
by means of moduli spaces ${\cal M}(O_a,O_{a'})$ of solutions of
Seiberg-Witten equations on 
the cylinder $Y\times \RR$ endowed with a metric $g_t+dt^2$ that
varies between $g_0+dt^2$ and $g_1+dt^2$ along the cylinder. The main
theorem, Theorem \ref{metrics} is presented at the beginning of
Section 6, followed by a ``model'' proof of the easiest case of
metrics and perturbation in the same chamber. The proof of the general 
statement of Theorem \ref{metrics} will only be given in Section 6.3,
after the necessary technical tools have been introduced.

The construction of the chain map, in the general case of Theorem
\ref{metrics}, requires a careful analysis of the
boundary structure of the moduli spaces ${\cal M}(O_a,O_{a'})$.
In Section 6.1 we present the properties of the moduli spaces ${\cal
M}(O_a,O_{a'})$. In Theorem \ref{seqcompact2}, we prove the existence
of a compactification for the moduli spaces ${\cal M}(O_a,O_{a'})$,
obtained by adding fibered products of lower dimensional moduli
spaces. This is the analogue of the results of Section 4.1 in the case 
of flow lines. We then give a transversality result, in Lemma
\ref{surjT} and Lemma \ref{noflow}, under the hypothesis that at least 
one of the critical orbits $O_a$ and $O_{a'}$ has a free $U(1)$
action, or that the relative Morse index $\mu(\theta_0)-\mu(\theta_1)$ 
is non-negative. In the case of the moduli space ${\cal M}(\theta_0,
\theta_1)$, with $\mu(\theta_0)-\mu(\theta_1)<0$, we show in Corollary 
\ref{reducibleflow} and in Lemma \ref{unique:reducible} that the
transversality result fails. The moduli space ${\cal
M}(\theta_0,\theta_1)$ consists in this case of a unique reducible
solution and the linearization has a non-trivial cokernel of dimension 
equal to the index $|\mu(\theta_0)-\mu(\theta_1)|$. 
We prove in Theorem \ref{glue:I} the gluing theorem for all the non-singular
boundary strata in $\M(O_a,O_{a'})^*$, namely those with trivial
Cokernels. Similarly, we prove in Theorem \ref{glue:H} the gluing
theorem for all the non-singular boundary strata in $\M^P(O_a,O_b)^*$.
                           
Section 6.2 contains the general theory needed in order to deal with
the gluing theorems in the presence of obstructions. It relies heavily 
on Taubes' technique developed in \cite{Ta} and \cite{Ta2}, together
with Donaldson's results in \cite{D2}. We recall in Lemma
\ref{spectral:lemma} the necessary eigenvalue
splitting for the Laplacians at the approximate solution, as proved in 
Lemma \ref{split:convergence}, Lemma \ref{surj2}, and Lemma
\ref{apprkercoker} of Section 4.2. In Proposition \ref{obstruction} we 
introduce the obstruction bundle with the canonical obstruction
section. In Lemma \ref{solveqbeta} we derive the fixed point argument
which we use to complete the proof of the non-obstructed gluing of
Theorem \ref{equivgluing} and to formulate the gluing with obstruction 
in Proposition \ref{obstruction}.
In the remaining of Section 6.2 we analyze the modified boundary
strata of $\M(O_a,O_{a'})^*$ and $\M^P(O_a,O_b)^*$. In Theorem
\ref{geomobstr} and  Theorem \ref{noglue} we identify the extra
boundary components in $\M(O_a,O_{a'})^*$ due to the zeroes of the
obstruction section, and we define in Proposition \ref{gluing:obstr}
the additional gluing maps. Similarly, in Theorem
\ref{boundary:P:obstr}, we identify the extra boundary components in
$\M^P(O_a,O_b)^*$ due to the zeroes of the obstruction sections, and
we introduce the additional gluing maps.

In Section 6.3, we finally complete the proof of the topological
invariance. The modification of the boundary structure discussed in
Section 6.2 prescribes correction terms for the maps $I$ and $H$ (but
not $J$), so that the identities $ID-DI=0$ and $id-JI=DH+HD$ can still 
be satisfied. We first discuss, in Lemma \ref{sec:obstr:euler} and Lemma
\ref{obstr:sec:euler}, some more properties of the zeroes of the
obstruction sections and some identities obtained by counting these
zeroes, which are useful in checking the identities satisfied by the
coefficients of the maps $I$ and $H$. We then introduce, in Definition 
\ref{def:I:H:change}, the modified maps $I$ and $H$, and we prove, in
Lemma \ref{check:ID=DI}, Lemma \ref{check:1}, Lemma \ref{check:2}, and 
Lemma \ref{check:3}, that the necessary identities hold, thus
completing the proof of Theorem \ref{metrics} on the topological
invariance of the equivariant Seiberg--Witten Floer homology.
It is clear, from the way the proof is structured, that the argument
breaks down for the non-equivariant Floer homology. While the
invariance within the same chamber is still verified (nothing changes
in that part of the argument), in the case of metrics and
perturbations in different chambers we need essentially the
contribution of the reducible points in order to construct the chain
map $I$ and chain homotopy $H$, as one can see form Definition
\ref{def:I:H:change}, and the Lemmata \ref{check:ID=DI},
\ref{check:1}, \ref{check:2}, and \ref{check:3}.

Section 7 deals with the wall crossing formula for the Casson-type
invariant of homology spheres obtained as the Euler characteristic of
the (non-equivariant) Seiberg-Witten Floer homology. In Section 7.1 we
derive relations 
between the equivariant and the non-equivariant Floer homologies. We
prove in Theorem \ref{isomorphism} that, when we have $b^1(Y)>0$, the
equivariant and the non-equivariant Floer homologies are
isomorphic. This is not a surprising result: in fact, it is
conceptually like considering the equivariant homology of a manifold
with a free group action, for which it is well known that one
recovers the homology of the quotient. In the case of 
rational homology spheres, we prove in Theorem \ref{exact} that there
is an exact sequence relating the 
equivariant Floer homology with the non-equivariant and with the
polynomial algebra $H_*(BU(1),\RR)$. The results of both Theorem
\ref{isomorphism} and \ref{exact} are derived by
considering filtrations of the complexes and the associated spectral
sequences, as in Lemma \ref{spectr:seq1} and Lemma \ref{spectr:seq2}.
We also give an explicit expression
of the connecting homomorphism in the exact sequence in Proposition
\ref{connect:homom:prop}. 
The exact sequence that relates the metric independent equivariant
Floer homology with the metric dependent non-equivariant one is the
tool we exploit in Section 7.2, in order to derive, in a purely
algebraic way, the wall crossing formula for the Euler
characteristic in Theorem \ref{wallcross2}. The argument is based on
the topological invariance of the equivariant Floer homology proved in 
Theorem \ref{metrics}, and the result of Proposition
\ref{wallcross1} on the counting of the ranks of the equivariant Floer 
groups. The result is generalized to multiple wall crossings, and
to the case of the $J$-invariant perturbations of \cite{Chen}, in
Propositions \ref{wallcross:n}, 
\ref{wallcross:general}, and in Corollary \ref{Jwallcross}.
In Section 7.3 we show how the same wall crossing formula can be
derived geometrically by considering the local model of the
parameterized moduli space of critical points ${\cal M}(g,\nu)$ along
with a deformation of the metric and perturbation $(g,\nu)$.

Given the length of this work, we
have thought it useful to add an index of notation at the end of the
paper. 

\vskip .3in

\subsection{Three dimensional Seiberg-Witten theory} 

Any three--manifold admits a $Spin$-structure. A choice of the metric
determines a (non-canonical) choice of a  ``trivial'' $Spin$-structure
with spinor bundle 
$S$. A $Spin_c$-structure is therefore obtained by twisting $S$ with a
line bundle $L$.
Suppose given a three-manifold $Y$ with a $Spin_c$ structure $S\otimes
L$. Consider the space ${\cal A}_k$ \index{${\cal A}_k$}
of pairs $(A,\psi)$, where $A$ is a
$U(1)$ connection on the line bundle $L$ and $\psi$ is a section of
$S\otimes L$, endowed with a fixed $L^2_k$-Sobolev completion.

The group $\tilde{\cal G}_{k+1}$ is \index{$\tilde{\cal G}_{k+1}$}
the gauge group of maps of $Y$ in $U(1)$ locally modeled on the Lie
algebra 
$$Lie(\tilde{\cal G}_{k+1})=L^2_{k+1}(\Omega^0(Y,i\RR)),$$
acting on pairs $(A,\psi)$ by
\[ \lambda : (A,\psi)\mapsto (A-2\lambda^{-1}d\lambda, \lambda\psi). \]
The group ${\cal G}_{k+1}$ is the subgroup of $\tilde{\cal G}_{k+1}$
of gauge transformations $\lambda$ satisfying the condition
\beq
\frac{i}{2\pi}\int_Y c_1(L)\wedge \lambda^{-1}d\lambda =0. 
\label{int=0}
\eeq
Here $c_1(L)$ represents the Chern class of the line bundle $L$ that
defines the $Spin_c$-structure on $Y$. Clearly, the group ${\cal G}$
coincides with the full gauge group $\tilde{\cal G}$ in the case when
$c_1(L)=0$ rationally, hence in particular for rational homology spheres. 

The functional \index{${\cal C}$}
\beq
{\cal C}(A,\psi)=\frac{-1}{2}\int_Y (A-A_0)\wedge (F_A+F_{A_0})
+\frac{1}{2}\int_Y <\psi,\partial_A\psi>dv
\label{functional}
\eeq
was first introduced by Kronheimer and Mrowka in the proof of the Thom
conjecture \cite{KM}. It
is defined on the space ${\cal A}$ of connections and sections and
it is invariant under the action of the identity component
${\cal G}_{k+1}$. Thus, it descends to a real-valued functional on the
space ${\cal B}_k={\cal A}_k/{\cal G}_{k+1}$.

The first order increment of this functional defines a 1-form on the 
$L^2_k$-tangent space $T{\cal A}$, 
\beq
\label{1forma}
{\cal F}\mid_{(A,\psi,\rho)}(\alpha,\phi)=\int_Y -\alpha\wedge 
(F_A -*\sigma(\psi,\psi)) +\int_Y <\phi,\partial_A\psi>. 
\eeq
Thus, the gradient flow of the functional (\ref{functional}), with
respect to the $L^2$-inner product, is given
by the paths of connections and sections $(A(t),\psi(t))$ that satisfy the
equations 
\begin{equation}
\frac{d}{dt}\psi=-\partial_A\psi,
\label{3SW1}
\end{equation}
and
\begin{equation}
\frac{d}{dt}A=-*F_A+\sigma(\psi,\psi),
\label{3SW2}
\end{equation}
where the 1-form $\sigma(\psi,\psi)$ is given in local 
coordinates by $\frac{1}{2}<e_i\psi,\psi>e^i$. \index{$\sigma(\psi,\psi)$}

These equations can be thought of as the three-dimensional reduction
of the Seiberg-Witten equations on four-manifolds introduced in
\cite{W} (see also \cite{Au1}, \cite{Au2}, \cite{Wa}, \cite{Ma}, 
\cite{Wa2}, \cite{GWa}). 
In fact, we can consider the four-manifold $Y\times \RR$
with a cylindrical metric $g+dt^2$, and with $Spin_c$ structure determined by
the pullback of $S\otimes L$ via the projection $\pi:Y\times \RR\to
Y$. Thus we have $S^+\otimes L\cong \pi^*(S\otimes L)$. Over $Y$ we have
$S^+\otimes L\cong S^-\otimes L$ under Clifford multiplication by
$dt$, and the identifications $\Omega^{2+}(Y\times \RR, i\RR)\cong
\pi^*(\Omega^1(Y,i\RR))$ and $\Omega^1(Y\times \RR, i\RR)\cong
\pi^*(\Omega^0(Y,i\RR)\oplus\Omega^1(Y,i\RR))$, that is
$\pi^*(\rho(t))=*\rho(t) +\rho(t)\wedge dt$ and
$\pi^*(f(t),\alpha(t))=\alpha(t)+f(t)dt$ on $Y\times \RR$.  

Consider a pair $(\AA,\Psi)$ on
$Y\times \RR$, where $\AA$ is a $U(1)$ connection on the determinant
line bundle of $S^+\otimes L$ and $\Psi$ is a section of the spinor
bundle $S^+\otimes L$. An element $(\AA,\Psi)$ is in a temporal gauge
if the $dt$-component of $\AA$ vanishes identically. 
Thus, a path $(A(t),\psi(t))$ corresponds to an element $(\AA,\Psi)$ in a
temporal gauge.

\begin{lem}
\label{rem43}
The Seiberg-Witten equations (\ref{3SW1}) and (\ref{3SW2}) on $Y\times
\RR$ are equivalent to
\begin{equation} \label{41}
F^+_{\AA}=\Psi\cdot\bar\Psi, 
\end{equation}
and
\begin{equation} \label{42}
D_{\AA}\Psi=0, 
\end{equation}
where $\Psi\cdot\bar\Psi$ is the self-dual two-form given in local
coordinates by $\Psi\cdot\bar\Psi=\langle e_ie_j\Psi,\Psi\rangle
e^i\wedge e^j$ and 
$D_{\AA}$ is the Dirac operator twisted with the connection $\AA$
acting on sections of $S^+\otimes L$. 
\end{lem}

\noindent\underline{Proof.} The Dirac operator $D_{\AA}: S^+\otimes L
\to S^- \otimes L$ on the four-manifold, 
twisted with the connection $\AA$, has the form
\[ D_{\AA}=\partial_t +\partial_{A(t)}, \]
where $\partial_A$ is the self-adjoint Dirac operator on $Y$ twisted with
the time dependent connection $A(t)$. For the curvature equation
(\ref{41}), write $F_{\AA}^+=\frac{1}{2}(F_{\AA}+*_4 F_{\AA})$. Since
$F_{\AA}^-$ acts trivially on the positive spinors, the action of
$F_{\AA}^+$ corresponds precisely to the action of $\frac{dA}{dt}+*_3
F_A$. Here we have introduced the notation $*_4$ and $*_3$ to
distinguish the Hodge $*$-operator on the 4-manifold $Y\times \RR$ and 
on the 3-manifold $Y$. In the following, we shall drop the subscript,
since it will be clear which $*$-operator is being used. Thus, given
the expression of the 2-form $\Psi\cdot\bar\Psi$ and of the one form
$\sigma(\psi(t),\psi(t))$, we can write equation (\ref{41}) as
\[ \frac{1}{2}(F_{it}+\epsilon^{itjk} F_{jk}) e^i\wedge dt
=\frac{1}{4}<e_ie_t\psi,\psi> e^i\wedge dt, \]
with an implicit sum over repeated indices and the symbol $\epsilon$
denoting the sign of the permutation $\{ itjk \}$. Upon applying 
the $*_3$ operator and using the identification of $S^+\otimes L$ 
and $S^-\otimes L$ on $Y\times \RR$ under Clifford multiplication by
$dt$, we can identify this equation with the equation (\ref{3SW2}).

\noindent $\diamond$

The Seiberg--Witten equations on the 4-manifold $Y\times\RR$ are
preserved by the action of the gauge group \index{$\tilde{\cal
G}_{Y\times \RR}$} 
$\tilde{\cal G}_{Y\times \RR}$ of maps of $Y\times \RR$ in
$U(1)$. Any element $(\AA,\Psi)$ can be transformed into a temporal
gauge by effect of a gauge transformation.
Suitable Sobolev completions of these spaces will be
introduced in Proposition \ref{sobolev}. For a general overview
of the Seiberg-Witten theory on four-manifold see \cite{D}, \cite{Ma2},
\cite{Mo}, \cite{Sa}. 

The critical points of the functional ${\cal C}$ are  pairs $(A,
\psi)$ that satisfy
\beq
\begin{array}{c}
         \partial_A\psi=0, \\
         {}*{}F_A=\sigma(\psi,\psi).
\end{array}
\label{extrem}
\eeq

Let ${\cal M}_k$ be the moduli space of solutions of (\ref{extrem}) in
${\cal B}_k$ and $\tilde{\cal M}_k$ be the moduli space of solutions
in $\tilde{\cal B}_k={\cal A}_k/\tilde{\cal G}_{k+1}$. \index{${\cal
M}$} 

\begin{lem}
\label{ellreg}
By elliptic regularity, ${\cal M}_k$ can be represented by smooth
elements.
\end{lem}

Lemma \ref{ellreg} follows from the Sobolev embedding theorems. In the
following we drop the subscript $k$ and just write ${\cal M}$
for the moduli space of critical points.

The deformation complex that determines the virtual dimension of
${\cal M}$ is given by 
\[ 0\to\Omega^0(Y, i\RR)\oplus\Omega^1(Y, i\RR)\oplus\Gamma (S
\otimes L)\stackrel{L}{\to}\Omega^0(Y, i\RR)\oplus\Omega^1(Y, i\RR)\oplus
\Gamma (S\otimes L)\to 0, \]
with an index zero Fredholm operator \index{$L\mid_{(A,\psi)}$}
\beq
L\mid_{(A,\psi)}(f,\alpha,\phi)=\left\{\begin{array}{l}
            T\mid_{(A,\psi)}(\alpha,\phi)+G\mid_{(A,\psi)}(f)\\
            G^*\mid_{(A,\psi)}(\alpha,\phi)
          \end{array}\right. 
\label{L}
\eeq
defined between the $L^2_k$ and the $L^2_{k-1}$-Sobolev completions of 
the spaces above, where the
operator 
\beq
G\mid_{(A,\psi)}(f)=(-df,f\psi)
\label{G}
\eeq
is the infinitesimal action of
the gauge group and $G^*$ is the adjoint with respect to the $L^2$
pairing.
The map $T$ is the linearization of the \index{$T\mid_{(A,\psi)}$}
equation (\ref{extrem}) at a pair $(A,\psi)$,
\beq
T\mid_{(A,\psi)}(\alpha,\phi)=\left\{
\begin{array}{c}
      *d\alpha -2iIm\sigma(\psi,\phi) \\
      \partial_{A}\phi +\alpha \psi.
\end{array}\right.
\label{T}
\eeq

The Hessian of the functional ${\cal C}$ is given by a quadratic form
in the increment $(\alpha,\phi)\in T{\cal A}$,
\beq
\label{2form}
\nabla{\cal F}\mid_{(A,\psi)}(\alpha,\phi)=<\alpha, *d\alpha 
-2iIm\sigma(\psi,\phi)> + <\phi,\partial_A\phi +\alpha \psi>. 
\eeq
This descends to the same operator $T$ on the $L^2_k$-tangent space
$T_{[A,\psi]}{\cal B}_k$, when $(A,\psi)$ is a solution of
(\ref{extrem}), since the condition 
\[ G^*_{(A,\psi)}( T_{(A,\psi)}(\alpha,\phi))=0 \]
is satisfied. The operator $T$ is essentially self-adjoint.

As in the case of Donaldson theory \cite{A}, the linearization is a first 
order elliptic operator, hence its spectrum is not in general bounded
{}from below, and this affects the definition of the index of critical
points, as we are going to discuss in Proposition \ref{relmorseind}.

\section{Non-equivariant Floer theory}

There are two disjoint situations in which the non-equivariant
Seiberg-Witten Floer homology can be defined. They have been
considered in \cite{Wa}, \cite{Ma}, \cite{GWa} and are
summarized in this section. 

We consider the functional ${\cal C}$ perturbed with a co-closed
1-form $\rho$, \index{${\cal C}_\rho$}
\beq
{\cal C}_\rho(A,\psi)={\cal C}(A,\psi)-2i\int_Y (A-A_0)\wedge *\rho.
\label{pertfunct}
\eeq
We have the corresponding perturbed critical point equations
\beq
\begin{array}{c}
         \partial_A\psi=0, \\
         {}*{}F_A=\sigma(\psi,\psi)+2i\rho.
\end{array}
\label{extremP}
\eeq

Notice that the first integral cohomology group of $Y$ counts the
homotopy classes of gauge transformations, namely $H^1(Y,\ZZ)\cong
\pi_0(\tilde{\cal G})$ under the identification 
\[ \lambda\to h=[\frac{i}{2\pi}\lambda^{-1}d\lambda ]\in H^1(Y,\ZZ). \]
We denote with $H$ the subgroup of the classes $h\in H^1(Y,\ZZ)$ that satisfy
\[ \langle c_1(L)\cup h, [Y]\rangle =\frac{i}{2\pi}\int_Y c_1(L)\wedge
\lambda^{-1}d\lambda =0.\]  

We have the following compactness result for the set of critical
points of the Seiberg-Witten functional.

\begin{prop}
The moduli space $\tilde{\cal M}_\rho$ of solutions of
(\ref{extremP}) modulo the action of the full gauge group $\tilde
{\cal G}$ is compact. 
The space ${\cal M}_\rho$ of solutions, modulo the subgroup
${\cal G}\subset\tilde{\cal G}$ of gauge transformations satisfying
(\ref{int=0}), consists of a copy of \index{${\cal G}$}
$\tilde{\cal M}_\rho$ for each class in $H^1(Y,\ZZ)/H$. ${\cal
  M}_\rho$ is compact iff $H^1(Y,\ZZ)/H$ is finite. \index{${\cal M}$} 
In particular, ${\cal M}_\rho$ is compact if $b^1(Y)=0$.
\label{compactM}
\end{prop}

Notice that the condition that $H^1(Y,\ZZ)/H$ is finite corresponds to
$c_1(L)=0$ rationally. 

\begin{prop}
If the manifold $Y$ has $b^1(Y)>0$, then,
for an open set of small perturbations, the perturbed 
equations (\ref{extremP}) do not admit reducible solutions. 
Moreover, by the Sard-Smale
theorem, for a generic choice of $\rho$ the
corresponding moduli space $\tilde{\cal M}_\rho$ is a smooth manifold
that is cut out transversely by the equations. $\tilde{\cal M}_\rho$
is compact and zero-dimensional, hence it consists of a finite set of points.
\label{b1>0}
\end{prop}

All the moduli spaces come with a natural orientation
defined by the determinant line bundle of the Fredholm linearization,
as in \cite{D3}.

When $Y$ is a homology sphere the perturbation
$\rho$ can be written as $\rho=*d\nu$ and the equations
(\ref{extremP}) admit one gauge class of reducible solutions $[\nu,
0]$. 

\begin{prop}
\label{nu}
If $Y$ is a homology sphere, then there is a unique gauge class of reducible
solutions $\theta=[\nu, 0]$ of (\ref{extremP}) with $\rho=*d\nu$. If
\index{$\theta=[\nu, 0]$}
the metric on $Y$ is such that $Ker(\partial^g_\nu)=0$, where
$\partial^g_\nu$ is the self-adjoint Dirac operator on $Y$
\index{$\partial^g_\nu$}
twisted with the $U(1)$-connection $\nu$, then the
analysis of the local Kuranishi model shows that the reducible point
is isolated and non-degenerate.
\end{prop}

Propositions \ref{compactM}, \ref{b1>0}, and \ref{nu} are proven in
\cite{GWa}, \cite{Ma}, and \cite{Wa} respectively. See also the more
recent \cite{Lim}. 

\subsection{Chamber structure} 

In this subsection we 
prove that $Ker(\partial^g_\nu)=0$ is satisfied
for generic metrics. The condition $Ker(\partial^g_\nu)\neq 0$ defines
a chamber structure in the space of metrics. The walls form a
stratified space with a codimension one top stratum and higher
codimensional strata, in the sense described in Theorem
\ref{wallstrata}. 

Notice that there is a more abstract approach \cite{Kosch} describing a
stratification of the space of Fredholm operators by their index and
kernel dimension. However, the result we are interested in does not
follow directly by simply applying the results of \cite{Kosch}. In fact, 
we know from the general result of Section 1 of \cite{Kosch} that 
the space ${\cal F}_0(S,Y)$ of index zero Fredholm operators, acting on
the bundle $S$ (or $S\otimes L$) over $Y$, has the structure of a real
analytic manifold, 
with a stratification given by the sets ${\cal F}_{0,n}(S,Y)$ of
Fredholm operators of index zero with the dimension of the kernel
equal to $n$. In our specific problem, however, we are considering the 
particular map from the space ${\cal M}et \times Z^1(Y,\RR)$ of
metrics and perturbations $(g,*d\nu)$ to the space ${\cal F}_0(S,Y)$
given by
\begin{equation}
\label{Diracmap}
\begin{array}{c}
{\cal X}: {\cal M}et \times Z^1(Y,\RR) \to {\cal F}_0(S,Y)
\\[2mm] (g,*d\nu)\mapsto \partial^g_\nu =: {\cal X}(g,\nu). \end{array}
\end{equation}
Thus, we need a specific result that shows how the image of the map
$\Xi$ lies in ${\cal F}_0(S,Y)$ with respect to the stratification of
\cite{Kosch}. This is the purpose of our Theorem \ref{codim:1} and
Theorem \ref{wallstrata}. 
Although it is quite possible that a ``proof by library search'' of
these results may be 
provided combining the circle of ideas in \cite{Kosch}, \cite{Maier},
and the more recent \cite{Feehan}, we prefer to give a direct proof
that covers our specific case. 
A note
should be added: during the long period between the initial submission
of this paper and the completion of its refereeing process, a
discussion of the chamber structure has also appeared in \cite{Lim}
and \cite{OkTel}. 

\begin{rem}
Suppose $Y$ is a homology sphere. There is a quaternion structure on the
spinor bundle $S$, $J: S\to \bar S$ locally given by $(z_1,z_2)\mapsto
(-\bar z_2,\bar z_1)$. The action of $J$ extends to the configuration
space ${\cal A}$ as $(A,\psi)\mapsto (A^*,J\psi)$. This means that $J$
acts on $\nu$ by $J\nu=-\nu$, hence
\[ J(\partial^g_\nu\psi)=\partial^g_{-\nu}J\psi. \]
In other words, the Dirac operator $\partial^g$ on a three-manifold $Y$
is quaternion linear, but the twisted Dirac
operator $\partial^g_\nu$ is only complex linear.
\end{rem}

A $J$-invariant perturbation of the Seiberg-Witten equations has been
constructed by W. Chen \cite{Chen}, where the 
Dirac operator is perturbed with a smooth real function $f$ on $Y$. This
choice of perturbation leads to a different chamber structure in the
space of metrics which lies inside one chamber of the more general
perturbation  
\[ \partial_A\psi=f \psi, \]
\[ *F_A = \sigma (\psi, \psi) + *d\nu. \]
In this setting, transversality
can be achieved by adding a perturbation by a function of the holonomy.
This choice of perturbations lead to a different wall crossing
formula. However, we are going to show in Section 7, Proposition
\ref{wallcross:general} and Corollary \ref{Jwallcross}, that the wall
crossing formula for the $J$-invariant perturbations can be derived
with the same technique that we employ in the case on perturbations
$\nu$. 

The following theorems discuss the chamber structure in the space of
metrics and perturbations.

\begin{thm}
Let $Met$ be the space of all Riemannian metrics
on a homology 3-sphere $Y$. Consider the twisted Dirac operator
$\partial_{\nu}^g$ associated with the chosen metric $g$ and
the connection $\nu$.
The condition $Ker\partial_{\nu}^g \neq 0$ determines a real codimension
one subset in the space of $Met \times Z^1(Y, i\RR)$.
\label{codim:1}
\end{thm}

\noindent\underline{Proof.} Let $g_0$ be a metric on $Y$ such that
$Ker\partial_{\nu_0}^{g_0} \neq 0$ for the connection $\nu_0$.
We can decompose the spinor space as ${\cal H}\oplus {\cal H}^\perp$, under the
$L^2$ inner product,
where ${\cal H}=Ker \partial_{\nu_0}^{g_0}$, equipped with a Hermitian
metric from the $Spin_c$ structure. Consider the
Dirac operator $\partial_{\nu}^{g}$ for $(g, \nu)$ sufficiently
close to $(g_0 , \nu_0)$ in the $C^\infty$-topology. Under the isometry
identification of the spinor spaces for $g_0$ and $g$,
the Dirac operator for $g$ can be considered to act on the spinor space
of $g_0$, still denoted by $\partial_\nu^g$.

\noindent\underline{\bf Claim 1 :}
$\partial_\nu^g$ acting on the spinor space $S_{g_0}$ is self-adjoint
if and only if the metrics $g_0$ and $g$ define the same volume element.

Suppose $dvol_g = f \ dvol_{g_0}$ for a positive function
$f$ on $Y$, then with a direct calculation we get
\[<\partial_\nu^g \psi, \phi >_{g_0} = < \psi, f\partial_\nu^g(f^{-1}
\phi) >_{g_0}. 
\]
Claim 1 is then immediate.  Denote by $Met^0$ the space of metrics which
have the same volume element as the metric $g_0$. 

\noindent\underline{\bf Claim 2 :}
If two metrics $g_1$ and $g_2$ are conformal,
that is, $g_1 =e^{2u} g_2$ for a real function on $Y$, then the
multiplication by $e^{-u}$ defines an isomorphism between
$Ker \partial_\nu^{g_1}$ and $Ker \partial_\nu^{g_2}$.

This is the consequence of the following relation: under the isometry
identification of the spinor spaces for $g_1$ and $g_2$, we have
(see \cite{Hi} or \cite{LaMi} Theorem II.5.24)
\[ \partial_\nu^{g_1} = e^{-u} \partial_\nu^{g_2} e^u.\]
Notice that, in the formula for the variation of the Dirac operator under
conformal changes in the metric, different conventions are used in the 
literature. Here we are following the convention and notation of
\cite{Hi}. 

Therefore, we only need to prove that the condition
$Ker\partial_{\nu}^g \neq 0$ determines a real codimension
one subset in the space of $Met^0 \times Z^1(Y, i\RR)$.

We want to reduce the problem of the existence of solutions of the equation
$\partial_{\nu}^g \psi = 0 $ to a finite dimensional problem on ${\cal H}$.

As a map from $Met^0 \times Z^1(Y, i\RR) \times \Gamma (S)$ to $\Gamma (S)$,
the linearization of the equation $\partial_{\nu}^g \psi = 0$ at
$(g_0, \nu_0,0)$ 
is invertible when restricted to ${\cal H}^\perp$.
Thus, the implicit function theorem provides a unique
map $q: {\cal U}\to {\cal H}^\perp$ defined on a neighborhood ${\cal U}$ of
$(g_0, \nu_0, 0)$ in $Met^0 \times  Z^1(Y, i\RR) \times {\cal H}$, such
that

\[
(1-\Pi) \partial_{\nu}^g (\phi+ q( g, \nu, \phi ) ) = 0
\]
for all $(g, \nu, \phi )  \in {\cal U}$, where $\Pi$ is the projection
onto ${\cal H}$.

Therefore,  the operator $ \partial_{\nu}^g$ has non-trivial kernel if
and only if the equation
\[
\Pi \partial_{\nu}^g (\phi+ q( g, \nu, \phi ) )=0
\]
admits solutions in ${\cal H}$. This is a finite dimensional
problem.  Define a map
\[
  L: (Met^0 \times  Z^1(Y, i\RR)) \cap {\cal U} \longrightarrow U({\cal H})
\]
 by setting
\[
L(g, \nu) (\phi) = \Pi \partial_{\nu}^g (\phi+ q( g, \nu, \phi ) )
\]

Direct calculation implies that
$L(g, \nu)$ is a Hermitian transformation of the space ${\cal H}$,
that $L(g, \nu) \in  U({\cal H})$, where $ U({\cal H})$ is the
Hermitian transformation group on ${\cal H}$.
The kernel of $ \partial_{\nu}^g$ is non-trivial if and only if the kernel of
$L(g, \nu)$ is non-trivial. The determinant is a real-valued function
on $U({\cal H})$. Thus, we have a real-valued function $ f(g, \nu)
=det(L(g, \nu))$ on the neighborhood of $(g_0, \nu_0)$ in $Met^0
\times  Z^1(Y, i\RR)$. Those $(g, \nu)$ with non-trivial kernel have
value $0$ for this function. 

Now we only need to check that the derivative of $f(g, \nu)$ at
$(g_0, \nu_0)$ is surjective, then the condition $Ker\partial_{\nu}^g
\neq 0$ determines a real codimension one subset in the space of $Met^0
\times Z^1(Y, i\RR)$ by Morse theory.
It can be checked by differentiating $f(g, \nu)$ at
$(g_0, \nu_0)$ along $(0, \alpha)$-direction, for $\alpha \in
\Omega^1(Y, i\RR)$. Since
\[
Df_{(g_0, \nu_0)}(0, \alpha) =Tr(\phi\mapsto \Pi ( \frac 12 \alpha . \phi )),
\]
which is non-zero for suitable choice of $\alpha$.

Theorem \ref{codim:1} now follows from Claim 2 and the fact that any metric
is conformally equivalent to a metric in $Met^0$. 

\noindent $\diamond$

Let ${\cal W}$ denote the wall in the space of metrics and
perturbations $Met\times Z^1(Y,i\RR)$,
\[ {\cal W}=\{ (g, *d\nu)| \ Ker \partial_{\nu}^g \neq 0 \}. \]
Motivated by Chen's work \cite{Chen}, we analyze more carefully the
structure of ${\cal W}$. This analysis will be useful in Section 7, in 
establishing the results on the wall crossing formulae.

\begin{thm}
${\cal W}$ is a stratified space with the top stratum ${\cal W}_1$
consisting of those $(g,\nu)$ with $Ker \partial_\nu^g \cong \CC$. In
general, the set of pairs $(g,\nu)$ with
$Ker \partial_\nu^g \cong \CC^n$ is a codimension $2n-1$ subset ${\cal
  W}_n$ in $Met \times Z^1(Y, i\RR)$.
\label{wallstrata}
\end{thm} 

\noindent\underline{Proof:}
As in the proof of Theorem \ref{codim:1}, we only need to prove the
result for $(g, \nu)$ in $Met^0 \times Z^1(Y, i\RR)$ (see Claim 1 and Claim 2
in the proof of Theorem \ref{codim:1}).
Consider a real Hilbert bundle ${\cal L}$ over
\[ Met^0 \times Z^1(Y, i\RR) \times (L^2_1(S) - \{ 0\})\]
whose fiber over $(g, \nu, \psi)$ is
\[
{\cal L}_{(g, \nu, \psi)} = \{ \phi \in L^2_1(S) |
Re \langle \phi, i\psi \rangle_g = 0 \}
\]
Define a section $\zeta$ of ${\cal L}$ by assigning to $(g, \nu, \psi)$
the element $\partial_\nu^g \psi$.

\noindent\underline{\bf Claim:} $\zeta$ is transverse to the zero section of
${\cal L}$.

We need to prove that the differential map of $\zeta$ is surjective at
zeroes of $\zeta$. Suppose that $(g_0, \nu_0, \psi_0)$ (with $\psi\neq
0$) satisfies $\partial_{\nu_0}^{g_0} \psi_0 = 0$. 
Differentiating $\zeta$ with respect to the 
directions tangent to $Z^1(Y, i\RR) \times (L^2_1(S) - \{ 0\})$ only, we see
that the differential map is
\[ \begin{array}{rl}
 {\cal D}\zeta : \qquad
\Omega^1_{L^2_2} (Y, i\RR) \oplus L^2_1(S) & \longrightarrow
\{ \phi \in L^2_1(S) |
Re \langle \phi, i\psi \rangle_g = 0 \}\\[2mm]
(\nu_1, \psi_1) & \mapsto \partial_{\nu_0}^{g_0} \psi_1 +
\displaystyle{\frac 12} 
\nu_1. \psi_0.
\end{array}
\]
If $\phi \in {\cal L}_{(g_0, \nu_0, \psi_0)}$ is orthogonal to the image
of ${\cal D}\zeta$, then $\phi$ satisfies:
\[
\left\{\begin{array}{l}
(1). \ Re \langle \phi, i\psi_0 \rangle_{g_0} = 0, \\[2mm]
(2). \ Re \langle \phi, \nu_1. \psi_0 \rangle_{g_0} = 0, \qquad \hbox{ for any
$\nu_1 \in \Omega^1(Y, i\RR)$.}\\[2mm]
(3). \ \partial_{\nu_0}^{g_0} \phi = 0.
\end{array}\right.\]
{}From the second equation we see that there exists a
function $f: Y \to \RR$ such that $\phi = if \psi_0$. Substitute this into
equation (3), using $\partial_{\nu_0}^{g_0} \psi_0 = 0$. We obtain $df = 0$
which implies $f = C$ is a constant function.
Then
\[ Re \langle iC \psi_0 ,  i\psi_0 \rangle_{g_0} = C |\psi_0|^2 \]
implies $C =0$. Therefore, $\phi = 0$, that means, ${\cal D}\zeta$ is
surjective at $(g_0, \nu_0, \psi_0)$. It is easy to see that
the index of ${\cal D}\zeta$ is the index of $\partial_{\nu_0}^{g_0}$, which
is 1, since $i\psi_0$ is orthogonal to the image of
$\partial_{\nu_0}^{g_0}$. 

{}From the Claim, $\zeta^{-1} (0)$ is a Banach manifold and the projection
\[
\Pi: \zeta^{-1} (0) \to Met^0 \times \Omega^1(Y, i\RR)
\]
is a Fredholm operator of index 1. Note that for any
$(g, \nu)  \in  Met^0 \times \Omega^1(Y, i\RR)$ we have
$\Pi^{-1}(g,\nu) = Ker \partial_\nu^g$. 

Moreover, at $(g_0,\nu_0,\psi_0)$ we have
\[\left\{
\begin{array}{l}  Ker(\Pi_*) = \{ \phi | \partial^{g_0}_{\nu_0}\phi=0\} \\[2mm]
dim Ker(\Pi_*)  - dim Coker (\Pi_*)  = 1.
\end{array} \right.
\]
Therefore, $dim Coker (\Pi_*) = dim  Ker(\Pi_*) - 1$. Then
the Theorem follows,  with the top stratum of codimension one
described in Theorem \ref{codim:1}, and the stratification given by
$$ {\cal W}_n= {\cal X}^{-1}( {\cal F}_{0,2n}(S,Y) \cap {\cal X}({\cal M}et
\times Z^1(Y,\RR))), $$
where ${\cal X}$ is the map of (\ref{Diracmap}). Thus, the structure of
stratified set on ${\cal W}$ is induced by the structure on
$Image({\cal X})$ inside the stratified set ${\cal F}_0(S,Y)$.

\noindent $\diamond$

Notice that on any three-manifold it is possible to find special
metrics for which the dimension of the space of harmonic spinors is
arbitrarily large. The result for $S^3$ was proved by Hitchin
\cite{Hi} and recently generalized to all manifolds of dimension 3 mod
4 by B\"ar \cite{Bar}.

\subsection{Perturbation of flow lines}

We introduce suitable moduli spaces of gradient flow lines connecting
critical points. We prove in the Section 2.3, Proposition
\ref{transverse}  that, generically, these
are smooth manifolds that are cut out transversely, hence with the
dimension prescribed by the index theorem. This property depends on an
accurate choice of a class of perturbations for the gradient flow
equations.

Consider the space of connections and spinor sections $(\AA,\Psi)$ on $Y\times
\RR$ with the product metric $g+dt^2$, topologized with the
weighted Sobolev norms \cite{Fl} \cite{LM}. Here we choose the weight
\index{$e_\delta(t)$}
$e_\delta(t)=e^{\tilde\delta t}$, where $\tilde\delta$ is a smooth
function with bounded derivatives, $\tilde\delta : \RR\to
[-\delta,\delta]$ for some fixed positive number $\delta$, such that
$\tilde\delta (t)\equiv -\delta$ for $t\leq -1$ and $\tilde\delta
(t)\equiv \delta$ for $t\geq 1$. \index{$\| f \|_{2,k,\delta}$}
The $L^2_{k,\delta}$ norm is defined as $\| f \|_{2,k,\delta}=\|
e_\delta f \|_{2,k}$, where the Sobolev norms are defined with respect 
to a fixed reference connection $\AA_0$ in $L^2_{k,loc}$. By the
Sobolev multiplication theorem of Proposition \ref{sobolev} below, the 
spaces $L^2_{k,\delta, \AA_0}$ and $L^2_{k,\delta, \AA_1}$ are equal
whenever we have $\AA_0 - \AA_1 =\alpha_0 + \alpha_1$, with $\alpha_0$ 
in $L^2_{\ell,\delta, \AA_0}$ and $\alpha_1$ is ${\cal C}^\ell$
bounded, with $\ell \geq k-1$, cf. \cite{MorMrow}.
The weight $e_\delta$ imposes an exponential decay
as asymptotic condition along the cylinder. A proof of the following
Proposition can be found, for instance, in Section 9 of \cite{Eich}.

\begin{prop}
Let $Y$ be a compact oriented three-manifold endowed with a fixed
Riemannian metric $g_0$. Consider the cylinder $Y\times \RR$ with the
metric $g_0+dt^2$.
The weighted Sobolev spaces $L^2_{k,\delta}$ on the manifold $Y\times
\RR$ satisfy the following Sobolev embeddings.

(i) The embedding $L^2_{k,\delta}\hookrightarrow L^2_{k-1,\delta}$ is
compact for all $k\geq 1$.

(ii) If $k> m+2$ we have a continuous embedding
$L^2_{k,\delta}\hookrightarrow {\cal C}^m$.

(iii) If $k> m+3$ the embedding $L^2_{k,\delta}\hookrightarrow {\cal
  C}^m$ is compact.

(iv) If $2<k'$ and $k\leq k'$ the multiplication map 
$L^2_{k,\delta}\otimes L^2_{k,\delta}\stackrel{m}{\to}
L^2_{k,2\delta}$ is continuous.

\noindent Consider a metric $g_t +dt^2$ on the cylinder $Y\times
\RR$ such that for a fixed $T$ we have $g_t\equiv g_0$ for $t\geq T$
and $g_t\equiv g_1$ for $t\leq -T$ and $g_t$ varies smoothly when
$t\in [-1,1]$. The same Sobolev embedding theorems hold for the
$L^2_{k,\delta}$ spaces on $(Y\times \RR,g_t +dt^2)$.
\label{sobolev}
\end{prop}

Choose smooth representatives $(A_0,\psi_0)$ and $(A_1,\psi_1)$ of $a$
and $b$ in ${\cal M}$. Choose a smooth path $(A(t),\psi(t))$ such that
for $t\leq 0$ it satisfies $(A(t),\psi(t))\equiv (A_0,\psi_0)$ and for
$t\geq 1$ it is $(A(t),\psi(t))\equiv (A_1,\psi_1)$.
The configuration space ${\cal A}_{k,\delta}(a,b)$ is given by
the space of pairs $(\AA,\Psi)$ on $Y\times \RR$ satisfying
\beq
(\AA,\Psi)\in (A(t),\psi(t))+L^2_{k,\delta} (\Omega^1(Y\times\RR)\oplus
\Gamma (S^+\otimes L)).
\label{confspab} \eeq

Consider the group ${\cal G}_{k+1,\delta}(a,b)$ of gauge transformations in
$\tilde{\cal G}_{Y\times \RR}$, locally modeled on
$L^2_{k+1,\delta}(\Omega^0(Y\times\RR, i\RR))$, that approach elements
$\lambda_{\pm\infty}$ in the stabilizers $G_a$ and $G_b$ of
$(A_0,\psi_0)$ and $(A_1,\psi_1)$ as $t\to
\pm\infty$. This gauge 
group acts on ${\cal A}_{k,\delta}(a,b)$ and we can form the quotient
${\cal B}_{k,\delta}(a,b)$. There is an action of $\RR$ by translations
on ${\cal B}_{k, \delta}(a, b)$.
 
We consider the perturbed gradient flow equations for a path
$(A(t),\psi(t))$,
\beq
\frac{d}{dt}\psi(t)=-\partial_{A(t)}\psi (t)
\label{3SW1P'}
\eeq
and
\beq
\frac{d}{dt}A(t)=\sigma(\psi(t),\psi(t))-*F_{A(t)}+2i\rho
+2q_{(\AA,\Psi)}(t).
\label{3SW2P'} \index{$q_{(\AA,\Psi)}$}
\eeq

Equations (\ref{3SW1P'}) and (\ref{3SW2P'}) can be rewritten in terms
of pairs $(\AA,\Psi)$ in the form
\beq
D_{\AA}\Psi=0
\label{4SW1P'}
\eeq
and
\beq
F_{\AA}^+=\Psi\cdot\bar\Psi +i\mu +P_{(\AA,\Psi)},
\label{4SW2P'}
\eeq
as proved in Lemma \ref{rem43}. \index{$P_{(\AA,\Psi)}$}

The perturbation $P=*q+q\wedge dt$ is a function of ${\cal B}(a,b)$
to $\pi^*(\Omega^1(Y,i\RR))\cong\Omega^{2+}(Y\times\RR,i\RR)$, such that the
corresponding equations in a temporal gauge
(\ref{3SW1P'}) and (\ref{3SW2P'}) are preserved under the action of $\RR$
by reparameterization of the path $(A(t),\psi(t))$. The class of such
perturbations is described as follows.

\begin{defin}
\label{calP}
The space of perturbations ${\cal P}$ is the space of maps
\index{${\cal P}$}
\[ P:{\cal B}_{k,\delta}(a,b)\to
L^2_{k,\delta}\Omega^{2+}(Y\times\RR,i\RR), \] 
that satisfy the following conditions.

\noindent (1) $P_{(\AA,\Psi)}= *q_{(\AA,\Psi)}(t)
+q_{(\AA,\Psi)}(t) \wedge dt$, where $q_{(\AA,\Psi)}(t)$
satisfies 
\[ q_{(\AA,\Psi)^T}(t)=q_{(\AA,\Psi)}(t+T), \]
for any $T\in \RR$, where $(\AA,\Psi)^T$ is the $T$-translate of
$(\AA,\Psi)$. 

\noindent (2) the $L^2_{k,\delta}$-norm of the perturbation
$P_{(\AA,\Psi)}$ is bounded uniformly with respect to $(\AA,\Psi)$;

\noindent (3) the linearization ${\cal D}P_{(\AA,\Psi)}$ is a
pseudodifferential operator of order $<1$.
In particular ${\cal D}P_{(\AA,\Psi)}$ is a compact
operator from the $L^2_{k,\delta}$ to the
$L^2_{k-1,\delta}$ tangent spaces. 

\noindent (4) for all $l\leq k-1$, we have  
\[ \| q_{(\AA,\Psi)}(t)\|_{L^2_l}\leq C_l \| \nabla {\cal
    C}_\rho (A(t),\psi(t)) \|_{L^2_l} \]
in the $L^2_l$-norm on $Y\times \{ t\}$, where
\beq
\nabla {\cal C}_\rho (A(t),\psi(t))=(-\partial_{A(t)}\psi (t),
\sigma(\psi(t),\psi(t))-*F_{A(t)}+2i\rho)
\eeq
is the gradient flow of the functional ${\cal C}_\rho$, and $0<C_l<1$.

\noindent (5) The inequality
\[ \| {\cal D}q_{(\AA,\Psi)}(\alpha(t),\phi(t)) \| \leq C_{(\AA,\Psi)}
\| \nabla{\cal C}_\rho (A(t),\psi(t)) \|\cdot \| (\alpha(t),\psi(t))
\| \]
holds for $t\geq T_0$. Here ${\cal D}q_{(\AA,\Psi)}$ is the
linearization of the perturbation $q_{(\AA,\Psi)}$.

\end{defin}

With a perturbation in the class ${\cal P}$ the equations
(\ref{3SW1P'}) and (\ref{3SW2P'}) are invariant 
with respect to the action of $\RR$ by translations along the gradient
flow lines, that is if $(A(t), \psi(t))$ is a solution of (\ref{3SW1P'}) and
(\ref{3SW2P'}), then $(A(t+T), \psi(t+T))$ is also a solution for any
$T\in \RR$. 

An example of 
perturbation with these properties has been constructed by
Fr\/oyshov \cite{Fr}. 

\begin{prop}
The class of perturbations introduced by Fr\/oyshov in \cite{Fr} is 
contained in our class ${\cal P}$. 
\label{froyshov}
\end{prop}

According to Fr\/oyshov's construction, for fixed smooth compactly
supported functions $\eta_1$, $\eta_2$, with 
$supp(\eta_1)\subset [-1,1]$ and
$\eta_2|_I(t)=t$ on an interval $I$ containing all the critical values
of ${\cal C}_\rho$, a function $h: {\cal B}_{k,\delta}(a,b)\to {\cal
  C}^m (\RR)$ is defined as 
\[ h_{(\AA,\Psi)}(T)= \int_{\RR}\eta_1(s-T) \eta_2
(\int_{\RR}\eta_1(t-s) {\cal C}_\rho(A(t),\psi(t))dt)ds. \]

Let $\Omega_\Xi^2(Y\times\RR)$ be the set of ${\cal C}^m$ 2-forms $\omega$ 
that are compactly supported in $Y\times\Xi$, where $\Xi$ is the
complement of a union of small intervals centered at the critical
values of ${\cal C}_\rho$. Fr\/oyshov's perturbation is constructed by setting
\[ P_{(\AA,\Psi)}=(h_{(\AA,\Psi)}^*(\omega))^+, \]
where $h_{(\AA,\Psi)}^*(\omega)$ is the pullback of $\omega$ along the
map $Id_Y\times h_{(\AA,\Psi)}:Y\times \RR\to Y\times \RR$.

As shown in \cite{Fr}, the function $h_{(\AA,\Psi)}$ is bounded with
all derivatives, uniformly with respect to 
$(\AA,\Psi)$. Moreover, by the choice of $\Xi$, the perturbation
$h_{(\AA,\Psi)}^*(\omega)$ is smooth and compactly supported, hence in
$L^2_{k,\delta}$.
  
Condition (1) holds, since the function $h_{(\AA,\Psi)}(t)$ satisfies 
\[ h_{(\AA,\Psi)}(t+\tau)=h_{(\AA,\Psi)^\tau}(t), \]
where $(\AA,\Psi)^\tau$ is the $\tau$-reparameterized solution
represented in a temporal gauge by $(A(t+\tau),\psi(t+\tau))$.
In fact,
\[ h_{(\AA,\Psi)^\tau}(T)=\int_{\RR}\eta_1(s-T) \eta_2
(\int_{\RR}\eta_1(t-s) {\cal C}_\rho(A(t+\tau),\psi(t+\tau))dt)ds= \] 
\[ =\int_{\RR}\eta_1(s-T) \eta_2(\int_{\RR}\eta_1(u-s-\tau) {\cal
  C}_\rho(A(u),\psi(u))du)ds= \]
\[ = \int_{\RR}\eta_1(v-T-\tau) \eta_2(\int_{\RR}\eta_1(u-v) {\cal
  C}_\rho(A(u),\psi(u))du)dv=h_{(\AA,\Psi)}(T+\tau). \]

Condition (2) holds: in fact, it is shown in \cite{Fr} that the function
$h_{(\AA,\Psi)}$ is bounded with all derivatives, uniformly with respect to 
$(\AA,\Psi)$. The Sobolev embeddings of Proposition \ref{sobolev}
provide the uniform bound in the $L^2_{k,\delta}$-norms.

Condition (3) and (5): we can write the function $h_{(\AA,\Psi)}$ with the
notation 
\[ h_{(\AA,\Psi)}=\eta_1 * (\eta_2 (\eta_1 * {\cal
C}_\rho(A,\psi))), \]
where $*$ denotes the convolution product on $\RR$.
We obtain the variation with respect to $(\AA,\Psi)$ of the form
\[ v_{(\AA,\Psi)}(\alpha,\phi)=\eta_1 * \left( \eta_2^\prime (\eta_1 * {\cal
C}_\rho(A,\psi)) \eta_1 * \langle \nabla{\cal C}_\rho (A,\psi),
(\alpha,\phi) \rangle \right). \]
Thus, as shown in Fr\/oyshov (\cite{Fr}, Prop.5), for
$\omega$ a ${\cal C}^m$ form, the linearization of the perturbation
$h_{(\AA,\Psi)}^*(\omega)$ at the point $(\omega,\AA,\Psi)$ is a
  bounded operator $K_{(\omega,\AA,\Psi)}:L^2_{k,\delta}\to {\cal C}^m$ with 
\[ supp\left(K_{(\omega,\AA,\Psi)}(\alpha,\Phi)\right)\subset
h^{-1}_{(\AA,\Psi)}(\Xi)\times Y. \] 
Thus Condition (3) follows, since $\eta_1$ and $\eta_2$ are smooth compactly
supported functions. In particular, if $\omega$ is a smooth form, the
linearization ${\cal D}P_{(\AA,\Psi)}$ is a smoothing operator. The
expression of the variation $v_{(\AA,\Psi)}$ also shows that Condition
(5) holds.

Condition (4) follows from an estimate on $\| dh_{(\AA,\Psi)} \|$,
\[ \| dh_{(\AA,\Psi)}(t) \|\leq C \| \nabla{\cal C}_\rho
(A(t),\psi(t)) \|, \] 
where the constant $C$ only depends on the asymptotic values $a$ and
$b$ of $(\AA,\Psi)$.
 
Other perturbations of the functional ${\cal C}$ can be used to achieve
transversality of the moduli space of flow lines. For instance see the 
discussion in \cite{Kronheimer}, and the class of perturbations
introduced in \cite{CMW}. The perturbations introduced in \cite{CMW}
have the advantage of being defined directly as perturbations of the
Chern--Simons--Dirac functional ${\cal C}$, instead of being
perturbations of the 4-dimensional Seiberg--Witten
equations, as in the case of the class considered here.

\subsection{Transversality of ${\cal M}(a,b)$}

Let ${\cal L}_{(\AA,\Psi)}$ be the linearization of equations
(\ref{4SW1P'}) and (\ref{4SW2P'}) on ${\cal B}_{k,\delta}(a,b)$.

The operator ${\cal L}$ is of the form 
\beq {\cal L}_{(\AA,\Psi,P)}(\alpha,\Phi)=\left\{\begin{array}{c}
D_{\AA}\Phi+ \alpha \Psi \\ d^+\alpha
-\frac{1}{2}Im(\Psi\cdot\bar\Phi)+{\cal D}P_{(\AA,\Psi)}(\alpha,\Phi)\\
G^*_{(\AA,\Psi)}(\alpha,\Phi)
\end{array}\right. \label{lin:P} \eeq
mapping
\[ L^2_{k,\delta}(\Omega^1(Y\times\RR,i\RR)\oplus \Gamma(S^+\otimes
L))\to L^2_{k-1,\delta} (\Omega^0(Y\times\RR,i\RR)\oplus \Omega^{2+}
(Y\times\RR,i\RR)). \]
The operator $G^*$ is the adjoint, in the $L^2_{0,\delta}$-pairing, of
the linearization of the gauge 
group action $G_{(\AA,\Psi)}(f)=(-df,f\Psi)$. 

As the following proposition shows,
the operator ${\cal L}_{(\AA,\Psi,P)}$ is obtained by adding the 
compact perturbation ${\cal D}P_{(\AA,\Psi)}$ to a Fredholm map
from $L^2_{k,\delta}$ to $L^2_{k-1,\delta}$, hence it is
Fredholm. Therefore \index{${\cal D}P_{(\AA,\Psi)}$}
we have a well defined relative Morse index of two critical points $a$
and $b$ in ${\cal M}$.

\begin{prop}
Suppose $a$ and $b$ are irreducible critical points for the functional
${\cal C}_\rho$. Let $\{ \lambda_a \}$ and $\{ \lambda_b \}$ be the
eigenvalues of the Hessian $T$ at the points $a$ and $b$.
Assume that the positive number $\delta$ satisfies the condition
$$\delta < \min \{ |\lambda_a|, |\lambda_b| \}.$$ 
Let $(\AA,\Psi)$ be a solution of (\ref{4SW1P'}) and (\ref{4SW2P'}) in
${\cal B}_{k,\delta}(a,b)$. Then the linearization ${\cal
  L}_{(\AA,\Psi)}$ is a Fredholm operator of index
\[ Ind({\cal L}_{(\AA,\Psi)})=\sigma(a,b). \]
The right hand side $\sigma(a,b)$ is the spectral flow of the operator
$\nabla{\cal F}$ along a path $(A(t),\psi(t))$ in ${\cal A}$
that corresponds to $(\AA,\Psi)$ under $\pi^*$. 
The quantity $\sigma(a,b)$ is independent of the path, hence
\[ \sigma(a,b)=\mu(a)-\mu(b) \] \index{$\mu(a)-\mu(b)$}
defines a relative Morse index of $a$ and $b$, where $\mu(a)$ is the
spectral flow of $\nabla{\cal F}$ on a path joining $a$ to a fixed
$[A_0,\psi_0]$ in ${\cal M}$.
\label{relmorseind}
\end{prop}

\begin{rem}
In the case with $Y$ a homology sphere, Proposition \ref{relmorseind}
holds for a reducible point $a=[\nu,0]$ under the assumption that, for
the chosen metric $g$ on $Y$, the condition
$Ker(\partial^g_\nu)=0$ is satisfied.
\end{rem}

We state and prove an analogue of
Proposition \ref{relmorseind} in the context of framed moduli spaces and
equivariant theory in Theorem \ref{relmorseind2}. We also prove, in
Theorem \ref{relmorseind2}, that the relative Morse index
of points in $\tilde{\cal M}$ is well defined. 

Consider the moduli space ${\cal M}(a,b)$ of solutions of the
equations (\ref{4SW1P'}) and (\ref{4SW2P'}) in ${\cal
  B}_{k,\delta}(a,b)$. 

\begin{prop}
When $a$ and $b$ are irreducible critical points of ${\cal C}_\rho$,
for a generic choice of the perturbation $P\in {\cal P}$, the moduli space 
${\cal M}(a,b)$ of gradient flow lines is a smooth oriented manifold, cut out
transversely by the equations, of dimension 
\[ \dim({\cal M}(a,b))=\mu(a)-\mu(b), \]
where $\mu(a)-\mu(b)$ is the relative Morse index of the critical points.
\label{transverse}
\end{prop}

\noindent\underline{Proof:} It is first necessary to know that there
are no reducible flow lines connecting the critical points $a$ and
$b$. This fact is an easy consequence of the definition of the
configuration space (\ref{confspab}), since the exponential weight in
the Sobolev norm forces elements in (\ref{confspab}) to decay at the
ends to the asymptotic values, which are irreducible by assumption.
On the convergence of flow lines to the endpoints $a$ and
$b$, see the results of Section 3.4.
Thus, ${\cal M}(a,b)$ lies entirely in the
irreducible component ${\cal B}^\prime_{k,\delta}(a,b)$, provided that 
at least one of the endpoints $a$ and $b$ is irreducible.
The statement then follows via the implicit function theorem, upon
showing that, for a 
generic choice of the perturbation $P$, the linearization ${\cal L}$
is surjective.  

Consider the operator 
\[ \hat{\cal L}_{(\AA,\Psi,P)}(\alpha,\Phi,p)={\cal
  L}_{(\AA,\Psi,P)}(\alpha,\Phi) +p_{(\AA,\Psi,P)}(\alpha,\Phi), \]
where we vary the perturbation by an element $p_{(\AA,\Psi,P)}$ of
$T_P{\cal P}$. This corresponds to varying the parameter
$\omega\in\Omega_\Xi^2(Y\times \RR)$ in Fr\/oyshov's class of
perturbations. 

The operator ${\cal L}$ is Fredholm, therefore $\hat{\cal L}$ has a
closed range. We show that $\hat{\cal L}$ is surjective by proving
that it has dense range.
 
Suppose given an element $(\beta,\xi,g)$ in
$L^2_{-k-1,-\delta}$ that is $L^2$-orthogonal 
to the range of the operator $\hat{\cal L}$. Here $\beta$ is an element in 
$\Omega^{2+}(Y\times\RR, i\RR))$, $\xi$ is a spinor, and $g$ is a
zero-form. The element $(\beta,\xi,g)$ is in the kernel of the adjoint
$\tilde{\cal L}^*$, which is an elliptic operator with
$L^2_{-k,-\delta}$ coefficients, thus $(\beta,\xi)$ lives in
$L^2_{-k,-\delta}$ by elliptic regularity. In fact, the perturbation
satisfies Condition (3) of Definition \ref{calP}.
If we consider the $L^2$-pairing of $L^2_{k,\delta}$ and $L^2_{-k,-\delta}$,
we get 
\[ \langle \beta, d^+\alpha -\frac{1}{2}Im(\Psi\cdot\bar\Phi))+{\cal
  D}P_{(\AA,\Psi)}(\alpha,\Phi)+p_{(\AA,\Psi,P)}(\alpha,\Phi) \rangle
+ \] 
\[ +\langle \xi, D_{\AA}\Phi+\alpha\Psi \rangle+ 
\langle g,G^*_{(\AA,\Psi)}(\alpha,\Phi)\rangle=0. \]

By varying $p\in {\cal P}$ we force $\beta\equiv 0$.
The remaining inner product 
\[ \langle \xi, D_{\AA}\Phi+i\alpha\Psi \rangle+
\langle g,G^*_{(\AA,\Psi)}(\alpha,\Phi)\rangle=0 \]
gives the following equations

(a) $(e_{-\delta}de_{\delta})g=\frac{1}{2}\xi\cdot\bar\Psi$ and

(b) $D_A\xi-g\Psi=0$.

\noindent We assume that $\Psi$ is not identically zero. Applying
$d^*$ to (a) and using (b) we obtain
$d^*(e_{-\delta}de_{\delta}g)+g|\Psi|^2=0$. Equivalently, we get
\[ (e_{\delta/2}d^*e_{-\delta/2}) (e_{-\delta/2}de_{\delta/2})
e_{\delta/2}g +|\Psi|^2 e_{\delta/2}g=0. \] 
The equation 
\[ \Delta_{\delta/2} e_{\delta/2}g + e_{\delta/2}g|\Psi|^2=0, \] 
with
\[ \Delta_{\delta}=e_{-\delta}\Delta e_{\delta}, \] 
implies that $g\equiv 0$, since $g$ decays at $\pm\infty$ and the
maximum principle applies. Then, by varying 
$\alpha$ alone in $\langle \xi, D_{\AA}\Phi+\alpha\Psi \rangle=0$,
we force $\xi$ to vanish on some arbitrary open
set. We obtain $(\beta,\xi)\equiv 0$.

Thus the operator $\tilde{\cal L}$ is surjective. This implies that
zero is a regular value for the map defined by the equations
(\ref{4SW1P'}) and (\ref{4SW2P'}). Therefore the moduli space  ${\cal
  M}od$ of triples $([\AA,\Psi],P)$ in ${\cal B}^\prime_{L^2_{k,\delta}}(a,b)
  \oplus {\cal P}$ that satisfy the equations is a 
smooth (infinite dimensional) manifold with tangent space
$Ker(\hat{\cal L})$. 

The projection $\Pi: {\cal M}od \to {\cal P}$ given by
$\Pi([\AA,\Psi],P)=P$ linearizes to a surjective Fredholm operator
${\cal D}\Pi:Ker(\tilde{\cal L})\to T_P{\cal P}$. The kernel of ${\cal
  D}\Pi$ is $Ker({\cal D}\Pi_{(\AA,\Psi,P)})=Ker({\cal L}_{(\AA,\Psi,P)})$.
The infinite dimensional Sard
theorem implies that the moduli space ${\cal M}(a,b)$, 
for a generic perturbation $P\in {\cal P}$, is
the inverse image under the projection map from ${\cal M}od$ to ${\cal
  P}$ of a regular value. Thus ${\cal M}(a,b)$ is a
smooth manifold which is cut out transversely by the equations.
Equivalently, the linearization ${\cal L}$ with a fixed generic $q$ is
surjective.  

The virtual dimension of the moduli space ${\cal M}(a,b)$ equals the
index of the Fredholm operator ${\cal L}$. According to Proposition
\ref{relmorseind}, this is the relative Morse index $\mu(a)-\mu(b)$.

The orientation of ${\cal M}(a,b)$ is given by a trivialization of the
determinant line bundle of the operator ${\cal L}$. This is obtained
given a choice of an orientation of 
\[ H^0_\delta(Y\times\RR)\oplus H^{2+}_\delta(Y\times\RR)\oplus
H^1_\delta(Y\times\RR), \]
the cohomology groups of $\delta$-decaying forms, as discussed in the
following Proposition.

\noindent $\diamond$ 

\begin{prop}
\label{orientation}
The manifold ${\cal M}(a,b)$ is oriented by a trivialization of the 
determinant line bundle of the operator ${\cal L}$. This is obtained
from an orientation of 
\[ H^0_\delta(Y\times\RR)\oplus H^{2+}_\delta(Y\times\RR)\oplus
H^1_\delta(Y\times\RR), \]
the cohomology groups of $\delta$-decaying forms.
\end{prop}

\noindent\underline{Proof:} Suppose given $[x]=[\AA,\Psi]\in {\cal
M}(a,b)$. Let $x=(A(t),\psi(t))$ be a temporal gauge representative
such that
\[ \lim_{t\to\infty}(A(t),\psi(t))=(A_b,\psi_b), \]
and
\[ \lim_{t\to -\infty}(A(t),\psi(t))=(A_a,\psi_a). \]
Consider the family of operators \index{${\cal L}_x$}
\[ {\cal L}_x(\alpha,\Phi)=\left\{ \begin{array}{l} d^+\alpha
-\frac{1}{2}Im(\Psi\cdot\bar\Phi) \\
D_{\AA}\Phi+\alpha\Psi \\
G^*_{(\AA,\Psi)}(\alpha,\Phi) \end{array}\right. \]
acting on the space of $L^2_{k,\delta}$-decaying 1-forms and spinor sections
on $Y\times \RR$.

An orientation of the moduli space ${\cal M}(a,b)$ is determines by a
trivialization of the determinant line bundle of the family of
operators ${\cal L}_x$. We can separate ${\cal L}_x$ in the first
order term and a perturbation,
\[ {\cal L}_x={\cal L}^1_x + {\cal L}^0_x, \]
with
\[ {\cal L}^1_x= d^+ + d^*_\delta + D_{\AA} \]
and 
\[ {\cal L}^0_x =  \left(\begin{array}{c} 
-\frac{1}{2}Im(\cdot\bar\Psi) \\ i<\Psi,\cdot> \end{array}\right). \]
The operators induced by ${\cal L}_x$ on the asymptotic ends are 
\[ L(\pm\infty)= L^1(\pm\infty)+L^0(\pm\infty), \]
where we have
\[ L^1(\infty)=\left( \begin{array}{ccc} *d & d^* & 0 \\ d & 0 & 0 \\ 0 &
0 & \partial_{A_b} \end{array} \right), \]
\[ L^1(-\infty)=\left( \begin{array}{ccc} *d & d^* & 0 \\ d & 0 & 0 \\ 0 &
0 & \partial_{A_a} \end{array} \right), \]
and
\[ L^0(\infty)= \left( \begin{array}{ccc} 0 & 0 & -2iIm\sigma(\cdot,
\psi_b) \\ 0 & 0 & -i<\cdot, \psi_b> \\ \cdot\psi_b & \psi_b & 0
\end{array} \right), \]
\[ L^0(-\infty)=\left( \begin{array}{ccc} 0 & 0 & -2iIm\sigma(\cdot,
\psi_a) \\ 0 & 0 & -i<\cdot, \psi_a> \\ \cdot\psi_a & \psi_a & 0
\end{array} \right), \]
acting on $\Omega^0(Y,i\RR)\oplus\Omega^1(Y,i\RR)\oplus\Gamma(S\otimes
L)$.

We can consider a
deformation ${\cal L}^\epsilon_x$ of the family ${\cal L}_x$ obtained
as in \cite{Mo} Section 6.6,
\[ {\cal L}^\epsilon_x={\cal L}^1_x +(1- \epsilon) {\cal L}^0_x. \]
The deformation changes the asymptotic operators in the form
\[ L^\epsilon(\pm\infty)= L^1(\pm\infty)+(1-\epsilon) L^0(\pm\infty), \]
We can guarantee that this is a deformation via Fredholm
operators provided that the weight $\delta$ is chosen such that
$\delta/2$ is not in the spectrum of $L^\epsilon(\pm\infty)$ for all
$\epsilon$ (\cite{LM} Theorem 6.2, and \cite{MMR} Lemma 8.3.1). That
is, if the spectrum of the operators $L^\epsilon(\pm\infty)$ is
uniformly bounded away from zero.

If this is the case, then
a trivialization of the determinant line of the family ${\cal
L}^\epsilon_x$ is obtained by a trivialization at $\epsilon=1$. This
induces a trivialization of ${\cal L}_x$. The trivialization at
$\epsilon=1$ is a trivialization of the determinant line of 
the operator $d^+ + d^*_\delta + D_{\AA}$.

The Dirac operator is complex linear and it preserves the orientation
induced by the complex structure on the spinor bundle  $S^+\otimes L$.
Thus a trivialization is obtained by an orientation of 
\[ H^0_\delta(Y\times\RR)\oplus H^{2+}_\delta(Y\times\RR)\oplus
H^1_\delta(Y\times\RR), \]
the cohomology groups of $\delta$-decaying forms, \cite{Mo},
\cite{MST}. 

However, the condition on the spectrum of $L^\epsilon(\pm\infty)$ may
not always be satisfied: the deformation
${\cal L}^\epsilon_x$ may not be through Fredholm operators. It is
still possible to obtain an orientation of ${\cal M}(a,b)$:
the following argument was suggested to us by L. Nicolaescu
\cite{Ni}. 

We can change the family ${\cal L}_x$ by a deformation such that
the new family ${\cal H}_x$ satisfies
\[ {\cal H}_x |_{Y\times [-1,1]} \equiv {\cal L}_x |_{Y\times [-1,1]}. \]
On $Y\times (-\infty,-2] \cup Y\times [2,\infty)$ it satisfies
\[ {\cal H}_x |_{Y\times (-\infty,-2]} = d^+ + d^*_\delta + D_{A_a} +
 \left(\begin{array}{c} 
-\frac{1}{2}Im(\cdot\bar\psi_a) \\ i<\psi_a,\cdot>
\end{array}\right), \]
\[ {\cal H}_x |_{Y\times [2,\infty)} = d^+ + d^*_\delta + D_{A_b} +
 \left(\begin{array}{c} 
-\frac{1}{2}Im(\cdot\bar\psi_b) \\ i<\psi_b,\cdot>
\end{array}\right). \]

A trivialization of the determinant line of ${\cal H}_x$ induces a
trivialization of the determinant of ${\cal L}_x$.

The index of the family $Ind({\cal H}_x)$ on $Y\times \RR$ equals the
index on $Y\times [-3,3]$ with APS boundary conditions \cite{APS},
\cite{Ni3}. 

Now we can consider the family of operators $\tilde{\cal H}_x$ that
satisfies 
\[ \tilde{\cal H}_x |_{Y\times [-3,3]}\equiv {\cal H}_x |_{Y\times
[-3,3]}. \]
On $Y\times (-\infty,-4]\cup Y\times [4,\infty)$ it satisfies
\[ \tilde {\cal H}_x |_{Y\times (-\infty,-4]} =d^+ + d^*_\delta +
D_{A_a}, \]
\[ \tilde {\cal H}_x |_{Y\times [4,\infty)} = d^+ + d^*_\delta +
D_{A_b}. \]

On the cylinder $Y\times [-4,-3]$ and $Y\times [3,4]$ the index of the
operator $\tilde {\cal H}_x$ can be obtained as
$Ind(\frac{\partial}{\partial t}+ H(t))$, where $H(t)$ is independent
of $x$ (but depends on the asymptotic values $a$ or $b$). Thus,
$Ind(\tilde {\cal H}_x)$ on $Y\times [-4,-3]$ and $Y\times [3,4]$ is
the spectral flow $SF(H(t))$ of $H(t)$.

The indices of ${\cal H}_x$ and $\tilde{\cal H}_x$ are related by the
excision formula
\[ Ind({\cal H}_x)-SF(H(t)) = Ind(\tilde {\cal H}_x), \]
generalizing the excision formula of \cite{Ni2}.

Thus the relative orientation of ${\cal H}_x$ and by $\tilde {\cal
H}_x$ is exactly  $(-1)^{SF(H)}$. Notice that \cite{Ni2}
provides examples where this spectral flow is computed explicitly and
is odd.

Finally we can introduce a deformation $\tilde{\cal H}^\epsilon_x$
with a homotopy that shrinks to zero
the spinor part. In this case the asymptotic operators $H(\pm\infty)$
remain constant, hence the deformation is through Fredholm
operators. This implies that the orientation determined by
$\tilde{\cal H}_x$ is the same as the one determined by $d^+ +
d^*_\delta + D_{\AA}$, that is by an orientation of
\[ H^0_\delta(Y\times\RR)\oplus H^{2+}_\delta(Y\times\RR)\oplus
H^1_\delta(Y\times\RR). \]

\subsection{Floer homology}

The Floer complex has generators 
\[ C_q =\{ b\in {\cal M}' \mid \mu(b)=q \}, \]
where ${\cal M}'$ is the irreducible part of the moduli space of
critical points. The Morse index $\mu(b)$ is computed with respect to
a fixed element $[A_0,\psi_0]$ in ${\cal M}$,
$[A_0,\psi_0]=[\nu,0]$ in the case of a homology sphere.

The boundary operator is given by 
\beq
\partial a=\sum_{b\mid \mu(a)-\mu(b)=1} \epsilon(a,b) b,
\label{boundary}
\eeq
where $\epsilon(a,b)$ is the algebraic sum over the paths joining $a$
and $b$ of the signs given by the orientation,
\[ \epsilon(a,b)=\sum_{\gamma\in\hat{\cal M}(a,b)} \epsilon_\gamma. \]

Here ${\cal M}(a,b)$ is the moduli space of flow lines on $Y\times\RR$
\index{${\cal M}(a,b)$} \index{$\hat {\cal M}(a,b)$}
with asymptotic values $a$ and $b$. A description of ${\cal M}(a,b)$
will be given in the next section. The space $\hat{\cal M}(a,b)$ is the  
quotient of ${\cal M}(a,b)$ by the action of $\RR$ by translations.
A compactness result for $\hat{\cal M}(a,b)$ is needed in order to
make sense of $\epsilon(a,b)$. This result will follow from the more
general result proved in Section 4 in the equivariant setup, in
Theorem \ref{seqcompact}, Theorem \ref{equivgluing}, and Proposition
\ref{equivgluing:k}. 

The property that $\partial\circ \partial=0$ relies on the gluing
formula of Lemma \ref{glue} that follows from an accurate analysis of
the properties of the gradient flow moduli space ${\cal M}(a,b)$. In
our setting, again, this result will follow from the more general
results in the equivariant context, see Theorem \ref{equivgluing}. 

In the case where $Y$ is a homology sphere,
we need to ensure that in the expression of $\partial^2$
there is no contribution coming from trajectories that break through
the unique reducible solution. In other words, no component of the
form $\hat{\cal M}(a,\theta)\times \hat{\cal M}(\theta,b)$ can appear in
the boundary of $\hat{\cal M}(a,b)$ for $\mu(a)-\mu(b)=2$. This has been
proved in \cite{Wa}. In fact, the following gluing formula holds.

\begin{lem}
Suppose given $a$, $b$ and $c$ in ${\cal M}'$, irreducible critical
points with $\mu(a)>\mu(b)>\mu(c)$. Then, for large enough $T$, there
is a local diffeomorphism 
\[ \hat{\cal M}(a,b)\times \hat{\cal M}(b,c)\times [T,\infty) \to
\hat{\cal M}(a,c). \]
If $Y$ is a homology sphere and $\theta=[\nu,0]$ is the unique reducible
critical point, there is a local diffeomorphism
\[ \hat{\cal M}(a,\theta)\times \hat{\cal M}(\theta,c)\times U(1)\times
[T,\infty) \to \hat{\cal M}(a,c). \]
In this case $U(1)$ is the stabilizer of the reducible solution
$\theta=[\nu,0]$. 
\label{glue}
\end{lem}

In Theorem \ref{equivgluing}, we prove the gluing formula in the
equivariant setup. A proof of Lemma \ref{glue} can be found in
\cite{Ma} and \cite{Wa}, and it follows from our equivariant result,
as discussed at the end of Section 4.2.
The result has an immediate corollary.

\begin{corol}
Suppose $Y$ is a homology sphere. Let $\theta$ be the unique reducible
solution. If $a$ and $c$ are irreducible critical points such that
$\mu(a)-\mu(c)=2$, then generically there will be no boundary strata
of the form $\hat{\cal M}(a,\theta)\times\hat{\cal M}(\theta,c)$. In fact for
dimensional reasons the moduli space $\hat{\cal M}(\theta,c)$ 
of gradient flow lines is generically empty if $\mu(\theta)-\mu(c)=1$. 
\label{empty1}
\end{corol}

The property that $\partial^2=0$ in the Floer complex follows then
{}from the fact that the matrix elements
\[ <\partial\partial a,c>=\sum_b <\partial a,b><\partial b,c>, \]
of the operator $\partial^2$ are the algebraic sum of the points of
the oriented zero-dimensional manifold
\[ \cup_{b\in {\cal M}'}\hat{\cal M}(a,b)\times\hat{\cal M}(b,c)=
\partial\hat{\cal M}(a,c). \]

In the case with non-trivial $b^1(Y)$, we construct the Floer homology
under the assumption that $c_1(L)\neq 0$ rationally and that the
perturbation is restricted to the trivial cohomology class
$[*\rho]=0$. 

\begin{rem}
\label{Ybadcase}
For a 3-manifold $Y$ with $b_1>0$ and a $Spin_c$-structure with
$c_1(L)=0$ rationally, there is no nice way at this stage
to formulate the corresponding monopole homology. In fact, the 
Chern--Simons--Dirac functional is $\RR$-valued on the configuration 
space as long as the perturbation term represents a trivial de Rham 
cohomology class. In this case, however, the condition $c_1(L)=0$
implies the existence of a reducible set of critical points that is a
torus $T^{b_1(Y)}$. These can be degenerate, in the sense of
Morse-Bott, even in the framed configuration space. Perturbing the
functional with a 1-form $\rho$ that is non-trivial cohomologically
can destroy this reducible set, but the functional would no longer be
$\RR$-valued. There is then no uniform energy bound on the space of
flow lines of a fixed virtual dimension. This creates a problem in the
compactification by broken trajectories (see Theorem \ref{seqcompact}
and Theorem \ref{equivgluing}). The right framework for this bad case
seems to be a Novikov type complex, where trajectories with the same
virtual dimension but with different energies are counted separately
as coefficients of a power series. Since this case has important
applications in the gluing formulae, we deal with it separately in
\cite{MW4}. 
\end{rem}

Under our assumptions, when we have $b^1(Y)>1$ and a nontrivial $c_1(L)$, 
the Floer homology groups can be proved to be
independent of the metric and of the perturbation. In fact a chain map
and a chain homotopy are constructed by a cobordism argument between
moduli spaces for two different metrics and perturbations. 
This result follows from the topological invariance of the equivariant 
Floer homology proved in Theorem \ref{metrics}, in Section 6, and the
equivalence of equivariant and non-equivariant Floer theories in the
case of manifolds with $b^1(Y)>1$ and a nontrivial $c_1(L)$, proved in 
Section 7, Theorem \ref{isomorphism}.

In the case with $b^1(Y)=1$ one
expects to find a dependence on the choice of the perturbation, see
\cite{Lim}, 
however, since we are only considering perturbations that are
cohomologically trivial, we obtain independence of the metric and
perturbation as in the $b^1(Y)>1$ case.
A similar dependence was detected in \cite{Au2} in the case of the invariant
of three-manifolds obtained by counting points in $\tilde{\cal M}$ with the
orientation. In \cite{Ma} it is proved in the case $b^1(Y)>0$ that
this invariant is in fact the Euler characteristic of the Floer homology.
The same invariant was introduced in \cite{CMWZ2}
following the Quantum Field Theory formulation of Seiberg-Witten
theory. 

In the case of a homology sphere the metric dependence problem is more
complicated. In fact due to the reducible solution a cobordism
argument does not work and more generally the construction of a chain
homotopy can fail due to the presence of moduli spaces of gradient
flow lines that connect the irreducibles to the reducible critical
point. The space of metrics and perturbations breaks into chambers
with codimension-one walls, so that the Floer groups are isomorphic
for metrics that belong to the same chamber and are in general
non-isomorphic when the metric crosses a wall. We shall discuss the
wall-crossing phenomenon in Section 7.

\section{Morse-Bott theory}

We are now going to introduce the equivariant Floer complex. This can
be constructed for all three-manifolds. Clearly in the case of an
integral homology sphere there will be no question of different
$Spin_c$-structures. In the case of a rational homology sphere, we
have finitely many choices of $Spin_c$-structures and there is a
reducible point corresponding to each of these. 
In the case of manifolds with $b^1(Y)>0$ there are infinitely many
possible choices of $Spin_c$-structures. This gives rise to a family
of Floer complexes corresponding to the different choices of the
$Spin_c$-structures. The relation between these requires further
investigation and is analyzed elsewhere (cf. \cite{CMW}). In all
the following we always assume to work with a fixed choice of the
$Spin_c$-structure. 

In order to consider reducible as well as irreducible generators,
we introduce a framed configuration space with a $U(1)$-action, where
the functional ${\cal C}_\rho$ is defined as a $U(1)$-invariant real
valued functional.
In order to apply the analogue of the finite dimensional equivariant
Morse theory \cite{AB2}, we need ${\cal C}$ to be a Morse-Bott
function. That is, we have to ensure that the Hessian is
non-degenerate on the normal bundle to the critical $U(1)$-orbits. 

\subsection{Framed moduli space}

\begin{defin}
Let $x_0$ be a fixed base point in $Y$.\index{${\cal B}^0$}
We define the space ${\cal B}^0$ to be the quotient of ${\cal A}$
with respect to the action of the subgroup ${\cal G}^0\subset {\cal
  G}$ of gauge transformations $\lambda$ that act as the \index{${\cal G}^0$}
identity on the fiber of $S\otimes L$ over $x_0$ and that satisfy the
condition (\ref{int=0}).
The space ${\cal B}^0$ is the framed configuration space. 
\label{framed}
\end{defin}

The action of the group ${\cal G}^0$ on ${\cal A}$ is free,
therefore the space ${\cal B}^0$ is an infinite dimensional
Banach manifold (using a fixed $L^2_k$ norm) that carries a residual $U(1)$
action. There is a fibration ${\cal B}^0\to{\cal B}$ over the unframed
configuration space with fiber $U(1)$.    
The solutions of the three dimensional Seiberg-Witten equations
(\ref{3SW1}) and (\ref{3SW2}) in ${\cal B}^0$ form the framed moduli
space ${\cal M}^0$, that is the critical
set of the functional (\ref{functional}) modulo based gauge
transformations. 

As in the case of Donaldson theory \cite{D2}, an equivalent
description of the framed configuration space can be given as the triples
$(A, \psi, \phi)$ with $(A,\psi)\in {\cal A}$ and $\phi$ a unit vector
in the fiber $S\otimes L |_{x_0}$. The full gauge group acts freely on
this space. Solutions of the Seiberg-Witten equations in this
configuration space modulo the full gauge group provide another model
of framed moduli space. This has been used in \cite{Au1} and
\cite{Au2}. We use the description given in definition \ref{framed}, 
since ${\cal M}^0$ has \index{${\cal M}^0$}
an explicit $U(1)$ action which allows us to work equivariantly.

Since the action of the base point preserving gauge transformations on
${\cal A}$ is free, the reducible solutions with $\psi\equiv 0$ now
have trivial stabilizer, hence they are smooth points in ${\cal M}^0$.
The reducible part of the unframed moduli space ${\cal M}$ 
corresponds exactly to the fixed point set for the
$U(1)$ action on ${\cal M}^0$. 

\begin{lem}
Consider the unperturbed equations (\ref{extrem}) in ${\cal B}^0$. Let
$[A,0]$ be a solution that is a fixed point of the $U(1)$-action. Then
the virtual tangent space of ${\cal M}^0$ at the point $[A,0]$ is
$H^1(Y,\RR)\oplus Ker(\partial_A)$. Moreover, the set of fixed points
in ${\cal M}^0$ is identified with the torus $H^1(Y,\RR)/H^1(Y,\ZZ)$
together with one point determined by the choice of the
$Spin_c$-structure, out of  
a discrete set given by the torsion part of
$H^1(Y,\ZZ)$.
\label{fixedpoints}
\end{lem}

\noindent\underline{Proof:} Fixed points are flat $U(1)$-connections
modulo gauge: these are the representations of $\pi_1(Y)$ into $U(1)$. 
The linearization $T$ at a point $[A,0]$ is of the form
$(-*d,\partial_A)$. For the element in the torsion part of
$H^1(Y,\ZZ)$ see Theorem \ref{Morse} below.

\noindent $\diamond$

We need to perturb the equations in some generic way in order to
have ${\cal M}^0$ cut out transversely. For the perturbed equations
(\ref{extremP}), the fixed point set is described by the equation
$*F_A=2i\rho$. 

\begin{lem}
\label{free:orbits}
If Y has non-trivial $b^1(Y)$ and the functional ${\cal C}$ is 
perturbed with a generic co-closed 1-form $\rho$, then the set of
critical orbits ${\cal M}^0$ contains no fixed point and is cut out
transversely by the equations. The Hessian of the perturbed functional
${\cal C}_\rho$ is non-degenerate in the directions normal to the
critical orbits.  
\end{lem}

\noindent\underline{Proof:} If the Chern class $c_1(L)$ is non-trivial,
choose a perturbation $\rho$ with $[*\rho]\neq i\pi c_1(L)$, or
perturb with a harmonic form if $c_1(L)=0$. This implies that there
are no solutions of the equation $F_A=2i*\rho$.

Consider the linearization $\hat L_{(A,\psi,\rho)}$ of the equations
(\ref{extremP}), where we allow the perturbation to vary,
\[ \hat L_{(A,\psi,\rho)}(\alpha,\phi,\eta)= L_{(A,\psi,\rho)}
(\alpha,\phi) -2i\eta. \]

The operator $\hat L$ has closed range, since $L$ is Fredholm.
We show that $\hat L$ is surjective. Let $(\beta,\xi,g)$ be an element
that is $L^2$-orthogonal to the range of $\hat L$. Then
$(\beta,\xi,g)$ is in the kernel of the adjoint, hence by elliptic
regularity we can consider the $L^2$ pairing of $L^2_k$ and $L^2_{-k}$,
\[ \langle \beta, -*d\alpha -df +2\sigma(\psi,\phi)- 2i\eta \rangle +
\langle \xi, \partial_A \phi+\alpha\psi+f\psi \rangle+ \]
\[ +\langle g,G^*(\alpha,\phi)\rangle=0. \]

The argument is analogous to the proof of Proposition
\ref{transverse}. By varying $\eta$ we force $\beta\equiv 0$.
The vanishing of 
\[ \langle \xi, \partial_A\phi+\alpha\psi+f\psi \rangle+ 
\langle g,G^*(\alpha,\phi)\rangle \]
gives an equation $\Delta g +1/2 g |\psi|^2=0$ which implies $g\equiv
0$ by the maximum principle. Then by varying $\phi$ and $\alpha$ we
get $\partial_A\xi=0$ and $\sigma(\xi,\psi)=0$. The latter is satisfied if 
$\xi$ is an imaginary multiple of $\psi$, $\xi=i\lambda\psi$, where
neither of the two vanishes. Both $\xi$ and $\psi$ are in the kernel of
$\partial_A$, thus if either of them vanishes on an open
set it has to vanish identically (and we know that $\psi$ is not
identically zero). If we have $\xi=i\lambda\psi$, we obtain that $\xi$
is identically zero as a consequence of the
vanishing of the inner product $\langle \xi, f\psi\rangle$ for
arbitrary smooth compactly supported functions $f$.

This is enough to show that for a generic perturbation $\rho$ the
moduli space ${\cal M}$ (and therefore also ${\cal M}^0$) is cut out
transversely, as in the analogous proof of Proposition \ref{transverse}.

\noindent $\diamond$ 

When $b^1(Y)=0$ the virtual tangent space at a solution of
$*F_A=2i*d\nu$ is identified with $Ker(\partial_A)$
and the perturbation $\rho=*d\nu$ is not enough to ensure that
${\cal M}^0$ is cut out transversely and that the fixed point set is
separated from the other components of ${\cal M}^0$.

\begin{thm}
Let Y be a rational homology sphere. Suppose we choose a 
perturbation by a co-closed 1-form $\rho$ and a generic metric. Then
the framed moduli space ${\cal M}^0$ consists of a disjoint union
of finitely many circles (corresponding to the irreducible part of the
unframed ${\cal M}$) and finitely many points (the
reducibles of ${\cal M}$). These correspond to different choices of
the $Spin_c$-structure: for a fixed $Spin_c$-structure we have a
unique fixed point of the $U(1)$-action.
Moreover, the Hessian of the functional
${\cal C}_\rho$ is non-degenerate in the normal directions to the
critical orbits.  
\label{Morse}
\end{thm}

\noindent\underline{Proof:}
The choice of a perturbation $\rho=*d\nu$ makes the fixed point set
into the finitely many solutions modulo gauge of
\beq
F_A=d\nu,
\label{reduced}
\eeq 
namely flat $U(1)$ connections modulo gauge. These are representations
of the $\pi_1(Y)$ into $U(1)$, hence they are identified with the finitely
many elements in the group $H^1(Y,\ZZ)$. Each of these elements
specifies one choice of the $Spin^c$--structure. Thus, for fixed
$Spin^c$--structure, there is a unique fixed point
$\theta=[A_0+\nu,0]$ of the
$U(1)$-action in the moduli space ${\cal M}^0$. If the metric on $Y$ satisfies
$Ker(\partial_A)=0$ at a solution $[A,0]$ of
(\ref{reduced}), then the virtual tangent space at this fixed point is
zero-dimensional, hence the Hessian of ${\cal C}_\rho$ is non-degenerate
at the fixed points. The condition is satisfied for a generic metric
because of Lemma \ref{codim:1}.

The irreducible component of ${\cal M}$ is a zero dimensional manifold
if the perturbation $\rho$ is generic. The moduli space ${\cal M}$ is
compact, since $b^1(Y)=0$, and therefore so is ${\cal M}^0$. However,
we have to show that the unique fixed point is separated from the other
components of ${\cal M}^0$. We want to show that no sequence of
irreducible solutions can converge to a reducible solution. We can
apply the same perturbative argument used in \cite{Wa}, based on the
local Kuranishi model. In fact if
$(A,\psi)$ is an irreducible solution which is sufficiently close to a
reducible solution $(A_0+\nu, 0)$, then we have an expansion
\[ A=A_0+\nu +\epsilon \alpha_1 +\epsilon^2 \alpha_2 + \cdots, \]
and
\[ \psi=\epsilon \psi_1 +\epsilon^2 \psi_2 + \epsilon^3 \psi_3
+\cdots. \]
Using the equation (\ref{extremP}) and the condition that
$Ker(\partial_{A_0+\nu})=0$ we 
get that $\psi_i\equiv 0$ for all $i$, in contradiction with the
assumption that $(A,\psi)$ is an irreducible solution.
Thus, for a fixed choice of the $Spin^c$--structure, the framed moduli
space ${\cal M}^0$ consists of finitely many circles 
and a unique point fixed by the $U(1)$ action.

\noindent $\diamond$

\begin{corol}
Under the choice of a generic perturbation (and of a generic metric in
the case with $b^1(Y)=0$), the functional ${\cal C}_\rho$ satisfies
the Morse-Bott property.
\end{corol}

\begin{rem}
The condition that $\partial_A$ satisfies $Ker(\partial_A)=0$ at  the
reducible critical point breaks the space of metrics and perturbations
into chambers, as discussed in Theorem \ref{wallstrata}. The Floer
groups can be expected to 
change when crossing a wall corresponding to metrics with non-trivial
$Ker(\partial_A)$. The problem of possible dependence of the metric in
this case was already addressed by Donaldson in \cite{D}. In Section
7, we prove that indeed the non-equivariant Floer groups do change
when crossing a wall. In fact, we derive explicit wall crossing
formulae for the Euler characteristic of the Floer homology.
\end{rem} 

\subsection{Gradient flow lines}

Denote by $O_a$ the critical orbit in ${\cal M}^0$ that corresponds to
a critical point $a \in {\cal M}$. From \index{$O_a$}
Proposition \ref{Morse} we know that ${\cal M}^0$ is a $U(1)$-fibration
over ${\cal M}$, where $O_a$  is a circle if $a$ is irreducible and
$O_a$ is the point $a$ itself otherwise. 

Let us introduce the configuration space ${\cal A}_{k,\delta}(O_a,
O_b)$. Consider the space \index{${\cal A}_{k,\delta}(O_a,O_b)$}
${\cal A}_{k,\delta}((A_0,\psi_0),(A_1,\psi_1))$ of elements 
$(\AA,\Psi)$ of the form
\beq (\AA,\Psi)\in (A(t),\psi(t))+L^2_{k,\delta}
(\Omega^1(Y\times\RR)\oplus \Gamma (S^+\otimes L)),
\label{confspOaOb1} \eeq
where $(A(t),\psi(t))$ is a smooth path such that
$(A(t),\psi(t))\equiv (A_0,\psi_0)$ for $t\leq 0$ and
$(A(t),\psi(t))\equiv (A_1,\psi_1)$ for $t\geq 1$, with ${\cal
  G}^0$-gauge classes $[A_0,\psi_0]\in O_a$ and $[A_1,\psi_1]\in
O_b$. 

The space ${\cal A}_{k,\delta}(O_a,O_b)$ is given by
\beq {\cal A}_{k,\delta}(O_a,O_b)=\bigcup_{\lambda\in U(1)} {\cal
A}_{k,\delta}(\lambda(A_0,\psi_0),\lambda
(A_1,\psi_1)). \label{confspOaOb2} \eeq 
The space is endowed with a $U(1)$-action.

Let ${\cal B}^0_{k\delta}(O_a,O_b)$ be the quotient of ${\cal
  A}_{k,\delta}(O_a, O_b)$ modulo the free action of the gauge group
\index{${\cal B}^0_{k\delta}(O_a,O_b)$}
${\cal G}^0_{k+1,\delta}(O_a,O_b)$ of based gauge transformations,
\index{${\cal G}^0_{k+1,\delta}(O_a,O_b)$}  
modeled on the $L^2_{k+1,\delta}$ completion of the Lie algebra, that
decay, as $t\to \pm\infty$, to elements $\lambda_{\pm\infty}$ 
in the stabilizers $G_a$ and $G_b$. 
The quotient has an induced $U(1)$-action and endpoint maps
\[ e^+_a : {\cal B}^0_{k\delta}(O_a,O_b) \to O_a \]
and 
\[ e^-_b : {\cal B}^0_{k\delta}(O_a,O_b) \to O_b. \]
\index{$e^+_a$, $e^-_b$}

We denote by ${\cal M}(O_a,O_b)$ the moduli space of solutions of
equations (\ref{4SW1P'}) and (\ref{4SW2P'}) in ${\cal
    B}^0_{k\delta}(O_a,O_b)$. These will be our moduli spaces of flow
lines connecting the orbits $O_a$ and $O_b$.
\index{${\cal M}(O_a,O_b)$}

\begin{rem}
Suppose one of $O_a$ and $O_b$ is not the fixed point $\theta$. Then
no reducible 
solution arises among the flow lines ${\cal M}(O_a,O_b)$. 
\label{noreducibles}
\end{rem}

In fact, either $(A_a,\psi_a)$ or $(A_b,\psi_b)$ has non-trivial
spinor. Thus  $\psi(t)$ has to be non-trivial, in fact, the
exponential weight on the Sobolev norms forces the elements in
(\ref{confspOaOb1}) to decay to the endpoints, at least one of which
has non-vanishing spinor. Thus, the space
${\cal A}_{k,\delta}(O_a,O_b)$ only contains irreducible lines. 
The action of ${\cal G}^0_{k+1,\delta}(O_a,O_b)$ is
free and ${\cal B}^0_{k,\delta}(O_a,O_b)$ is a manifold with a free
$U(1)$-action. 

\subsection{Relative Morse Index}

We can rephrase Proposition \ref{relmorseind} in the case of framed
moduli spaces as follows.

\begin{thm}
Suppose $O_a$ and $O_b$ are critical orbits of ${\cal C}_\rho$ in
${\cal B}^0$. If $Y$ is a rational homology sphere, assume that the
metric and perturbation are chosen generically so that ${\cal C}_\rho$
satisfies the
Morse-Bott condition. Let $\{ \lambda_a \}$ and $\{ \lambda_b \}$ be
the eigenvalues of the Hessian operators $Q_a$, $Q_b$ (that is the
Hessian in ${\cal B}^0$ restricted to the directions orthogonal to the
$U(1)$-orbits). Choose a weight $\delta$ satisfying $0<\delta<\min\{
|\lambda_a |, |\lambda_b | \}$. Then the linearization ${\cal
  L}_{(\AA,\Psi)}$ at a solution $[\AA,\Psi]\in{\cal
    B}^0_{k,\delta}(O_a,O_b)$ of (\ref{4SW1P'}) and (\ref{4SW2P'}) is
a Fredholm operator from $L^2_{k,\delta}$ to $L^2_{k-1,\delta}$. The
virtual dimension of the moduli space ${\cal M}(O_a,O_b)$ is given by
the index of ${\cal L}_{(\AA,\Psi)}$ and is obtained as
\[ \dim {\cal M}(O_a,O_b)=\sigma(O_a,O_b)+1 -\dim G_a, \]
where $\sigma(O_a,O_b)$ is the spectral flow of the
operator $\nabla{\cal F}$ on a path $(A(t),\psi(t))$ in ${\cal A}$
corresponding to $(\AA,\Psi)$. The quantity
$\sigma(O_a,O_b)$ is independent of the path in ${\cal B}^0_k$; by
additivity of the spectral flow, it can be written as
\[  \sigma(O_a,O_b)=\mu(O_a)-\mu(O_b), \]
where $\mu(O_a)$ is the flow of $\nabla{\cal F}$ on a path connecting
the orbit $O_a$ to a fixed orbit in ${\cal M}^0$.
\label{relmorseind2}
\end{thm}

\noindent\underline{Proof:}
The fact that the linearization ${\cal L}_{(\AA,\Psi)}$  on
the spaces ${\cal B}_{k,\delta}(a,b)$ and on ${\cal
B}_{k,\delta}(O_a,O_b)$ is Fredholm follows from Theorem 6.2 of
\cite{LM}, or Lemma 8.3.1 of \cite{MMR} (cf. the previous discussion
of this point in Proposition \ref{orientation}), provided
that the operator $T$ has trivial kernel at the points $a$ and $b$ and
the weight $\delta$ is smaller than the least eigenvalue of $T$, see
also \cite{T}. 

We shall write $Ind_e ({\cal L}_{(\AA,\Psi)})$ for the index formula
in ${\cal B}^0_{k\delta}(O_a,O_b)$. The subscript denotes the fact
that we are computing \index{$Ind_e$}
$$ Ind_e ({\cal L}_{(\AA,\Psi)})= \dim Ker_e({\cal L}_{(\AA,\Psi)})
-\dim Coker({\cal L}_{(\AA,\Psi)}), $$
where $Ker_e({\cal L}_{(\AA,\Psi)})$ is the extended kernel, which
consists of solutions $(\alpha, \Phi)$ of ${\cal
L}_{(\AA,\Psi)}(\alpha, \Phi)=0$, \index{$Ker_e$} with $(\alpha,
\Phi)-(\alpha_-,\phi_-)$ in $L^2_{k,\delta}$, where $(\alpha_-,\phi_-)$
is a solution of $L_{A_a,\psi_a}(\alpha_-,\phi_-)=0$, namely a tangent
vector to the orbit $O_a$ in ${\cal T} {\cal B}^0$. 
Notice that this $Ker_e({\cal
L}_{(\AA,\Psi)})$ is the correct space representing the tangent space
of ${\cal M}(O_a,O_b)$ in ${\cal B}^0_{k\delta}(O_a,O_b)$. The
cokernel is simply given by $L^2_{k-1,\delta}$ solutions of ${\cal
L}^*_{(\AA,\Psi)}(f,\beta,\xi)=0$, as in the setting of \cite{APS} I,
pg. 58. The formula for the index of ${\cal L}_{(\AA,\Psi)}$  on
${\cal B}^0_{k,\delta}(a,b)$,
\[ Ind_e({\cal L}_{(\AA,\Psi)})=\sigma(O_a,O_b) +1 -\dim G_a, \]
then follows from the splitting of the index that will be proven in
Corollary \ref{indexsplit}, together with the additivity of the
spectral flow, and \cite{APS} III pg.95.

Thus we obtain the virtual dimension of ${\cal M}(O_a,O_b)$,
\[ \dim {\cal M}(O_a,O_b)=\sigma(O_a,O_b)+1-\dim G_a, \]
where the $+1$ contribution depends upon the presence
of the $U(1)$-action and $\dim
G_a$ is the dimension of the stabilizer of the point $a$.

The relative Morse index is well defined. To see this, we have to
examine the spectral flow of the operator $\nabla{\cal
  F}_{(A(t),\psi(t))}$ on a path $(A(t),\psi(t))$ in ${\cal
  A}$ with endpoints $(A,\psi)$ and $\lambda(A,\psi)$, where $\lambda$
is a gauge transformation satisfying
\[ \frac{i}{2\pi}\int_Y c_1(L)\wedge \lambda^{-1}d\lambda =0. \] 
According to \cite{APS} III pg.95, and \cite{GWa}, this spectral flow is
just the index of ${\cal L}$ on the manifold $Y\times S^1$.
The index on a closed four-manifold is
\[ Ind({\cal L})=c_1(L)^2-\frac{2\chi+3\sigma}{4} \]
(cf. Corollary 4.6.2 of \cite{Mo}, or Theorem 2.3.8 of \cite{Ma2}), and
this quantity vanishes in the case of a manifold of the form $Y\times S^1$.

\noindent $\diamond$

Given the relative Morse index, we can define the Morse index of a critical 
orbit up to fixing arbitrarily the index of a particular solution. In the 
case with $b^1(Y)>0$ there is no canonical choice, hence the grading of
the Floer complex is only defined up to an integer.
When $b^1(Y)=0$ we can remove this ambiguity by fixing the trivial solution
$\theta=[A_0+\nu,0]$ to have index zero.

Notice that, in the case with $b^1(Y)>0$,  $c_1(L)\neq 0$, and
$[*\rho]=0$, the relative Morse index would be defined only up to a 
periodicity if we considered solutions modulo the full gauge group
$\tilde{\cal G}$. This is related to the fact that the functional
${\cal C}$ is well defined on ${\cal A}/{\cal G}$, but is only defined
as a circle valued functional on ${\cal A}/\tilde{\cal G}$. The
relation between $\ZZ$-graded and $\ZZ_l$-graded Seiberg-Witten Floer
theory will be sketched briefly in Section 4.1, in Remark
\ref{correction}. 

\subsection{Decay estimate}

In this subsection we introduce some analytic properties of the
functional ${\cal C}_\rho$. We show that finite energy solutions of
the flow equations necessarily converge to asymptotic values that lie
on some critical orbit. Moreover, we show that, if $O_a$ or $O_b$ is
an irreducible orbit, then the flow lines in the  moduli space ${\cal 
M}(O_a, O_b)$ decay exponentially towards the endpoints. Notice that
in this case the moduli space ${\cal M}(O_a, O_b)$
only contains irreducible flow lines. 

We give the following preliminary definition.
\begin{defin}
A smooth path $(A(t),\psi(t))$ in ${\cal A}$ is of finite energy if
the integral
\beq
\int_{-\infty}^{\infty} \| \nabla {\cal C}_\rho(A(t),\psi(t))
\|^2_{L^2} dt <\infty  
\label{finenergy}
\eeq
is finite.
\end{defin}

Notice that any solution of (\ref{3SW1P'}) and (\ref{3SW2P'}) with
asymptotic values in $O_a$ and $O_b$ is of finite energy, in fact
in this case the total variation of the functional ${\cal C}_\rho$
along the path $(A(t),\psi(t))$ is finite and
(\ref{finenergy}) satisfies
\[ \int_{-\infty}^{\infty} \| \nabla {\cal C}_\rho
(A(t),\psi(t)) \|_{L^2}^2 dt\leq C ({\cal C}_\rho(a)-{\cal C}_\rho
(b)), \] 
because of the assumptions on the perturbation $q_{(\AA,\Psi)}$.
Finite energy solutions of the flow equations have nice properties:
they necessarily decay to asymptotic values that are critical points
of ${\cal C}_\rho$ as we prove in Corollary \ref{finenergy2} and in Theorem
\ref{decay}. We 
begin by introducing some analytic properties of the functional ${\cal
  C}_\rho$ (see also \cite{Fr}, \cite{MST}, \cite{Wa}). 

\begin{lem}
Let ${\cal M}_\rho$ be the moduli space of critical points of ${\cal
  C}_\rho$, with $\rho$ a sufficiently small perturbation.
For any $\epsilon >0$ there is a $\lambda >0$ such that, if the
$L^2_1$-distance of a point $[A,\psi]$ of ${\cal B}$ to all the
points in ${\cal M}_\rho$ is at least $\epsilon$, then 
\[ \| \nabla{\cal C}_\rho(A,\psi) \|_{L^2} >\lambda. \]
\label{PalaisSmale}
\end{lem}

\noindent\underline{Proof:} For a sequence $[A_i,\psi_i]$ of
elements of ${\cal B}$ with a distance at least $\epsilon$ from all
the critical points, such that
\[ \| \nabla{\cal C}_\rho(A_i,\psi_i) \|_{L^2}\to 0, \]
as $i\to\infty$, we would have
\[ \| *F_{A_i}-\sigma(\psi_i,\psi_i)-i\rho \| + \|
\partial_{A_i}\psi_i \| \to 0. \]
Thus, there is a constant $C$ such that
\[ \int_Y |*F_{A_i}-\sigma(\psi_i,\psi_i)-i\rho |^2 +
|\partial_{A_i}\psi_i |^2 dv <C. \]
If the perturbation $\rho$ is sufficiently small, the Weizenb\"ock
formula implies that
\[ \int_Y |F_{A_i}|^2 + |\sigma(\psi_i,\psi_i)|^2
+\frac{\kappa}{2}|\psi_i|^2 + 2|\nabla_{A_i}\psi_i|^2 dv <C. \]
Thus we have a uniform bound on the norms $\| \psi_i \|_{L^4}$, $\|
F_{A_i} \|_{L^2}$, and $\|\nabla_{A_i}\psi_i\|_{L^2}$. An elliptic
estimate shows that there is a subsequence that converges in the
$L^2_1$ norm to a solution of the critical point equations
(\ref{extremP}), and this contradicts the assumption.

\noindent $\diamond$

\begin{corol}
Let $(A(t),\psi(t))$ be a smooth finite energy solution of equations
(\ref{3SW1P'}) and (\ref{3SW2P'}) with a smooth perturbation $q$.
Then there exist critical 
points $a$ and $b$ of ${\cal C}_\rho$, such that the 
$\lim_{t\to\pm\infty}(A(t),\psi(t))$ are in the gauge classes of $a$
and $b$.
\label{finenergy2}
\end{corol}

\noindent\underline{Proof:}
The finite energy condition (\ref{finenergy}) implies that
\[ \| \nabla{\cal C}_\rho(A(t),\psi(t)) \|\to 0 \]
as $t\to\pm\infty$. 
The Palais-Smale condition of Lemma \ref{PalaisSmale} implies that
there exist $T$ large, such that
for $|t| > T$,  $(A(t), \psi(t))$ lies in a very small
$\epsilon$-neighborhood of critical points of ${\cal C}_\rho$.

\noindent $\diamond$

Now we prove the exponential decay property.

\begin{thm}
Let $O_a$ and $O_b$ be non-degenerate critical orbits in ${\cal
  B}_0$. There exists a weight $\delta >0$ such that the following holds. 
Suppose given any solution
$[\AA,\Psi]$ of (\ref{4SW1P'}) and (\ref{4SW2P'}) that is represented
by a smooth pair $(A(t),\psi(t))$ in a temporal gauge, with asymptotic
values $(A_a,\psi_a)$ and $(A_b,\psi_b)$ in $O_a$ and $O_b$. Then there
exists a constant $K$ such that, for $t$
outside an interval $[-T,T]$, the distance in any fixed
${\cal C}^l$-topology of $(A(t),\psi(t))$ from the
endpoints is
\[ dist_{{\cal C}^l}((A(t),\psi(t)),(A_i,\psi_i))< K \exp(-\delta
|t|), \]
with $i=a$ if $t< -T$ and $i=b$ if $t>T$. 
\label{decay}
\end{thm}

\noindent\underline{Proof:}
The proof consists of a few steps. Let us consider the decaying as
$t\to \infty$; the other case with $t\to -\infty$ is analogous. 

For simplicity of notation we shall prove the Theorem in the case of
flow lines in ${\cal B}$, with the action of the gauge group ${\cal
G}$ of gauge maps that satisfy (\ref{int=0}), and perturbation $\rho$
satisfying $[*\rho]=0$ in cohomology. This ensures that the functional
${\cal C}_\rho$ is $\RR$-valued. All the claims and the proofs extend
directly to the case of the based space with the
${\cal G}^0$-action and a non-degenerate critical orbit $O_b$ in the
quotient space ${\cal B}_0$. In this case, the distance of $(A,\Psi)$
from the orbit $O_b$ is the minimal distance from points on the orbit.

\begin{lem} 
\label{near:critical}
Suppose $b$ is a  non-degenerate critical point of ${\cal C}_\rho$. Then
there exists a constant $C_b$ such that if the
$L^2$-distance from $[A, \psi]$ to $b$ is sufficiently small,
then we have the following estimate  of the $L^2$-distance from $[A,
\psi]$ to $b$: 
\[ dist_{L^2}( [A, \psi], b )\leq  C_b \|\nabla {\cal C}_\rho (A, \psi)
\|_{L^2}. \]  
\end{lem}

\noindent\underline{Proof:} Consider the Hessian operator $T_b$ acting
as an unbounded operator on the space of $L^2$ connections and sections.
Since $b$ is a non-degenerate critical
point, we have that
\[
h_b  = \max \{ \displaystyle{\frac {1}{|\lambda_i|}} | \
\hbox{$\lambda_i$ is a eigenvalue of the Hessian operator $T_b$ at $b$}\}
\]
exists and is bounded. We know that
 $[A, \psi] \mapsto \nabla {\cal C}_\rho (A, \psi)$ defines a
$L^2$-tangent section, which is smooth and transverse  to zero at
$b$.
Thus, for any $\epsilon>0$,  we may choose a small neighborhood $U_b$ of $b$
which may be identified with a small
neighborhood of $0$ in the $L^2$-tangent space of ${\cal B}$ at $b$,
such that for all $(A,\psi)=b+(\alpha,\psi)$ in $U_b$  we have
\[ \| \nabla {\cal C}_\rho (A,\psi) - T_b (\alpha,\phi) \|_{L^2} \leq
\epsilon, \]
where we write $\nabla {\cal C}_\rho (A,\psi)$ as the sum of the
linear and a non-linear term,
\[ \nabla {\cal C}_\rho (A,\psi)=T_b (\alpha,\phi)+N(\alpha,\phi),
\]
with $N(\alpha,\phi)=(\sigma(\phi,\phi),\alpha\cdot\phi)$.

When the neighborhood $U_b$ is small enough, we can ensure that 
\[ \| N(\alpha,\phi) \|\leq 1/2 \|\nabla {\cal C}_\rho (A,\psi) \|,
\]
so that we have
\[ \| T_b(\alpha,\phi) \|\leq \| \nabla {\cal C}_\rho (A,\psi) \| + \|
N(\alpha,\phi) \| \]
\[ \leq 3/2 \| \nabla {\cal C}_\rho (A,\psi) \|. \]
 
Thus, we get the following estimate:
\[
\begin{array}{lll}
&& dist_{L^2}( [A, \psi], b ) \\[2mm]
&=&\|(\alpha,\phi)\|_{L^2}\\[2mm]
&\leq& h_b \|(T_b(\alpha,\phi))\|_{L^2}
\\[2mm]
&\leq&\displaystyle{\frac 32}h_b\|(\nabla{\cal C}_\rho ([A,
\psi]))\|_{L^2(Y)} \\[2mm]
&\leq& \displaystyle{\frac 32}
h_b \|\nabla {\cal C}_\rho ([A, \psi])\|_{L^2}. \end{array}
\]
The Lemma follows upon
choosing the constant $C_b$ with
\[  C_b >  \displaystyle{\frac 32} h_b. \]

\noindent{\bf Claim 1:} 
Let $(A(t),\psi(t))$ be a representative of a path in
${\cal B}$, such that  
\[ \lim_{t\to \infty} (A(t),\psi(t))=(A_b,\psi_b) \]
where $b=[A_b,\psi_b]$ is a non-degenerate critical point in ${\cal
M}_\rho$. Then there is a $T_0>>0$ and a constant $K_b$, such that the
inequality 
\[ |{\cal C}_\rho(A(t),\psi(t))-{\cal C}_\rho(A_b,\psi_b) |\leq K_b \|
\nabla {\cal C}_\rho (A(t),\psi(t)) \|_{L^2}^2 \]
holds for all $t\geq T_0$.

\noindent\underline{Proof of Claim 1:} 

Choose $T_0$ such that for all $t\geq T_0$ the path $(A(t),\psi(t))$
lies in a neighborhood $U_b$ of $b$ for which the result of Lemma
\ref{near:critical} holds. Notice that such a $T_0$ depends in general
on which path $(A(t),\psi(t))$ is considered.

With the notation of (\ref{1forma}) and (\ref{2form}), we have
\[ {\cal C}_\rho(A,\psi)\sim{\cal C}_\rho(A_b,\psi_b)+{\cal
F}_b(\alpha,\phi)+\nabla{\cal F}_b(\alpha,\phi), \]
where $(A,\psi)=b+(\alpha,\phi)$ in $U_b$. Since $b$ is a critical
point, we have ${\cal F}_b\equiv 0$.

Thus we get
\[ |{\cal C}_\rho(A,\psi)-{\cal C}_\rho(A_b,\psi_b)| \leq | \langle
T_b(\alpha,\phi), (\alpha,\phi)\rangle |. \]
Now applying Lemma \ref{near:critical} we obtain the estimate
\[ |{\cal C}_\rho(A,\psi)-{\cal C}_\rho(A_b,\psi_b)| \leq \|
T_b(\alpha,\phi) \| \cdot \| (\alpha,\phi) \| \]
\[ \leq  \frac{3}{2}  C_b \| \nabla {\cal C}_\rho (A,\psi) \|_{L^2}^2. \]
This completes the proof of Claim 1, with the constant $K_b\geq 3C_b /2$.

\noindent{\bf Claim 2:}
For a finite energy solution $(A(t),\psi(t))$ of (\ref{3SW1P'}) and
(\ref{3SW2P'}), the inequality 
\[ \frac{1}{2} \int_t^\infty \|\nabla {\cal C}_\rho (A(s),\psi(s))
\|_{L^2}^2 ds \leq {\cal C}_\rho (A(t),\psi(t))- \]
\[ -{\cal C}_\rho(A_b,\psi_b) \leq \frac{3}{2} \int_t^\infty \|\nabla
{\cal C}_\rho (A(s),\psi(s)) \|_{L^2}^2 ds \]
holds for large $t$. 

\noindent\underline{Proof of Claim 2:}
Without loss of generality we can assume that the perturbation in
${\cal P}$ satisfies Condition (4) of \ref{calP} with $C_0< 1/2$,
so that
\[ \|q_{(\AA,\Psi)}(t)\|_{L^2}< \frac{1}{2}\| \nabla{\cal C}_\rho
(A(t),\psi(t)) \|_{L^2}. \]
Thus, we can replace the equality
\[ {\cal C}_\rho (A(t),\psi(t))-{\cal C}_\rho
((A_b,\psi_b))= 
\int_t^\infty -\frac{d}{ds}{\cal C}_\rho (A(s),\psi(s))
ds \]
\[ =\int_t^\infty -<\frac{d}{ds} (A(s),\psi(s)),\nabla {\cal C}_\rho
(A(s),\psi(s))> ds = \]
\[ =\int_t^\infty \|\nabla {\cal C}_\rho (A(s),\psi(s)) \|^2 ds, \] 
that holds for solutions of the unperturbed equations with the
inequality of Claim 2 for solutions of the perturbed equations.

\noindent{\bf Claim 3:} Let $(A(t)),\psi(t))$ be a finite energy
solution of the equations (\ref{3SW1P'}) and (\ref{3SW2P'}). The quantity 
\[ E(t)=\int_t^\infty \|\nabla {\cal C}_\rho (A(s),\psi(s)) \|_{L^2}^2 ds \]
decays exponentially as $t\to \infty$.

\noindent\underline{Proof of Claim 3:}
In fact, the inequality of Claim 2 gives the first inequality in the following
estimate: 
\[ E(t)\leq 2 ({\cal C}_\rho (A(t),\psi(t))-{\cal
  C}_\rho(A_b,\psi_b)) \leq \]
\[  \leq K_b \| \nabla {\cal C}_\rho (A(t),\psi(t))
\|^2  = -K_b \frac{d}{dt} E(t). \]
The second inequality follows from Claim 1.

\noindent{\bf Claim 4:} For large $t$ we have the inequality 
\[ dist_{L^2}((A(t),\psi(t)),(A_b,\psi_b))\leq K (\int_{t-1}^\infty
\| \nabla {\cal C}_\rho (A(s),\psi(s)) \|_{L^2}^2 ds)^{1/2}, \]
when $x(t)=(A(t),\psi(t))$ is a finite energy solution of (\ref{3SW1P'})
and (\ref{3SW2P'}). 

\noindent\underline{Proof of Claim 4:}
We can prove that the following inequality holds true for
$t\geq T_0$:
\[ \int_{t}^\infty \| \frac{d}{ds}x(s) \|_{L^2}ds \leq
\frac{2K_b}{1-C_1} ({\cal 
C}_\rho(x(t))- {\cal C}_\rho (b))^{1/2}, \]
where $T_0$ is such that for all $t\geq T_0$ the perturbation
$q_{(\AA,\Psi)}$ satisfies the inequality  
\[ \| q_{(A(t),\psi(t))}(t)\|_{L^2} \leq C_1 \|
\nabla {\cal C}_\rho (A(t),\psi(t)) \|_{L^2}, \]
with $0< C_1 <1$, as in property (4) of definition \ref{calP}.
Moreover, $T_0$ is such that Lemma \ref{near:critical} holds. 
The proof of this inequality follows \cite{Si} Lemma 3.1 pg.542. We
have 
\[ \| \frac{d}{dt}x(t) \|^2 + \| \nabla {\cal C}_\rho (x(t)) \|^2 + 2
\langle \frac{d}{dt}x(t), \nabla{\cal C}_\rho(x(t)) \rangle = \]
\[ = \| q(t) \|^2\leq C_1^2 \| \frac{d}{dt}x(t) \|^2 \leq C_1 \|
\frac{d}{dt}x(t) \|^2. \]
Thus we get
\[ -\langle \frac{d}{dt}x(t),\nabla{\cal C}_\rho(x(t))\rangle \geq
\frac{(1-C_1)}{2} (\| \frac{d}{dt}x(t) \|^2 + \| \nabla {\cal C}_\rho (x(t))
\|^2) \]
\[ \geq (1-C_1) \| \frac{d}{dt}x(t) \| \cdot \| \nabla{\cal C}_\rho
(x(t)) \|, \]
for all $t\in [T_0,\infty)$. That is,
\[ -\frac{d}{dt}{\cal C}_\rho (x(t))\geq (1-C_1) \| \frac{d}{dt}x(t)
\| \cdot \| \nabla{\cal C}_\rho(x(t)) \|. \] 
We obtain the inequality
\[ -\frac{d}{dt} ({\cal C}_\rho(x(t))-{\cal C}_\rho(b))^{1/2} \geq
\frac{(1-C_1)}{2} ({\cal C}_\rho(x(t))-{\cal C}_\rho(b))^{-1/2}\|
\frac{d}{dt}x(t) \| \cdot \| \nabla{\cal C}_\rho(x(t)) \| \]
\[ \geq \frac{(1-C_1)}{2 K_b} \| \frac{d}{dt}x(t) \|. \]
The last inequality follows from Claim 1.   
We can now integrate both sides to obtain
\[ \int_t^\infty \|\frac{d}{ds}x(s) \|_{L^2}ds \leq \frac{2K_b}{(1-C_1)}
\int_t^\infty -\frac{d}{ds} ({\cal C}_\rho(x(s))-{\cal
C}_\rho(b))^{1/2}ds \]
\[ =\frac{2K_b}{(1-C_1)} ({\cal C}_\rho(x(t))-{\cal C}_\rho(b))^{1/2}. \]

Thus by Claim 2 we have
\[ dist_{L^2}((A(t),\psi(t)),(A_b,\psi_b))\leq K (\int_{t}^\infty
\| \nabla {\cal C}_\rho (A(s),\psi(s)) \|_{L^2}^2 ds)^{1/2}. \]

The exponential decay of $E(t)$ proves the claim of the Theorem for
the case of $L^2$-topology. Smooth estimates then follow by a
bootstrapping argument and elliptic regularity. This completes the
proof of the Theorem.

\noindent $\diamond$

Analogous exponential decay estimates have been proven in 
\cite{MST}, cf. \cite{Au1} and \cite{Wa}. 

\subsection{Transversality of ${\cal M}(O_a,O_b)$}

We have the following transversality result for the moduli spaces
${\cal M}(O_a,O_b)$. 

\begin{prop}
Suppose given two orbits $O_a$ and $O_b$ in ${\cal M}^0$.
For a generic choice of the perturbation $P\in{\cal P}$, the space
${\cal M}(O_a, O_b)$ is a smooth $\RR \times U(1)$-manifold
with dimension given by $\mu(O_a)- \mu(O_b)+1-\dim G_a$.
${\cal M}(O_a, O_b)$ is non-empty only if $\mu(O_a)- \mu(O_b) \ge 1$.
There are endpoint maps
\[ e_a^+: {\cal M}(O_a, O_b) \to O_a \]
and
\[ e_b^-:  {\cal M}(O_a, O_b) \to O_b \]
that are smooth $U(1)$-equivariant maps.
\label{flowlines}
\end{prop}

The computation of the virtual dimension follows from Theorem
\ref{relmorseind2}. The transversality statement follows from
Proposition \ref{transverse}, together with the Remark
\ref{noreducibles}, to the effect that the smooth manifolds ${\cal
M}(a,b)$ of Proposition \ref{transverse} are just the quotient of
${\cal M}(O_a, O_b)$ with respect to the free $U(1)$-action.
The properties of the endpoint maps follow from the discussion on the
convergence to the endpoints of Section 3.4.

We have the analogue of Corollary \ref{empty1} in this case.

\begin{corol}
Suppose $\theta$ is the reducible critical point in ${\cal M}$ and  $b$ any
point in ${\cal M}$.
Then, after a generic perturbation, ${\cal M}(\theta, O_b)$ is a smooth
$\RR \times U(1)$-manifold with dimension $\mu(\theta) -\mu(O_b)$. Moreover, 
${\cal M}(\theta, O_b)$ is non-empty only if $\mu(\theta)  -\mu(O_b) \ge 2$.
\label{empty2}
\end{corol}

\section{Boundary structure}

The purpose of this part of the work is to study the compactification
of the moduli spaces $\hat{\cal M}(O_a,O_b)$. 
\index{$\hat{\cal M}(O_a,O_b)$} The first step consists
of showing the existence of a compactification, obtained by adding a
certain set of broken trajectories, namely trajectories that break
through other critical orbits of intermediate relative Morse index. 
We shall introduce the notation $\hat{\cal M}(O_a,O_b)^*$ for the
\index{$\hat{\cal M}(O_a,O_b)^*$}
compactified moduli spaces. Proving the existence of the
compactification by broken trajectories 
is dealt with in Section 4.1. Section 4.2 then deals with the
gluing construction, and Section 4.3 analyzes the fine structure of
the compactification.

 It is perhaps useful to recall why in Floer theory it is necessary to 
develop the full gluing construction, as in Section 4.2.
The results of Section 4.1 are sufficient to prove 
the existence of a
compactification for the moduli spaces of flow lines, obtained by adding
broken trajectories. The points of the ideal boundary are not completely 
identified, and the purpose of the gluing construction of Section 4.2 is
precisely to identify all the broken trajectories that appear in the
boundary. Recall also that, in the Floer complex, the boundary operator
is obtained by counting elements in the oriented moduli spaces 
$\hat{\cal M}(a,b)^*$ for critical points $a$, $b$ of relative index one.
The statement that this is a boundary operator, namely that $\partial 
\circ \partial =0$, depends precisely upon having a complete
description of the ideal boundary in the compactified
moduli spaces $\hat{\cal M}(a,c)^*$ for critical points $a$ and $c$ of
relative index two. In other words, proving that all the broken trajectories
$$ \cup_{\mu(a)>\mu(b)>\mu(c)} \hat{\cal M}(a,b)^*\times \hat{\cal
M}(b,c)^* $$
appear in the compactification  $\hat{\cal M}(a,c)^*$ is necessary in
order to have a chain complex.
In the equivariant setting we are considering here, there are
components of the boundary operator $D$ which are obtained by pullback
and push-forward of forms on compactified moduli spaces    
$\hat{\cal M}(O_a,O_b)^*$, for orbits of relative index one, and on
$\hat{\cal M}(O_a,O_c)^*$ for orbits of relative index two (see the
explicit description of the boundary operator given in Section 5, in
(\ref{equivD}). The proof that $D\circ D=0$ relies upon a version
of Stokes' theorem for manifolds with corners (as in \cite{AB2}).
For this reason, it is necessary to show that the compactification of
the moduli spaces of flow lines of Section 4.1 has the structure of a
smooth manifold with corners. This is proved in the case of the
codimension one boundary in Section 4.2, and the generalization to the 
strata of higher codimension is considered in Section 4.3.

Thus, in Floer theory, a very elaborate analysis of the
compactification of the moduli spaces of flow lines is required, in
order to construct the Floer complex. This situation is essentially
different from other issues of compactification in gauge theory. For
instance, in the case of Donaldson invariants \cite{DK}, one only
needs to show the existence of the Uhlenbeck compactification of the
moduli space of anti-self-dual connections. The definition of the
invariants and the proof of the diffeomorphism invariance do not
require to show that every ideal anti-self-dual connection actually
appears in the Uhlenbeck compactification, nor they require the
existence of a fine structure (such as that of smooth manifold with
corners) on the compactification. 

To start our analysis, we  fix a unique way to identify the space of
unparameterized trajectories \index{$\hat{\cal M}(O_a,O_b)$}
\[ \hat{\cal M}(O_a,O_b)={\cal M}(O_a,O_b)/\RR \]
with a subset of the space of parameterized
trajectories ${\cal M}(O_a,O_b)$. It is sufficient to choose the 
parameterization $x(t)$ of $\hat x\in \hat{\cal M}(O_a,O_b)$ which
satisfies the equal energy condition
\begin{equation}
\int_{-\infty}^0 \| \nabla {\cal C}_\rho (A(t),\psi(t))
\|^2_{L^2}dt=\int_0^\infty\| \nabla {\cal C}_\rho (A(t),\psi(t))
\|^2_{L^2}dt. 
\label{equal:energy}
\end{equation} 

This lifting of $\hat{\cal M}(O_a,O_b)$ to ${\cal M}(O_a,O_b)$ is
unique. In fact, in the family of gradient flows $\{
[\AA, \Psi]^T, T\in \RR\}$ that represent the class $x\in \hat{\cal
M}(O_a,O_b)$ there is a unique element which satisfies the equal
energy condition (\ref{equal:energy}). The lifting is often referred
\index{${\cal M}^{bal}(O_a,O_b)$}
to as ${\cal M}^{bal}(O_a,O_b)\subset {\cal M}(O_a,O_b)$, the balanced 
moduli space. We have $\hat{\cal M}(O_a,O_b)\cong {\cal
M}^{bal}(O_a,O_b)$. In the following, in order to avoid exceeding use 
of different notation, we shall always write $\hat{\cal M}(O_a,O_b)$,
unless we need to make explicit use of the equal energy condition
(\ref{equal:energy}), in which case we may recall that $\hat{\cal
M}(O_a,O_b)$ is realized by the balanced moduli space ${\cal
M}^{bal}(O_a,O_b)\subset {\cal M}(O_a,O_b)$, as in the proof of
Theorem \ref{equivgluing} in Section 4.2. 

\subsection{Convergence theorem}

The following theorem describes convergence in the moduli space $\hat{\cal
  M}(O_a,O_b)$, and proves the existence of a compactification, which
we denote $\hat{\cal M}(O_a,O_b)^*$, obtained by
adding broken trajectories.

\begin{thm}
Consider the moduli space $\hat{\cal M}(O_a,O_b)$, with
$\mu(O_a)-\mu(O_b)=k+1$, $k\geq 0$.
The space $\hat{\cal M}(O_a,O_b)$ is precompact. Namely, any sequence
$[\hat x_i]$ of elements in $\hat{\cal M}(O_a,O_b)$ has a subsequence which
either converges in norm to another solution $[\hat x]\in\hat{\cal
M}(O_a,O_b)$, or converges to a broken trajectory. This means that
there are critical orbits $O_{c_1},\ldots, O_{c_k}$, with
$\mu(O_a)\geq \mu(O_{c_1})>\ldots >\mu(O_{c_k})\geq \mu(O_b)$, 
trajectories $[y_j]\in {\cal M}(O_{c_j},O_{c_{j+1}})$, and a sequence of
real numbers $T_{i_k,j}\in \RR$ such that the sequence of parameterized
trajectories $x_{i_k}^{ T_{i_k,j}}$ converges smoothly on compact sets
\[ x_{i_k}^{ T_{i_k,j}}\to y_j. \]
Here $x_{i_k}^{ T_{i_k,j}}$ denotes the lifting of $\hat x_{i_k}$ to
the space of parameterized trajectories specified by the condition
$x_{i_k}^{T_{i_k,j}}(0)=x_{i_k}(T_{i_k,j})$, where $x_{i_k}$ is the
equal energy lift of $\hat x_{i_k}$.
\label{seqcompact}
\end{thm}

\noindent\underline{Proof:} The proof of Theorem \ref{seqcompact}
consists of several steps. 
Suppose that the perturbation $P_{(\AA,\Psi)}$ is smooth. Given
a sequence of unparameterized trajectories $[\hat x_i]$ in $\hat{\cal
M}(O_a,O_b)$, choose representatives $x_i$ of the corresponding equal
energy lift in ${\cal M}(O_a,O_b)$. By elliptic regularity these can
be represented by smooth solutions.

Let $(A_i(\pm\infty),\psi_i(\pm\infty))$ be the asymptotic values of
the elements $x_i$ at the ends of the cylinder $Y\times \RR$. These
represent elements $[A_i(\pm\infty),\psi_i(\pm\infty)]$ in $O_a$ and
$O_b$ respectively.

\noindent{\bf Claim 1:} There is a subsequence $\{ x_{i_k} \}$ that
converges smoothly on compact sets to a solution $y$ of the perturbed
flow equations. 

\noindent\underline{Proof of Claim 1:} We first show that there is a
subsequence that converges uniformly on compact sets. By the use of
a bootstrapping argument is then possible to improve the convergence
to ${\cal C}^\infty$ on compact sets.

It is useful to recall the following weak version of Arzela'-Ascoli.

\begin{prop}
Let $K$ be a compact subset of $\RR$ and $(X, \|\cdot\|)$ a normed
metric space. $S$ is a subset of continuous functions from $K$ to $X$.
Suppose the following conditions are
satisfied 

\noindent (i) pointwise bound:
\[ \sup_{f\in S} \| f(t) \| <\infty \]
for all $t\in K$;

\noindent (ii) local equicontinuity: given $\epsilon >0$,
for all $t$ in $K$ there exists a
neighborhood $V_{t}$ such that, for all $\tau\in V_{t}$
\[ \| f(\tau)-f(t) \| < \epsilon \]
uniformly in $f\in S$.

Then the set $S$ is uniformly bounded, that is
\[ \sup_{t\in K, f\in S} \| f(t) \| =M <\infty. \]
\label{Ascoli}
\end{prop}

We show that the sequence $x_i$ satisfies properties (i) and (ii) of
Proposition \ref{Ascoli}. 

\noindent\underline{Equicontinuity:} Choose $t\leq t'$. We have
\[ dist_{L^2_1}(x_i(t),x_i(t'))\leq K \int_t^{t'} \| \nabla {\cal
C}_\rho (x_i(s))\|_{L^2_1}ds, \]
where we denote by $dist_{L^2_1}(x(t),x(t'))$ the quantity
$$ dist_{L^2_1}(x(t),x(t'))= | \int_t^{t'} \| \frac{d}{ds} x(s)
\|_{L^2_1}ds |. $$
Since the perturbation satisfies condition (4) of Definition
\ref{calP}, we obtain
\[ dist_{L^2_1}(x_i(t),x_i(t'))\leq K |t-t'|^{1/2} |\int_t^{t'} \| \nabla {\cal
C}_\rho (x_i(s))\|^2_{L^2_1}ds |^{1/2} \leq \]
\[ \tilde K |t-t'|^{1/2} |\int_{t-1}^{t'+1} \| \nabla {\cal
C}_\rho (x_i(s))\|^2_{L^2}ds |^{1/2}\leq \tilde K E |t-t'|^{1/2}. \]
Here the first step comes from the H\"older inequality and the second step
from Lemma 6.14 of \cite{MST}.

\noindent\underline{Pointwise bound:} Fix $t_0$. We can assume that 
\[ \| \nabla {\cal C}_\rho (x_i(t_0)) \|_{L^2}> \lambda \]
for all $i>>1$, where $\lambda$ is the constant of Lemma
\ref{PalaisSmale}. In fact, if there is a subsequence $x_{i_k}$ such
that this condition is not satisfied, then by Lemma \ref{PalaisSmale}
the elements $x_{i_k}(t_0)$ lie in an $\epsilon$-neighborhood of a
critical point hence their norms are bounded.

We can therefore choose a sequence $t_i$ of real numbers such that 
\[ \| \nabla {\cal C}_\rho (x_i(t_i)) \|_{L^2}= \lambda \]
and 
\[ \| \nabla {\cal C}_\rho (x_i(t)) \|_{L^2}> \lambda \]
for $t\in (t_i,t_0]$. Lemma \ref{PalaisSmale} then implies that 
there is an element $(A_a,\psi_a)$ with $[A_a,\psi_a]\in O_a$ and
there are gauge transformations $\lambda_i$ such that 
\[ dist_{L^2_1}(\lambda_i x_i(t_i),(A_a,\psi_a))\leq \epsilon \]
for all $i$, and therefore we can write
\[ dist_{L^2_1}(\lambda_i x_i(t_0),(A_a,\psi_a))\leq dist_{L^2_1}(\lambda_i
x_i(t_i),(A_a,\psi_a))+ 
dist_{L^2_1}(\lambda_i x_i(t_i),\lambda_i x_i(t_0))  \]
\[ \leq \epsilon + K \int_{t_i}^{t_0} \| \nabla {\cal C}_\rho
(\lambda_i x_i(t))\|_{L^2_1} dt. \]
Using the inequality
\[ \| \nabla {\cal C}_\rho
(\lambda_i x_i(t))\|_{L^2_1} \leq \frac{\| \nabla {\cal C}_\rho
(\lambda_i x_i(t))\|^2_{L^2_1} }{\| \nabla {\cal C}_\rho
(\lambda_i x_i(t))\|_{L^2}}\leq \frac{\| \nabla {\cal C}_\rho
(\lambda_i x_i(t))\|^2_{L^2_1}}{\lambda}, \]
for $t\in [t_i,t_0]$, and Lemma 6.14 of \cite{MST}, we obtain
\[ dist_{L^2_1}(\lambda_i x_i(t_0),(A_a,\psi_a))\leq \epsilon + \frac{
\tilde K}{\lambda} 
\int_{t_i-1}^{t_0+1} \| \nabla {\cal C}_\rho
(\lambda_i x_i(t))\|^2_{L^2} \leq \epsilon + \tilde K E /\lambda. \]

Now Proposition \ref{Ascoli} implies that, given any compact set
$K\subset \RR$, we have
\begin{equation}
\sup_{t\in K} \| \lambda_i x_i(t) \|_{L^2_1(Y\times \{t\})}=M <\infty 
\label{uniform}
\end{equation}
uniformly in $i$. 

For simplicity of notation we refer in the following to the gauge
transformed sequence $\lambda_i x_i$ simply as $x_i$.

From the uniform estimate (\ref{uniform}) and from the equations
(\ref{3SW1P'}) and (\ref{3SW2P'}) we obtain a uniform bound of the
$L^2$ norms of $\frac{d}{dt} x_i(t)$ on $Y\times K$:
\[ \| \frac{d}{dt} x_i(t) \|_{L^2(Y\times \{ t \}) }\leq C \| x_i(t)
\|_{L^2_1(Y\times \{t \})}, \]
hence the left hand side is bounded uniformly with respect to $t\in K$
and $i$. We use the fact that the perturbations 
$P_{(\AA_i,\Psi_i)}$ are uniformly bounded with respect to
$(\AA_i,\Psi_i)$ according to condition (2) of Definition \ref{calP}. 

In the following we shall always consider the set $K$ to be some large
interval $[-T,T]$. The previous estimates 
provide a uniform bound of the norms 
$\| x_i \|_{L^2_1(Y\times [-T,T])}$. In fact, we have
\[ \| x_i \|^2_{L^2_1(Y\times [-T,T])}\leq 2T M + 2T CM. \]

In order to bound the higher Sobolev norms we need the following gauge
fixing condition. 
Up to gauge transformations $\lambda_i$ in the 
group $\tilde {\cal G}_{k+1}(Y\times [-T,T])$, it is possible to make
$\AA_i-\AA_0$ into a sequence of co-closed 1-forms.

We have elliptic estimates for the domain $Y\times [-T,T]$,
\[ \|\AA_i-\AA_0\|_{L^2_k(Y\times [-T,T])} \leq C(\|(d^*
+d^+)(\AA_i-\AA_0)\|_{L^2_{k-1}(Y\times [-T',T'])} +\] 
\[ + \|\AA_i-\AA_0\|_{L^2_{k-1}(Y\times [-T,T])}) \]
and 
\[ \| \Psi_i \|_{L^2_k(Y\times [-T,T])}\leq C \left( \| \nabla_A
\Psi_i \|_{L^2_{k-1}(Y\times [-T',T'])}  
  +\| \Psi_i \|_{L^2_{k-1}(Y\times [-T,T])} \right), \]
where $[-T',T']$ is any strictly smaller interval contained in $[-T,T]$.
We also have an estimate given by the curvature equation:
\[ \| F_{\AA_i}^+\|_{L^2_k(Y\times [-T,T])} \leq
\|\Psi_i\cdot\bar\Psi_i \|_{L^2_{k}(Y\times [-T,T])} + \| i\mu +
P_{(\AA_i,\Psi_i)} \|_{L^2_{k}(Y\times [-T,T])}. \] 
These provide a uniform bound on the $L^2_k$-norms of the $(\AA_i,\Psi_i)$
on a smaller $Y\times [-T', T']$.  
By the Sobolev embeddings, this implies that on a smaller $Y\times
[-T'', T'']$ there is a subsequence $x_{i_k}=(\AA_{i_k},\Psi_{i_k})$
that converges uniformly with all derivatives.

This completes the proof of Claim 1: there is a subsequence $x_{i_k}$
that converges smoothly on compact sets to a solution $y=(\AA,\Psi)$
of (\ref{3SW1P'}) and (\ref{3SW2P'}) on $Y\times \RR$. 

The limit $y$ is of finite energy, hence it defines an element of some
moduli space ${\cal M}(O_c,O_d)$ according to Corollary
\ref{finenergy2}. 

For simplicity of notation, assume that the sequence $x_i$ itself
converges to $y$ smoothly on compact sets.

\noindent{\bf Claim 2:} If the limit $y$ is an element of ${\cal
M}(O_a,O_b)$ then the convergence $x_i\to y$ is strong in the
$L^2_{k,\delta}$ norm. 

\noindent\underline{Proof of Claim 2:} There exists a $T_0$ such that
for all $t\leq -T_0$ we have
\[ dist_{L^2_1}(y(t),(\tilde A_a,\tilde \psi_a))\leq \epsilon/2, \]
for some element $(\tilde A_a,\tilde \psi_a)$ with $[\tilde A_a,\tilde
\psi_a]\in O_a$.
For $i>>1$ we also have
\[ dist_{L^2_1}(x_i(-T_0),y(-T_0))\leq \epsilon/2, \]
hence 
\[ dist_{L^2_1}(x_i(-T_0),(\tilde A_a,\tilde \psi_a))\leq \epsilon. \]
This implies that $(\tilde A_a,\tilde \psi_a)$ is the same as the
element $(A_a,\psi_a)$. Moreover,
the exponential decay property ensures that we
have 
\begin{equation}
dist_{L^2_1}(x_i(t),(A_a,\psi_a))\leq \epsilon 
\label{uniformnear}
\end{equation} 
for all $t\leq -T_0$. This, together with the uniform convergence on
compact sets, implies that the convergence is uniform on all $Y\times
\RR$. Thus, we have a uniform exponential decay towards
the endpoints, hence convergence in the $L^2_{1,\delta}$ norm.

Notice that in this case we can ensure that there is a unique $T_0$
such that (\ref{uniformnear}) is satisfied for all $t\leq -T_0$,
whereas, in the proof of the exponential decay, the interval $(-\infty,
T_0]$ depends on the solution $(A(t),\psi(t))$. 

To improve the convergence to $L^2_{k,\delta}$ we need to choose
suitable gauge transformations.

\begin{lem}
There exist gauge transformations $\lambda_i$ in
${\cal G}^0_{2,\delta}(O_a,O_b)$ such that the forms $\lambda_i
(\AA_i-\AA_0)$ are co-closed, that is 
\[ d^*_{\delta} \lambda_i(\AA_i-\AA_0)=0. \]
\end{lem}

\noindent\underline{Proof of Lemma:} Let $\alpha_i =\AA_i-\AA_0$. The
element $e_{\delta/2}d^*e_{-\delta}\alpha_i$ is $L^2$-orthogonal to
the kernel of the Laplacian $\Delta_{\delta/2}$. Thus we can define
the elements 
\[ g_i= e_{-\delta/2} \Delta_{\delta/2}^{-1}
(e_{\delta/2}d^*e_{-\delta}\alpha_i), \]
with 
\[ \Delta_{\delta/2}^{-1}:
Ker(\Delta_{\delta/2})^{\perp}_{L^2_{1,\delta}}
\to Ker(\Delta_{\delta/2})^{\perp}_{L^2_{2,\delta}}. \]
These are elements in the $L^2_{2,\delta}$ completion of the Lie
algebra hence they define gauge transformations $\lambda_i\in {\cal
G}^0_{2,\delta}(O_a,O_b)$. 
We have
\[ \lambda_i(\AA_i-\AA_0)=\beta_i=\alpha_i-dg_i. \]
Thus, we can compute $d^*_{\delta}\beta_i$ as
\[ d^*_{\delta}\beta_i=d^*_\delta \alpha_i - d^*_\delta d
e_{-\delta/2} \Delta_{\delta/2}^{-1}
(e_{\delta/2}d^*e_{-\delta}\alpha_i) \]
\[ = d^*_\delta \alpha_i - e_{\delta/2} (d^*_{\delta/2}d_{\delta/2})
\Delta^{-1}_{\delta/2} e_{\delta/2}d^*e_\delta \alpha_i \]
\[ = d^*_\delta \alpha_i -d^*_\delta \alpha_i =0. \]

This completes the proof of the Lemma.

Now we can use the elliptic estimates to improve the convergence to
$L^2_{k,\delta}$.  This completes the proof of Claim 2.

Thus, we see that
broken trajectories arise when the limit $y$
is in some ${\cal M}(O_c,O_d)$ where either $O_c\neq O_a$ or $O_d\neq
O_b$. These are constrained by the condition 
\[ {\cal C}_\rho(a)\geq {\cal C}_\rho(c)\geq {\cal C}_\rho(d)\geq
{\cal C}_\rho(b). \] 
Without loss of generality, assume that $y$ is in ${\cal
M}(O_a,O_d)$ and ${\cal C}_\rho(a)>{\cal C}_\rho(d)> {\cal
C}_\rho(b)$. 

There exists a sequence of real numbers $T_i$ such that ${\cal
C}_\rho(x_i(T_i))=\alpha$, with $\alpha$ satisfying ${\cal
C}_\rho(d)>\alpha> {\cal C}_\rho(b)$. 
Consider the reparameterization $x_i^{ T_i}(t)=x_i(t+T_i)$. This is
another possible lifting of $\hat x_i$ to the parameterized
trajectories (which does not satisfy the equal energy condition). We
can apply the result of Claim 1 to the sequence $x_i^{T_i}$ and we
obtain a subsequence that converges smoothly on compact sets to a
solution $\tilde y$ of the equations (\ref{3SW1P'}) and (\ref{3SW2P'})
on $Y\times \RR$. Now the asymptotic values $\tilde y(\pm\infty)$ are
constrained by
\[ {\cal C}_\rho(d)\geq{\cal C}_\rho(\tilde y(-\infty))\geq {\cal
C}_\rho(\tilde y(\infty))\geq {\cal C}_\rho(b). \]

The process can be repeated: at each step we either get strong
convergence or another broken trajectory. If we start with relative
index $k+1$, the process ends at most at the 
order $k$, in fact the condition
\[ \mu(O_a)\geq \mu(O_c)>\mu(O_d)\geq \mu(O_b) \]
has to be satisfied for dimensional reasons.

\noindent $\diamond$

\begin{corol}
Notice that, in particular, if $\mu(O_a)-\mu(O_b)=1$ with $O_a$ a free 
orbit, or if $\mu(O_a)-\mu(O_b)=2$, with $O_a=\theta$ the unique fixed
point, the previous argument shows that $\hat{\cal M}(O_a,O_b)$ is
already compact, $\hat{\cal M}(O_a,O_b)=\hat{\cal M}(O_a,O_b)^*$.
\end{corol}

A brief notational remark about the case $b^1(Y)>0$:

\begin{rem}
In the case of a three-manifold with $b^1(Y)>0$, 
the Seiberg-Witten Floer homology considered in this paper differs
from the construction of \cite{Ma}, where the connected component of
the full gauge group $\tilde {\cal G}(Y)$ is considered. The relation
is illustrated in the following.
\label{correction}
\end{rem}

With our notation, we have $\tilde {\cal G}(Y)$, the full gauge group
on $Y$, and ${\cal G}(Y)$, the group of gauge transformations on $Y$
satisfying the condition (\ref{int=0}). We have another possible
choice of an even smaller subgroup, the identity component ${\cal
G}_c(Y)$ of the gauge group. These all have subgroups of
\index{${\cal G}^0$} 
base point preserving gauge transformations, $\tilde {\cal G}^0(Y)$,
${\cal G}^0(Y)$, ${\cal G}^0_c(Y)$.

Thus, we have configuration spaces $\tilde {\cal B}^0$, ${\cal B}^0$,
and ${\cal B}^0_c$. The map ${\cal B}^0_c \to \tilde{\cal B}^0$ has
fiber $H^1(Y,\ZZ)$, and the map ${\cal B}^0 \to \tilde {\cal B}^0$ has
fiber the group $H$ \index{${\cal B}^0$}
\[ H=\{ h\in H^1(Y,\ZZ) | \langle c_1(L)\cup h, [Y]\rangle =0\}. \]
The group $H$ is the kernel of the morphism
\[ \xi: H^1(Y,\ZZ)\to \ZZ, \]
\[ \xi(h)=\langle c_1(L)\cup h, [Y]\rangle. \]
The quotient $H^1(Y,\ZZ)/H$ is either trivial or $\ZZ$. It is
trivial only when $c_1(L)=0$ rationally, that is in the case of Remark 
\ref{Ybadcase} which will be considered separately in \cite{MW4}.
Thus we have the following diagram

$$
\spreaddiagramrows{-1pc}
\spreaddiagramcolumns{-1pc}
\diagram
& {\cal B}^0_c(Y)={\cal A}(Y)/{\cal G}^0_c(Y) \dlto_{H \quad}
\ddto^{H^1(Y,Z)} \\  
{\cal B}^0(Y)={\cal A}(Y)/{\cal G}^0(Y) \drto^{\qquad Z=H^1(Y,Z)/H} & \\
& \tilde {\cal B}^0(Y)={\cal A}(Y)/\tilde{\cal G}^0(Y) \\
\enddiagram
$$ 

\vspace{.2in}

As we are about to see, the choice of $ \tilde {\cal B}^0(Y)$ gives
rise to a periodic Floer homology of periodicity
\[ l =g.c.d \{ \langle c_1(L)\cup h, [Y] \rangle, h\in H^1(Y,\ZZ)
\}. \]

The choice of ${\cal B}^0(Y)$, that was suggested to us by R.G. Wang,
represents the minimal cover of $\tilde {\cal B}^0(Y)$ that
gives rise to a $\ZZ$-graded Floer homology. The choice of ${\cal
B}^0_c(Y)$ also gives rise to a $\ZZ$-graded complex.

In \cite{Ma} the Floer homology is constructed using the group
${\cal G}_c(Y)$. The relation to our construction here is as follows.
Let $\{ O_a^h \}_{h\in H}$ be the family of distinct critical orbits
in ${\cal B}^0_c$ that correspond to a critical orbit $O_a$ in ${\cal
B}^0$. Notice that in ${\cal B}^0_c$ the relative Morse index satisfies
$$ \mu(O_a)-\mu(O_a^h)=\int_Y c_1(L)\wedge h=0. $$ 

We have the configuration space of flow lines
\[ {\cal A}_{k,\delta}(O_a^h,O_b^{h'})=\bigcup_{\lambda\in U(1)} {\cal
A}_{k,\delta} (\lambda (A_a,\psi_a),\lambda (A_b,\psi_b)), \]
with $(A_a,\psi_a)$ and $(A_b,\psi_b)$ representatives of elements 
$[A_a,\psi_a]$ and $[A_b,\psi_b]$ in $O_a^h$ and $O_b^{h'}$
respectively, in the space ${\cal B}^0_c(Y)$. The space
${\cal A}_{k,\delta}(O_a^h,O_b^{h'})$ is acted upon by
the gauge group of gauge transformations in $\tilde
{\cal G}(Y\times \RR)$ that decay to elements of $G_a$ and $G_b$
respectively. We can form the quotient space ${\cal B}^0_c
(O_a^h,O_b^{h'})$. 

We can consider the moduli spaces of flow lines ${\cal
M}(O_a^h,O_b^{h'})$ in the space ${\cal B}^0_c(O_a,O_b)$. In the
analysis of the boundary structure of ${\cal M}(O_a^h,O_b^{h'})$,
there may be a question of whether a sequence of elements $[x_i]$
of ${\cal M}(O_a^h,O_b^{h'})$ can converge to a solution $y$ that
lives in another ${\cal M}(O_a^h,O_b^{h'})$ (for different $h$ and
$h'$). Claim 2 of Theorem \ref{seqcompact} rules out this
possibility. Thus, a nice compactification of the moduli spaces 
${\cal M}(O_a^h,O_b^{h'})$ can be achieved in this setting as well,
for all manifolds with $b^1(Y)>0$ and $c_1(L)\neq 0$. There is,
therefore, a well defined associated Seiberg-Witten Floer homology
which is finer than the one we are considering here.  

On the other hand we can also consider critical orbits $\tilde O_a$ in
the configuration space $\tilde {\cal B}^0(Y)$. In this case the
relative Morse index is not well defined, in fact the spectral flow
between two gauge equivalent orbits satisfies
\[ \mu(\lambda \tilde O_a)-\mu(\tilde O_a)=\int_Y c_1(L)\wedge h(\lambda), \]
with $h(\lambda)=\lambda^{-1}d\lambda$. Thus the grading of the Floer
complex is defined only up to a periodicity of
\[ l=g.c.d \{ \langle c_1(L)\cup h(\lambda), [Y]\rangle \}. \] 

Therefore in the configuration space $\tilde {\cal B}(\tilde
O_a,\tilde O_b)$ obtained as a quotient of 
\[ {\cal A}_{k,\delta}(\tilde O_a,\tilde O_b)=\bigcup_{\lambda\in U(1)} {\cal
A}_{k,\delta} (\lambda (A_a,\psi_a),\lambda (A_b,\psi_b)), \]
with $(A_a,\psi_a)$ and $(A_b,\psi_b)$ representatives of elements 
$[A_a,\psi_a]$ and $[A_b,\psi_b]$ in $\tilde O_a$ and $\tilde O_b$
respectively, in the space $\tilde {\cal B}^0(Y)$. The space
${\cal A}_{k,\delta}(\tilde O_a,\tilde O_b)$ is acted upon by
the gauge group ${\cal G}_{k+1,\delta}(\tilde O_a, \tilde O_b)$ of
gauge transformations in $\tilde {\cal G}(Y\times \RR)$ that decay to
elements of $G_a$ and $G_b$. 

Now the space ${\cal M}(\tilde O_a, \tilde O_b)$ has infinitely many
components of virtual dimensions
\[ \mu(O_a)-\mu(O_b)+ k l, \]
with $k\in\ZZ$. Each component has uniformly bounded energy, hence it
has a compactification by boundary strata of broken trajectories. 

The Floer homology defined this way, namely with generators in
$\tilde{\cal M}$ and with the boundary operator defined by counting
the components of minimal energy in ${\cal M}(\tilde O_a, \tilde
O_b)$, is finitely generated and $\ZZ_l$-graded. There are various
other ways of defining the Floer homology for $b^1(Y)>0$. For
instance, a way of obtaining a $\ZZ$--graded Floer homology which is
also finitely generated (again for $c_1(L)\neq 0$) is by a filtration
as in \cite{FintStern} or \cite{Weiping}. The comparison between
these different versions will be discussed in \cite{MW4}. For the
purpose of this paper we shall only consider the 
Seiberg-Witten Floer homology defined by the choice of the gauge group
${\cal G}$.

\subsection{Gluing theorem}

Theorem \ref{seqcompact} proves that lower dimensional moduli spaces
appear naturally in 
the compactification of the spaces $\hat{\cal M}(O_a,O_b)$. In the rest
of this section we describe a gluing formula, thus proving that all
the broken trajectories in the ideal boundary which break through one
intermediate critical point actually occur as codimension one strata
in $\hat{\cal M}(O_a,O_b)^*$. The case of multiply broken trajectories 
is analyzed in Section 4.3.

In order to show that the gluing construction is well defined over 
the gauge classes, we need a preliminary discussion of slices of the 
gauge action.

\begin{defin}
Let $\Gamma=(\AA,\Psi)$ be an element in
${\cal A}_{k,\delta}(O_a,O_b)$. The slice ${\cal S}_\Gamma$ at
$\Gamma$ for the
action of the gauge group ${\cal G}^0_{k+1,\delta}(Y\times \RR)$, is
the set of elements $(\AA,\Psi)+(\alpha,\phi)$ in the configuration
space ${\cal A}_{k,\delta}(O_a,O_b)$ with the following properties.
The element $(\alpha,\phi)$ satisfies 
\[  (\alpha,\phi)\in Ker(G^*_{(\AA,\Psi)}), \] 
Where $G^*_{(\AA,\Psi)}$ is the adjoint with respect to the
$L^2_\delta$ norm of the linearization of the gauge group action.
There is a $T_0>>0$, such that, on $Y\times [T_0,\infty)$ we have
\[ (\alpha,\phi)=(\tilde\alpha,\tilde\phi)+ (e_b^-)^*(A_1,\psi_1), \]
and on $Y\times (-\infty, -T_0]$ we have
\[ (\alpha,\phi)=(\tilde\alpha,\tilde\phi)+
(e_a^+)^*(A_0,\psi_0), \]
with $(\tilde\alpha,\tilde\phi)$ is in $L^2_{k,\delta}$ on $Y\times
\RR$. The elements $(A_0,\psi_0)$ and $(A_1,\psi_1)$ are
representatives in the configuration space ${\cal A}$ over $Y$ of
elements in the orbits $O_a$ and $O_b$.
\label{slices}
\end{defin}

The following Lemma shows that there is a uniform
choice of slices for the gluing construction, cf. Lemma 2.1.4 and
2.1.5 of \cite{MorMrow}.

\begin{lem}
Given $\Gamma=(\AA,\Psi)$ in the space ${\cal A}_{k,\delta}(O_a,O_b)$,
consider $\mu_{(\AA,\Psi)}$, the first positive eigenvalue of the
operator $G^*_{(\AA,\Psi)}G_{(\AA,\Psi)}$. There is a ball of radius
$r(\mu)$ around $\Gamma$ in the configuration space ${\cal
A}_{k,\delta}(O_a,O_b)$ such that the intersection of the ball with
the slice ${\cal S}_{\Gamma}$ embeds into ${\cal
B}^0_{k,\delta}(O_a,O_b)$.

Moreover, there is a radius $r$ such that we have an embedding 
\[ {\cal S}_{(\AA,\Psi)}\cap B_{r}(\AA,\Psi) \hookrightarrow {\cal
B}^0_{k,\delta}(O_a,O_b) \] 
for all $[\AA,\Psi]$ in ${\cal M}(O_a,O_b)$. \index{${\cal S}_{\Gamma_{ab}}$} 
\label{uniformslice}
\end{lem} 

\noindent\underline{Proof:} The claim that the radius can be taken to
be $r=r(\mu)$ follows from Lemma 10.3.1 of \cite{MMR} (cf. similar
arguments in Lemma 2.1.5 of \cite{MorMrow}). 
We need to show that there is a uniform
lower bound for $\mu_{(\AA,\Psi)}$ as $[\AA,\Psi]$ varies in
${\cal M}(O_a,O_b)$. Consider a sequence of elements $x_i$ in ${\cal
M}(O_a,O_b)$. Theorem \ref{seqcompact} implies that a subsequence
converges smoothly on compact sets to a limit $y$ that defines an
element in some ${\cal M}(O_c,O_d)$. Since we are considering the
adjoint operator $G^*_{(\AA,\Psi)}$ with respect to the
$L^2_\delta$ norm, for a generic choice of metric and
perturbation, we have $\mu_y >0$, as discussed in Lemma
\ref{near:critical}. We also have $\mu_{x_i}\to \mu_y$, and this 
provides the uniform bound. In fact, according to Theorem
\ref{seqcompact}, only finitely many different components ${\cal
M}(O_c,O_d)$ can appear in the boundary of ${\cal M}(O_a,O_b)$.

\noindent $\diamond$

The following result is also useful in constructing a
gluing map at the level of gauge classes, cf. Lemma 2.1.8 of
\cite{MorMrow}. 

\begin{lem}
\label{slicelemma}
Consider a compact subset $K\subset {\cal M}(O_a,O_b)$. Then there is
a radius $r$, such that, for every $[\Gamma]\in K$, \index{$V_{\Gamma}$}
there is a $U(1)$-invariant submanifold $V_{\Gamma}\subset {\cal
A}_{k,\delta}(O_a,O_b)$ with the following properties.
Suppose given any solution $(\AA,\Psi)$ of equations (\ref{4SW1P'}) and
(\ref{4SW2P'}) that is in the ball of radius $r$ around $\Gamma$. Then
$V_{\Gamma}$ contains an element gauge equivalent to $(\AA,\Psi)$.
Moreover, the quotient by the action of the gauge group
${\cal G}^0_{k+1,\delta}(Y\times \RR)$ induces a map
\[ V_{\Gamma} \to {\cal M}(O_a,O_b) \]
that is a diffeomorphism on a neighborhood of $[\Gamma]$.
Moreover, there is a $T_0$ such that, for all $T\geq T_0$ we can find
families  of gauge transformations 
\[ \{ \lambda^\pm_x \mid x\in V_\Gamma \} \]
such that the elements in the set \index{$W^\pm_\Gamma$}
\[ W^\pm_\Gamma =\{ \lambda^\pm_x(x) \mid x\in V_\Gamma \} \]
are in a temporal gauge on $(-\infty, -T]$ in the case of $W^-_\Gamma$
and on $[T,\infty )$ in the case of $W^+_\Gamma$. The set
$W^\pm_\Gamma$ has the same properties of $V_\Gamma$, namely the
quotient with respect to the gauge action induces a diffeomorphism
\[ W^\pm_\Gamma \to {\cal M}(O_a,O_b) \]
on a neighborhood of $\Gamma$.
\end{lem}

\noindent\underline{Proof:} It is sufficient to take $V_\Gamma$ to be the
intersection of the space of solutions to (\ref{4SW1P'}) and
(\ref{4SW2P'}) in ${\cal A}_{k,\delta}(O_a,O_b)$ with the slice ${\cal
S}_{\Gamma}$. According to Lemma \ref{uniformslice}, there is a radius
$r_1$ such that the intersection of this set with the ball of radius
$r_1$ around $\Gamma$ embeds in ${\cal M}(O_a,O_b)$. 
In general, elements of $V_\Gamma$ are not in a temporal gauge,
however, following Corollary 3.1.7 of \cite{MMR}, there exists a
constant $C$ independent of $x\in V_\Gamma$ 
such that we obtain $\| \lambda^\pm_x(x) \|\leq C \| x \|$, and the
elements in the set  
\[ W^\pm_\Gamma =\{ \lambda^\pm_x(x) \mid x\in V_\Gamma \} \]
are in a temporal gauge on $(-\infty, -T]$ and on $[T,\infty )$,
respectively. 
This follows from the exponential decay of solutions towards the
endpoints. Upon rescaling the radius $r$ to $r/C$ the properties
stated for $V_\Gamma$ hold for $W^\pm_\Gamma$. The constant $C$
depends only on the set $K$, cf. Lemma 2.1.8 and Proposition 2.6.4 of
\cite{MorMrow}. 

\noindent $\diamond$

Now we can introduce the gluing construction. This identifies the
codimension one boundary in the compactification ${\cal
M}(O_a,O_b)^*$. 

\begin{thm}
Suppose given $O_a$, $O_b$ and $O_c$ in ${\cal M}^0$  with
$\mu(O_a)>\mu(O_b)>\mu(O_c)$.  
Assume that $b$ is irreducible. Then, given a compact set 
\[ K\subset {\cal M}(O_a, O_b) \times _{O_b} {\cal M}(O_b, O_c), \]
there are a lower bound $T_0(K)>0$ and a smooth map \index{$\#$}
\[ \#: K\times [T_0,\infty)\to {\cal M}(O_a, O_c) \]
\[ ((\AA_1,\Psi_1),(\AA_2,\Psi_2),T)\mapsto
(\AA_1\#_T\AA_2,\Psi_1\#_T\Psi_2), \] 
such that $\#_T$ is an embedding for all $T>T_0(K)$.
The gluing map $\#$ induces a smooth embedding \index{$\hat \#$}
\[
\hat\#:  \hat K\times [T_0,\infty) \to \hat{\cal M}(O_a, O_c),
\]
where
\[  \hat K\subset \hat{\cal M}(O_a, O_b) \times _{O_b} \hat{\cal
  M}(O_b, O_c). \]
This local diffeomorphism is compatible with the orientation, where the
orientation on $K\times [T_0,\infty)$ is induced by the product
orientation on 
$$ {\cal M}(O_a, O_b) \times \RR \times  {\cal M}(O_b, O_c). $$
If $b$ is reducible, then $O_b$ is just a point and there is a similar
orientation preserving local diffeomorphism 
\[
\hat\#: \hat{\cal M}(O_a, b) \times \hat{\cal M}(b,  O_c)
\times [T_0,\infty) \to \hat{\cal M}(O_a, O_c).
\]
Moreover, any sequence of trajectories in $\hat{\cal M}(O_a,O_c)$
converging to a broken trajectory in $\hat{\cal
M}(O_a,O_b)\times_{O_b} \hat{\cal M}(O_b,O_c)$ lies eventually in the
image of the gluing map.
\label{equivgluing}
\end{thm}

In order to prove Theorem \ref{equivgluing},
we first define a pre-gluing map $\#_T^0$. This provides an
approximate solution via the following construction.  
Consider  classes $[x]=[\AA_1,\Psi_1]$ and $[y]=[\AA_2,\Psi_2]$ in
${\cal M}(O_a,O_b)$ and ${\cal M}(O_b,O_c)$ respectively. 
Let
\[ (A_1(t),\psi_1(t))\in {\cal A}_{k,\delta}(O_a,O_b), ~~
(A_2(t),\psi_2(t))\in {\cal A}_{k,\delta}(O_b,O_c) \] 
be temporal gauge representatives of $[x]$ and $[y]$.

Choose slices \index{${\cal S}_{\Gamma_{ab}}$}
\beq {\cal S}_{\Gamma_{ab}}\subset {\cal A}_{k,\delta}(O_a,O_b)
\label{slice1} \eeq
and 
\beq {\cal S}_{\Gamma_{bc}}\subset {\cal A}_{k,\delta}(O_b,O_c),
\label{slice2} \eeq
determined respectively by the elements 
\[ \Gamma_{ab} \in {\cal
A}_{k,\delta}((A_a,\psi_a),(A_b,\psi_b)) \]
and 
\[ \Gamma_{bc} \in {\cal
A}_{k,\delta}((A_b,\psi_b),(A_c,\psi_c)). \] 
We choose them so that there is a representative $x$ in the ball of
radius $r$ in ${\cal A}_{k,\delta}(O_a,O_b)$ centered at $\Gamma_{ab}$
and a representative $y$ in the ball of radius $r$ in ${\cal
A}_{k,\delta}(O_b,O_c)$ centered at $\Gamma_{bc}$. 
According to Lemma \ref{slicelemma},
there are gauge transformations $\lambda_1^+$ and $\lambda_2^-$ such that
we have
\beq \lambda_1^+ x \in W^+_{\Gamma_{ab}}, ~~~ \lambda_2^- y \in
W^-_{\Gamma_{bc}}. \label{slice:cond} \eeq  
Thus, for $|t|>T_0$ and for all $T\geq T_0$, 
we can write $\lambda_1^+ x(t)=\lambda_b
(A_b,\psi_b)+(\alpha_1(t),\phi_1(t))$ and 
$\lambda_2^- y(t)=\lambda_b
(A_b,\psi_b)+(\alpha_2(t),\phi_2(t))$. Here we have $\lambda_b\in
{\cal G}(Y)$ and $\lambda_b (A_b,\psi_b)$ is a representative of an element
in the orbit $O_b$ that satisfies
\[ \lim_{t\to\infty}(A_1(t),\psi_1(t))= \lambda_b
(A_b,\psi_b)=\lim_{t\to -\infty} (A_2(t),\psi_2(t)). \]
The elements $(\alpha_i,\phi_i)$ have exponentially decaying
behavior, as in Proposition \ref{decay}.

We construct an approximate solution $x\#^0_T y=(\AA_1\#_T^0 \AA_2,
\Psi_1\#_T^0 \Psi_2)$ of the form \index{$x\#^0_T y$}
\beq x\#^0_T y =\left\{ \begin{array}{lr}\lambda_1^+
(\AA_1^{2T},\Psi_1^{2T})&t\leq -1\\ 
\lambda_b (A_b,\psi_b)+\rho^-(t)(\alpha_1(t+2T),\phi_1(t+2T))+& \\
\rho^+(t)(\alpha_2(t-2T),\phi_2(t-2T))& -1\leq t\leq 1 \\ 
\lambda_2^- (\AA_2^{-2T},\Psi_2^{-2T})&t\geq 1.
\end{array}\right. \label{pre:map} \eeq
Here $\rho^\pm(t)$ are smooth cutoff functions with bounded
derivative, such that $\rho^-(t)$ is equal
to one for $t\leq -1$ and to zero for $t\geq 0$ and $\rho^+(t)$ is
equal to zero for $t\leq 0$ and to one for $t\geq 1$. 

\begin{lem}
\label{inducedgluing}
The pre-gluing map $\#_T^0$ descends to a map at the level of gauge
classes. Namely, it is well defined with respect to the choice of
slices of the gauge action.
\end{lem}

\noindent\underline{Proof:} Consider a different choice of slices 
\[ {\cal S}_{\tilde\Gamma_{ab}}\subset {\cal A}_{k,\delta}(O_a,O_b) \]
determined by 
\[ \tilde\Gamma_{ab} \in {\cal A}_{k,\delta} ((\tilde
A_a,\tilde \psi_a),(\tilde A_b,\tilde \psi_b)) \]
and
\[ {\cal S}_{\tilde\Gamma_{bc}}\subset {\cal A}_{k,\delta}(O_b,O_c) \]
determined by 
\[ \tilde\Gamma_{bc} \in {\cal A}_{k,\delta}((\tilde
A_b,\tilde \psi_b),(\tilde A_c,\tilde \psi_c)). \]
Suppose that $x$  and $y$ are now in balls of radius $r'\leq r$
centered at 
$\tilde\Gamma_{ab}$ and $\tilde\Gamma_{bc}$ respectively. Then,
according to Lemma \ref{slicelemma}, there are gauge
transformations $\tilde \lambda_1^+$ and $\tilde\lambda_2^-$ that
conjugate $x$ and $y$ so that we get
\[ \tilde \lambda_1^+ x \in W^+_{\tilde\Gamma_{ab}}, ~~~
\tilde\lambda_2^- y \in W^-_{\tilde\Gamma_{bc}}. \] 
Consider the pre-glued element of the form
\[ \tilde x\#^0_T \tilde y=\left\{ \begin{array}{lr}\tilde\lambda_1^+
(\AA_1^{2T},\Psi_1^{2T})&t\leq -1\\ 
\tilde \lambda_b (\tilde A_b,\tilde
\psi_b)+\rho^-(t)(\tilde\alpha_1(t+2T),\tilde\phi_1(t+2T))+& \\  
\rho^+(t)(\tilde\alpha_2(t-2T),\tilde\phi_2(t-2T))& -1\leq t\leq 1 \\ 
\tilde \lambda_2^- (\AA_2^{-2T},\Psi_2^{-2T}) &t\geq 1.
\end{array}\right. \]
The temporal gauge representatives are defined up to an ambiguity
given by the action gauge transformations that are constant in the
$\RR$ direction, that is by the action of ${\cal G}(Y)$. Thus, 
on the interval $[-1,1]$, the
gauge transformations $\tilde\lambda_i^\pm \circ {\lambda_i^\pm}^{-1}$
are constant with respect to $t$,
\[ \tilde\lambda_i^\pm \circ {\lambda_i^\pm}^{-1}|_{Y\times [-1,1]} =
\tilde \lambda_b \lambda_b^{-1}. \] 
Thus $\tilde\lambda_1^+ \circ {\lambda_1^+}^{-1}$ and
$\tilde\lambda_2^- \circ {\lambda_2^-}^{-1}$ 
define a gauge transformation $\lambda$ in ${\cal G}(Y\times \RR)$
(depending on $T$) such that  
\[ \tilde x\#^0_T \tilde y= \lambda( x\#^0_T y ) \]
is satisfied. Thus, Lemma \ref{uniformslice} and Lemma
\ref{slicelemma} imply that the 
pre-gluing map induces a well defined  map $\#_T^0$ on compact subsets
of broken trajectories,
\[ \#_T^0: K\subset \hat{\cal M}(O_a,O_b)\times_{O_b}\hat{\cal M}(O_b,O_c)\to
\hat{\cal B}^0_{k,\delta}(O_a,O_c). \]

\noindent $\diamond$

The proof of Theorem \ref{equivgluing} will consist of showing that the
pre-glued approximate solution $x\#^0_T y$ can be perturbed to an
actual solution of the flow equations by a small perturbation.

We need a preliminary discussion of some properties of the
linearization of the Seiberg--Witten equations at the pre-glued solution.

Consider the weight function $e_{\delta^T}(t)=e^{\delta^T(t)}$, with
$\delta^T(t)$ a function that satisfies
$e_{\delta^T}(t)=e^{\delta (t+T)}$ for $t < 0$ and
$e_{\delta^T}(t)=e^{\delta (t-T)}$ for $t >0$. Here
$e_\delta(t)=e^{\delta(t)}$ is the usual weight function introduced in
Proposition \ref{sobolev}. We define the norm $L^2_{k,\delta(T)}$ as
$\| f \|_{L^2_{k,\delta(T)}}=\| f e_{\delta^T} \|_{L^2_k}$. We can
choose the weight function $e_{\delta^T}(t)$ as above, so that the
weighted norms satisfy the estimates
\beq C_1 e^{\delta T} \| f \|_{L^2_{k,\delta}} \leq  \| f
\|_{L^2_{k,\delta(T)}} \leq C_2 e^{\delta T} \| f \|_{L^2_{k,\delta}}, 
\label{equivnorms} \eeq
with constants $C_1$ and $C_2$ independent of $T$.

Suppose given classes $[x]=[\AA_1,\Psi_1]$ and $[y]=[\AA_2,\Psi_2]$ in
${\cal M}(O_a,O_b)$ and ${\cal M}(O_b,O_c)$ respectively. We can
consider the operator \index{${\cal L}_x$}
\beq {\cal L}_x(\alpha,\Phi)=\left\{ \begin{array}{l} d^+\alpha
-\frac{1}{2}Im(\Psi_1\cdot\bar\Phi)+{\cal
D}P_{(\AA_1,\Psi_1)}(\alpha,\Phi)\\
D_{\AA_1}\Phi+\alpha\Psi_1 \\
G^*_{\Gamma_{ab}}(\alpha,\Phi) \end{array}\right. \label{lin:x} \eeq
acting on the space of $L^2_{1,\delta}$-decaying 1-forms and spinor
sections on $Y\times \RR$.
Analogously we have 
\beq {\cal L}_y(\alpha,\Phi)= \left\{ \begin{array}{l} d^+\alpha
-\frac{1}{2}Im(\Psi_2\cdot\bar\Phi)+{\cal
D}P_{(\AA_2,\Psi_2)}(\alpha,\Phi)\\
D_{\AA_2}\Phi+\alpha\Psi_2 \\
G^*_{\Gamma_{bc}}(\alpha,\Phi) \end{array}\right. \label{lin:y} \eeq
acting on $L^2_{1,\delta}$ 1-forms and spinor sections on $Y\times
\RR$. We also consider \index{${\cal L}_{x\#_T^0 y}$}
\beq {\cal L}_{x\#_T^0 y}(\alpha,\Phi)=\left\{ \begin{array}{l} d^+\alpha
-\frac{1}{2}Im(\Psi_1\#^0_T\Psi_2\cdot\bar\Phi)+{\cal
D}P_{(\AA_1\#^0_T \AA_2,\Psi_1\#^0_T\Psi_2)}(\alpha,\Phi)\\
D_{\AA_1\#^0_T\AA_2}\Phi+\alpha\cdot (\Psi_1\#^0_T\Psi_2) \\
G^*_{\Gamma_{ac}}(\alpha,\Phi) \end{array}\right. \label{lin:T} \eeq 
acting on $L^2_{1,\delta(T)}$ 1-forms and spinor sections on $Y\times
\RR$. 
Moreover, if $(A_b,\psi_b)$ is a representative of a point on $O_b$,
with $O_b$ a free $U(1)$-orbit, consider the operator \index{${\cal L}_b$}
\beq {\cal L}_b(\alpha,\Phi)=\left\{ \begin{array}{l} d^+\alpha
-\frac{1}{2}Im(\psi_b\cdot \bar\Phi)
\\ D_{A_b}\Phi +\alpha\psi_b \\
G^*_{(A_b,\psi_b)}(\alpha,\Phi) \end{array}\right. \label{lin:birr} \eeq
acting on $L^2_1$ 1-forms and spinor sections on $Y\times \RR$.
If $O_b\equiv b$ is fixed by the $U(1)$-action, consider the operator
\beq {\cal L}_b(\alpha,\Phi)=\left\{ \begin{array}{l} d^+\alpha 
\\ D_{A_b}\Phi \\
G^*_{(A_b,\psi_b)}(\alpha,\Phi) \end{array}\right.  \label{lin:bred} \eeq
acting on $L^2_1$ 1-forms and spinor sections on $Y\times \RR$.

The central technique in the gluing construction is the study of the
behavior of the eigenspaces of small eigenvalues of the Laplacians 
${\cal L}_{x\#_T^0 y}^*{\cal L}_{x\#_T^0 y}$ and ${\cal L}_{x\#_T^0
y}{\cal L}_{x\#_T^0 y}^*$ of the linearization
of the Seiberg--Witten equations at the approximate solutions. The
general philosophy can be summarized as follows: we are trying to
paste together solutions of a system of non-linear elliptic
equations. If the equations were linear, we would simply encounter the 
obstruction to inverting the operator ${\cal L}_{x\#_T^0 y}^*$. This
obstruction is represented by $Coker({\cal L}_{x\#_T^0 y})$. However,
the equations being non-linear, the presence of eigenvectors with
small eigenvalues of ${\cal L}_{x\#_T^0
y}{\cal L}_{x\#_T^0 y}^*$ also represents obstructions. (cf. the
discussion of this phenomenon in \cite{Ta}, pg. 169). For this reason, 
we introduce the notation ``approximate cokernel'', $ApprCoker({\cal
L}_{x\#_T^0 y})$ (cf. Definition \ref{defapprker}) to denote this
eigenspace of small eigenvalues. The notation may be slightly
confusing, as it may suggest a space of approximate eigenvectors,
whereas it simply denotes the eigenspace of small eigenvalues, as
explained in Definition \ref{defapprker}. We trust that the reader
will not be confused by our, perhaps unorthodox, choice of notation. 

So the central technical issue in the gluing construction becomes
relating the small eigenvalues eigenspace of the Laplacian ${\cal L}_{x\#_T^0
y}{\cal L}_{x\#_T^0 y}^*$ (we call this space $ApprCoker({\cal
L}_{x\#_T^0 y})$) to the kernels of the Laplacians ${\cal L}_x {\cal
L}_x^*$, ${\cal L}_y {\cal L}_y^*$, and ${\cal L}_b {\cal L}_b^*$, and
similarly for the other  
Laplacians ${\cal L}_x^* {\cal
L}_x$ and ${\cal L}_y^* {\cal L}_y$ and the small eigenvalues
eigenspace of ${\cal L}_{x\#_T^0
y}^*{\cal L}_{x\#_T^0 y}$ (which we call $ApprKer({\cal L}_{x\#_T^0
y})$). 

First of all we can identify explicitly the kernel and cokernel of the
operator ${\cal L}_b$ as in (\ref{lin:birr}) and (\ref{lin:bred}).

\begin{lem}\label{ker:coker:Lb}
The operator ${\cal L}_b$ defined as in
(\ref{lin:bred}), that is, for $O_b\equiv b$ the fixed point of the
$U(1)$-action, has trivial $L^2_{1,\delta}$--kernel and 
trivial $L^2_{0,\delta}$--cokernel. The operator ${\cal L}_b$
defined as in (\ref{lin:birr}), with $O_b$ a free $U(1)$-orbit, also has
trivial kernel in $L^2_{1,\delta}$, and trivial cokernel in
$L^2_{0,\delta}$. 
\end{lem}

\noindent\underline{Proof:} We prove the case of
(\ref{lin:bred}), where the critical orbit $O_b\equiv b$ is the unique 
fixed point of the $U(1)$ action $b=\theta$. The irreducible case,
with a free $U(1)$ orbit $O_b$ follows by an analogous argument. 

Let $(\alpha,\Phi)$ satisfy 
$$ \left\{ \begin{array}{l} d^+\alpha=0 \\ D_{A_b}\Phi=0 \\
G^*_{A_b,0}(\alpha,\Phi)=0, \end{array} 
\right. $$
where the operator $G^*_{A_b,0}$ is the adjoint of the infinitesimal gauge
group action with respect to the $L^2_{0,\delta}$ inner product.
In this case, this is simply given by the operator $e_{\delta}
d^* e_{-\delta}$. The Dirac equation $D_{A_b}\Phi=0$ can
be rewritten in terms of the 3-dimensional Dirac operator,
$(\partial_t + \partial_{A_b})\Phi =0$. Upon expanding the spinor on a 
basis of eigenvectors for the 3-dimensional Dirac operator
$\partial_{A_b}$, we obtain $\Phi=\sum a_k(t) \phi_k$ and the
equation has solutions $\Phi(t)=\sum_k a_k(0) e^{\lambda_k t}
\phi_k$. Since none of these solutions is in $L^2_{1,\delta}$, we have
obtained $\Phi\equiv 0$. Now we consider the curvature equation. This
gives $d^+ \alpha=0$ and $e_{\delta} d^* e_{-\delta} \alpha =0$. These 
two conditions imply that $e_{-\delta/2}\alpha$ is a harmonic
representative in $H^{2+}_{\delta/2}(Y\times \RR)$, with respect to the
Laplacian $\Delta_{\delta/2}=e_{\delta/2}\Delta e_{-\delta/2}$. For a
rational homology sphere $Y$ (the only case where we have to consider the
presence of the reducible point $\theta$), this is necessarily
trivial. The argument for the cokernel is similar. Consider an element
in the cokernel, namely, suppose given a 
triple $(\beta,f,\xi)$ of a self-dual 2-form, a 0-form, and a spinor
section of the bundle of negative spinors over $Y\times \RR$, such
that we have
$$ \langle \beta, d^+\alpha \rangle + \langle f, e_{\delta}
d^* e_{-\delta} \alpha \rangle + \langle \xi, D_{A_b}\Phi \rangle =0$$
in the $L^2$ pairing. We then have $D_{A_b}\xi =0$ which implies
$\xi=0$ by the previous argument, and $d^* \omega=e_{-\delta} d
e_{\delta}f$. This equation implies $d^* \omega =0$ and $e_{-\delta} d 
e_{\delta}f=0$. 
The first equation then implies that $\omega$ is a harmonic
representative in $H^{2+}(Y\times \RR)$. Again, this is trivial since
the only case where the reducible $\theta$ appears is when $Y$ is a
rational homology sphere.
The remaining equation gives $f=c e_{-\delta}$. The condition
$f(x_0)=0$ at the base point $x_0=(y_0,t_0)$ implies $c= 0$, hence
$f\equiv 0$.

\noindent $\diamond$

The following Lemma, together with  Lemma \ref{surj2} and 
Lemma \ref{apprkercoker}, contains the key argument that will be used
in Proposition \ref{apprkercoker2} to relate the spaces
$ApprCoker({\cal L}_{x\#_T^0 y})$ and $ApprKer({\cal L}_{x\#_T^0 y})$
of eigenvectors of small eigenvalues to the kernels
and cokernels of ${\cal L}_x$ and ${\cal L}_y$ and in Corollary
\ref{indexsplit} for the formula on the splitting of the index. 

\begin{lem}
Let $x$ and $y$ be elements in ${\cal M}(O_a,O_b)$ and ${\cal
M}(O_b,O_c)$ respectively. 
Suppose given a sequence $\xi_k$ of $L^2_{1,\delta}$ 1-forms and
spinor sections on $Y\times \RR$. Assume these elements live in the
tangent space, in the $L^2_{1,\delta}$--norm, to the fiber
$(e_a^+)^{-1}(x_a)$ of the endpoint map
$$ e_a^+ : {\cal S}_{\Gamma_{ac}} \to O_a, $$ 
with $x_a$ the asymptotic value at $t\to -\infty$ of $x\in {\cal
M}(O_a,O_b)$, and for a given choice of slices ${\cal
S}_{\Gamma_{ab}}$ and ${\cal 
S}_{\Gamma_{bc}}$ as in (\ref{slice1}) and (\ref{slice2}), and 
$$ {\cal S}_{\Gamma_{ac}} \subset {\cal A}_{1,\delta}(O_a,O_c). $$
Suppose given a corresponding sequence of gluing parameters $T_k$ with 
$T_k \to \infty$. We
assume that $\xi_k$ and $T_k$ satisfy
$e^{\delta T_k}\| \xi_k \|_{L^2_{1,\delta}}=1$ and $\| {\cal
L}_{x\#^0_{T_k} y} \xi_k \|_{L^2_{0,\delta(T_k)}}\to 0$, with ${\cal
L}_{x\#^0_{T_k} y}$ as in (\ref{lin:T}).  
We have $L^2_1$ convergence
$e^{\delta T_k}\xi_k |_{Y\times [-T_k,T_k]}\to 0$. 
Then there exist
$L^2_{1,\delta}$-elements $u$ and 
$v$ in $Ker({\cal L}_x)$ and $Ker({\cal L}_y)$ respectively, with
${\cal L}_x$ and ${\cal L}_y$ as in (\ref{lin:x}) and (\ref{lin:y}),
such that we have convergence 
$$\| \rho_{1-T_k}^- \xi_k^{-T_k}- u \|_{L^2_{1,\delta}}\to 0 $$ 
and
$$ \| \rho_{-1+T_k}^+\xi_k^{T_k} - v \|_{L^2_{1,\delta}}\to 0. $$ 
Here we have $\rho^+_{-1+T_k}(t)=\rho^+(t-1+T_k)$ and
$\rho^-_{1-T_k}(t)=\rho^-(t+1-T_k)$. 
\label{split:convergence}
\end{lem}

\noindent\underline{Proof:} We first want to prove that 
the supports of the $\xi_k$ become more and more
concentrated at the asymptotic ends as $k\to\infty$.
Consider the operator ${\cal L}_b$ as in (\ref{lin:birr}) or
(\ref{lin:bred}). 
If $\zeta:Y\times \RR\to[0,1]$ is a smooth function which is equal to 1 on
$Y\times [-1/2,1/2]$ and equal to zero outside $Y\times
(-1,1)$. Moreover suppose that, for any fixed 
$t$, the function $\zeta$ is constant on $Y\times \{t\}$. Let
$\zeta_k(t)=\zeta(\frac{t}{T_k})$.
Then, using the pointwise estimate $| \zeta_k |\leq 1$, we have  
\[ \| {\cal L}_b ~\zeta_ke^{\delta T_k}\xi_k \|_{L^2_{0,\delta}} \leq \|
\zeta_k^\prime e^{\delta T_k}\xi_k \|_{L^2_{0,\delta}} + \| 
\zeta_k e^{\delta T_k}{\cal L}_b \xi_k \|_{L^2_{0,\delta}}  \]
\[ \leq\frac{1}{T_k} \max |\zeta^\prime | +\|\zeta_k ({\cal L}_{x\#^0_{T_k}
y}-{\cal L}_b)e^{\delta T_k} \xi_k \|_{L^2_{0,\delta}} +\| e^{\delta
T_k} {\cal L}_{x\#^0_{T_k} y} \xi_k \|_{L^2_{0,\delta}}  \] 
\[ \leq \frac{1}{T_k} \max |\zeta^\prime |+ \sup_{t\in[-T_k,T_k]}\|
{\cal L}_{x\#^0_{T_k} y}-{\cal L}_b \|~ \|e^{\delta
T_k}\xi_k\|_{L^2_{1,\delta}} +  
\| {\cal L}_{x\#^0_{T_k} y}\xi_k \|_{L^2_{0,\delta(T_k)}}. \] 
The first and last term tend to zero as $k\to\infty$. The remaining
term 
$$ \sup_{t\in[-T_k,T_k]} \| {\cal L}_{x\#^0_{T_k} y}-{\cal L}_b \|~
\|e^{\delta T_k}\xi_k\|_{L^2_{1,\delta}} $$  
is bounded by 
\[ \sup_{t\in[-1,1]}\| {\cal L}_{x\#^0_{T_k} y}-{\cal L}_b\|  
+ \sup_{t\in[-T_k,-1]}\| {\cal L}_{x^{2T_k}}-{\cal L}_b\|
+\sup_{t\in[1,T_k]}\| {\cal L}_{y^{-2T_k}}-{\cal L}_b\|. \]
All these terms tend to zero because of the exponential decay to the
critical point $b$ of the trajectories $(A_1(t),\psi_1(t))$ and
$(A_2(t),\psi_2(t))$, and because of condition (5) of Definition
\ref{calP}. 
Thus, we have $\| {\cal L}_b ~\zeta_k e^{\delta T_k}\xi_k
\|_{L^2_{0,\delta}}\to 0$ as $k\to\infty$. 

Lemma \ref{ker:coker:Lb} shows that the condition
$\| {\cal L}_b ~\zeta_k e^{\delta T_k}\xi_k \|_{L^2_{0,\delta}}\to 0$
implies the convergence $\zeta_k e^{\delta T_k}\xi_k\to 0$ in the
$L^2_{1,\delta}$-norm, and therefore we have 
$$\| e^{\delta T_k} \xi_k
\|_{L^2_1(Y\times [-T_k,T_k])}\to 0.$$  
In other words, the $\xi_k$ get localized at the asymptotic ends
exponentially fast.

This result allows us to rephrase the convergence condition $\| {\cal
  L}_{x\#^0_{T_k} y} \xi_k \|\to 0$ in terms of the Fredholm operators
${\cal L}_{x}$ and ${\cal L}_{y}$: 
\[ \| {\cal L}_{x} (\rho^-_{1-T_k}\xi_k^{-T_k})
\|_{L^2_{0,\delta}}\leq \|{\rho^\prime}^-_{1-T_k} 
\xi_k^{-T_k}\|_{L^2_{0,\delta}} 
+ \| \rho^- {\cal L}_{x\#^0_{T_k} y} \xi_k
\|_{L^2_{0,\delta(T_k)}} \leq \] 
\[ C \| {\rho^\prime}^-_1 e^{\delta T_k} \xi_k
\|_{L^2_{0,\delta}} +
 \| \rho^- {\cal L}_{x\#^0_{T_k} y} \xi_k
\|_{L^2_{0,\delta(T_k)}} \leq \] 
\[   C \| e^{\delta T_k}\xi_k\|_{L^2_1(Y\times [-1,1])} + \| {\cal
L}_{x\#^0_{T_k} y} \xi_k \|_{L^2_{0,\delta(T_k)}} \to 0, \]
where $\rho^-_{1-T_k}(t)=\rho(t+1-T_k)$ and $\xi_k^{-T_k}(t)=\xi_k(t-T_k)$.

The linearization ${\cal L}_x$ is a Fredholm operator. The sequence 
$\rho^-_{1-T_k}\xi_k^{-T_k}$ is uniformly bounded in the
$L^2_{1,\delta}$ norm: in fact, by the estimate (\ref{equivnorms}) we
have  
\[ \| \rho^-_{1-T_k}\xi_k^{-T_k} \|_{L^2_{1,\delta}}\leq C \| 
\rho^-_1 e^{\delta T_k} \xi_k \|_{L^2_{1,\delta}} \leq C. \]
This implies the existence of an $L^2_{1,\delta}$-element $u$ in
$Ker({\cal L}_x)$ such that, upon passing to a subsequence, we obtain
the $L^2_{1,\delta}$ convergence 
$$ \| \rho^-_{1-T_{k'}}\xi_{k'}^{-T_{k'}}- u \|_{L^2_{1,\delta}} \to 0. $$
This implies the convergence
\beq \| \rho^-_1 e^{\delta T_{k'}}\xi_{k'}-
\rho^- e^{\delta T_{k'}}u^{T_{k'}}\|_{L^2_{1,\delta}}\to
0. \label{slide:conv:1} \eeq
Similarly we obtain an $L^2_{1,\delta}$-element $v$ in $Ker({\cal
L}_{y})$  and a subsequence such that 
$$ \| \rho^+_{-1+T_{k'}} \xi_{k'}^{T_{k'}}- v \|_{L^2_{1,\delta}} \to 0. $$
and
\beq \| \rho^+_{-1}e^{\delta T_{k'}}\xi_{k'} -\rho^+ e^{\delta T_{k'}}
v^{-T_{k'}}\|_{L^2_{1,\delta}} \to 0. \label{slide:conv:2} \eeq
Here we have $\rho^+_{-1}(t)=\rho^+(t-1)$ and
$\rho^-_1(t)=\rho^-(t+1)$.

\noindent $\diamond$

Thus, by Lemma \ref{split:convergence} we have identified the
eigenspace of small eigenvalues of ${\cal L}_{x\#^0_T y}^* {\cal
L}_{x\#^0_T y}$ with a subspace of $Ker({\cal L}_x)\times Ker({\cal
L}_y)$. Lemma \ref{apprkercoker} will prove the reverse inclusion. The 
next Lemma proves a stronger result, namely that the operator ${\cal
L}_{x\#^0_T y}^* {\cal L}_{x\#^0_T y}$ is uniformly invertible, in the 
rescaled norm $L^2_{1,\delta(T)}$, on the complement of the eigenspace
of small eigenvalues. This is the essential result on the eigenvalue
splitting which we are going to use in Section 6, in Lemma
\ref{spectral:lemma}. 

Let us introduce the following notation. Let $F_{\#_T}$ be the
\index{$F_{\#_T}$}
linearization of the pre-gluing map (\ref{pre:map}). Thus, $F_{\#_T}$
is a map
\beq F_{\#_T}|_{(x,y)}: Ker({\cal L}_x) \times Ker({\cal L}_y) \to 
{\cal T}_1 ((e_a^+)^{-1}(x_a)) 
\subset {\cal T}_1 ({\cal B}^0_{k,\delta}(O_a,O_b)). \label{lin:pre-glue} \eeq
Here $x_a$ is the asymptotic value as $t\to -\infty$ of $x$ and
$e_a^+$ is the endpoint map $e_a^+: {\cal S}_{\Gamma_{ac}} \to O_a$.
Since the pre-glued solution (\ref{pre:map}) involves translations, the 
map $F_{\#_T}$ in (\ref{lin:pre-glue}) has operator norm bounded by $C_T\sim C
e^{-\delta T}$, namely, we have
\beq \| F_{\#_T}|_{(x,y)}(u,v) \|_{1,\delta} \leq C e^{-\delta T} \|
(u,v) \|_{1,\delta}. \label{estFT} \eeq
We may also consider the same map acting as
\beq F_{\#_T}|_{(x,y)}: Ker({\cal L}_x) \times Ker({\cal L}_y) \to
{\cal T}_{1,T} ((e_a^+)^{-1}(x_a))\subset {\cal
T}_{1,T} ({\cal B}^0_{k,\delta(T)}(O_a,O_b)), \label{lin:pre-glue:T} \eeq
with ${\cal T}_{1,T}$ the virtual tangent space in the
$L^2_{1,\delta(T)}$ norm. In this case $F_{\#_T}$ is bounded in the
operator norm uniformly in $T$, given the estimate (\ref{equivnorms})
for the rescaled norms.

\begin{lem}
\label{surj2}
There exist $T_0$ such that 
\[ \| {\cal L}_{x\#^0_T y} \xi \|_{L^2_{0,\delta(T)}}\geq C_T \|\xi
\|_{L^2_{1,\delta}}, \] 
for all $T\geq T_0$ and for all $\xi$ in the orthogonal complement
\[ {\cal T}^\perp=\left( F_{\#_T} (Ker({\cal L}_x) \times Ker({\cal L}_y))
\right)^\perp \]
in the space of $L^2_{1,\delta}$-connections and sections. With this
choice of norms, the constant $C_T$ grows like $C e^{\delta T}$. If we 
consider the rescaled norms $L^2_{\ell,\delta(T)}$ or the original
norms $L^2_{\ell,\delta}$ on both the source and target space, the
constant is bounded uniformly in $T$, given the estimate (\ref{equivnorms})
for the rescaled norms.
\end{lem}

\noindent\underline{Proof:}
Suppose there is a sequence of parameters $T_k$ and elements $\xi_k$ in the
orthogonal complement
\[ \xi_k \in\left( F_{\#_{T_k}} (Ker({\cal L}_x) \times Ker({\cal L}_y))
\right)^\perp \]
in the space of $L^2_{1,\delta}$-connections and sections,   
such that $e^{\delta T_k} \| \xi_k \|_{L^2_{1,\delta}}=1$ and $\| {\cal
L}_{x\#^0_{T_k} y} \xi_k \|_{L^2_{0,\delta(T_k)}}\to 0$. We use the
estimates of Lemma \ref{split:convergence} to derive a 
contradiction with the assumption $\xi_k\in{\cal T}^\perp$.   
We have
\[ 1=\lim_k e^{2\delta T_k} \|\xi_k\|^2_{L^2_{1,\delta}}=\lim_k
\langle \rho^-_1 e^{\delta T_k}
\xi_k,e^{\delta T_k}\xi_k\rangle_{L^2_{1,\delta}} + \langle
\rho^+_{-1}e^{\delta T_k}\xi_k,e^{\delta
T_k}\xi_k\rangle_{L^2_{1,\delta}}, \]  
since the remaining term satisfies 
\[ \lim_k\langle (1-\rho^-_1-\rho^+_{-1})e^{\delta T_k}\xi_k,e^{\delta
T_k}\xi_k\rangle=0. \]
This follows from Lemma \ref{split:convergence}, since 
$(1-\rho^-_1-\rho^+_{-1})$ is supported in $[-2,2]$, and we
are considering the case where $Ker({\cal L}_b)=0$ in $L^2_{1,\delta}$. 

Moreover, we know from Lemma \ref{split:convergence} that, upon
passing to a subsequence we have convergence as in
(\ref{slide:conv:1}) and (\ref{slide:conv:2}).
This implies that, upon passing to
a subsequence, we obtain
\[ \lim_k \langle \rho^-_1 e^{\delta T_k}\xi_k, e^{\delta T_k} \xi_k
\rangle_{L^2_{1,\delta}} =\lim_k e^{2\delta T_k} \langle \rho^-
u^{T_k}, \xi_k \rangle_{L^2_{1,\delta}}, \]
and correspondingly
\[ \lim_k \langle \rho^+_{-1}e^{\delta T_k}\xi_k,e^{\delta
T_k}\xi_k\rangle_{L^2_{1,\delta}}=\lim_k e^{2\delta T_k} \langle
\rho^+ v^{-T_k}, \xi_k \rangle_{L^2_{1,\delta}}. \]
The sum of these two terms therefore gives
\[ \lim_k e^{2\delta T_k} \langle F_{\#_{T_k}}(u,v), \xi_k
\rangle_{L^2_{1,\delta}}. \]
But this term is zero since we are assuming that the elements $\xi_k$
are in the orthogonal complement
\[ \xi_k \in \left( F_{\#_{T_k}} (Ker({\cal L}_x) \times Ker({\cal L}_y))
\right)^\perp. \]

\noindent $\diamond$

\begin{rem}
The results of Lemma \ref{split:convergence} and Lemma \ref{surj2} can 
be reformulated and proved, with minor modifications, for the operator 
${\cal L}_{x\#^0_T y}^*$ and elements 
$$(u,v)\in Ker({\cal L}_x^*)\times Ker({\cal
L}_y^*). $$ 
Namely, suppose given a sequence $\xi_k$ in 
${\cal T}_0({\cal S}_{\Gamma_{ac}})$, bounded uniformly in the
$L^2_{0,\delta(T_k)}$-norms, such that ${\cal L}_{x\#^0_{T_k} y}^* \xi_k
\to 0$ as $k\to\infty$. Then we have convergence $\| {\cal L}_b^*
\zeta_k \xi_k \|_{L^2_{1,\delta}} \to 0$, which implies
$L^2_{0,\delta}$--convergence of the  
elements $\zeta_k \xi_k$ to zero.
Moreover, there exist elements $u\in Coker({\cal
L}_x)$ and $v\in Coker({\cal L}_y)$ satisfying
$$ \| \rho_{1-T_k}^- \xi_k^{-T_k} -u \|_{0,\delta} \to 0 $$
and
$$ \| \rho_{-1+T_k}^+ \xi_k^{T_k} -v \|_{0,\delta} \to 0. $$
Thus, for all $T\geq T_0$, there is an estimate
$$ \| {\cal L}_{x\#^0_T y}^* \xi \|_{1,\delta} \geq C \| \xi
\|_{0,\delta}, $$
for all $\xi$ in
$$ (F_{\#_T}(Coker({\cal L}_x)\times
Coker({\cal L}_y)))^\perp \subset  
{\cal T}_0({\cal S}_{\Gamma_{ac}}). $$
The constant $C$ is independent of $T\geq
T_0$. \label{split:convergence:coker} 
\end{rem}

Thus, in Lemma \ref{surj2} and Remark \ref{split:convergence:coker},
we have obtained invertibility of the linearization of the
pre-glued solution on the orthogonal complement of the linearization of
the pre-gluing map. The norm of the inverse operator grows
exponentially with the gluing parameter $T$. This is an effect of
using the weighted Sobolev norms that are not translation
invariant. If we use the rescaled norm $\delta(T)$, or the original
$\delta$--norm on both the source and target space, we obtain uniform
invertibility, with the norm of the inverse operator uniformly bounded
in $T$. 

Consider the two Laplacians of the gluing construction, namely, the
self-adjoint operators \index{$H^0_{x\#_T^0 y}$} \index{$H^1_{x\#_T^0 y}$}
\[ H^0_{x\#_T^0 y}= {\cal L}_{x\#^0_T y}{\cal L}^*_{x\#^0_T y} \]
and
\[ H^1_{x\#_T^0 y}= {\cal L}^*_{x\#^0_T y} {\cal L}_{x\#^0_T y} \]
acting on $L^2_{0,\delta(T)}$ (or $L^2_{0,\delta}$) and on
$L^2_{1,\delta(T)}$ (or $L^2_{1,\delta}$), respectively.

\begin{lem}
\label{apprkercoker}
The operator $H^1_{x\#_T^0 y}$ has at least
\[ \dim Ker ({\cal L}_x) + \dim Ker ({\cal L}_y) \]
independent eigenvectors with small eigenvalues $\mu_T$ satisfying $\mu_T\to
0$ as $T\to \infty$. All these eigenvectors have fast decaying
eigenvalues, namely $\mu_T$ satisfying $\mu_T \leq Ce^{-\delta T}$.
The operator $H^0_{x\#_T^0 y}$ has at least
\[ \dim Coker ({\cal L}_x) + \dim Coker ({\cal L}_y) \]
independent eigenvectors with small eigenvalues, namely with
eigenvalues $\mu_T$ satisfying the property that $\mu_T\to
0$ as $T\to \infty$. Again, these have fast
decaying eigenvalues, namely with $\mu_T$
satisfying $\mu_T \leq Ce^{-\delta T}$. 
\end{lem}

\noindent\underline{Proof:} Let us prove the case of $H^1_{x\#_T^0 y}$ 
first. 
Suppose given an element $u$ in the $L^2_{1,\delta}$--kernel
$Ker({\cal L}_x)$. Consider the 
element $\xi_T=\rho^- u^{T}$. We obtain the following estimate
\[ \| {\cal L}_{x\#^0_T y}(\xi_T) \|_{L^2_{0,\delta(T)}} \leq \| {\cal
L}_x u \|_{L^2_{0,\delta}}  + \| {\cal L}_{x\#^0_T y} e^{\delta T} \xi_T
\|_{L^2_{1,\delta}(Y\times [-1,1])}. \] 
The first term on the right vanishes, and the second term can be
estimated as
\[  \| {\cal L}_{x\#^0_T y} e^{\delta T}\xi_T
\|_{L^2_{1,\delta}(Y\times [-1,1])} \leq C \|
e_\delta u \|_{L^2_1(Y\times [T-1,T+1])}, \]
by the effect of time translations on the norms.
This last term is certainly bounded by 
$$ C\cdot \| \xi_T \|_{L^2_{1,\delta}},   $$ 
Thus, using the relation \ref{equivnorms}, we obtain
$$ \| {\cal L}_{x\#^0_T y}(\xi_T) \|_{L^2_{0,\delta}} \leq C
e^{-\delta T} \| \xi_T \|_{L^2_{1,\delta}}. $$
The same argument can be repeated starting with elements 
$v\in Ker({\cal L}_y)$ and constructing $\xi_T=\rho^+ v^{-T}$. 
Notice moreover that if $u$ and $w$ are two independent vectors in
$Ker({\cal L}_x)$, then the corresponding elements $\xi^u_T=\rho^-
u^{T}$ and $\xi^w_T=\rho^- w^{T}$ are
independent vectors in $ApprCoker ({\cal L}_{x\#^0_T y})$, since
we have
\[ \langle \xi^u_T,\xi^w_T \rangle_{L^2_{1,\delta(T)}} \to \langle
u,w\rangle_{L^2_{1,\delta}}. \]
This proves the statement for $H^1_{x\#_T^0 y}$.
The argument for $H^0_{x\#_T^0 y}$ is analogous. 

\noindent $\diamond$

Thus we can give the following definition.

\begin{defin}
We define the approximate kernel and cokernel \index{$ApprKer$}
\index{$ApprCoker$} 
\[ ApprKer ({\cal L}_{x\#^0_T y}) ~~~  ApprCoker({\cal L}_{x\#^0_T
y})\]
of the linearization at
the pre-glued solution as the span of the eigenvectors of
$H^1_{x\#_T^0 y}$ and $H^0_{x\#_T^0 y}$ respectively, with small
eigenvalues, namely with 
eigenvalues $\mu_T$ satisfying $\mu_T\to 0$ as $T\to \infty$.
\label{defapprker}
\end{defin}

Lemma \ref{split:convergence} and Lemma \ref{surj2}, together with
Remark \ref{split:convergence:coker}, and
Lemma \ref{apprkercoker} give the following result. 

\begin{prop}
We have isomorphisms
\beq ApprKer ({\cal L}_{x\#^0_T y})\stackrel{\cong}{\to} Ker({\cal L}_x)
\times Ker({\cal L}_y) \label{apprkers:isom} \eeq
and
\beq ApprCoker({\cal L}_{x\#^0_T y})\stackrel{\cong}{\to} Coker({\cal
L}_x)\times Coker({\cal L}_y),
\label{apprcokers:isom} \eeq
given by the projections. 
Moreover, we obtain constraints on the dimensions of the actual
kernel and cokernel of the linearization at the pre-glued solution:
\[ \dim Ker ({\cal L}_{x\#^0_T y})\leq \dim Ker({\cal L}_x) +
\dim Ker({\cal L}_y) \]
and 
\[ \dim Coker ({\cal L}_{x\#^0_T y})\leq \dim Coker({\cal L}_x) + \dim
Coker({\cal L}_y). \] 
\label{apprkercoker2}
\end{prop}

\noindent\underline{Proof:} We have proved in Lemma
\ref{split:convergence} and Lemma \ref{surj2} that
$$ \dim ApprKer ({\cal L}_{x\#^0_T y})\leq \dim Ker({\cal L}_x) + \dim 
Ker({\cal L}_y), $$
and in Lemma \ref{apprkercoker} we have proved the reverse estimate,
namely that
$$ \dim ApprKer ({\cal L}_{x\#^0_T y})\geq \dim Ker({\cal L}_x) + \dim 
Ker({\cal L}_y). $$
This, together with the result of Lemma \ref{surj2} shows that the
linearization of the gluing map $F_{\#}$ gives the isomorphism of the
Proposition. The result for the cokernels is analogous, by Remark
\ref{split:convergence:coker} and Lemma \ref{apprkercoker}.

\noindent $\diamond$

\begin{corol}
A direct consequence of Proposition \ref{apprkercoker2} 
is the splitting of the index:
\beq \label{ind:split:formula}
 Ind({\cal L}_{x\#^0_T y})=Ind({\cal L}_x)+ Ind({\cal L}_y). \eeq
This can be rephrased as
\[ Ind_e({\cal L}_{x\#^0_T y})=Ind_e({\cal
L}_x)+ Ind_e( {\cal L}_y)-1, \]
in the case of gluing across a free orbit $O_b$, 
or as
\[ Ind_e({\cal L}_{x\#^0_T y})=Ind_e({\cal
L}_x)+ Ind_e( {\cal L}_y) \]
in the case of the fixed point $b=\theta$,
with $Ind_e$ as in Theorem \ref{relmorseind2},
satisfying
$$ Ind_e({\cal L}_{\AA,\Psi})=\mu(O_a)-\mu(O_b) +1 -\dim(G_a). $$ 
\label{indexsplit}
\end{corol}

\noindent\underline{Proof:}
The index of the linearization ${\cal L}_{x\#^0_T y}$ can be computed
as 
$$ Ind({\cal L}_{x\#^0_T y})=\dim(ApprKer ({\cal L}_{x\#^0_T y})) -
\dim (ApprCoker({\cal L}_{x\#^0_T y})). $$
This can be seen, for instance, by the
supertrace formula for the index. 
Proposition \ref{apprkercoker2} then gives the result. 
The statement for $Ind_e$ follows from the relation 
$$ Ker({\cal L}_x)\cong Ker_e({\cal L}_x)/\RR \ \hbox{and} \ Ker({\cal
L}_y)\cong Ker_e({\cal L}_y)/\RR, $$  
with $Ker_e$ the space of extended $L^2_{1,\delta}$ solutions as
described in Theorem \ref{relmorseind2}.

\noindent $\diamond$

All the previous discussion makes no assumption on the cokernels of the 
linearizations ${\cal L}_x$ and ${\cal L}_y$, so it applies equally to 
the unobstructed case analyzed in the rest of this section (the case
of the moduli spaces $\hat{\cal M}(O_a,O_b)$ of flow lines), and to the
obstructed case of the moduli spaces ${\cal M}(O_a,O_{a'})$ and ${\cal 
M}^P(O_a,O_b)$ introduced in Section 6, Theorem \ref{metrics}. The
obstructed case is discussed in Lemma \ref{solveqbeta} and Proposition
\ref{obstruction}. 

For the remaining of this
section we assume that $Coker({\cal L}_x)=0$ and $Coker({\cal
L}_y)=0$. When these
conditions are satisfied, we say that our gluing theory is
unobstructed. In this
case, we obtain the following result.

\begin{prop}
Let $K$ be a compact set 
\[ K\subset \hat{\cal M}(O_a,O_b)\times_{O_b}\hat{\cal M}(O_b,O_c). \]
There exist a bound $T_0(K)> 0$ such that, for all $T\geq T_0(K)$ and
for all broken trajectories 
\[ \left((A_1(t),\psi_1(t)),(A_2(t),\psi_2(t))\right)\in K, \] 
the Fredholm operator ${\cal L}_{x\#^0_T y}$ 
\[ {\cal L}_{x\#^0_T y}: {{\cal T}_1}_{x\#^0_T y}\to {{\cal
T}_0}_{x\#^0_T y} \] 
is surjective.
We are also assuming that $Coker({\cal L}_x)=0$ and $Coker({\cal
L}_y)=0$. Then the composition of the pre-gluing map $\#^0_T$ with the
orthogonal projection on the actual kernel $Ker ({\cal L}_{x\#^0_T y})$
gives an isomorphism 
\beq Ker({\cal L}_x)\times Ker({\cal L}_y)
    \stackrel{\cong}{\to} Ker ({\cal L}_{x\#^0_T
y}). \label{kers:isom} \eeq
\label{approxsurjective2}
\end{prop}

\noindent\underline{Proof:} We
know that ${\cal L}_{x\#^0_T y}$ is Fredholm of index $\mu(O_a)-\mu(O_c)$. We
also know that $\dim Ker_e({\cal L}_x)=\mu(O_a)-\mu(O_b)+\dim(O_a)$
  and $\dim Ker_e({\cal L}_y)=\mu(O_b)-\mu(O_c)+\dim(O_b)$.

We need to know that for any pair $x=(\AA_1,\Psi_1)$ and $y=(\AA_2,\Psi_2)$
there is a bound $T_0=T(x,y)$ such that ${\cal
  L}_{x\#^0_T y}$ is surjective for $T\geq T_0$. The compactness of $K$ will
ensure that there is a uniform such bound $T(K)$.

Lemma \ref{surj2} provides the crucial step in the argument: the
linearization ${\cal L}_{x\#^0_T y}$ has a bounded right inverse when
restricted to the orthogonal complement of $F_\# (Ker({\cal L}_x) \times
Ker({\cal L}_y))$.

Now Proposition \ref{approxsurjective2} follows form the splitting of
the index \ref{indexsplit}, since we obtain the
estimate 
\[ \dim Ker({\cal L}_{x\#^0_T y} )\leq \dim Ker ({\cal L}_x) + \dim
Ker({\cal L}_y) \]
\[ = Ind({\cal L}_x)+Ind({\cal L}_y) = Ind({\cal L}_{x\#^0_T y})\leq
\dim Ker({\cal L}_{x\#^0_T y}), \]
hence
\[ \dim Ker({\cal L}_{x\#^0_T y} )= \dim Ker ({\cal L}_x)
+ \dim Ker({\cal L}_y). \]

\noindent $\diamond$

In terms of the extended kernels, the isomorphism (\ref{kers:isom})
gives rise to the isomorphism
\beq \label{kers:e:isom} Ker_e({\cal L}_{x\#^0_T y})\cong Ker_e({\cal
L}_x)\times_{\RR} Ker_e({\cal L}_y). \eeq
The space on the right hand side is the tangent space to the fibered
product 
$$ {\cal M}(O_a,O_b)\times_{O_b} {\cal M}(O_b,O_c). $$

Now we can proceed to give a proof of Theorem \ref{equivgluing}.

\noindent\underline{Proof of Theorem \ref{equivgluing}:}
Given the approximate solution obtained via the pre-gluing map
$\#_T^0$, we prove
that this can be perturbed to an actual solution. 
Similar arguments are presented in \cite{DK}, \cite{Schwartz}.

In order to deform a pre-glued approximate solution to an actual
solution we need to construct a projection of pre-glued solutions onto
${\cal M}(O_a,O_c)$. Thus we need a normal bundle to ${\cal
M}(O_a,O_c)$ inside ${\cal B}^0_{k\delta} (O_a,O_c)$. We proceed
according to the following construction. \index{${\cal U}(O_a,O_c)$}
Let ${\cal U}(O_a,O_c)$ be the image in ${\cal B}^0(O_a,O_c)$
of the pre-gluing map $\#^0_T$ defined on ${\cal
M}(O_a,O_b)\times_{O_b}{\cal M}(O_b,O_c)$. Consider the Hilbert
bundles ${\cal T}_1$ and ${\cal T}_0$ given by the $L^2_{1,\delta}$ and 
$L^2_{0,\delta}$ tangent bundles of ${\cal A}^0(O_a,O_c)$.
Choose a slice 
\[ {\cal S}_{\Gamma_{ac}}\subset {\cal A}(O_a,O_c), \]
such that the pre-glued element $x\#^0_T y$ lies in a ball of radius
$r$ in ${\cal A}(O_a,O_c)$ centered at $\Gamma_{ac}$.
The flow equation
\beq
\label{f}
f(x\#^0_T y)=\left\{\begin{array}{l}
D_{\AA_1\#_T^0\AA_2}(\Psi_1\#^0_T\Psi_2)\\
F_{\AA_1\#_T^0\AA_2}^+
-(\Psi_1\#^0_T\Psi_2)\cdot\overline{(\Psi_1\#^0_T\Psi_2)} 
-i\mu -P_{(\AA_1\#_T^0\AA_2,\Psi_1\#^0_T\Psi_2)}\\
G^*_{\Gamma_{ac}}(\AA_1\#_T^0\AA_2,\Psi_1\#^0_T\Psi_2) \end{array}\right.
\eeq
induces a bundle map $f: {\cal T}_1 \to {\cal T}_0$. The linearization
${\cal L}_{x\#^0_T y}$ is the fiber derivative of $f$.

We assume here that the linearizations ${\cal L}_x$ and ${\cal L}_y$ are
surjective. Thus, consider the space 
\[ {\cal K}=\bigcup_{K\times [T_0,\infty)} Ker_e({\cal
  L}_x)\times_{\RR} Ker_e({\cal L}_y). \]
The image of ${\cal K}$ under the linearization $F_{\#}$ of the
pre-gluing map defines a sub-bundle of ${\cal T}_1$. We
consider the space ${\cal T}^\perp_{x\#^0_T y}$ of elements of the
tangent space ${\cal T}_1$ at the point $x\#^0_T y$ that are
orthogonal to \index{${\cal T}_{x\#^0_T y}$}
\[ {\cal T}_{x\#^0_T y}=F_{\#} (Ker_e({\cal L}_x)\times_{\RR}
Ker_e({\cal L}_y)). \]
The space ${\cal T}^\perp$ gives the normal bundle of the gluing
construction. 
In the case of the fixed point $b=\theta$, we have
$$ {\cal K}=\cup_{K\times [T_0,\infty)} Ker_e({\cal L}_x)\times
Ker({\cal L}_y) $$
and the normal bundle of the gluing construction given by
$$  {\cal T}^\perp_{x\#^0_T y}= (F_{\#_T}( Ker_e({\cal L}_x)\times
Ker({\cal L}_y)))^\perp. $$ 

Now we want to define the actual gluing map $\#$ that provides a
solution of the flow equations in ${\cal M}(O_a,O_c)$.
This means that we want to obtain a section $\sigma$ of ${\cal T}_1$ such
that the image under the bundle homomorphism given by the flow
equation is zero in ${\cal T}_0$. \index{$\sigma(x,y,T)$}
Moreover, we want this element 
$$ \sigma(x,y,T)=\sigma(\AA_1,\Psi_1,\AA_2,\Psi_2,T)$$ to
converge to zero sufficiently rapidly as $T\to\infty$, so that the
glued solution will converge to the broken trajectory in the limit
$T\to\infty$. 

The perturbation of the approximate solution to an actual solution can
be obtained as a fixed point theorem in Banach spaces, via the following
contraction principle. 

\begin{rem}
Suppose given a
smooth map $f:E\to F$ between Banach spaces of the form
\[ f(x)=f(0)+Df(0)x+N(x), \]
with $Ker(Df(0))$ finite dimensional, with a right inverse
$Df(0)\circ G=Id_F$, and with the nonlinear part $N(x)$ satisfying the
estimate 
\beq
\label{nonlinestimate}
\| GN(x)-GN(y) \|\leq C(\|x\|+\|y\| )\|x-y\| 
\eeq
for some constant $C>0$ and $x$ and $y$ in a small neighborhood
$B_{\epsilon(C)}(0)$. Then, with the initial condition $\| G(f(0))
\|\leq \epsilon/2$, there is 
a unique zero $x_0$ of the map $f$ in  $B_\epsilon(0)\cap G(F)$. This
satisfies $\| x_0 \|\leq \epsilon$.
\label{fixedpoint}
\end{rem}

The map $f$ is given in our case by the flow equation, viewed as a
bundle homomorphism ${\cal T}_1\mapsto {\cal T}_0$.

\begin{prop} Lemma \ref{surj2} and Proposition \ref{approxsurjective2}
imply that the 
linearization ${\cal L}_{x\#^0_T y}$ is uniformly invertible on the
orthogonal complement of the actual kernel
\[ {\cal T}^\perp_{x\#^0_T y} = Ker_e({\cal L}_{x\#^0_T y})^\perp \cong
(F_{\#_T}(Ker_e({\cal L}_x)\times_{\RR} Ker_e({\cal L}_y)))^\perp. \]
Consider the right inverse map of ${\cal L}$ restricted to ${\cal
T}^\perp$,  
\[ G: {\cal T}_0\to {\cal T}^\perp. \]
There is a $T(K)$ and a constant $C>0$ independent of $T$, such that
we have
\[ \|G_\chi\xi\|_{L^2_{1,\delta}}\leq C \|\xi \|_{L^2_{0,\delta}} \]
for $\chi=(x,y,T)\in K\times [T(K),\infty)$. 
Moreover, we can write $f$ as a
sum of a linear and a non-linear term, where the linear term is ${\cal
  L}$ and the nonlinear term is
\[ N_{(A(t),\psi(t))}(\alpha,\phi)=(\sigma(\phi,\phi)+{\cal
  N}q_{(A(t),\psi(t))}(\alpha,\phi),\alpha\cdot\phi). \]
Here we write the perturbation $q$
of equation (\ref{3SW2P'}) as sum of a linear and a non-linear term, 
$2q={\cal D}q+{\cal N}q$. 
Then the conditions of Remark \ref{fixedpoint} are satisfied.
This provides the existence of a unique correction term 
\[ \sigma(\AA_1,\Psi_1,\AA_2,\Psi_2,T)\in B_\epsilon(0)\cap {\cal
  T}^\perp \] 
satisfying $f(\sigma)=0$. This element $\sigma$ is smooth and it
decays to zero when $T$ is very large, as proved in Section 6, using
the estimate (\ref{gluedecay}). This means that the glued solution
converges to the broken trajectory $(x,y)$ for $T\to\infty$.
The gluing
map is given by
\begin{equation} \label{glue:map}
(A_1\#_T A_2,\psi_1\#_T\psi_2)=(A_1\#_T^0 A_2,\psi_1\#_T^0\psi_2) +
\sigma(\AA_1,\Psi_1,\AA_2,\Psi_2,T), \end{equation}
with a rate of decay  
\beq \| \sigma(\AA_1,\Psi_1,\AA_2,\Psi_2,T)
\|_{L^2_{1,\delta}}\leq C e^{-\delta T} \label{error:estimate} \eeq 
as $T\to \infty$.
\label{glue=fixedpoint}
\end{prop}

In Section 6, in Lemma \ref{solveqbeta} and Proposition
\ref{obstruction}, we shall consider a similar fixed point problem in
the presence of obstructions coming form non-vanishing cokernels. The
proof of Proposition \ref{glue=fixedpoint} is given in Section 6,
after Lemma \ref{solveqbeta}, since it follows from the more
general case discussed there. The rate of decay  $\|
\sigma(\AA_1,\Psi_1,\AA_2,\Psi_2,T) \|_{1,\delta}\leq C 
e^{-\delta T}$ as $T\to \infty$ is derived in Lemma \ref{solveqbeta}
in Section 6.2.

\begin{rem}
\label{rem:isom}
Notice that, if we have a nontrivial cokernel of ${\cal L}_x$ or ${\cal
L}_y$, then not only ${\cal T}^\perp$ is not be a bundle, but in
general ${\cal L}_{x\#^0_T y}$ may not be invertible on all of
the orthogonal complement of the actual kernel. We still know that
${\cal L}_{x\#^0_T y}$ is uniformly invertible on the orthogonal
complement of the approximate kernel, but in this case the  
kernel is only a proper subspace of the approximate kernel.
Thus there may be
elements $\chi=(x,y,T)$ in $K\times [T_0,\infty)$ that cannot be
perturbed to an actual solution. 
In fact, as seen in the proof of Lemma \ref{surj2},
the condition that the composition $\tilde F_{\#_T}$ of the
linearization $F_{\#_T}$ of the pre-gluing map $\#^0_T$ with the 
orthogonal projection onto $Ker ({\cal L}_{x\#^0_T y})$ gives an isomorphism
\[ Ker_e({\cal L}_x)\times_{\RR} Ker_e({\cal L}_y)
    \stackrel{\cong}{\to} Ker_e({\cal L}_{x\#^0_T y}) \]
is equivalent to the condition that
\[ Coker({\cal L}_x)\times Coker({\cal
L}_y)\stackrel{\cong}{\to} 
Coker ({\cal L}_{x\#^0_T y}). \] 
However, in the general case, where we may have non-trivial
$Coker({\cal L}_x)$ and $Coker({\cal L}_y)$ we expect to have the
isomorphism only at 
the level of the approximate kernel and cokernel, whereas the actual
kernel and cokernel satisfy the weaker condition
\[ \dim Ker({\cal L}_{x\#^0_T y}) =\dim Ker({\cal L}_x) +
\dim Ker({\cal L}_y) -k, \]
for some $k\geq 0$, and 
\[ \dim Coker({\cal L}_{x\#^0_T y}) =\dim Coker({\cal L}_x) +
\dim Coker({\cal L}_y) -k, \]
according to the computation of Lemma \ref{surj2}.
This case occurs for instance in the case discussed in Section 6, 
Theorem \ref{noglue} and Theorem \ref{boundary:P:obstr}.
\end{rem}

To continue with the proof of Theorem \ref{equivgluing}, we now look at 
the induced map
$$ \hat \# : \hat K \times [T_0,\infty) \to \hat{\cal M}(O_a,O_c). $$
Recall that we have an
identification of the moduli spaces $\hat{\cal M}(O_a,O_b)$
with the balanced moduli spaces ${\cal M}^{bal}(O_a,O_b)$ of classes
in ${\cal M}(O_a,O_b)$ that satisfy the equal energy condition
(\ref{equal:energy}). Thus, we can define the induced gluing map
$\hat \#$ by restricting the gluing map $\#$ of (\ref{glue:map}) 
to the subspace
$$ {\cal M}^{bal}(O_a,O_b)\times_{O_b} {\cal M}^{bal}(O_b,O_c) \subset 
{\cal M}(O_a,O_b)\times_{O_b} {\cal M}(O_b,O_c), $$
and then composing the image $x\#_T y \in {\cal M}(O_a,O_c)$ with a
time translation that determines an element $\widehat{x\#_T y}\in {\cal
M}^{bal}(O_a,O_c)$. We define $x\hat\#_T y=\widehat{x\#_T y}$.

The identification $\hat{\cal M}(O_a,O_b)\cong {\cal
M}^{bal}(O_a,O_b)$ determines an identification of the tangent spaces 
\begin{equation}
\label{tangent:id}
{\cal T}_{x} {\cal M}(O_a,O_b) \cong {\cal T}_{x} \hat{\cal
M}(O_a,O_b) \oplus \RR \cdot U_x 
\end{equation}
for $x\in {\cal M}^{bal}(O_a,O_b)$. Here we have  
$$ U_x = \Pi_{{\cal S}} \nabla {\cal C}_\rho (\lambda_1^+ x (t)), $$
where the gauge transformation $\lambda^+_1$ is as in
(\ref{slice:cond}) and $\Pi_{{\cal S}}$ is the projection on the
tangent space to the slice ${\cal T}({\cal S}_{\Gamma_{ab}})$. 
Here ${\cal T}$ denotes the
$L^2_{1,\delta}$ tangent space.

In order to prove that $\hat\#$ is a local diffeomorphism, it suffices
to show that the linearization $\tilde F_{\hat\#_T}$ is an
isomorphism, and then prove injectivity of $\#_T$.
Recall that we have tangent spaces
$$ {\cal T}_{x} {\cal M}(O_a,O_b)\cong Ker_e({\cal L}_x) \ \hbox{and} \
{\cal T}_{y} {\cal M}(O_a,O_b)\cong Ker_e({\cal L}_y). $$
We have a decomposition as in (\ref{tangent:id}) with an element $U_y$ 
defined analogously. Moreover, we have the isomorphism of 
Proposition \ref{approxsurjective2}, cf. Remark \ref{rem:isom},
which gives an isomorphism
\beq \label{isom:xy} \tilde F_{\#_T} :
{\cal T}_{(x,y)}({\cal M}(O_a,O_b)\times_{O_b}{\cal M}(O_b,O_c))
\stackrel{\cong}{\to} {\cal T}_{x \#_T y} {\cal M}(O_a,O_c), \eeq
where $\tilde F_{\#_T}$ is the composition of the linearization
$F_{\#_T}$ of the pre-gluing map with the projection onto the kernel
$Ker_e({\cal L}_{x\#_T^0 y})$. 

We can consider the subspace
\beq\label{subsp:xy}
{\cal T}_{(x,y)}(\hat{\cal M}(O_a,O_b)\times_{O_b}\hat{\cal
M}(O_b,O_c)) \oplus \RR\cdot (U_x,-U_y) \eeq
inside the space
$$ \begin{array} {rl} {\cal T}_{(x,y)}({\cal M}(O_a,O_b)\times_{O_b}{\cal
M}(O_b,O_c))\cong & {\cal T}_{(x,y)}(\hat{\cal
M}(O_a,O_b)\times_{O_b}\hat{\cal 
M}(O_b,O_c)) \\[2mm]
& \RR\cdot U_x \oplus \RR\cdot U_y. \end{array} $$
Here we have $x$ and $y$ in ${\cal M}^{bal}(O_a,O_b)$ and ${\cal
M}^{bal}(O_b,O_c))$, respectively.

\noindent{\bf Claim:} For all $T\geq T_0$, sufficiently large, the
restriction of the isomorphism 
(\ref{isom:xy}) to the subspace (\ref{subsp:xy}) gives the desired
isomorphism 
$$ \tilde F_{\hat\#_T}: {\cal T}_{(x,y)}(\hat{\cal
M}(O_a,O_b)\times_{O_b}\hat{\cal 
M}(O_b,O_c)) \oplus \RR\cdot (U_x, -U_y) \stackrel{\cong}{\to}{\cal
T}_{x \hat\#_T y}\hat{\cal M}(O_a,O_c). $$ 

\noindent\underline{Proof of Claim:}  Suppose there exist a sequence
$T_n\to \infty$ such that there exist $(\xi_n, \eta_n, \tau_n)$ in
the space (\ref{subsp:xy}), with the property that 
$$ F_{\#_{T_n}}(\xi_n,\eta_n, \tau_n) \in \RR \cdot U_{T_n}, $$
where $U_{T_n}$ is defined as
$$ U_{T_n}=\Pi_{{\cal S}} \nabla{\cal C}_\rho (\lambda_n (x\hat\#_{T_n}
y) (t) ), $$
with gauge transformations $\lambda_n$ as in (\ref{slice:cond}) and
the projection $\Pi_{{\cal S}}$ on the slice ${\cal
S}_{\Gamma_{ac}}$. We are using the identification
$$ {\cal T}_{x \#_{T_n} y} {\cal M}(O_a,O_c)\cong {\cal T}_{x \#_{T_n}
y} \hat {\cal M}(O_a,O_c) \oplus \RR \cdot U_{T_n}. $$

We have 
$$ F_{\#_{T_n}}(\xi_n,\eta_n, \tau_n)= \tilde\tau_n U_{T_n}, $$
for some $\tilde \tau_n \in \RR$. We can normalize the elements
$(\xi_n, \eta_n, \tau_n)$ so that $\tilde \tau_n =1$ for all $n$.
By the result of Proposition
\ref{glue=fixedpoint}, this implies that the pre-glued elements
satisfy
$$ \| F_{\#^0_{T_n}}(\xi_n,\eta_n) + \tau_n
F_{\#^0_{T_n}}(U_x,-U_y)-F_{\#^0_{T_n}}(U_x,U_y)
\|_{L^2_{1,\delta(T_n)}} \to 0 $$ 
as $n\to \infty$. We are using the fact that the error term 
in the gluing map (\ref{glue:map}) decays to zero as $T_n\to \infty$
as in Proposition \ref{glue=fixedpoint} and in (\ref{gluedecay}).
This then implies the following convergence
$$ \| \xi_n +(\tau_n -1)U_x \|_{L^2_{1,\delta}(Y\times
(-\infty, T_n-1])}\to 0 $$ and 
$$ \| \eta_n -(1+\tau_n)U_y
\|_{L^2_{1,\delta}(Y\times [+1-T_n,\infty)) }\to 0, $$
where we use the norm estimate for the map (\ref{lin:pre-glue:T}). 

This gives a contradiction, according to the
decomposition (\ref{tangent:id}), which prescribes that
$\xi_n \perp U_x$ and $\eta_n \perp U_y$, cf. Proposition 2.56
of \cite{Schwartz}.  

We can estimate explicitly the norm of the isomorphism
$\tilde F_{\hat\#_T}$. This is bounded by the product of the norm of the
linearization $F_{\#_T}$ at the pre-glued solution (\ref{lin:pre-glue}), 
the norm of the isomorphism 
$Ker({\cal L}_{x\#_T^0 y}) \cong Ker({\cal L}_{x\#_T y})$, and the norm of 
the isomorphism $Ker({\cal L}_{x\#_T y}) \cong Ker({\cal L}_{x\hat\#_T y})$
given by time translation.  The first norm is bounded
uniformly in $T$, if we use the rescaled norms $L^2_{1,\delta(T_n)}$
on the tangent space of ${\cal M}(O_a,O_c)$. The second norm is
bounded by a constant $C_2=C_2(K)$ over the compact set  
$K$ because of the decay of the error term in the gluing map
(\ref{glue:map}), cf. (\ref{gluedecay}). The third norm is bounded by
$C e^{\delta \tau(x\#_T y)}$, where 
$\tau(x\#_T y)$ is the unique time shift that maps $x\#_T y$ to
the element $\widehat{x\#_T y}$ satisfying the equal energy condition
(\ref{equal:energy}). Thus, on the compact set $K$ this norm is also
bounded by a term $C e^{\delta\cdot \tau(K)}$. Summarizing, we have
obtained the estimate $\| \tilde F_{\hat\#_T} \|\leq C $, on the
operator norm of $\tilde F_{\hat\#_T}$, uniformly in $T$, for
$$ \tilde F_{\hat\#_T}: {\cal T}_{1,\delta}({\cal M}(O_a,O_b)) \times
{\cal T}_{1,\delta}({\cal M}(O_b,O_c)) \to {\cal T}_{1,\delta(T)}
{\cal M}(O_a, O_c), $$
or $\| \tilde F_{\hat\#_T} \|\leq C e^{-\delta T}$, for
$$ \tilde F_{\hat\#_T}: {\cal T}_{1,\delta}({\cal M}(O_a,O_b)) \times
{\cal T}_{1,\delta}({\cal M}(O_b,O_c)) \to {\cal T}_{1,\delta}{\cal
M}(O_a, O_c). $$ 

We still have to prove the injectivity of the map $\#_T$ for all
sufficiently large $T\geq T_0$. Suppose there exist a sequence $T_k
\to \infty$ and elements $(x_{1,k},y_{1,k})\neq (x_{2,k},y_{2,k})$ in 
${\cal M}^{bal}(O_a,O_b)\times_{O_b}{\cal M}^{bal}(O_b,O_c)$,
such that
$$ x_{1,k}\hat\#_{T_k}y_{1,k} = x_{2,k}\hat\#_{T_k}y_{2,k} $$
for all $k$. 
Then, there exist a based $L^2_{2,\delta}$ gauge transformation
$\lambda$ on $Y\times \RR$ such that we have
$$ \lim_{k\to\infty} \| x_{1,k}\#^0_{T_k}y_{1,k} - \lambda
(x_{2,k}\#^0_{T_k}y_{2,k} ) \|_{L^2_{1,\delta}} = 0. $$
We are using here the fact that the time shifts
agree as $T_k \to \infty$. Upon passing to a subsequence, we know that 
there are limits $x_i =\lim_k x_{i,k}$ and $y_i=\lim_k y_{i,k}$, for
$i=1,2$, by the assumed compactness of $\hat K$. We obtain 
$$ x_1 \hat\#_{T_k} y_1 = x_1 \hat\#_{T_k} y_2, $$
for all $T_k$, but for $T_k\to\infty$ we have convergence of
$x_i \hat\#_{T_k} y_i$ to the broken trajectory $(x_i,y_i)$. We obtain 
$(x_1,y_1)=(x_2,y_2)=(x,y)$. Now consider the elements
$(u_{i,k},v_{i,k})=(x_{i,k}-x, y_{i,k}-y)$, for $i=1,2$. They satisfy 
$(u_{1,k},v_{1,k})\neq (u_{2,k},v_{2,k})$ and 
$$ F_{\#_{T_k}}|_{(x,y)} (u_{1,k},v_{1,k})=F_{\#_{T_k}}|_{(x,y)}
(u_{2,k},v_{2,k}), $$
for all $k$, which contradicts the fact that the map $F_{\#}$ is a
local diffeomorphism, as we just proved. 

Now a comment about orientations. Notice that the isomorphism of
Proposition \ref{approxsurjective2} induces an isomorphism 
$$ \Lambda^{max} Ker_e({\cal L}_x)\otimes \Lambda^{max}Ker_e({\cal L}_y)
\to \Lambda^{max} Ker_e({\cal L}_{x\#_T y}). $$
Thus, in the unobstructed case we are considering, where all the
cokernels are trivial, we obtain that the gluing map $\#_T$ is
compatible with the orientations. Using the decompositions
(\ref{tangent:id}) we obtain that the gluing map $\hat \#_T$ is also
compatible with the induced orientations on the moduli spaces
$\hat{\cal M}(O_a,O_b)$, cf. Proposition 3.6 of \cite{Schwartz}.
Notice that the orientation on
$$ (\hat{\cal M}(O_a,O_b)\times_{O_b} \hat{\cal M}(O_b,O_c))\times \RR $$
induced, under the decomposition (\ref{tangent:id}),  by the product
orientation on 
$$ {\cal M}(O_a,O_b)\times_{O_b}{\cal M}(O_b,O_c) $$
agrees with the pullback of the product 
orientation of
$$ \hat{\cal M}(O_a,O_b)\times \RR\times \hat{\cal
M}(O_b,O_c) $$ 
under the map $(x,y,T)\to (x,T,y)$.

We now come to the last statement in Theorem \ref{equivgluing}. We
prove that any sequence of trajectories in ${\cal
M}(O_a,O_c)$ converging to a broken trajectory in ${\cal
M}(O_a,O_b)\times_{O_b} {\cal M}(O_b,O_c)$
lies eventually in the image of the gluing map.

This requires a preliminary statement about the endpoint maps 
$$e_b^- : {\cal M}(O_a,O_b)\to O_b \ \hbox{ and } \ e_b^+ : {\cal
M}(O_b,O_c)\to O_b. $$

\noindent{\bf Claim:}
Suppose given a sequence $(x_n,y_n)$ in 
$$ K_1\times K_2 \subset {\cal M}(O_a,O_b)\times {\cal M}(O_b,O_c), $$
with $K_1$ and $K_2$ compact sets in the $L^2_{2,\delta}$-topology. 
Suppose that the endpoints 
$e_b^-(x_n)$ and $e_b^+(y_n)$ converge to the same element $x_b$ on
the critical orbit $O_b$, in the $L^2_{1}$-topology on ${\cal B}^0$. 
Then there is a subsequence $(x_n',y_n')$ converging to an element
$(x',y')$ in the fibered product ${\cal M}(O_a,O_b)\times_{O_b} {\cal
M}(O_b,O_c)$. 

\noindent\underline{Proof of Claim:}
The compactness of $K_1\times K_2$ ensures the existence of a
convergent subsequence $(x_n',y_n')$ with limit $(x',y')$. Moreover,
we can estimate the distance
$$ \| e_b^-(x') - e_b^+(y') \|_{L^2_1} \leq  \| e_b^-(x_n)- e_b^-(x')
\|_{L^2_1} +$$ 
$$ \| e_b^+(y_n)- e_b^+(y')\|_{L^2_1} +
\| e_b^-(x_n)- e_b^+(y_n)\|_{L^2_1}. $$ 
The last term goes to zero by hypothesis, and the first and second
term on the right hand side go to zero by the continuity of the
endpoint maps. Thus, we have $e_b^-(x')=e_b^+(y')=x_b$.

Now we return to the statement on the range of the gluing map. This is 
the analogue in our setting of the method of continuity used in
Section 7.3 of \cite{DK}, and of the arguments of
Lemma 4.5.1 and Section 4.6 of \cite{MorMrow}

Suppose given a sequence $X_k$ of elements in ${\cal M}(O_a,O_c)$
which converges smoothly on compact sets (as in Theorem
\ref{seqcompact}) to a broken trajectory $(x,y)$ in the fibered
product ${\cal M}(O_a,O_b)\times_{O_b} {\cal M}(O_b,O_c)$. 

For a sequence of gluing parameters $T_k \to \infty$, let $(x_k,y_k)$
be the projection onto the slices 
${\cal S}_{\Gamma_{ab}} \times {\cal S}_{\Gamma_{bc}}$, 
$$ {\cal S}_{\Gamma_{ab}}\subset {\cal A}_{k,\delta}(O_a,O_b) \
\hbox{and} \   {\cal S}_{\Gamma_{bc}}\subset {\cal
A}_{k,\delta}(O_b,O_c), $$
of the cut off elements $(\rho^-_{1-T_k} X_k^{-T_k}, \rho^+_{-1+T_k}
X_k^{T_k})$.   
The elements $x_k$ and $y_k$ are no longer solutions of the
Seiberg--Witten equations. However, we can choose the sequence $T_k
\to \infty$ such that the elements 
$x_k$ and $y_k$ lie respectively within neighborhoods of radius
$\epsilon /4$ of the solutions $x$ and $y$.
We choose $\epsilon$ so that it satisfies 
$$ \epsilon \leq \min \{ \epsilon_{ac}, \epsilon_{ab},
\epsilon_{bc} \}. $$
Here we have $\epsilon_{ab}$ and $\epsilon_{bc}$ determined as in 
Lemma \ref{solveqbeta}, for the 
contraction principle for the fixed point problem
(\ref{contraction:map}) in ${\cal B}_{2,\delta}(O_a,O_b)$ and in
${\cal B}_{2,\delta}(O_b,O_c)$, respectively. Similarly,  
the constant $\epsilon_{ac}$ is the constant for the contraction
principle of Lemma \ref{solveqbeta} in ${\cal B}_{2,\delta}(
O_a,O_c)$. (cf. Proposition \ref{glue=fixedpoint} and Remark
\ref{fixedpoint}.) 
By applying the result of Proposition
\ref{glue=fixedpoint}, there are then unique elements $\tilde x_k$ and
$\tilde y_k$ in a $\epsilon /2$--neighborhood of $x_k$ and $y_k$,
respectively, 
which satisfy the equations. Notice that these elements $(\tilde x_k,
\tilde y_k)$ will in general have $e_b^- (\tilde x_k)\neq e_b^+
(\tilde y_k)$, hence they do not define an element in the fibered
product. However, we can estimate that the distance between the
endpoints goes to zero. In fact, we have
$$ \begin{array}{c}
\| e_b^- (\tilde x_k)-e_b^+(\tilde y_k) \|_{L^2_1} \leq \\[2mm]
C(dist_{L^2_1}(\lambda_k^- \tilde x_k (t), e_b^- (\tilde x_k)) +
dist_{L^2_1}(\lambda_k^+ \tilde y_k (t), e_b^+(\tilde y_k)))\\[2mm] 
+ C (\|  \tilde x_k - x \|_{L^2_{2,\delta}} +\|  \tilde y_k -y
\|_{L^2_{2,\delta}} )\\[2mm]
+ C(dist_{L^2_1} (\lambda^- x(t) , e_b^- (x)) + dist_{L^2_1}(\lambda^+ y
(t), e_b^+(y))). \end{array} $$
The first two terms and last two terms on the right hand side decay
exponentially like $C e^{-\delta |t|}$ by the results of 
Section 3.4, Theorem \ref{decay}. The gauge elements $\lambda_k^-$,
$\lambda_k^+$, $\lambda^-$, and $\lambda^+$ are as in
(\ref{slice:cond}), cf. Lemma \ref{slicelemma}. 
The remaining term is bounded by $\epsilon$.

Since in the statement of Theorem \ref{equivgluing} we are only
interested in the codimension one boundary, we may as well assume
that, upon passing to a subsequence, the sequence $(\tilde x_k, \tilde 
y_k)$ converges in the strong topology to a pair $(x',y')$. By the
result of the previous Claim, we know that this element $(x',y')$ is
in the fibered product ${\cal M}(O_a,O_b)\times_{O_b} {\cal
M}(O_b,O_c)$. 

Now consider the pre-glued solutions $x' \#_{T_k}^0 y'$. By our
construction, the original elements $X_k$ lie within neighborhoods of
radius $\epsilon$ of $x' \#_{T_k}^0 y'$ in the
$L^2_{2,\delta}$--norms. By  Lemma 
\ref{solveqbeta}, and Remark \ref{fixedpoint} there is a unique
zero of the map $f$ in this neighborhood, that is, a unique solution
of the Seiberg--Witten equations, obtained as a small deformation of
the approximate solution $x' \#_{T_k}^0 y'$. This implies the desired
equality $X_k= x' \#_{T_k} y'$.

This completes the proof of Theorem \ref{equivgluing}.

\noindent $\diamond$

Before discussing the fine structure of the compactification, we can
add a brief comment about the gluing result in the non-equivariant
setting, as stated in Lemma \ref{glue}. The difference in the
non-equivariant setting is the presence of the extra $U(1)$--gluing
parameter, in gluing across the reducible $\theta$. Namely, the gluing 
map is of the form
\[  \# :\hat{\cal M}(a,\theta)\times \hat{\cal M}(\theta,c)\times U(1)\times
[T_0,\infty) \to \hat{\cal M}(a,c). \]

The reason for the presence of the extra $U(1)$ gluing parameter is
the fact that the reducible point $\theta$ is not a smooth point in
the non-equivariant moduli space, hence, in order to formulate the
pre-gluing and gluing construction, it is necessary to lift the
pre-gluing to the framed moduli space, in the proximity of the
reducible point. More explicitly, we define the pre-gluing map as
{\tiny 
\beq \label{pre-glue:non-eq}
x\#_{T,u}^0 y=\left\{ \begin{array}{lr}
\lambda^+_1 (\AA_1^{2T},\Psi_1^{2T}) & t\leq -2 \\[2mm]
\lambda_0 (A_0+\nu,0) + \exp(iu (t+2))
(\alpha_1(t+2T),\phi_1(t+2T))  & -2\leq t\leq -1 \\[2mm]
\lambda_0 (A_0+\nu,0) +\exp(iu)  \cdot & -1\leq t\leq 1 \\
\left( \rho^-(t)(\alpha_1(t+2T),\phi_1(t+2T))+
\rho^+(t)(\alpha_2(t-2T),\phi_2(t-2T)) \right) & \\[2mm]
\lambda \lambda_2^- (\AA_2^{-2T},\Psi_2^{-2T}) & t\geq 1
\end{array}\right. \eeq  
}
with the extra gluing parameter $\lambda=\exp(iu)$ in $U(1)$.

This pre-gluing map gives a choice of a lift of
$$ \hat{\cal M}(a,\theta)\times \hat{\cal M}(\theta,c)\times U(1) $$
to $\hat{\cal M}(O_a,\theta)\times \hat{\cal M}(\theta,O_c)$, then the 
gluing construction works as in the equivariant case and composition
with the projection 
onto the quotient by the $U(1)$ action then provides the resulting
gluing map with values in $\hat{\cal M}(a,c)$.

\subsection{Multiple gluing theorem and corner structure}

We now generalize the result of Theorem \ref{equivgluing} of the
previous subsection to the case of multiple gluings of broken
trajectories in boundary strata of higher codimension. The purpose is
to identify the fine structure of the compactification, namely to show 
that the spaces $\hat{\cal M}(O_a,O_b)$ of flow lines compactify to a
${\cal C}^\infty$ manifold with corners, in the sense of
\cite{Melrose}. We shall follow the notation $\hat{\cal M}(O_a,O_b)^*$ 
to distinguish the compactification from the original space $\hat{\cal
M}(O_a,O_b)$. 

\begin{thm}
\label{corners}
The compactification $\hat{\cal M}(O_a,O_b)^*$ has the structure of a
smooth manifold with corners, with codimension $k$ boundary faces of
the form
\beq\label{k-strata} 
\bigcup_{c_1, \cdots c_k} \hat{\cal M}(O_a, O_{c_1})^* \times _{O_{c_1}}
\hat{\cal M}(O_{c_1}, O_{c_2})^*\times \cdots \times_{O_{c_k}}
\hat{\cal M}(O_{c_k}, O_b)^*.
\eeq
Here the union is over all possible sequences of critical orbits
$O_{c_1}, \cdots, O_{c_k}$ with decreasing indices.  
\end{thm}

We also have the following.

\begin{corol}
The endpoint maps $e_a^+$ and $e_b^-$ and their derivatives extend
continuously over the boundary and on the boundary they coincide with
$e_a^+$ and $e_b^-$ on $\hat{\cal M}(O_a, O_{c_1})$ and
$\hat{\cal M}(O_{c_k}, O_b)$ respectively. Thus, the maps $e_a^+$ and
$e_b^-$ are fibrations with compact fibers in the category of smooth
manifolds with corners.
\label{strata}
\end{corol}

\noindent\underline{Proof of Theorem \ref{corners}:}
We proceed as follows. First we prove that $\hat{\cal M}(O_a,O_b)^*$
is a $t$-manifold (has a ${\cal C}^\infty$ structure with corners) in
the sense of Definition 1.6.1 of \cite{Melrose}. This amounts to
showing that the multiple gluing maps define a compatible set of
charts in the sense of Section 1.6 of \cite{Melrose} on open sets in
$\hat{\cal M}(O_a,O_b)$ near the boundary.
We then show that this $t$-manifold has the structure of a smooth
manifold with corners, in the sense of Section 1.8 of \cite{Melrose},
by showing that the boundary faces satisfy the condition (1.8.7) of
\cite{Melrose}.  

Consider broken trajectories $(x_0, x_1, \ldots, x_k)$
in a compact subset
\[ \hat K \subset
 \hat{\cal M}(O_a, O_{c_1}) \times _{O_{c_1}}\hat{\cal M}(O_{c_1},
O_{c_2})\times \cdots \times_{O_{c_k}} 
\hat{\cal M}(O_{c_k}, O_b),
\]
where $O_{c_1}, \ldots, O_{c_k}$ are the critical points with 
decreasing indices 
\[
\mu(O_a)>\mu(O_{c_1}) > \mu(O_{c_2}) > \cdots > \mu(O_{c_k}) > \mu(O_b).
\]

We introduce multiple gluing maps
$$ 
\hat\#_{i_1,\ldots,i_k} :  \hat K \times [T_0,\infty)^k 
\to \hat{\cal M}(O_a,O_b). 
$$
Here $( T_1,\ldots, T_k )\in [T_0,\infty)^k$ are gluing parameters,
and the indices $\{ i_1,\ldots,i_k \}$
are a permutation of the set $\{ 1,\ldots , k\}$ which specifies in
which order the multiple gluing is performed as a sequence of $k$
gluings as in Theorem \ref{equivgluing}. For instance, we have
$$ \hat\#^{T_1,\ldots, T_k}_{1,2,\ldots,k} (x_0, x_1, \ldots, x_k)=(\cdots 
(((x_0\hat\#_{T_1} x_1)\hat\#^{T_2} x_2)\hat\#_{T_3}x_4) \cdots
)\hat\#_{T_k}x_k, $$
or
$$ \hat\#^{T_1,\ldots, T_k}_{2,1,3,\ldots,k}(x_0, x_1, \ldots,
x_k)=(\cdots(x_0\hat\#_{T_1}(x_1\hat\#^{T_2} x_2))\hat\#_{T_3}x_3)
\cdots)\hat\#_{T_k}x_k, $$
etc.

The gluing construction is non-canonical, in the sense that the
identification of the normal bundle ${\cal T}^\perp$ in the pre-gluing
construction, as in Lemma \ref{surj2}, is 
dependent on the order of the gluing, hence there is no obvious
associativity law for multiple gluings. However, we are going to show
that the maps $\hat\#^{T_1,\ldots, T_k}_{i_1,\ldots,i_k}$, or rather their 
inverses, to be consistent with the notation of (1.6.1) of
\cite{Melrose}, define a system of charts of a ${\cal C}^\infty$
structure with corners, or $t$-manifold.
We prove the following result, which ensures that the
$\hat\#^{T_1,\ldots, T_k}_{i_1,\ldots,i_k}$ are local
diffeomorphisms as needed.

\begin{prop}
Suppose given $O_a$, $O_{c_1}, \cdots O_{c_k}$
 and $O_b$ in ${\cal M}^0$  with decreasing indices. 
Then, given a compact set $\hat K$ in
\[
\hat{\cal M}(O_a, O_{c_1}) \times _{O_{c_1}}
\hat{\cal M}(O_{c_1}, O_{c_2})\times \cdots \times_{O_{c_k}}
\hat{\cal M}(O_{c_k}, O_b),
\]
there is a lower bound $T_0(K)>0$ such that the gluing maps
$\hat\#^{T_1,\ldots, T_k}_{i_1,\ldots,i_k}$ define local
diffeomorphisms \index{$\hat \#$}
\[
\hat\#_{i_1,\ldots,i_k} : \qquad
\hat K \times [T_0(K), \infty)^k \longrightarrow \hat{\cal M}(O_a, O_b).
\]
where the orientation on the left hand side is the one induced by the
product orientation on
\[
 \Bigl(\bigl((\hat{\cal M}(O_a, O_{c_1}) \times \RR \times
\hat{\cal M}(O_{c_1}, O_{c_2}))\times \RR\bigr)
\times \cdots \times \RR \Bigr)  \times
\hat{\cal M}(O_{c_k}, O_b).
\]
\label{equivgluing:k}
\end{prop}

\noindent\underline{Proof:}
We want to generalize the analogous statement proved in Theorem
\ref{equivgluing} for the case of one intermediate critical point. 

First, we need to show the analogue of the error estimate
(\ref{error:estimate}), which shows that the glued solutions 
$\#^{T_1, \ldots, T_k}_{i_1,\ldots, i_k} (x_0,\ldots,x_k)$
converge to broken trajectories in the limit when a certain subset $\{
T_{i_1}, \ldots T_{i_m} \}$ of the gluing parameters satisfies
$T_{i_\ell}\to \infty$. 

We then need to show that a sequence of solutions in $\hat {\cal
M}(O_a,O_b)$ that converges to a broken trajectory $(x_0,\ldots, x_k)$ 
in the boundary
$$ \hat{\cal M}(O_a,O_{c_1})\times_{O_{c_1}} \cdots
\times_{O_{c_k}}\hat{\cal M}(O_{c_k}, O_b) $$
lies eventually in the range of a gluing map $\#^{T_1, \ldots,
T_k}_{i_1,\ldots, i_k}$. We first prove that, when the gluing
parameters $\{ T_1, \ldots,T_k \}$ satisfy $T_i \to \infty$, for all
$i=1,\ldots, k$, we have an estimate \index{$\sigma(x,y,T)$}
\beq \| \sigma_k(x_0,\ldots,x_k,T_1, \ldots,T_k) \|_{L^2_{1,\delta}}
\leq C e^{-\delta T}, \label{gluedecay:k} \eeq
with $T=\min \{ T_1, \ldots,T_k \}$. This is the analogue of the
estimate (\ref{gluedecay}) for multiple gluing. Here $\sigma_k$ is the
error term
$$ \sigma_k(x_0,\ldots,x_k,T_1, \ldots,T_k) = \#^{T_1,
\ldots,T_k}_{i_1,\ldots, i_k}(x_0,\ldots,x_k) - {\#^0}^{T_1,
\ldots,T_k}_{i_1,\ldots, i_k}(x_0,\ldots,x_k), $$
with ${\#^0}^{T_1,\ldots,T_k}_{i_1,\ldots, i_k}$ the pre-gluing map,
namely the composite of $k$ pre-gluing maps as in (\ref{pre:map}), in
the order specified by $( i_1,\ldots, i_k )$.

We prove the estimate (\ref{gluedecay:k}) for the case of the gluing map
$\#^{T_1, \ldots, T_k}_{1,\ldots, k}$. The same procedure works with
a permutation of the order of gluing. 
We proceed by induction on $k$. The case $k=1$ follows from Theorem
\ref{equivgluing}. Assume as induction hypothesis that the gluing map 
$\#^{T_1, \ldots, T_{k-1}}_{1,\ldots, k-1}$ can be written as
$$ \#^{T_1, \ldots, T_{k-1}}_{1,\ldots, k-1}(x_0,\ldots, x_{k-1})= $$
$$ {\#^0}^{T_1, \ldots, T_{k-1}}_{1,\ldots, k-1}(x_0,\ldots, x_{k-1}) +
\sigma_{k-1}(x_0,\ldots, x_{k-1},T_1, \ldots, T_{k-1}), $$
with ${\#^0}^{T_1, \ldots, T_{k-1}}_{1,\ldots, k-1}$ the pre-gluing
map (i.e. the composition of $k-1$ pre-gluing maps (\ref{pre:map}), and
with the error term $\sigma_{k-1}$ satisfying the estimate
$$ \| \sigma_{k-1}(x_0,\ldots, x_{k-1},T_1, \ldots, T_{k-1})
\|_{L^2_{1,\delta}} \leq C
e^{-\delta T}, $$
with $T=\min \{ T_1, \ldots,T_{k_1} \}$.
We can then form the pre-glued solution $y\#^0_{T_k} x_k$, with
$$ y=\#^{T_1, \ldots, T_{k-1}}_{1,\ldots, k-1}(x_0,\ldots,
x_{k-1}). $$
By Proposition \ref{glue=fixedpoint} this approximate solution can be
deformed to an actual solution $y\#_{T_k} x_k$, which by definition is 
the same as
$$ y\#_{T_k} x_k=\#^{T_1, \ldots, T_k}_{1,\ldots, k}(x_0,\ldots,
x_k). $$
We need to estimate the difference $y\#_{T_k} x_k - y_0 \#^0_{T_k}
x_k$, where $y_0$ is the approximate solution
$$ y_0 ={\#^0}^{T_1, \ldots, T_{k-1}}_{1,\ldots, k-1}(x_0,\ldots,
x_{k-1}). $$
We have
\beq \begin{array}{c}
 \| y\#_{T_k} x_k - y_0 \#^0_{T_k}x_k \|_{L^2_{1,\delta}}
\leq \| y\#_{T_k} x_k - y\#^0_{T_k}x_k \|_{L^2_{1,\delta}}
\\[2mm]
 + \| y\#^0_{T_k}x_k -y_0 \#^0_{T_k}x_k
\|_{L^2_{1,\delta}}. \end{array} \label{estim:induction} \eeq  
Now we can estimate the second term on the right hand side by
\[ \| y\#^0_{T_k}x_k -y_0 \#^0_{T_k}x_k
\|_{L^2_{1,\delta}} \leq \| F_{\#_{T_k}} \| \cdot \| y- y_0
\|_{L^2_{1,\delta}}. \]
We know that (\ref{lin:pre-glue}) has norm bounded by $C
e^{-\delta T_k}$. Moreover, the quantity $\| y- y_0 \|$ is
given by 
$$ \| \sigma_{k-1}(x_0,\ldots,x_{k-1},T_1,\ldots,T_{k-1})
\|_{L^2_{1,\delta}}. $$
Using the induction hypothesis, we obtain
$$ \| y\#^0_{T_k}x_k -y_0 \#^0_{T_k}x_k \|_{L^2_{1,\delta}} \leq C
e^{-\delta T}. $$
The first term in the right hand side of (\ref{estim:induction}), on
the other hand, is given by 
$$ \| y\#_{T_k} x_k - y\#^0_{T_k}x_k \|_{L^2_{1,\delta}}
= \| \sigma(y,x_k, T_k) \|_{L^2_{1,\delta}}. $$
By the estimate (\ref{gluedecay}) (cf. Proposition \ref{glue=fixedpoint},
Lemma \ref{solveqbeta}, Remark \ref{fixedpoint}, and Proposition
\ref{glue=fixedpoint}), this satisfies the estimate
$$ \| \sigma(y,x_k, T_k) \|_{L^2_{1,\delta}} \leq C
e^{-\delta T_k}. $$

Thus, we have the next step of the induction, namely we have obtained 
$$ \| \sigma_k(x_0,\ldots, x_k, T_1,\ldots, T_k)
\|_{L^2_{1,\delta}} 
\leq C e^{-\delta \cdot \min \{ T, T_k \} }, $$
for 
$$ \sigma_k(x_0,\ldots, x_k, T_1,\ldots, T_k)=y\#_{T_k} x_k - y_0
\#^0_{T_k}x_k.  $$

Now consider a sequence $X_\ell$ in
${\cal M}(O_a,O_b)$ that converges smoothly on compact sets, in the
sense of Theorem \ref{equivgluing}, to a broken trajectory
$(x_0,\ldots, x_k)$ in the fibered product
$$ {\cal M}(O_a,O_{c_1})\times_{O_{c_1}} \cdots \times_{O_{c_k}}{\cal
M}(O_{c_k},O_b). $$
We need to show that these $X_\ell$ 
belong eventually to the range of at least one of the multiple gluing maps
$\#_{i_1,\cdots,i_k}$. We proceed as in the analogous argument in
Theorem \ref{equivgluing}, for a single gluing. Convergence in the
sense of Theorem \ref{equivgluing} implies that there exist time
shifts $\{ T_{1,\ell}, \ldots, T_{k,\ell} \}$ with $T_{i,\ell}\to
\infty$ for $i=1,\ldots k$, and combinatorial data which describe the
order in which successive limits on compact sets and
reparameterizations are considered in Theorem \ref{equivgluing}. These 
data determine which gluing maps $\#_{i_1,\cdots,i_k}$ have the points
$X_\ell$ in the range. These combinatorial data can be described as
the set of binary trees with $k+1$ leaves, one root, and $k-1$
intermediate nodes. To one such tree we associate a cutoff element
\beq (\eta_0^{T_{1,\ell},\ldots, T_{k,\ell}} (X_\ell),\ldots,
\eta_\ell^{T_{1,\ell},\ldots, T_{k,\ell}}(X_\ell))
\label{multiple:cutoff} \eeq 
where the $\eta_i^{T_{1,\ell},\ldots, T_{k,\ell}}$ act as a composition of 
time shifts and cutoff functions, as specified by the combinatorics of 
the tree. 
\begin{figure}[ht]
\epsfig{file=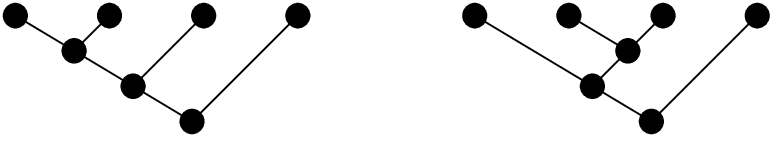,angle=270}
\end{figure}
For instance, the two trees of the figure correspond respectively to the
elements 
$$ \begin{array}{l} (\rho^-_{1-T_{1,\ell}}(\rho^-_{1-T_{2,\ell}}
(\rho^-_{1-T_{3,\ell}}X_\ell^{-T_{3,\ell}} )^{-T_{2,\ell}} )^{
-T_{1,\ell}}, \\[2mm]
\rho^+_{-1+T_{1,\ell}} (\rho^-_{1-T_{2,\ell}}
(\rho^-_{1-T_{3,\ell}}X_\ell^{-T_{3,\ell}})^{-T_{2,\ell}}
)^{T_{1,\ell}}, \\[2mm] 
\rho^+_{-1+T_{2,\ell}}(\rho^-_{1-T_{3,\ell}} X_\ell^{-T_{3,\ell}}
)^{T_{2,\ell}}, \\[2mm] 
\rho^+_{-1+T_{3,\ell}} X_\ell^{T_{3,\ell}} ) \end{array} $$  
and
$$ \begin{array}{l}
(\rho^-_{1-T_{1,\ell}} (\rho^-_{1-T_{3,\ell}} X_\ell^{-T_{3,\ell}}
)^{-T_{1,\ell}}, \\[2mm] 
\rho^-_{1-T_{1,\ell}} (\rho^+_{-1+T_{2,\ell}} (\rho^-_{1-T_{3,\ell}}
X_\ell^{-T_{3,\ell}} )^{T_{2,\ell}} )^{-T_{1,\ell}}, \\[2mm]
 \rho^+_{-1+T_{1,\ell}} (\rho^+_{-1+T_{2,\ell}} (\rho^-_{1-T_{3,\ell}} 
X_\ell^{-T_{3,\ell}} )^{T_{2,\ell}} )^{T_{1,\ell}}, \\[2mm]
\rho^+_{-1+T_{3,\ell}} X_\ell^{T_{3,\ell}} ). \end{array} $$ 
These correspond to cases where a sequence of solutions $X_\ell$ in,
say, ${\cal M}(O_a,O_b)$ converges smoothly on compact sets,
after the different reparameterizations specified above, to elements
in ${\cal M}(O_a,O_{c_1})$, ${\cal M}(O_{c_1},O_{c_2})$, ${\cal
M}(O_{c_2},O_{c_3})$, and ${\cal M}(O_{c_3},O_b)$, respectively.

We denote by $(x_{1,\ell},\ldots,x_{k,\ell})$ the projection onto the
slices of the cutoff elements (\ref{multiple:cutoff}). The elements
$x_{i,\ell}$ are no longer solutions of the Seiberg--Witten equations, 
however, by hypothesis, for $\ell$ large enough, they are contained in
$\epsilon/2^{k+1}$-neighborhoods of the solutions $x_i$, for
$i=0,\ldots,k$. Pick $\epsilon$ satisfying 
$$ \epsilon\leq \max_{0\leq j \leq k+1} \{ \epsilon_{c_j c_{j+1}} \}, $$
where the $\epsilon_{c_j c_{j+1}}$ is the constant for the
contraction principle, as in Lemma \ref{solveqbeta} in ${\cal
B}_{2,\delta}(O_{c_j},O_{c_{j+1}})$, and we take $c_0=a$ and
$c_{k+1}=b$. Then, proceeding as in the proof of Theorem
\ref{equivgluing}, we find elements $\tilde x_{i,\ell}$, for
$i=0,\ldots k$ that satisfy the Seiberg--Witten equations and are
contained in $\epsilon/2^k$-neighborhoods of the elements $x_i$. The
elements $(\tilde x_{0,\ell},\ldots,\tilde x_{k,\ell})$ are not
necessarily in the fibered product, however, by the same argument used 
in the proof of Theorem \ref{equivgluing}, we obtain convergence of
the elements $(\tilde x_{0,\ell},\ldots,\tilde x_{k,\ell})$ to a limit 
$(x_0',\ldots,x_k')$ in the fibered product
$$ {\cal M}(O_a,O_{c_1})\times_{O_{c_1}} \cdots \times_{O_{c_k}}{\cal
M}(O_{c_k},O_b). $$
By construction, the approximate solution
$$ {\#^0}^{T_{1,\ell},\ldots, T_{k,\ell}}_{i_1,\ldots i_k}
(x_0',\ldots x_k'), $$
with the order $(i_1,\ldots i_k)$ of gluing specified by the
combinatorial data, is contained in an $\epsilon$-neighborhood of
$X_\ell$ in the $L^2_{2,\delta}$-norm. By the contraction argument of
Lemma \ref{solveqbeta} 
together with the first part of the proof of this Proposition, we know 
that we must then have
$$ X_\ell =\#^{T_{1,\ell},\ldots, T_{k,\ell}}_{i_1,\ldots
i_k}(x_0',\ldots x_k'). $$  
This completes the proof. 

\noindent $\diamond$

Thus, the changes of coordinates
$$ \hat\#^{T_1,\ldots, T_k}_{i_1,\ldots,i_k}\circ
(\hat\#^{\tau_1,\ldots, \tau_\ell}_{j_1,\ldots,j_\ell} )^{-1} $$
are local diffeomorphisms between open subsets of 
$$ \RR^n_k =\RR^{n-k}\times [T_0, \infty)^k $$
and
$$ \RR^n_\ell=\RR^{n-\ell}\times [T_0, \infty)^\ell, $$
with $n=\dim \hat{\cal M}(O_a,O_b)$, as prescribed in Section 1.6 of
\cite{Melrose}. 

Notice, in particular, that given two multiple gluing maps
$\#^{T_1,\ldots, T_k}_{i_1,\ldots, i_k}$ and $\#^{T_1,\ldots,
T_k}_{j_1,\ldots j_k}$, for two different orders of gluing
$(i_1,\ldots, i_k)\neq (j_1,\ldots j_k)$, we obtain an estimate
\beq \| \#^{T_1,\ldots, T_k}_{i_1,\ldots, i_k}(x_0,\ldots, x_k) -
\#^{T_1,\ldots, T_k}_{j_1,\ldots j_k}(x_0,\ldots, x_k)
\|_{L^2_{1,\delta}} \leq C e^{-\delta T}, \label{closecharts} \eeq
for $T=\min \{ T_1,\ldots, T_k \}$. The proof of this estimate is
completely analogous to the argument used in the proof of Proposition
\ref{equivgluing:k}. For instance, we can estimate
\[ \| (x_0\#_{T_1} x_1)\#_{T_2} x_2 - x_0\#_{T_1}(x_1\#_{T_2} x_2) \|
\leq \]
\[  \| (x_0\#_{T_1} x_1)\#_{T_2} x_2 - (x_0\#_{T_1} x_1)\#^0_{T_2} x_2
\| + \| (x_0\#_{T_1} x_1)\#^0_{T_2} x_2 - x_0\#^0_{T_1} x_1\#^0_{T_2}
x_2 \| \]
\[ + \| x_0\#^0_{T_1} x_1\#^0_{T_2} x_2 - x_0\#^0_{T_1}(x_1\#_{T_2}
x_2) \| + \| x_0\#^0_{T_1}(x_1\#_{T_2}x_2) - x_0\#_{T_1}(x_1\#_{T_2}
x_2) \|  \]
\[ \leq C e^{-\delta T_2} + C e^{-\delta T_1 + T_2} + C e^{-\delta T_1
+ T_2} + C e^{-\delta T_1}. \]
The estimate (\ref{closecharts}) implies that, when all the $T_i$ go
to infinity, the glued elements $\#^{T_1,\ldots, T_k}_{i_1,\ldots,
i_k}(x_0,\ldots, x_k)$ and $\#^{T_1,\ldots, T_k}_{j_1,\ldots
j_k}(x_0,\ldots, x_k)$ end up in the same coordinate patch. This may
not be the case if only some of the $T_i$ tend to infinity and other
remain bounded.

We have shown, as a result of Proposition
\ref{equivgluing:k}, and by Lemma 1.7.1 of \cite{Melrose}, that the strata
$$ \hat{\cal M}(O_a, O_{c_1})^* \times_{O_{c_1}}\times \cdots
\times_{O_{c_k}} \hat{\cal M}(O_{c_k}, O_b)^* $$
inherit the structure of smooth $t$-sub-manifolds of codimension $k$ in 
$\hat{\cal M}(O_a, O_b)^*$.

The codimension one boundary strata described in Theorem 4.8,
compactified to smooth $t$-manifolds, define the boundary
hypersurfaces of $\hat{\cal M}(O_a, O_b)^*$, as in (1.8.3) and (1.8.4) 
of \cite{Melrose}.

Let us recall that a ${\cal C}^\infty$ manifold with
corners is a $t$-manifold where all the boundary faces are 
${\cal C}^\infty$ embedded sub-manifolds. According to Section 1.8 of
\cite{Melrose}, to ensure that this is the case, namely that the
compactification $\hat{\cal M}(O_a, 
O_b)^*$ has the structure of a ${\cal C}^\infty$ manifold with
corners, one only needs to check the following fact.

\begin{lem}
\label{corners:structure}
Every compactified codimension $k$ boundary face
$$ \hat{\cal M}(O_a, O_{c_1})^* \times _{O_{c_1}}\times \cdots
\times_{O_{c_k}} \hat{\cal M}(O_{c_k}, O_b)^* $$
of $\hat{\cal M}(O_a, O_b)^*$ is a component of (precisely) one
intersection of $k$ boundary hypersurfaces 
$$ H_{i_1} \cap \ldots \cap H_{i_k}. $$ 
\end{lem}

\noindent\underline{Proof of Lemma \ref{corners:structure}:} All the
boundary hypersurfaces are identified by Theorem
\ref{equivgluing}. Thus, we see that the Lemma is verified by
setting
$$ \begin{array}{lcl}
H_{i_1}& = & \hat{\cal M}(O_a,O_{c_1})^* \times_{O_{c_1}} \hat{\cal
M}(O_{c_1}, O_b)^* \\[2mm]
H_{i_2}& = & \hat{\cal M}(O_a,O_{c_2})^*\times_{O_{c_2}} \hat{\cal
M}(O_{c_2}, O_b)^* \\ 
\vdots & & \vdots \\
H_{i_{k-1}}& = & \hat{\cal M}(O_a,O_{c_{k-1}})^*\times_{O_{c_{k-1}}}
\hat{\cal 
M}(O_{c_{k-1}}, O_b)^* \\[2mm]
H_{i_k}& = & \hat{\cal M}(O_a,O_{c_k})^* \times_{O_{c_k}} \hat{\cal
M}(O_{c_k}, O_b)^* 
\end{array} $$

\noindent $\diamond$

This completes the proof of Theorem \ref{corners}.

\noindent $\diamond$

Corollary \ref{strata} now follows, by considering the restriction of
the asymptotic value maps $e_a^+$ and $e_b^-$ to the range of the
gluing maps. By the convergence property, when a subset of the gluing
parameters goes to infinity, we obtain that the asymptotic value maps
restrict to the corresponding maps on the boundary strata.

\section{Equivariant Homology}

Let us recall briefly the construction of the de Rham model for
$U(1)$-equivariant cohomology (and homology) on a finite dimensional
manifold with a $U(1)$ action. The main reference is \cite{AtBo}.

Let $W$ be the Weil algebra of the Lie algebra $i\RR$ of $U(1)$. This
is a free commutative graded algebra in one generator $\theta$ of
degree 1 and one generator $\Omega$ of degree 2, with differential
$\delta$ that satisfies 
\[ \begin{array}{cc}\delta\theta=\Omega&\delta\Omega=0.
   \end{array} \]

Let $M$ be a manifold with a $U(1)$ action. Consider the complex 
\[ C^* = W\otimes \Omega^*(M), \]
with differential
\[ d_{U(1)}=d-\Omega c(T). \]
Here $T$ is the vector field on $M$ generated by the infinitesimal
$U(1)$ action and $c$ is the unique derivation in $W$ that satisfies 
\[ \begin{array}{cc}c(\theta)=1&c(\Omega)=0.
   \end{array} \]

We choose a sub-complex of $C^*$ by taking the cochains on which
$c+c(T)$ and ${\cal L}(T)$ vanish, where ${\cal L}(T)$ is the Lie
derivative. Let us call this sub-complex $\Omega^*_{U(1)}(M)$.
An alternative description of the complex $\Omega^*_{U(1)}(M)$ is
\[ \Omega^*_{U(1)}(M)=\RR [\Omega]\otimes \Omega_0^*(M), \]
where $\Omega_0^*(M)$ are de Rham forms that are annihilated by the Lie
differentiation ${\cal L}(T)$.
The cohomology $H^*(\Omega^*_{U(1)}(M), d_{U(1)})$ is isomorphic to
the equivariant cohomology with real coefficients,
\[ H^*(\Omega^*_{U(1)}(M), d_{U(1)})\cong H^*_{U(1)}(M,\RR). \]

In order to compute equivariant homology with real coefficients a
de Rham complex can be constructed as in \cite{AB3} by considering
differential forms graded by the dimension of $M$ minus the degree, 
\[ {\Omega_*}_{U(1)}(M)=\RR [\Omega]\otimes \Omega_0^{\dim(M)-*}(M), \]
with boundary operator $\partial_{U(1)}=d-c(\theta) \otimes c(T)$. 
\index{$\partial_{U(1)}$}  
Austin and Braam \cite{AB3} proved that 
this complex computes the same homology as the equivariant complex of
currents introduced by Duflo and Vergne \cite{DV} and studied by Kumar
and Vergne \cite{KV}. 
With this understood, we can consider the complex
\index{${\Omega_*}_{U(1)}(O_a)$} 
${\Omega_*}_{U(1)}(O_a)$ associated to each critical orbit $O_a$ in ${\cal
  M}_0$. This will be a copy of the polynomial algebra
$\RR [\Omega]$ for the fixed point and a complex of $\RR [\Omega]$-modules 
with de Rham forms in degree zero and one in the case of orbits that 
come from irreducibles.

We can give a more explicit description of the boundary operator
$\partial_{U(1)}$ in this case. Let $O_a$ be a critical orbit with a
free $U(1)$ action. Then, the generators of $\Omega_0^{1-*}(O_a)$ are
a 1-form $\eta_a$ that generates $H^1(S^1)$, in degree zero, and a
zero-form $1_a$ (the constant function equal to one on $O_a$), in
degree one. Thus we have
$$ \partial_{U(1)}(\Omega^n\otimes 1_a)=0 $$
and
$$ \partial_{U(1)}(\Omega^n\otimes \eta_a)=-\Omega^{n-1} \otimes
1_a. $$ 

\subsection{The equivariant complex}

We can form the bigraded complex that computes equivariant Floer
homology as in \cite{AB}, defined by
\beq
{C_k}_{U(1)}(Y)=\bigoplus_{\mu(O_a)=i,i+j=k}{\Omega_j}_{U(1)}(O_a), 
\label{homology}
\eeq
with differentials \index{$D$}
\beq 
D_{a,b} \eta =\left\{\begin{array}{lr}
\partial_{U(1)}\eta&O_a=O_b\\
(-1)^{1-r(\eta)}(e_b^-)_*(e_a^+)^*\eta& \mu(O_a)>\mu(O_b)\\
0&otherwise.\end{array}\right.
\label{boundary2}
\eeq

Here $\eta$ is an equivariant differential form on the orbit $O_a$,
that is, an element of ${\Omega_*}_{U(1)}(O_a)$. The number $r(\eta)$
is the de Rham degree of $\eta$, that is, the maximum degree of the
elements of $\Omega_0^*(O_a)$ that appear in the expression of
$\eta$. Recall that ${\Omega_*}_{U(1)}(O_a)$ is the dual of the de Rham
complex, and the forms are graded by $\dim(O_a)-*=1-*$.

The analogous complex that computes equivariant Floer cohomology is
given by
\beq
{C^k}_{U(1)}(Y)=\bigoplus_{\mu(O_a)=i,i+j=k}{\Omega^j}_{U(1)}(O_a),
\label{cohomology}
\eeq
with coboundaries 
\beq 
\delta_{a,b} \eta =\left\{\begin{array}{lr}
d_{U(1)}\eta&O_a=O_b\\
(-1)^{r(\eta)}(e_a^+)_*(e_b^-)^*\eta& \mu(O_a)>\mu(O_b)\\
0&otherwise.\end{array}\right.
\label{coboundary}
\eeq
Here $\eta$ is an equivariant form, that is an element of $O_b$ and
$r(\eta)$ is the de Rham degree of $\eta$, i.e. the maximum degree
of the elements of $\Omega_0^*(O_b)$ that appear in the expression of
$\eta$. 

The boundary map is well defined, since the endpoint maps are
compatible with the boundary strata as stated in \ref{strata}. 
Notice that, as pointed out in \cite{AB}, for dimensional reasons the
maps $D_{a,b}$ are trivial whenever $\mu(O_a)\geq \mu(O_b)+3$. A more
explicit description of the boundary map will be given in the following.

\subsection{Equivariant Floer Homology}

We can prove that the composites $D^2$ and $\delta^2$ are zero. 
Moreover we can see that there is a duality at the level of forms that
induces a duality between equivariant Floer homology and cohomology as in
\cite{AB} and \cite{AB3}.

\begin{thm}\label{DD=0}
The composite maps $D^2$ and $\delta^2$ are zero, this means
that the identities 
\[ D^2_{a,c}=\sum_b D_{a,b}D_{b,c}=0, \]
\[ \delta^2_{a,c}=\sum_b \delta_{b,c}\delta_{a,b}=0, \]
hold, with $O_b$ that ranges among critical orbits satisfying 
$\mu(O_a)\geq\mu(O_b)\geq\mu(O_c)$.
Moreover there is a pairing $<,>$ of forms in ${C_*}_{U(1)}$ and in 
${C^*}_{U(1)}$ that satisfies
\[ <D \eta,\gamma>=<\eta,\delta \gamma>. \]
\end{thm}

\noindent\underline{Proof:} The statement is true for $O_a=O_c$.
Given critical orbits with $\mu(O_a)>\mu(O_c)$, we have the expression
\[ D^2_{a,c}\eta =(-1)^{r(\eta)}\sum_{\{ O_b
  |\mu(O_a)>\mu(O_b)>\mu(O_c)\} } (-1)^{ r((e_b^-)_*(e_a^+)^* \eta)} \]
\[ (e_c^-)_*(e^+_b)^*(e_b^-)_*(e_a^+)^* \eta + 
 - (e_c^-)_*(e_a^+)^* \partial_{U(1)}\eta +\partial_{U(1)}
(e_c^-)_*(e_a^+)^*\eta. \]
This expression vanishes since the 
co-dimension 1 boundary of $\hat{\cal M}(O_a,O_b)^*$ is the union 
of the components
$\hat{\cal M}(O_a,O_b)^*\times_{O_b}\hat{\cal M}(O_b,O_c)^*$, as in
\ref{equivgluing}. Thus, if we use Stokes theorem applied to the fibration
with boundary $e^-_c: \hat{\cal M}(O_a,O_c)^*\to O_c$, we obtain
\cite{AB} that
\[ \partial_{U(1)}(e_c^-)_*\gamma=(e_c^-)_*\partial_{U(1)}\gamma -
(-1)^{r(\gamma)-\dim F}(e_c^-)_{\partial,*}\gamma. \]
Here $\gamma$ is an equivariant form on $\hat{\cal M}(O_a,O_c)^*$ and
the map $(e_c^-)_{\partial,*}$ is the push-forward map induced by the
restriction of the bundle to the boundary of each fiber $F$. If we choose
$\gamma$ of the form $\gamma=(e_a^+)^*\eta$ where $\eta$ is an
equivariant form on the orbit $O_a$, then the map
$(e_c^-)_{\partial,*}$ can be written as
\[ (e_c^-)_{\partial,*}\gamma=\sum_b
(e_c^-)_*(e_b^+)^*(e_b^-)_*(e_a^+)^* \eta, \]
for all the $b$ that satisfy $\mu(O_a)>\mu(O_b)>\mu(O_c)$, as shown in
the diagram.  Notice that the sign
$(-1)^{r((e_b^-)_*(e_a^+)^* \eta)}$ is exactly the sign
$(-1)^{r((e_a^+)^*\eta)-\dim F}$.

$$
\spreaddiagramrows{-1pc}
\spreaddiagramcolumns{-1pc}
\diagram
&  & \hat{\cal M}(O_a,O_c)^* \xto[2,-2] \xto[2,2]  \\
&  & \hat{\cal M}(O_a,O_b)^*\times_{O_b}\hat{\cal M}(O_b,O_c)^* \dlto
\drto \uto \\ 
O_a   &  \hat{\cal M}(O_a,O_b)^* \lto \rto & O_b  & \hat{\cal
M}(O_b,O_c)^* \lto \rto & O_c \\
\enddiagram
$$ 

The statement about the pairing follows from the identity
\[ \int_{O_b}(e_b^-)_*(e^+_a)^*\eta\wedge\gamma= \int_{\hat{\cal
    M}(O_a,O_b)^*} (e^+_a)^*\eta\wedge (e^-_b)^*\gamma. \]
 
\noindent $\diamond$

\begin{defin}
We define the equivariant Floer homology and cohomology to be 
\[ {SWF_*}_{U(1)}=H_*({C_*}_{U(1)},D) \]
and \index{${SWF_*}_{U(1)}$}
\[ SWF^*_{U(1)}=H^*(C^*_{U(1)},\delta). \]
\label{floergroups}
\end{defin}

We now give a more detailed description of the equivariant Floer
complex, which will be useful in Section 6 and Section 7.

\begin{prop}
Let $\Omega^n \otimes
\eta_a$ and $\Omega^n\otimes 1_a$ be the generators of the equivariant 
Floer complex.
The only possibly non-trivial coefficients of the boundary map $D$ are
those of the form
\[ \langle  \Omega^{n-1}\otimes 1_a, D(\Omega^n\otimes \eta_a)
\rangle, \]
or of the form
\[ \langle \Omega^n\otimes 1_b, D (\Omega^n\otimes 1_a) \rangle, \]
\[ \langle \Omega^n\otimes \eta_b, D (\Omega^n\otimes \eta_a)
\rangle, \] 
when $\mu(O_a)-\mu(O_b)=1$, and 
\[ \langle \Omega^n\otimes 1_c, D (\Omega^n\otimes \eta_a)
\rangle, \]  
when $\mu(O_a)-\mu(O_c)=2$. In the case of the critical orbit $\theta$,
we obtain the boundary component 
\[ \langle \Omega^n\otimes 1_c, D(\Omega^n\otimes \theta) \rangle \]
when $\mu(\theta)-\mu(O_c)=2$, where we consider the one-dimensional
space $\hat{\cal M}(\theta, O_c)$ fibering over $O_c$ with
zero-dimensional fiber. We also have the component
\[ \langle \Omega^n\otimes \theta, D(\Omega^n\otimes \eta_a) \rangle,
\]
when $\mu(O_a)-\mu(\theta)=1$, where the space $\hat{\cal
M}(O_a,\theta)$ is one-dimensional. In this case the coefficient is
obtained by integrating the 1-form $\eta_a$ over the 1-dimensional
fiber of the endpoint map $e_\theta: \hat{\cal M}(O_a,\theta)\to \theta$. 
\label{boundaryterms}
\end{prop}

\noindent\underline{Proof:} 
For a given orbit $O_a$ the complex of equivariant forms is given by 
$$ \Omega_{*, U(1)}(O_a), $$
with the total grading 
$$\Omega_{j, U(1)}(O_a)=\bigoplus_{2k +l = j}\RR \Omega^k \otimes
\Omega_0^{1-l}(O_a).$$  
Here $\Omega$ is of degree 2 and $l$ is the grading in the dual
de Rham complex where forms are graded by $\dim(O_a)-*$.  

Thus we have
\[ \Omega_{0, U(1)}(O_a) = \RR < 1\otimes \eta_a >, \]
\[ \Omega_{1, U(1)}(O_a) = \RR < 1\otimes 1_a >, \]
\[  \Omega_{2, U(1)}(O_a) = \RR < \Omega\otimes \eta_a >, \]
\[ \Omega_{3, U(1)}(O_a) = \RR < \Omega\otimes 1_a >, \]
and so on, where $\eta_a$ is the one form that generates $H^1(S^1)$ and $1_a$
is the constant function equal to 1 on $O_a$.

The differential is given by $\partial_G=d-\tilde\Omega c(T)$, where
$\tilde \Omega $ is the element in the dual of the Lie algebra $i\RR$
such that $<\tilde \Omega, \Omega>=1$ under the trace pairing, and
$c(T)$ is the contraction with the vector field $T$ generated by the
infinitesimal action of $U(1)$.

As we have already discussed at the beginning of Section 5, the
differential $\partial_{U(1)}$ of the equivariant 
complex $\Omega_{*, U(1)}(O_a)$ acts as 
$$ D (\Omega^n \otimes 1_a)=0 $$
and
$$ D (\Omega^n \otimes \eta_a)= - \Omega^{n-1} 1\otimes 1_a. $$

The homology of this complex (for a fixed irreducible orbit $O_a$) is
therefore one copy of $\RR$ (corresponding to the generator $1\otimes
\eta_a$) in degree zero and zero in all other degrees. This is just
the usual result that the equivariant homology of $O_a$ with a free
action of $U(1)$ is isomorphic to the ordinary homology (with real
coefficients) of the quotient, that is of a point.

Now let us consider the bigraded complex where several
orbits are considered. The bigraded complex is of the form
\[ \bigoplus_{\mu(O_a)+j=*} \Omega_{j,U(1)}(O_a). \]

The equivariant boundary operator can be written explicitly in
components of the form 
\beq
D: \begin{array}{ccc}
\Omega^n \otimes 1_a &\mapsto & -n_{ab}\Omega^n \otimes 1_b \\
                        &        & \\
\Omega^n \otimes \eta_a &\mapsto & (n_{ab}\Omega^n \otimes \eta_b)\oplus
                                   (m_{ac}\Omega^n \otimes 1_c) \\
                        &        & \oplus(-\Omega^{n-1}\otimes 1_a).
\end{array}
\label{equivD}
\eeq

Here $n_{ab}=\# \hat{\cal M}(O_a,O_b)^*$, where the relative index is
$\mu(O_a)-\mu(O_b)=1$, so 
that $-(e_b^-)_* (e_a^+)^*1_a = -n_{ab}1_b$ and $(e_b^-)_*
(e_a^+)^*\eta_a = n_{ab}\eta_b$. The coefficient $m_{ac}$ comes from the
integration of the one-form $\eta_a$ over the one-dimensional fiber of
the moduli space ${\cal M}(O_a,O_c)^*$ with $\mu(O_a)-\mu(O_c)=2$, namely
$(e_c^-)_* (e_a^+)^*\eta_a = m_{ac}1_c$. The components of $D$ are
represented in the following diagram.

$$
\diagram
& \vdots \dto & \vdots \dto & \vdots \dto \\
\cdots \rto & \Omega^k\otimes 1_a \rto^{-n_{ab}} \dto_0 &
\Omega^k\otimes 1_b \rto^{-n_{bc}} \dto^0 & \Omega^k\otimes 1_c \rto
\dto^0 & \cdots\\
\cdots \rto & \Omega^k\otimes \eta_a \rto^{\phantom{mm} n_{ab}} \dto_{-1} 
\urrto^{m_{ac}\qquad}|{\phantom{mmm}} &
\Omega^k\otimes \eta_b \rto^{n_{bc}} \dto^{-1} & \Omega^k\otimes \eta_c \rto
\dto^{-1} & \cdots\\
\cdots \rto & \Omega^{k-1}\otimes 1_a \rto^{-n_{ab}} \dto   &
\Omega^{k-1}\otimes 1_b \rto^{-n_{bc}} \dto   & \Omega^{k-1}\otimes 1_c \rto
\dto   & \cdots\\
& \vdots      & \vdots      & \vdots
\enddiagram
$$

\begin{rem}
In view of the result of Proposition \ref{boundaryterms}, we notice
that the full strength of Theorem \ref{corners} is in fact more than
what is strictly necessary in order to prove that the boundary $D$ of
the equivariant Floer complex satisfies $D^2=0$. In fact, knowing that
the boundary operator only 
depends on the moduli spaces $\hat{\cal M}(O_a,O_b)$ of dimension at
most two, makes it only necessary, strictly speaking, to know the
properties of the compactification for moduli spaces $\hat{\cal
M}(O_b,O_c)$ of dimension at most three. 
\end{rem}

Clearly, it does not make any substantial difference to derive the
properties of the fine structure of the compactification $\hat{\cal
M}(O_a,O_b)^*$  in this reduced case or in the general case proved in
Theorem \ref{corners}. However, both here and in Section 6, where we 
analyze the case of obstructed gluing, it seems useful to identify 
explicitly what is the minimum requirement on the fine structure of
the compactification that is needed in the arguments, cf. Remark
\ref{Icomponents} and Remark \ref{Hcomponents}.

\subsection{Intermezzo: the coefficient $m_{ac}$} 
 
We now give a topological description of the coefficient $m_{ac}$ 
that appears in the component 
$$ \langle D(\Omega^n\otimes \eta_a),\Omega^n\otimes 1_{c} \rangle 
=m_{ac} $$ 
of the boundary operator, for $\mu(O_a)-\mu(O_c)=2$. 
This coefficient is defined as the pullback and pushforward of the 
1-form $\eta_a$ along the $U(1)$-equivariant endpoint fibrations 
$$ e_a^+ : \hat\M(O_a,O_c)^* \to O_a $$ and 
$$ e_c^- : \hat\M(O_a,O_c)^* \to O_c, $$ 
that is, 
$$ m_{ac} = (e_c^-)_* (e_a^+)^* \eta_a. $$ 
 
The 2-dimensional $U(1)$-manifold $\hat\M(O_a,O_c)^*$ fibers over the circles  
$O_a$ and $O_c$ with 1-dimensional fibers. Topologically, 
$\hat\M(O_a,O_c)^*$ consists of a collection of finitely many 
cylinders and finitely many tori.  
 
\begin{lem} 
The pullback 1-form $(e_a^+)^* \eta_a$ defines a Cartan connection on 
the principal $U(1)$-bundle $\hat\M(O_a,O_c)\to \hat\M(a,c)$. 
\label{connection:1form} 
\end{lem} 
 
\noindent\underline{Proof:} 
We want to obtain a splitting of the exact sequence 
$$ 0 \to {\cal T}(U(1)) \to {\cal T}\hat\M(O_a,O_c) \to {\cal 
T}\hat\M(O_a,O_c)/{\cal T}(U(1))\cong {\cal T}\hat\M(a,c)\to 0. $$ 
Since the endpoint map $e_a^+$ is $U(1)$-equivariant, and $\eta_a$ 
is the generator of $H^1(O_a)$, the kernel $Ker((e_a^+)^* \eta_a)$ 
defines a horizontal subspace of ${\cal T}\hat\M(O_a,O_c)$. 
 
\noindent $\diamond$ 
 
\begin{lem} 
Over each component of $\hat\M(O_a,O_c)$, the pushforward 
$$(e_c^-)_* (e_a^+)^* \eta_a$$ computes the winding number $W(\gamma)$ of the  
1-dimensional $\gamma=(e_c^-)^{-1} (x_c)$, for a point
$x_c\in O_c$, around the fiber    
of the $U(1)$-fibration $\pi: \hat\M(O_a,O_c)\to \hat\M(a,c)$. 
\label{wind:number} 
\end{lem} 
 
\noindent\underline{Proof:} Consider a connected component $\hat X$ of 
$\hat\M(a,c)$. For $\mu(O_a)-\mu(O_c)=2$, the component
$$\pi^{-1}(\hat X)\subset \hat\M(O_a,O_c)$$ is a cylinder or a torus.  
Choosing a point $x_c \in O_c$, one can  
identify the fiber $$(e_c^-)^{-1} (x_c)\cap \pi^{-1}(\hat X)$$ with a lift 
$\gamma$ of the path $\hat X \subset  \hat\M(a,c)$ to 
$\hat\M(O_a,O_c)$. In fact, the endpoint fibration $e_c^-$ is 
compatible with the boundary strata 
$\hat\M(O_a,O_b)\times_{O_b}\hat\M(O_b,O_c)$. 
Since the   
1-form $\eta_a$ is the generator of $H^1(O_a)$, the pullback 
$(e_a^+)^* \eta_a$ integrated along the path  
$$ \gamma= (e_c^-)^{-1} (x_c)\cap 
\pi^{-1}(\hat X)$$  
measures precisely the number of times this path winds around the 
fibers of the map $\pi$, 
$$ \int_{(e_c^-)^{-1} (x_c)\cap \pi^{-1}(\hat X)}  (e_a^+)^* \eta_a = 
W(\gamma). $$  
 
\noindent $\diamond$

Notice that in $\hat\M(O_a,O_c)$ we have a choice of a horizontal 
direction $\ell(x_a)$,  
given by the integral lines of $Ker((e_a^+)^*\eta_a)$, for any initial 
condition $x_a\in O_a$, and a direction 
given by the path $(e_c^-)^{-1}(x_c)$ for fixed $x_c$ in $O_c$. 
The horizontal line $\ell(x_a)$ can be identified with the fiber of 
the other endpoint map $(e_a^+)^{-1}(x_a)$. 
Both the lines $(e_c^-)^{-1}(x_c)$ and $(e_a^+)^{-1}(x_a)$ represent 
a homeomorphic lift of $\hat\M(a,c)$ to $\hat\M(O_a,O_c)$, however, in 
general, an  
element $[\AA,\Psi]\in (e_c^-)^{-1}(x_c)$ has a limit at $t\to-\infty$  
which is some $\lambda x_a$ in $O_a$, not equal to $x_a$, and 
similarly an element $[\AA,\Psi]\in (e_a^+)^{-1}(x_a)$ has a limit at 
$t\to +\infty$ which is some $\lambda x_c$ in $O_c$ not equal to 
$x_c$. Different elements in 
$(e_a^+)^{-1}(x_a)$ have different $\lambda x_c$ limits in $O_c$ and 
similarly for the other endpoint fibration. 
The intersection between the two lines then consists of finitely 
many points. By our identification of $\ell(x_a)=(e_a^+)^{-1}(x_a)$ 
with the horizontal direction, this number is equal to the winding 
number of the path $(e_c^-)^{-1}(x_c)$. This intersection number 
counts the solutions $[\AA,\Psi]\in \hat\M(O_a,O_c)$ that actually 
have prescribed endpoints $x_a$ and $x_c$ at both ends.

Now consider the $U(1)$-bundle $\M(O_a,O_c)\to \M(a,c)$ over the 
parameterized moduli space. Connected components $X$  of $\M(a,c)$ are  
infinite cylinders or strips of the form $X\cong [0,1]\times \RR$.  
Consider a choice of a base point $(y_0,t_0)$ on $Y\times  
\RR$, and a complex line $\ell_{y_0}$ in the fiber $W_{y_0}$ of the 
spinor bundle $W=S\otimes L$ over $y_0\in Y$. We choose $\ell_{y_0}$ 
so that it does not contain the spinor part $\psi$ of any irreducible 
critical point. Since there are only finitely many critical points we 
can guarantee such choice exists. We consider the line bundle 
\beq {\cal L}_{ac} =\M(O_a,O_c)\times_{U(1)} (W_{y_0}/\ell_{y_0}) \to 
\M(a,c) \label{Lac} \eeq 
with a section given by 
\beq s([\AA,\Psi])=([\AA,\Psi],\Psi(y_0,t_0)). \label{sec:ac} \eeq 
 
\begin{lem} 
For a generic choice of $(y_0,t_0)$ the section $s$ of (\ref{sec:ac}) 
has no zeroes on the boundary strata $\hat\M(a,b)\times \hat 
\M(b,c)$. This determines a  
trivialization of ${\cal L}_{ac}$ away from a compact set in $\M(a,c)$. 
The line bundle ${\cal L}_{ac}$ over $\M(a,c)$, with the 
trivialization $\varphi$ specified above, 
has relative Euler class satisfying 
$$ e({\cal L}_{ac},\varphi)=m_{ac}. $$ 
\label{rel:eul:class} 
\end{lem} 
 
\noindent\underline{Proof:} The section $s$ induces a trivialization 
of the line bundle ${\cal L}_{ac}$ away from a compact set in 
$\M(a,c)$. In fact, for  $T\to \infty$, we have  
$$ s((\AA,\Psi)^{-T}) \to s((A_a,\psi_a))\neq 0 $$ 
and 
$$ s((\AA,\Psi)^{T})\to s((A_c,\psi_c))\neq 0. $$ 
This corresponds to trivializing the fibration $\hat{\cal 
M}(O_a,O_c)\times \{ T \}$ with the horizontal lines 
$\ell(x_a,\lambda)=(e_a^+)^{-1}(\lambda x_a)$, for $T\leq -T_0$ and with 
the lines $\ell(x_c,\lambda)= (e_c^-)^{-1}(\lambda x_c)$ for $T\geq 
T_0$. The obstruction to extending the trivialization over all of 
$\M(a,c)$ is then measured precisely by the winding number 
$m_{ac}$. This is therefore the relative Euler class of ${\cal 
L}_{ac}$, which can also be computed as 
$$ e({\cal L}_{ac},\varphi)=\# s^{-1}(0). $$ 
 
\noindent $\diamond$ 
 
In particular, this implies that the coefficient $m_{ac}$ can be 
computed by counting zeroes of any transverse section  
$$ s: \M(a,c) \to {\cal L}_{ac} $$ 
which is non-vanishing away from a compact set, and induces the same 
trivialization $\varphi$, 
$$ m_{ac} =\# s^{-1}(0). $$ 
 
We discuss some identities satisfied by the coefficients $m_{ac}$. 
 
\begin{rem} 
Let $a$ and $d$ be two irreducible critical points with 
$\mu(a)-\mu(d)=3$. Assume that all the critical points $c_i$ with  
$\mu(a)> \mu(c_i)> \mu(d)$ are irreducible. Then we have the identity  
$$ \sum_{c_1: \mu(a)-\mu(c_1)=1} n_{a,c_1} m_{c_1,d} - 
\sum_{c_2:\mu(c_2) -\mu(d)=1} m_{a,c_2}n_{c_2,d}=0. $$ 
When the reducible critical point $\theta$ 
has index between $(\mu(a), \mu(c)) $, with 
$\mu(a)=1$ and $\mu(d)=-2$, we have the identity 
\beq 
\sum_{c_1:  \mu(c_1)=0} n_{a,c_1} m_{c_1,d} 
 + n_{a\theta} n_{\theta d}  
-\sum_{c_2:\mu(c_2) =- 1} m_{a,c_2}n_{c_2,d} =0.
\label{id:ad} 
\eeq 
\end{rem} 
 
This identity follows directly from the definition of the invariant 
$m_{ac}$ as pullback and pushforward of the 1-form $\eta_a$, using 
Stokes' theorem for fibrations with boundary on the compactified 
moduli space $\hat\M(O_a,O_d)^*$. The sign comes 
from the orientations in the gluing theorem. 
 However, the identity (\ref{id:ad}) also has a 
topological interpretation in terms of relative Euler 
classes. In fact, the counting above is the counting of the boundary 
points of the zero set of a strata transverse section of the line bundle 
${\cal L}_{ad}$ over the unparameterized moduli space $\M(a,d)$. For 
such section, $s^{-1}(0)$ is a 1-dimensional manifold with boundary.  
The counting above corresponds to a section that has 
$$ \partial s^{-1}(0) = \bigcup_{c_1:\mu(a)-\mu(c_1)=1} \hat\M(a,c_1) 
\times \bigl(s^{-1}(0)\cap \M(c_1,d)\bigr) $$  
$$\cup \bigcup_{c_2: \mu(c_2) 
-\mu(d)=1} \bigl(- s^{-1}(0)\cap \M(a,c_2)\bigr) \times \hat\M(c_2,d). $$ 
The other identity can be proved by the similar methods, with an 
extra $U(1)$ gluing parameter when gluing along 
$\theta$, as in Lemma \ref{glue}. 
 
Now consider the case of the 2-dimensional moduli spaces 
$\hat\M(O_a,\theta)$, for $\mu(O_a)-\mu(\theta)=2$,  
and $\hat\M (\theta, O_c)$ for $\mu(\theta) -\mu(O_c) =3$, 
with $\theta$ the 
unique reducible (the $U(1)$-fixed point in ${\cal M}^0$). 
 
The moduli spaces $\hat\M(O_a,\theta)$, with 
$\mu(O_a)-\mu(\theta)=2$, do not contribute components to the boundary  
operator $D$ of the Floer complex, because we have 
$$ (e_\theta)_* (e_a^+)^* \eta_a =0, $$ 
since the fibers of the endpoint map 
$$ e_\theta : \hat\M(O_a,\theta)\to \theta $$ 
are 2-dimensional.  
 
It is still true that $(e_a^+)^* \eta_a$ 
defines a connection on the principal $U(1)$-bundle over 
$\M(a,\theta)$. Moreover, a component $\hat X\subset \hat\M(a,\theta)$ 
is again topologically a circle or a line segment. Given a choice of a 
lift $\gamma$ of $\hat X$ inside the component $\pi^{-1}(\hat 
X)\subset \hat\M(O_a,\theta)$, we can still compute its winding number 
and follow the same construction given above. However, now 
the choice of the lift $\gamma$ is no longer determined 
by the endpoint fibration  
$$ e_\theta: \hat\M(O_a,\theta) \to \theta. $$ 
Different possible choices of the lift $\gamma$ can have 
different winding numbers $W(\gamma)$.  
 
Similarly, we can still construct the bundle ${\cal 
L}_{a,\theta}$ over $\M(a,\theta)$. However, the section 
(\ref{sec:ac}), in this case, does not give a trivialization away from  
a compact set. Similarly, the fibers of the endpoint map $e_\theta$ do  
not provide a trivialization at the $T\to \infty$ end, since they do 
not determine a lift of $\M(a,\theta)$.  
 
However, we claim that it is still possible to associate an invariant
$m_{a\theta}$ to the moduli space $\M(a,\theta)$, which can be
interpreted topologically as a relative Euler class. We briefly
illustrate here why this should be the case, then in Section 6.3, in Remark
\ref{sec:obstr:induction}, we give a definition of the
invariant $m_{a\theta}$ as counting  
of the zeroes of a particular section of ${\cal 
L}_{a,\theta}$, the {\em obstruction section}.

The choice of the base point $(y_0, t_0)$ determines a trivialization
of the $U(1)$-fibration ${\cal B}^0$ over the configuration space
${\cal B}$, in a small neighborhood of each irreducible critical point. 
If we think of parameterized flow lines between 
two critical points $a$ and $c$ as embedded in the irreducible part
${\cal B}^*$ of ${\cal B}$, then the $U(1)$-fibration   
$\M (O_a, O_c)$  has an induced trivialization away from a compact set 
in $\M (a, c)$, which agrees with the trivialization $\varphi$ given by the
cross section (\ref{sec:ac}) of the associated line bundle 
${\cal L}_{ac}$.  We can choose a small neighborhood $U_\theta$  
of the reducible critical point $\theta$ such that 
$U_\theta ^* = U_\theta \cap {\cal B}^*$ is 1-connected, hence 
there is also a unique trivialization of ${\cal B}^0$ over 
$U_\theta^*$. We assume that this trivialization  
is given by a constant cross section. Thus,
there is an induced trivialization  
away from a compact set in $\M (a, \theta)$ 
(where $\mu(O_a) -\mu (\theta) =2$) or 
$\M (\theta, O_c)$ (where $\mu(\theta) - \mu (O_c) =3$).
Using this trivialization, we can define  
$m_{a\theta}$ and $m_{\theta c}$ to be the relative 
Euler class of the corresponding associated 
line bundles  ${\cal L}_{a\theta}$ and ${\cal L}_{\theta c}$, 
respectively. 

We return to the description of the invariant $m_{ac}$ 
for the case of $O_a$ and $O_c$ free orbits of relative index two, 
in order to illustrate the point of view we shall adopt in Remark
\ref{sec:obstr:induction} in describing the invariant $m_{a\theta}$.

In general, we can construct the relative Euler class which 
is a cycle defined  by the zeros of the transverse section 
as follows.  Over the unparameterized moduli space $\hat\M(a,c)$,
consider the line bundle   
$$ {\cal L}_{ac}= \hat\M(O_a,O_c)\times_{U(1)} \CC, $$ 
associated to the principal $U(1)$-fibration  
$$ \pi: \hat\M(O_a,O_c) \to \hat\M(a,c). $$ 
The pullback under the quotient map $\M(a,c)\to \hat\M(a,c)$ induces a 
pullback line bundle, which we still denote $ {\cal L}_{ac}$, over the  
parameterized moduli space $\M(a,c)$. 
 
\begin{lem} 
Over the boundary strata 
$$ \hat\M(a,b)\times \hat\M(b,c) $$ 
of $\hat\M(a,c)$, the line bundle ${\cal L}_{ac}$ has the form 
$$ {\cal L}_{ac}\cong \pi_1^* {\cal L}_{ab} \otimes \pi_2^* {\cal 
L}_{bc}. $$   
\label{tens:line} 
\end{lem} 
 
\noindent\underline{Proof:} The result follows by showing that, for all 
critical orbits $O_p$ and $O_q$ with $\mu(O_p)> \mu(O_q)$, the principal 
$U(1)$-bundles  
$$ \pi: \hat\M(O_p,O_q) \to \hat\M(p,q) $$ 
are compatible with the gluing maps, namely they define a line bundle 
$$ {\cal L}_{pq} = \hat\M(O_p,O_q)^* \times_{U(1)} \CC \to 
\hat\M(p,q)^* $$ 
in the category of manifolds with corners. The compatibility of the 
$U(1)$-fibration with the gluing maps follows from our proof of 
Theorem \ref{equivgluing}. 
 
\noindent $\diamond$ 
 
We now specify a choice of a class of transverse sections of ${\cal 
L}_{ac}$ over  $\M(a,c)$. 
We prescribe a choice of  non-zero  
complex numbers $s_a$ for each critical point $a$.  
For each pair of critical points $a, b$ of relative index $1$,
we choose a path of non-zero complex numbers 
 $s_{ab}: \RR \to \CC-\{ 0\} $ connecting $s_a$ and $s_b$, 
 which is contractible   in $\CC-\{ 0\}$ and is constant
outside a compact set in $\RR$. Over   
the 1-dimensional $\M(a,b)$ we 
define nowhere vanishing sections of ${\cal L}_{ab}$ as pullbacks of 
the $s_{ab}$ on $\hat\M(a,b)$. We still denote these sections over 
$\M(a,b)$ by $s_{ab}$.  
Over the 2-dimensional $\M(a,c)$ we consider the class of all 
transverse sections $s_{ac}$ of ${\cal L}_{ac}$ that satisfy 
$$ s_{ac}= \pi^*_1 s_{ab}\otimes\pi^*_2  s_{bc} $$ 
over the product submanifolds 
$$ \hat \M(a,b)\times \hat \M(b,c) \times [0, \epsilon) \times \RR $$ 
where $[0, \epsilon)$ is the gluing parameter,  
and such that, for a sufficiently large $T$,  
$s_{ac}$ is the constant section 
$s_a$ over  
$$\hat \M (a, c) \times   (-\infty, -T]$$ 
and the constant section $s_b$ 
over 
$$ 
\hat \M (a, c) \times   [T, \infty). $$  
Clearly, any such section $s_{ac}$ 
gives the same trivialization $\varphi$ of ${\cal L}_{ac}$ 
away from a compact set  
in $\M (a, c)$. 
We can therefore compute the relative Euler 
class by counting the zeros of this transverse section, that is, 
$$ e({\cal L}_{ac},\varphi)=\# s^{-1}_{ac} (0). $$ 
 
For higher dimensional moduli spaces $\M(a,c)$, we proceed 
analogously, by inductively defining the class of transverse 
sections of ${\cal L}_{ac}$ which are compatible with the boundary 
strata of $\M (a, c)$, that is, we require that 
$$ s_{ac}= \pi^*_1s_{ab}\otimes \pi^*_2 s_{bc} $$ 
over all the product submanifolds 
$$ \hat \M(a,b)\times \hat \M(b,c) \times [0, \epsilon) \times \RR $$  
and that $s_{ac}$ agrees with  the constant section 
$s_a$ over  
$$\hat \M (a, c) \times   (-\infty, -T]$$ 
and with the constant section $s_b$ 
over 
$$ 
\hat \M (a, c) \times   [T, \infty). $$

\begin{lem} 
Consider the line bundle ${\cal L}_{ac}$ over a parameterized moduli 
space $\M(a,c)$, where $\mu(a)-\mu(c) =3$ and $a$ is  irreducible 
or where $a$ is reducible and $ \mu(c) =-4$, 
 with the cross section $ s_{ac}$ as above. Then the zero set 
$s_{ac}^{-1}(0)$ is a co-dimension 2 submanifold of $\M(a,c)$ which 
can be compactified by adding lower dimensional   
boundary strata such that the codimension one boundary strata are given by 
$$ \partial s_{ac}^{-1}(0) = \bigcup_{ \mu(b)-\mu(c)=1}  
s_{ab}^{-1}(0) \times \hat\M(b,c) $$ 
$$\cup \bigcup_{  \mu(b)-\mu(c) =2} \hat\M(a,b)\times 
s_{bc}^{-1}(0),$$ 
for $\mu(a) \neq 1$;  
$$ \partial s_{ac}^{-1}(0) = \bigcup_{ \mu(b)-\mu(c)=1}  
s_{ab}^{-1}(0) \times \hat\M(b,c) $$ 
$$\cup \bigcup_{ \mu(a)-\mu(b)=1} \hat\M(a,b)\times 
s_{bc}^{-1}(0)$$  
$$ \cup    s^{-1}_{ac} (0) \cap   \bigl( \hat\M (a, \theta) 
\times \M (\theta, c) \times U(1)\bigr)$$ 
for $\mu(a) =1$ and $\mu(c) =-2$. 
\label{zeroes:s:ac} 
\end{lem} 
 
\noindent\underline{Proof:} The result follows from the choice of the 
class of sections $s_{ac}$, together with the result of Lemma 
\ref{tens:line}. 
 
\noindent $\diamond$ 
 
If we compare the results of Lemma \ref{tens:line} and Lemma 
\ref{zeroes:s:ac}, with Lemma \ref{rel:eul:class}, we see that  
our previous definitions of the invariant  
$$ m_{ac} =e({\cal L}_{ac},\varphi)= \# s_{ac}^{-1}(0) $$  
are a particular case of the procedure illustrated here.  
In particular, the identities  
of Remark \ref{id:ad} are then the direct consequences of Lemma 
\ref{zeroes:s:ac}.

\section{Topological invariance} 
 
In the present section we show that the equivariant Floer groups, 
defined as in Definition \ref{floergroups}, are topological 
invariants, in the sense that   
they are independent of the choice of the metric and of the perturbation. 
 
Throughout all this discussion,  
recall that the Floer groups depend on the choice of a 
$Spin_c$-structure. We always assume to work with a fixed $Spin_c$-structure, 
so we never mention explicitly this dependence.  
 
The easiest case, for the proof of topological invariance, is under the 
assumption that $b^1(Y)\geq  
2$. In fact, in this case there are no reducible solutions, hence a 
cobordism argument can be used to construct a morphism of the chain 
complexes that gives an isomorphism of the cohomology groups.  
In the case when $b^1(Y)=1$, the metric and 
the perturbation cannot be chosen arbitrarily. This was already 
noticed in \cite{Au2} in the case of the 
invariant associated to the moduli space ${\cal M}$ (which is the 
Euler characteristic of the Floer homology, \cite{Ma}). However,  
given two generic metrics, it is possible to find a path of metrics joining 
them and an open set of sufficiently small perturbations $\rho$ such that 
the invariant does not change. The dependence which is introduced in this case 
is rather mild and it amounts to a choice of a homology class in 
$H^1(Y;\ZZ)$. Thus, the $b^1(Y)\geq 1$ case will follow by exactly the  
same argument that we present below in our simple model case of same 
chamber metrics and perturbations. 
 
The case of a homology sphere \cite{Wa} is rather different. In fact 
in this case there is an explicit dependence on the metric (\cite{D}, 
\cite{KM2}). The interesting fact about the equivariant Floer homology 
is that it is instead metric independent. This means that, 
in this case, though we still have a condition on the kernel of the 
Dirac operator  
($Ker(\partial_A)$ trivial at the reducible solution) that breaks the 
space of metrics and perturbations into chambers with co-dimension  
one walls, it is however  
possible to construct a chain homomorphism 
between the Floer complexes that correspond to different metrics, and 
a chain homotopy. That is, for the equivariant Floer homology the wall  
crossing consists only of a global grade shift. 
 
The main theorem of this section is the following.  
 
\begin{thm} 
\label{metrics} 
Let $Y$ be a homology sphere or a rational homology sphere. 
Suppose given two metrics $g_0$ and $g_1$ on $Y$ and perturbations 
$\nu_0$ and $\nu_1$. Then there exist a chain homomorphism $I$ 
between the equivariant Floer complexes 
$({C_*}_{U(1)(g_0,\nu_0)},D)$ and   
$({C_*}_{U(1)(g_1,\nu_1)},D)$, defined by considering a generic path 
$(g_t,\nu_t)$ and moduli spaces ${\cal M}(O_a,O_{a'})$ of solutions of 
the flow equations on $(Y\times \RR, g_t +dt^2)$.  
This chain homomorphism induces an isomorphism in 
homology up to an index shift given by the spectral flow of the  
Dirac operator $\partial^{g_t}_{\nu_t}$.  
The analogous construction works for cohomology groups. 
\end{thm} 
 
We first give a ``model proof'' in the much simpler case 
of two metrics and perturbations $(g_0,\nu_0)$ and 
$(g_1,\nu_1)$ that are in the {\em same} chamber. 
The proof in this case is obviously much simpler than 
the general statement of Theorem \ref{metrics}. In fact, the 
interesting part of the statement is only for {\em different} chambers.  
The reason why we present explicitly the proof of this case 
is because, once we know that all the 
$SWF_{*,U(1)}(Y,(g,\nu))$, for $(g,\nu)$ in the same 
chamber, are isomorphic, we can prove the general case 
just by analyzing the case of two {\em sufficiently 
close} metrics and perturbations $(g_0,\nu_0)$ and 
$(g_1,\nu_1)$ in two different chambers, across a  
codimension one wall. Moreover, having the model case 
presented explicitly, we can discuss the case of different 
chambers by pointing out the various points in the proof where the 
original argument has to be modified. We prove the general statement
of Theorem \ref{metrics} in Section 6.3.
 
\noindent\underline{Proof, Part I: the same chamber case:}  
Choose a path of metrics and perturbations $(g_t, \nu_t)$ with $t\in 
[0,1]$ that interpolates  
between $(g_0,\nu_0)$ and $(g_1,\nu_1)$. Consider the manifold 
$Y\times \RR$ endowed  
with the metric $g_t+dt^2$ on $Y\times \{ t \}$ extended as the 
product metric outside $Y\times [0,1]$.  
 
Consider on $Y\times \RR$ the perturbed ``gradient flow'' 
equations with respect to the metric $g_t$. Denote $\{ O_a, O_b, \cdots \}$ 
 the critical orbits for the metric $g_0$ and  $\{ O_{a'}, O_{b'}, 
 \cdots \}$ the critical orbits for the metric $g_1$.  
 
We adopt the convention of fixing the absolute grading of the Floer 
complex to be the relative grading with respect to the reducible 
solution $\mu(O_a)-\mu(\theta_0)$, and $\mu(O_{a'})-\mu(\theta_1)$, 
where $\theta_0=[\nu_0,0]$ and $\theta_1=[\nu_1,0]$. We perform a 
shift of grading in the complex ${C_k}_{U(1)}(Y, g_1,\nu_1)$ by 
setting $\mu(\theta_1)=-SF(\partial^{g_t}_{\nu_t})$, where 
$SF(\partial^{g_t}_{\nu_t})$ is the spectral flow
\index{$SF(\partial^{g_t}_{\nu_t})$} 
of the Dirac operator along the path of reducible solutions 
$[\nu_t,0]$. In the same chamber case this shift is irrelevant, since 
the spectral flow is trivial, but it will be relevant in the general 
case we prove in Section 6.3.  
 
We have moduli spaces ${\cal M}(O_a,O_{a'})$ of solutions modulo gauge 
\index{${\cal M}(O_a,O_{a'})$} 
transformations of suitably perturbed Seiberg--Witten equations 
on the manifold $(Y\times \RR, g_t+dt^2)$, as specified in (\ref{4SW1P't})  
and (\ref{4SW2P't}).  
 
We need an analogue of Proposition \ref{transverse} which ensures 
that, under a generic choice of the perturbation, the moduli spaces 
${\cal M}(O_a,O_{a'})$ are smooth manifolds of dimension 
$\mu(O_a)-\mu(O_{a'})+\dim(O_a)$, where  
$\mu(O_{a'})$ denotes the shifted grading. In Lemma \ref{surjT}, Lemma  
\ref{noflow}, Corollary \ref{reducibleflow}, and Lemma   
\ref{unique:reducible} we deal with these transversality issues.  
 
We construct a degree zero homomorphism of the Floer 
complexes \index{$I$} 
\[ I: {C_k}_{U(1)}(Y, g_0, \nu_0)\to {C_k}_{U(1)}(Y, g_1, 
\nu_1). \]  
We define the homomorphism $I$ as  
\beq \label{I:def} 
I_{a,a'}\eta=\left\{ \begin{array}{lr} 
 (e^-_{a'})_*(e_a^+)^* \eta & \mu(O_a)-\mu(O_{a'})\geq 0 \\ 
0 & \mu(O_a)-\mu(O_{a'})<0, \end{array}\right. \eeq 
with $\eta \in \Omega_{k-\mu(O_a), U(1)}(O_a)$. The maps $e^-_{a'}$ 
and $e_a^+$  \index{$e^+_a$, $e^-_b$}
are the end point maps for ${\cal M}(O_a,O_{a'})$. Clearly, such $I$ 
commutes with $\Omega c(T)$. 
Requiring that the expression (\ref{I:def}) is well defined, that is, 
that we have a well defined push-forward map $(e^-_{a'})_*$, implies 
checking some properties of the compactification of the moduli spaces 
${\cal M}(O_a,O_{a'})$. For the same chamber case, the analysis of 
the boundary structure of $\M(O_a,O_{a'})$ follows 
closely the model developed in Section 4 for the 
compactification of the moduli spaces of flow lines. 
However, when adapting this argument to the general 
case of different chambers, this compactification will 
require a much more delicate analysis. 
 
The next step is then to prove that the map $I$ is a chain 
homomorphism. We want the relation  $I\circ D-D\circ I=0$ to be 
satisfied. This corresponds to the expression  
\beq \begin{array}{ll} 
 &  (ID-DI)(\eta)_{a, b'}\\[1mm] 
=& \left\{ \begin{array}{l}
      (e_{b'}^-)_* (e_a^+)^*\partial_{U(1)}\eta - 
      \partial_{U(1)}(e_{b'}^-)_* (e_a^+)^*\eta   \\[2mm]
\hbox{for } O_b=O_a \hbox{ and } \mu(O_b)\geq \mu(O_{b'})\\[1mm]
\hbox{or }  O_{a'}=O_{b'}\hbox{ and } \mu(O_a)\geq \mu(O_{b'})  \\[3mm]  
    -(-1)^{r(\eta)} \sum_{\{ b| \mu(O_b) <\mu (O_a)  \} } 
    (e_{b'}^-)_*(e_b^+)^* (e_b^-)_* (e_a^+)^*\eta  \\[1mm] 
    +(-1)^{r((e^-_{a'})_*(e_a^+)^*\eta)} 
    \sum_{\{a'|\mu(O_{a'})>\mu(O_{b'})\}}  
    (e_{b'}^-)_*(e^+_{a'})^*(e^-_{a'})_*(e_a^+)^*\eta, \\[2mm]
\hbox{for } \mu(O_a)> \mu(O_b)\geq \mu(O_{b'}) \hbox{ and } 
\mu(O_a)\geq \mu(O_{a'}) >\mu(O_{b'})\end{array} \right.
\end{array} 
\label{I} 
\eeq 
where the first line in the right hand side of (\ref{I}) 
corresponds to $O_b=O_a$ and $\mu(O_b)\geq \mu(O_{b'})$  
or $O_{a'}=O_{b'}$ and $\mu(O_a)\geq \mu(O_{b'})$, and the  
other terms correspond to $\mu(O_a)> \mu(O_b)\geq  
\mu(O_{b'})$ and $\mu(O_a)\geq \mu(O_{a'}) >\mu(O_{b'})$.  
 
In order to show that the right hand side is zero, we consider 
components of the form 
\beq 
\begin{array}{l} 
\bigcup_{\{a'|\mu(O_{a'}) >\mu(O_{b'})\}} ({\cal M}(O_a,O_{a'})^* 
\times_{O_{a'}}  
\hat {\cal M}(O_{a'}, O_{b'})^*) \\[1mm] 
\bigcup_{\{ b| \mu(O_b) <\mu (O_a)  \}}(-\hat  {\cal M}(O_a,O_b)^* 
\times_{O_b}   
{\cal M}(O_b,  O_{b'})^*), 
\end{array} 
\label{boundary:t} 
\eeq 
We want to show that these 
are precisely all the components of the ideal boundary that appear in 
the actual boundary of the compactification 
${\cal M}(O_a,O_{b'})^*$. We need Stokes' theorem to 
apply, hence \index{${\cal M}(O_a,O_{a'})^*$}
we need to prove that the compactification $\M(O_a, 
O_{a'})^*$ with boundary strata (\ref{boundary:t}) has the structure 
of a smooth manifold with corners. 
We also need 
an analogue of Corollary \ref{strata}, which shows that the endpoint 
maps   \index{$e^+_a$, $e^-_b$}
\[ e_a^+ : {\cal M}(O_a,O_a')^*\to O_a \] 
and  
\[ e_{a'}^-: {\cal M}(O_a,O_a')^*\to O_{a'} \] 
are fibrations compatible with the boundary strata and with compact 
fibers, as in Theorem \ref{glue:I}.
Thus, using Stokes' theorem again,  
we write the expression $(e_a^+)^*\eta $ 
in terms of the push-forward map on the co-dimension one boundary 
strata. This proves that (\ref{I}) vanishes identically, as expected.  
Notice that the sign $(-1)^{r((e^-_{a'})_*(e_a^+)^*\eta)}$ is just 
$(-1)^{r(\eta)}$.  
 
Now consider another path of metrics and perturbations $(\tilde g_t, 
\tilde \nu_t)$ that connects $(g_1,\nu_1)$ 
to $(g_0,\nu_0)$. We construct the corresponding homomorphism  
\[ J: {C_k}_{U(1)}(Y, g_1, \nu_1)\to {C_k}_{U(1)}(Y, g_0, 
\nu_0), \] 
as in the previous case, \index{$J$}
\beq \label{J:def} 
J_{a',a}\eta=\left\{ \begin{array}{lr} 
 (e^-_{a})_*(e_{a'}^+)^* \eta & \mu(O_{a'})-\mu(O_{a})\geq 0 \\ 
0 & \mu(O_{a'})-\mu(O_{a})<0, \end{array}\right. \eeq 
with $\eta \in \Omega_{k-\mu(O_{a'}), U(1)}(O_{a'})$. The maps $e^+_{a'}$ 
and $e_a^-$  are the end point maps for the moduli space ${\cal 
M}(O_{a'},O_{a})$ of solutions on $(Y\times \RR, \tilde g_t + dt^2)$.

The statement of the theorem now follows if we show that there is a 
chain homotopy $H$ such that 
\beq 
id_k -(JI)_k =D_{k+1}H_k+ H_{k-1} D_k.  
\label{chainhomotopy} 
\eeq 
 
In order to define $H$ let us consider the manifold $Y\times \RR$ 
endowed with the metric $\gamma_1$ which is  
\beq \label{gamma1} 
\gamma_1=\left\{ \begin{array}{lr} g_0+dt^2 & t< -2\\ 
g_{t+2}+dt^2 & t\in [-2,-1]\\ 
g_1+dt^2 & t\in [-1,1] \\ 
\tilde g_{2-t}+dt^2 & t\in [1,2]\\ 
g_0+dt^2 & t>2. 
\end{array}\right. \eeq 
 
Consider a path of metrics $\gamma_\sigma$ with $\sigma\in 
[0,1]$ that connects $\gamma_0 =g_0+dt^2$ to $\gamma_1$, such that for 
all $\sigma$ the metric $\gamma_\sigma$ is the product metric 
$g_0+dt^2$ outside a fixed large interval $[-T,T]$. 
 
Let ${\cal M}^{P}(O_a, O_b)$ be the parameterized  moduli space of 
\index{${\cal M}^{P}(O_a, O_b)$}
$(A(t),\psi(t),\sigma)$, solutions of the perturbed gradient flow 
equations with respect to the metric $\gamma_\sigma$, modulo gauge 
transformations. 
Lemma \ref{surjT}, with minor modifications to accommodate the presence 
of parameters, can be employed to show that , under a generic choice 
of the perturbation, the  
moduli spaces ${\cal M}^{P}(O_a, O_b)$ are smooth manifolds of dimension 
$\mu(O_a)-\mu(O_b) +2-\dim G_a $, cut out transversely by the equations. 
Denote by $\tilde e_a^+$ and $\tilde e_b^-$ the end point maps  
{}from ${\cal M}^{P}(O_a, O_b)$ to $O_a$ and $O_b$ respectively.  
 
Now we can define the degree-one map $H$ to be 
\[ H:  {C_k}_{U(1)}(Y, g_0, \nu_0)\to {C_{k+1}}_{U(1)}(Y, g_0, \nu_0) 
\] \index{$H$}
\beq \label{H:def} 
H_{a, b}: \quad  \eta \to (-1)^{1-r(\eta)} 
  (\tilde e_b^-)_*( \tilde e_a^+)^*\eta, \eeq 
with $\eta \in \Omega_{k-\mu(O_a), U(1)} (O_a)$. 
 
Again we need to verify some properties of the compactified moduli 
spaces ${\cal  
M}^{P}(O_a, O_b)^*$, which ensure that the 
map $H$ is well defined and has the desired properties.  
 
The identity (\ref{chainhomotopy}) which proves that $H$ is a chain 
homotopy can be rewritten as the following two identities: 
 
\beq\begin{array}{ll} 
& id_{a,a}-\sum_{a'} I_{a,a'} J_{a',a}\\[1mm] 
 = & \sum_{\{ b |\mu(O_b)\le\mu (O_a)\} } 
 D_{a,b}H_{b,a} +  \sum_{\{ b | \mu(O_b)\ge\mu(O_a) \}}  
H_{a,b}D_{b,a}\end{array}\label{aa}\eeq 
and, for $a \ne b$,  
\beq\begin{array}{ll} 
& -\sum_{a'} I_{a,a'} J_{a',b}\\[1mm] 
 = & \sum_{\{ c |\mu(O_c) \le \mu (O_a)\} } 
 D_{a,c}H_{c, b}+  \sum_{\{ c |\mu(O_c) \ge \mu (O_a)\} } 
H_{a,c}D_{c, b}.\end{array}\label{ab}\eeq 
 
The identities (\ref{aa}) and (\ref{ab}) can be proved by applying  
Stokes theorem again, in the way we 
discussed already. To this purpose, consider the components 
\beq \label{strataP:1} 
\begin{array}{c} 
\partial^{(1)} {\cal M}^P(O_a, O_a  )^* = \bigcup_{a'}({\cal M}(O_a, 
O_{a'})^* \times  
_{  O_{a'}} {\cal M}(O_{a'}, O_a)^* ) \cup \{ -O_a \} 
\\[1mm] 
\bigcup_{\{b| \mu(O_b)\ge\mu(O_a)\}} ((-1)^{dim  \hat {\cal M}(O_b, O_a) } 
 {\cal M}^P(O_a, O_b)^* 
\times_{  O_b}  
\hat {\cal M}(O_b, O_a)^*)\\[1mm] 
 \bigcup_{\{b| \mu(O_b)\le\mu(O_a)\}} ((-1)^{dim{\cal M}^P(O_b, O_a)} 
 \hat {\cal M}(O_a, O_b)^* 
\times_{  O_b}  
{\cal M}^P(O_b, O_a)^*)\end{array} 
\eeq 
and, for $a \ne b$, 
\beq \label{strataP:2} 
\begin{array}{c} 
\partial^{(1)} {\cal M }^P(O_a, O_b )^* = \bigcup_{a'}({\cal M}(O_a, 
O_{a'})^* \times  
_{  O_{a'}} {\cal M}(O_{a'}, O_b)^* )\\[1mm] 
 \bigcup_{\{c | \mu(O_c)\ge\mu(O_a)\}} ((-1)^{dim \hat{\cal M}(O_c,  O_b)} 
 {\cal M}^P(O_a, O_c)^* 
\times_{  O_c}  
\hat {\cal M}(O_c,  O_b)^*)\\[1mm] 
 \bigcup_{\{c| \mu(O_c) \le\mu(O_a)\}} ((-1)^{dim {\cal M}^P(O_c, , O_a)} 
\hat {\cal M}(O_a, O_c)^* 
\times_{  O_c}   
{\cal M}^P(O_c, , O_a)^*).\end{array} 
\eeq 
We need to know that these components 
are precisely the actual co-dimension one boundary in the 
compactifications ${\cal M}^P(O_a,O_a)^*$  
and  ${\cal M}^P(O_a, O_b)^*$ respectively. 
If we prove that the compactification $\M^P(O_a,O_b)^*$ has the 
structure of a smooth   
manifold with corners, with the actual boundary strata 
in the codimension one boundary consisting  
precisely of the components of (\ref{strataP:1}) and 
(\ref{strataP:2}), then the argument is complete. 
 
Again, since we 
are dealing only with the model case of same chamber metrics and 
perturbations, the necessary properties  
of the compactification $\M^P(O_a,O_b)^*$ follow from 
the analysis introduced in Section 4, as in the case of 
moduli spaces of flow lines. This is another crucial point that 
requires a much more accurate analysis in the general 
case of different chambers. 
The signs in the formulae (\ref{boundary:t}), (\ref{strataP:1}), and 
(\ref{strataP:2}) denote the difference of orientation, as in
Theorem \ref{glue:I} and Theorem \ref{glue:H}. 
 
\noindent $\diamond$ 
 
In particular, we can describe more explicitly the components of the 
maps $I$ and $H$. This identifies what are, in fact, the minimal 
requirements on the compactifications of the moduli spaces 
$\M(O_a,O_{a'})$ and $\M^P(O_a,O_b)$ that are necessary for the 
argument given above to work. 
 
Let $\eta_a$ be the one-form that generates $H^1(S^1)$ and $1_a$ be 
the constant function equal to 1 on $O_a$, so that we can write the 
generators of the equivariant complex in the form $\Omega^n \otimes 
\eta_a$ and $\Omega^n\otimes 1_a$. We have the following explicit 
description of the maps $I$ and $H$. 
 
\begin{rem} 
\label{Icomponents} 
The components of the map $I$ vanish whenever we have 
$\mu(O_a)-\mu(O_{a'})\geq 2$. Moreover, the only possibly non-trivial 
components of the chain map $I$ are  \index{$n_{aa'}$}
\[ n_{aa'}=\langle \Omega^n\otimes \eta_{a'}, I(\Omega^n\otimes \eta_a) 
\rangle, \] 
and 
\[ n_{aa'}=\langle \Omega^n\otimes 1_{a'}, I(\Omega^n\otimes 1_a) \rangle, \] 
when $\mu(O_a)-\mu(O_{a'})=0$, with $O_a$ and $O_{a'}$ free 
$U(1)$-orbits. In this case we need to know that the  
moduli space ${\cal M}(O_a,O_{a'})$ has a nice compactification ${\cal 
M}(O_a,O_{a'})^*$ analogous to the compactification of the moduli 
spaces of flow lines analyzed in Sections 4.2 and 4.3. We then have 
$$  n_{aa'}=\# \left({\cal M}(a,a')={\cal M}(O_a,O_{a'})/U(1)\right). $$ 
We also have \index{$m_{ab'}$}
\[ m_{ab'}=\langle \Omega^n\otimes 1_{b'}, I(\Omega^n\otimes \eta_a) \rangle, 
\] 
when $\mu(O_a)-\mu(O_{b'})=1$ and with $O_a$ and $O_{b'}$ free 
$U(1)$-orbits. In this case we need to know that the 
compactification ${\cal M}(O_a,O_{b'})^*$ is obtained by adding the 
co-dimension one boundary strata (\ref{broktraj}) and satisfies the 
analogue of Theorem \ref{equivgluing}, hence the push-forward map 
$(e^-_{b'})_*$ is well defined. Moreover, in the case of the 
critical points  
$\theta_0$ and $\theta_1$ we have \index{$n_{a\theta_1}$}
\[ n_{a,\theta_1}=\langle \Omega^n\otimes \theta_1, I(\Omega^n\otimes \eta_a) 
\rangle, \]  
when $\mu(O_a)-\mu(\theta_1)=0$. The coefficient in this case is obtained by 
integrating along the 1-dimensional fibers of the map  
$ e_{\theta_1}: {\cal M}(O_a,\theta_1)\to \theta_1$. Thus, we need to 
know that the moduli space ${\cal M}(O_a,\theta_1)$ is compact. This 
gives $$ n_{a,\theta_1}=\# \M(a,\theta_1). $$ 
We also have the component \index{$n_{\theta_0 b'}$}
\[ n_{\theta_0,b'}=\langle \Omega^n\otimes 1_{b'}, I(\Omega^n\otimes 
\theta_0) \rangle,  
\] 
when $\mu(\theta_0)-\mu(O_{b'})=1$. Here we have a 1-dimensional 
space ${\cal M}(\theta_0,O_{b'})$ fibering over $O_{b'}$ with 
zero-dimensional fiber, and again we need to know that we have 
${\cal M}(\theta_0,O_{b'})^*={\cal M}(\theta_0,O_{b'})$ in this case, 
which gives 
$$ n_{\theta_0,b'}=\# \M(\theta_0,b'). $$ 
\end{rem}

\begin{rem} 
\label{Hcomponents} 
The components of the map $H$ vanish whenever we have 
$\mu(O_a)-\mu(O_{b})\geq 1$. Moreover, the only possibly non-trivial 
components of the chain homotopy $H$ are \index{$n^P_{ab}$} 
\[ n^P_{ab}=\langle \Omega^n\otimes \eta_{b}, H(\Omega^n\otimes \eta_a) 
\rangle, \] 
and 
\[ n^P_{ab}=\langle \Omega^n\otimes 1_{b}, H(\Omega^n\otimes 1_a) \rangle, \] 
when $\mu(O_a)-\mu(O_b)=-1$, with $O_a$ and $O_b$ free 
$U(1)$-orbits. In this case we need to know that the one-dimensional  
moduli space ${\cal M}^P(O_a,O_b)$ is compact. This gives 
$$ n^P_{ab}=\# \left( \M^P(a,b)={\cal M}^P(O_a,O_b)/U(1) \right). $$ 
We also have \index{$m^P_{ac}$}
\[ m^P_{ac}=\langle \Omega^n\otimes 1_{c}, H((\Omega^n\otimes \eta_a) 
\rangle, \] 
when $\mu(O_a)-\mu(O_c)=0$ and with $O_a$ and $O_c$ free 
$U(1)$-orbits.  
In this case we need to know that the 
compactification ${\cal M}^P(O_a,O_c)^*$ of the 2-dimensional moduli 
space ${\cal M}^P(O_a,O_c)$ is obtained by adding the 
co-dimension one boundary strata (\ref{strataP:2}) and satisfies the 
analogue of Theorem \ref{equivgluing}, hence the push-forward map 
$(\tilde e_c^-)_*$ used in the definition of $H$ is well defined. 
 Moreover, in the case of the critical point $\theta_0$
\index{$n^P_{a,\theta_0}$} 
\[ n^P_{a,\theta_0}=\langle \Omega^n\otimes \theta_0, 
H(\Omega^n\otimes \eta_a) \rangle, \]  
when $\mu(O_a)-\mu(\theta_0)=-1$. In this case the moduli space ${\cal 
M}^P(O_a,\theta_0)$ is one-dimensional. 
Thus, in this case, the coefficient of $H$ is obtained by 
integrating along the 1-dimensional fibers of the map  
$ e_{\theta_0}: {\cal M}^P(O_a,\theta_0)^*\to \theta_0$. Thus, we need  
to know that the 1-dimensional moduli space ${\cal M}^P(O_a,\theta_0)$ 
has a nice compactification ${\cal M}^P(O_a,\theta_0)^*$. This gives 
$$ n^P_{a,\theta_0}=\# \M^P(a,\theta_0). $$ 
We also have the component \index{$n^P_{\theta_0,b}$}
\[ n^P_{\theta_0,b}=\langle \Omega^n\otimes 1_{b}, H(\Omega^n\otimes 
\theta_0) \rangle,  
\] 
when $\mu(\theta_0)-\mu(O_{b})=0$. Here we have a one-dimensional 
space ${\cal M}^P(\theta_0,O_{b})$ fibering over $O_{b}$ with 
zero-dimensional fiber, and again we need to know that ${\cal 
M}^P(\theta_0,O_b)^*$ is a nice compactification. 
This gives 
$$ n^P_{\theta_0,b}=\# \M^P(\theta_0,b). $$ 
\end{rem}

This gives a fair description of the model argument 
in the easy case of metric and perturbations in the  
same chamber, and of the minimal requirements on the 
compactifications, in order to extend the argument to the more general  
case. Now we can begin to analyze the general case.

\subsection{The boundary structure of ${\cal M}(O_a,O_{a'})$} 
 
In this section we analyze some properties of the moduli spaces 
${\cal M}(O_a,O_{a'})$ that 
are needed in the proof of Theorem \ref{metrics}. Here we make no 
assumption about the metrics and perturbations $(g_0,\nu_0)$, 
$(g_1,\nu_1)$. We shall soon restrict our attention to the case of two  
sufficiently close metrics and perturbations on opposite sides of a 
codimension one wall, which is the only case we need to complete the 
proof of Theorem \ref{metrics}. 
 
Let us first give a more precise account of the construction of the 
moduli spaces ${\cal M}(O_a,O_{a'})$. 
On $(Y\times\RR, g_t +dt^2)$ we consider the equations  
\beq 
D_{\AA}\Psi+ \rho\Psi=0 
\label{4SW1P't} 
\eeq 
and 
\beq 
F_{\AA}^+=\Psi\cdot\bar\Psi +i\mu +P_{(\AA,\Psi)}. 
\label{4SW2P't} 
\eeq 
Here the form $\mu$ is determined by the path $\nu_t$ of 
perturbations, and the additional 
perturbation $P$ satisfies conditions (2)-(5) of Definition \ref{calP} 
and the modified version of condition (1): 
 
\noindent (1') the perturbation $P_{(\AA,\Psi)}$ restricted to the 
interval $(1,\infty)$ can be written as $P_{(\AA,\Psi)}= *q_{(\AA,\Psi)}(t)  
+q_{(\AA,\Psi)}(t) \wedge dt$ with respect to the constant metric 
$g_1$, where $q_{(\AA,\Psi)}(t)$ satisfies  
\[ q_{(\AA,\Psi)|_{(1,\infty)}^T}(t)=q_{(\AA,\Psi)}(t+T), \] 
for any $T>0$ and for any $t>1$. Correspondingly, when restricted to 
the interval $(-\infty,0)$ the perturbation $q_{(\AA,\Psi)}(t)$, 
defined with respect to the constant metric $g_0$, 
satisfies 
\[ q_{(\AA,\Psi)|_{(-\infty,0)}^T}(t)=q_{(\AA,\Psi)}(t+T), \] 
for all $T<0$ and $t<0$.   
 
In the Dirac equation (\ref{4SW1P't}), we add the perturbation term 
given by a 1-form $\rho$ on $Y\times \RR$ that is supported on some 
compact set $Y\times [-T_0,T_1]$, with $T_0>0$ and $T_1>1$.  
Notice that, according to condition (1'), the equations 
(\ref{4SW1P't}) and (\ref{4SW2P't}) are translation invariant outside 
the interval $[-T_0,T_1]$. 
The moduli space ${\cal M}(O_a,O_{a'})$ is defined as the set of gauge  
equivalence classes of solutions of (\ref{4SW1P't}) and 
(\ref{4SW2P't}) in ${\cal B}^0_{k,\delta}(O_a,O_{a'})$.

Now, we first proceed as in Section 4.1 and prove the existence of a 
compactification for ${\cal M}(O_a,O_{a'})$. Theorem 
\ref{seqcompact2} (the analogue of Theorem \ref{seqcompact}) shows the  
existence of a compactification and identifies the ideal boundary with 
fibered products of lower dimensional moduli spaces, as in 
(\ref{broktraj}) and (\ref{broktraj:k}). This convergence argument 
follows very closely the argument for Theorem \ref{seqcompact}, 
whereas the gluing arguments will require a more 
sophisticated analysis.

\begin{thm} 
\label{seqcompact2} 
Suppose either $O_a$ or $O_{a'}$ is an irreducible critical orbit. 
Consider the moduli space ${\cal M}(O_a,O_{a'})$, with 
$\mu(O_a)-\mu(O_{a'})+\dim(O_a) \geq 1$. 
The codimension-one boundary of ${\cal M}(O_a,O_{a'})$ is a subset of 
the space  
\beq 
\label{broktraj} 
\bigcup_{O_b}\hat{\cal M}(O_a,O_b)\times_{O_b}{\cal M}(O_b,O_{a'}) 
\cup \bigcup_{O_{b'}} {\cal M}(O_a,O_{b'})\times_{O_{b'}}\hat{\cal 
M}(O_{b'},O_{a'}). \eeq 
Similarly, the components of higher codimension in the ideal boundary 
consist of fibered products of the form 
\beq 
\label{broktraj:k} \begin{array}{c} 
\bigcup_{c_1,\cdots, c_k, b} \hat{\cal 
M}(O_a,O_{c_1})\times_{O_{c_1}}\cdots \\[2mm]\times_{O_{c_k}} \hat{\cal 
M}(O_{c_k},O_b)\times_{O_b} {\cal M}(O_b,O_{a'}) \\[2mm] 
\bigcup_{b',c_1',\cdots, c_\ell '} {\cal 
M}(O_a, O_{b'})\times_{O_{b'}}\hat{\cal 
M}(O_{b'},O_{c_1'})\times_{O_{c_1'}} \cdots \\[2mm]\times_{O_{c_\ell '}} 
\hat{\cal M}(O_{ c_\ell '},O_{a'}). \end{array} 
\eeq 
\end{thm} 
 
\noindent\underline{Proof:} 
As in the proof of Theorem \ref{seqcompact}, we want to show that  
the space  is precompact. Namely, any sequence 
$x_i$ of elements in ${\cal M}(O_a,O_{a'})$ has a subsequence which 
either converges in norm to another solution $x\in{\cal 
M}(O_a,O_{a'})$, or converges to a broken trajectory.  
 
Given a sequence $x_i$ in ${\cal M}(O_a,O_{a'})$, we can follow the  
steps of the argument given in Theorem \ref{seqcompact}.  
 
\noindent{\bf Claim 1:} There is a subsequence $\{ x_{i_k} \}$ that 
converges smoothly on compact sets to a finite energy solution $y$ of 
the perturbed equations (\ref{4SW1P't}) and (\ref{4SW2P't}) on 
$(Y\times\RR, g_t +dt^2)$.  
 
The proof of Claim 1 can be divided in the two cases of a 
compact set contained in the complement of $Y\times [0,1]$, where the 
metric is constant, and of the set $Y\times [0,1]$ where the metric 
changes. The first case is a direct application of the proof given in 
Theorem \ref{seqcompact}. In the second case, we obtain convergence as 
in \cite{KM2}, Lemma 3.19. We can find a sequence of 
gauge transformations $\lambda_i$ such that the forms 
$\lambda_i(\AA_i-\AA_0)$ are co-closed and annihilate normal vectors 
at the boundary. Lemma 4 of \cite{KM} then shows that the sequence 
$\lambda_i(\AA_i,\Psi_i)$ has a subsequence that converges smoothly on 
$Y\times [0,1]$. 
  
Notice that the notion of finite energy, as 
well as the results of Corollary \ref{finenergy2} and Theorem 
\ref{decay}, extend to the case of solutions on $(Y\times\RR, g_t 
+dt^2)$, since the metric is constant outside the compact set $Y\times 
[0,1]$. This implies that the argument given in Theorem 
\ref{seqcompact} proves the following claim as well. 
 
\noindent{\bf Claim 2:} If the limit $y$ is an element of ${\cal 
M}(O_a,O_{a'})$ then the convergence $x_i\to y$ is strong in the 
$L^2_{k,\delta}$ norm.  
 
Thus we need to check whether there are obstructions to the 
convergence in norm, that is, whether broken trajectories arise as limits. 
 
Suppose the element $y$ is a limit smoothly on compact sets, but is not a 
limit in norm. Again, since the result of Corollary \ref{finenergy2} 
holds, the solution $y$ defines an element in some moduli space ${\cal 
M}(O_b,O_{b'})$.  
 
Since the metric is translation invariant only outside the interval 
$[0,1]$, consider the restriction of the sequence $x_i$ and of the 
limit $y$ to the interval $(1,\infty)$ with the constant metric $g_1$. 
We can choose a value $\alpha$ satisfying 
\[ {\cal C}_\rho(O_{b'})> \alpha > {\cal C}_\rho (O_{a'}). \] 
We can find times $T_i>1$ such that ${\cal C}_\rho(x_i(T_i))=\alpha$. 
The restriction to $(1,\infty)$ of the reparameterized sequence 
$x_i^{T_i}$ is again a sequence of solutions of the flow equations 
with the constant metric $g_1$ on $Y\times (1,\infty)$. A subsequence 
converges  
smoothly on compact sets to a solution $\tilde y$. This defines an 
element in some ${\cal M}(O_{c'},O_{d'})$, subject to the constraint  
\[ {\cal C}_\rho(O_{b'})\geq {\cal C}_\rho(O_{c'}) > {\cal 
C}_\rho(O_{d'}) \geq {\cal C}_\rho (O_{a'}). \] 
 
Notice that the proof of Theorem \ref{seqcompact2} is essentially 
along the same lines as the proof of the analogous result, Theorem 
\ref{seqcompact}, for the moduli spaces of flow lines $\hat{\cal 
M}(O_a,O_b)$.  
If we could guarantee that the cokernels of the linearizations are 
always trivial, then we could also extend the results of Section 4.2 
and 4.3 to this case, showing that the same properties of the 
compactification hold true. However, the essential new feature which 
appears in dealing with the moduli spaces ${\cal M}(O_a,O_{a'})$ is 
precisely that not all the cokernels are trivial. Thus, we need a much 
more accurate analysis of the fine structure of the compactification. 
 
The strategy we follow is to restrict our attention to just the cases 
which are needed for the proof of Theorem \ref{seqcompact2}, and show 
that, in those cases, we obtain the desired properties of the 
compactification. The purpose of the rest of Section 6 is to 
deal with these issues in detail. 
 
In order to identify where the problem of cokernels may arise, we 
study the transversality issue for the moduli spaces ${\cal 
M}(O_a,O_a')$. The following Lemma shows that the appearance of 
non-trivial cokernels is confined to the case of the moduli space 
${\cal M}(\theta_0,\theta_1)$ with 
$\mu(\theta_0)-\mu(\theta_1) < 0$. 
 
\begin{lem} 
Suppose given $(g_0,\nu_0)$, $(g_1,\nu_1)$ and a generic path 
$(g_t,\nu_t)$ connecting them. Consider the manifold $Y\times \RR$ 
with the metric $g_t +dt^2$, constant outside the interval 
$[0,1]$. Suppose given $x$, a solution of equations (\ref{4SW1P't}) 
and (\ref{4SW2P't}), with asymptotic values on the critical orbits 
$O_a$ and $O_{a'}$. Consider the linearization ${\cal L}_x$ of 
equations (\ref{4SW1P't}) and (\ref{4SW2P't}). Under a generic choice of the 
metric and of the perturbations $P$ and $\rho$, we have the following 
result. If the solution 
$x=(\AA,\Psi)$ satisfies $\Psi\neq 0$, then 
the linearization 
${\cal L}_x$ is surjective. If we have 
$O_a=\theta_0$ and $O_{a'}=\theta_1$ and the solution $x$ is 
of the form $x=(\AA,0)$, then the linearization ${\cal L}_x$ is 
surjective, provided that the relative Morse index satisfies 
$\mu(\theta_0)-\mu(\theta_1)\geq 0$.  
\label{surjT} 
\end{lem} 
 
The result of Lemma \ref{surjT} can be rephrased as the following 
transversality statement for the moduli spaces ${\cal M}(O_a,O_a')$. 
 
\begin{lem} 
Consider the following cases: 
 
\noindent (a) at least one of the asymptotic critical 
orbits $O_a$ and $O_{a'}$ has a free $U(1)$-action.  
 
\noindent (b) we have $O_a=\theta_0$ and $O_{a'}=\theta_1$ under the 
condition that $\mu(\theta_0)-\mu(\theta_1)\geq 0$. 
 
Then, for a generic choice of the metric and perturbations, 
the moduli space ${\cal M}(O_a,O_{a'})$ is a smooth 
oriented manifold 
that is cut out transversely by the equations, of dimension 
$\mu(O_a)-\mu(O_{a'})+\dim(O_a)$.  
In particular, in case (a) the moduli space is generically empty when 
the virtual dimension $\mu(O_a)-\mu(O_{a'})+\dim(O_a)$ is negative. 
\label{noflow} 
\end{lem}   
 
\noindent\underline{Proof of Lemma \ref{surjT} and \ref{noflow}:} No 
solution with  
$\Psi\equiv 0$ arises under the hypothesis of case (a). Thus the 
argument of Proposition \ref{transverse} and Proposition 
\ref{flowlines}, applied to the linearization   
\[ {\cal L}_{(\AA,\Psi,P,\rho)}(\alpha,\Phi)=\left\{\begin{array}{c} 
D_{\AA}\Phi+ \alpha \Psi + \rho \Phi \\ d^+\alpha 
-\frac{1}{2}Im(\Psi\cdot\bar\Phi)+{\cal D}P_{(\AA,\Psi)}(\alpha,\Phi)\\ 
G^*_{(\AA,\Psi)}(\alpha,\Phi) 
\end{array}\right. \] 
gives the desired transversality result. 
In case (b), we can find  
solutions in ${\cal M}(\theta_0,\theta_1)$ with $\Psi\neq 0$, as well 
as a solution with $\Psi\equiv 0$. The linearization is surjective at 
solutions with $\Psi\neq 0$ because of the argument of Proposition 
\ref{transverse} and Proposition \ref{flowlines}. Consider the 
linearization ${\cal L}_x$ at a  
solution of the form $x=(\AA,0)$. We have 
\[ {\cal L}_{(\AA,0,P,\rho)}(\alpha,\Phi)=\left\{\begin{array}{c} 
D_{\AA}\Phi+ \rho \Phi \\ d^+\alpha 
+{\cal D}P_{(\AA,0)}(\alpha,\Phi)\\ 
G^*_{(\AA,0)}(\alpha,\Phi) 
\end{array}\right. \] 
 
Consider the operator 
\[ \hat{\cal L}_{(\AA,0,P,\rho)}(\alpha,\Phi,p,\eta)={\cal 
  L}_{(\AA,0,P)}(\alpha,\Phi) +\left(\begin{array}{c} 
\eta \Phi \\ 
p_{(\AA,0,P)}(\alpha,\Phi)\end{array}\right), \] 
where we vary the perturbation of the curvature equation by an element 
$p_{(\AA,0,P)}$ of $T_P{\cal P}$, and the perturbation of the Dirac 
equation by a compactly supported 1-form $\eta$. 
We follow the analogous argument of Proposition 
\ref{transverse}. Suppose $(\beta,\xi)$ is an element orthogonal to 
the range of $\hat{\cal L}_{(\AA,0,P,\rho)}$. We have terms of the form 
\[ \langle \beta, d^+\alpha +{\cal 
D}P_{(\AA,0)}(\alpha,\Phi)+p_{(\AA,0,P)}(\alpha,\Phi) \rangle  
 \]  
\[ +\langle \xi, D_{\AA}\Phi+\eta\cdot\Phi \rangle=0. \] 
By varying $p\in {\cal P}$ we force $\beta\equiv 0$. If 
$\mu(\theta_0)-\mu(\theta_1)=0$, then the condition 
\[ \langle \xi,  D_{\AA}\Phi+\eta\cdot \Phi \rangle =0 \] 
implies that $\xi$ is in the kernel of $D_{\AA-\eta}$, but for a 
generic choice of $\eta$ this operator has trivial kernel. If 
$\mu(\theta_0)-\mu(\theta_1)>0$, the argument is similar:  
choose $\Phi$ to be a non-trivial element in $Ker(D_{\AA})$. This is 
possible since we  
are assuming that $\mu(\theta_0)-\mu(\theta_1)>0$, hence that $\dim 
Ker({\cal L}_x)>0$. We show that, by varying $\eta$ we force $\xi$ to 
vanish. In fact, suppose  
\[ \langle \xi, \eta\cdot \Phi \rangle =0 \] 
is satisfied for all possible compactly supported 1-forms $\eta$, and 
for a non-trivial $\xi$. Then we can follow the argument of \cite{Mo} 
6.2.1: choose some small open set where both $\Phi$ and $\xi$ are 
non-trivial and approximately constant. The isomorphism 
\[ \RR^4\otimes_{\RR} \CC \to Hom_{\CC} (S^+, S^-) \] 
induced by Clifford multiplication implies that there exist $\eta\in 
\RR^4$ such that 
\[ Re( \langle \xi, \eta\cdot \Phi \rangle) >0 \] 
at some point in the small open set $U$. By a trivialization of 
$T^*(Y\times R)$ on $U$ the vector $\eta$ can be extended to a 1-form 
$\eta$ on $U$. This 1-form can be extended to a  compactly supported 
1-form in a neighborhood of $U$, so that we get 
\[ \int_{Y\times\RR} Re( \langle \xi, \eta\cdot \Phi \rangle) >0, \] 
which contradicts the initial assumption. 
 
Notice that, since we consider the gauge group  
${\cal G}^0(O_a,O_{a'})$ of gauge transformations that decay to 
elements $\lambda(\pm\infty)$ in the stabilizers $G_a$ and $G_{a'}$, a 
constant element $g$ is necessarily trivial, hence the element 
$x=(\AA,0)$ is a manifold point in ${\cal B}^0(O_a,O_{a'})$, fixed by 
the $U(1)$-action.  
Thus, the operator ${\cal L}_{(\AA,\Psi,P,\rho)}$ is surjective, 
for a generic choice of the perturbations $P$ and $\rho$. 
This proves 
that in case (a) or case (b) the moduli space ${\cal 
M}(O_a,O_{a'})$ is cut out transversely by the equations and of 
dimension equal to the virtual dimension prescribed by the index 
theorem.  
The orientation of the moduli spaces ${\cal M}(O_a,O_{a'})$ is 
obtained as in Proposition \ref{orientation}, from an 
orientation of the homology groups  
\[ H^0_\delta(Y\times\RR)\oplus H^{2+}_\delta(Y\times\RR)\oplus 
H^1_\delta(Y\times\RR). \] 
 
\noindent $\diamond$ 
 
Notice that, since we are dealing with the case of homology spheres, 
we have only finitely many critical orbits $O_a$ and $O_{a'}$. Thus we 
can assume that there is a fixed compact set $Y\times [-T_0,T_1]$, 
with $T_0\geq 0$ and $T_1\geq 1$, such that transversality of all the 
moduli spaces of  
case (a) and (b) of Lemma \ref{noflow} can be achieved by a 
perturbation $\rho$ supported on $Y\times [-T_0,T_1]$, together with 
the perturbation $P$ of the curvature equation.  
Notice that outside $[-T_0,T_1]$ the metric is translation invariant  
and the equations coincide with the gradient flow equations. 
For simplicity of notation in the following we shall assume that 
$[-T_0,T_1]$ is just the interval $[0,1]$. 
 
We have the following result. 
 
\begin{corol} 
\label{reducibleflow} 
If either $O_b$ or $O_{b'}$ is a critical orbit with a free 
$U(1)$-action, then there are no reducible flow  
lines in the moduli space ${\cal M}(O_b,O_{b'})$, and the moduli space 
is generically empty if the virtual dimension is negative.  
 
Consider instead the case $O_a=\theta_0$ and $O_{a'}=\theta_1$ with 
$\mu(\theta_0)-\mu(\theta_1)<0$. In this case no solution with 
$\Psi\neq 0$ appears in the moduli space ${\cal M}(\theta_0,\theta_1)$.  
\end{corol} 
 
Both statements are consequences of Lemma 
\ref{surjT}. In fact, if there were any solution, 
then by argument of Lemma \ref{noflow} the linearization ought to be 
surjective for a generic choice of the perturbations, and this is 
incompatible with the value of the index. 
 
This result shows that reducible solutions in ${\cal M}(O_a,O_{a'})$ 
may appear only in the moduli spaces ${\cal M}(\theta_0, 
\theta_1)$, with $\mu(\theta_0)-\mu(\theta_1)<0$, which is in fact the 
only case not covered by the transversality result of Lemma \ref{surjT}. 
 
We can give a more precise description of this moduli spaces as 
follows. 
 
\begin{lem} 
\label{unique:reducible} 
Consider a choice of metrics and perturbations $(g_0,\nu_0)$ and 
$(g_1,\nu_1)$ in different chambers, with a negative  
spectral flow $SF(\partial^{g_t}_{\nu_t})<0$. In this case, the 
moduli space ${\cal M}(\theta_0, 
\theta_1)$ consists of a unique reducible class $x=[\AA,0]$. In 
particular, ${\cal M}(\theta_0,\theta_1)$ is non-empty, though the 
virtual dimension $\mu(\theta_0)-\mu(\theta_1)<0$ is negative. 
Thus, we have 
\[ \dim Coker({\cal L}_x)= | \mu(\theta_0)-\mu(\theta_1) |= 
-SF(\partial^{g_t}_{\nu_t}). \] 
\end{lem} 
 
\noindent\underline{Proof:} 
Lemma \ref{noflow} shows that ${\cal M}(\theta_0, 
\theta_1)$ can only contain reducible solutions. In particular, 
$x=[\AA,0]$ is obtained as the unique solution of the 
equations 
\[  d^+ \AA =\mu \ \hbox{ and } \ d^* \AA=0. \] 
Here $d^*$ is the operator on the four-manifold $Y\times \RR$, with 
$*_4$ the Hodge operator with respect to the metric $g_t + dt^2$, and 
the perturbation form $\mu$ is of the form  
$\mu=d\nu_t + *_3 d\nu_t \wedge dt$, with the $*_3$ operator with 
respect to the metric $g_t$.  
 
This solution of the four-dimensional Seiberg-Witten 
equations gives an element $x=[\AA,0]$ which is a fixed point of the 
$U(1)$-action on the configuration space 
${\cal B}^0 (\theta_0,\theta_1)$. Notice that, though the point
$x=[\AA,0]$ is a smooth point in ${\cal B}^0 (\theta_0,\theta_1)$
(fixed by the $U(1)$-action), it is a singular point in ${\cal M}(\theta_0, 
\theta_1)$, due to the lack of transversality, that is, the presence
of a non-trivial Cokernel.
 
\noindent $\diamond$ 
 
We see, from the previous discussion, that one first problem we 
encounter in extending the proof of the easy case of Theorem 
\ref{metrics} to the harder case is the compactification of the moduli  
spaces $\M(O_a,O_{a'})$. In fact, consider the case of the 
compactification of a moduli space ${\cal M}(\theta_0, O_{a'})$, under  
the assumption that $\mu(\theta_0)-\mu(\theta_1)<0$. 
 
For instance, in the case $\mu(\theta_0)-\mu(\theta_1)=-2$, 
by the result of Theorem \ref{seqcompact2}, we know that the  
product  
\beq 
\label{singboundary} 
{\cal M}(\theta_0,\theta_1)\times \hat{\cal M}(\theta_1,O_{a'})  
\eeq 
is a possible component of the ideal boundary. 
The counting of the virtual dimensions is correct as prescribed by 
Theorem \ref{seqcompact2} for components of co-dimension one in the 
ideal boundary. In fact, we have 
\[ \hbox{virtdim}({\cal M}(\theta_0,\theta_1)\times  
\hat{\cal M}(\theta_1,O_{a'})) = 
 \mu(\theta_0)-\mu(\theta_1)+\mu(\theta_1)-\mu(O_{a'})-1 \] 
\[ =\mu(\theta_0)-\mu(O_{a'})-1= \dim({\cal M}(\theta_0,O_{a'}))-1.\]  
However, from the result of Lemma \ref{unique:reducible}, we know that  
the actual dimension of ${\cal M}(\theta_0,\theta_1)$ is zero, since 
it consists of the unique solution $x=(\AA,0)$, hence the actual dimension 
of (\ref{singboundary}) is equal to 
$\mu(\theta_1)-\mu(O_{a'})$. This gives 
\[ \dim ({\cal M}(\theta_0,\theta_1)\times \hat{\cal 
M}(\theta_1,O_{a'})) = \mu(\theta_1)-\mu(O_{a'})-1 \] 
\[ =\mu(\theta_0)+2 -\mu(O_{a'})-1= 
\dim ({\cal M}(\theta_0,O_{a'}))+1. \] 
 
Clearly, this implies that, if an entire component of the form 
(\ref{singboundary}) actually   
appears in the compactification (that is, if they can be spliced 
together and deformed to actual solutions as described in Section 
4.2), then the resulting compactification ${\cal  
M}(\theta_0, O_{a'})^*$ does not have the structure of smooth manifold  
with corners described in Section 4.3. This could cause problems in 
applying the arguments based on Stokes' theorem in the proof of 
Theorem \ref{metrics}. Section 6.2 deals with this problem. The main 
result is stated in Theorem \ref{noglue}. 
 
We now proceed with general properties of the compactification of the 
moduli spaces $\M(O_a,O_{a'})$ and $\M^P(O_a,O_b)$ which we can derive  
in general, regardless of the behavior of the singular 
components of the ideal boundary. 
 
We have the following gluing theorem. 
 
\begin{thm} 
\label{glue:I} 
For $O_a$ and $O_{a'}$ irreducible, given a compact set 
$$ K\subset \hat{\cal M}(O_a,O_b)\times_{O_b} {\cal 
M}(O_b,O_{a'}) $$ 
or 
$$ K\subset {\cal M}(O_a,O_{b'})\times_{O_{b'}} \hat{\cal 
M}(O_{b'},O_{a'}), $$ 
there is a lower bound $T_0(K)>0$ and 
gluing maps 
\beq\label{glue:irred:I} \# : K\times [T,\infty)  
\to {\cal M}(O_a,O_{a'}), \eeq 
such that $\#_T$ is a smooth embedding for all $T\geq T_0(K)$. 
The gluing maps have the property that 
any sequence of solutions in ${\cal M}(O_a,O_{a'})$ that converges to 
an element in the boundary is eventually contained in the range of the gluing  
map. The map $\#$ is orientation preserving with respect to the 
orientation given by  
$$ (\hat{\cal M}(O_a,O_b)\times \RR^+)\times_{O_b} {\cal 
M}(O_b,O_{a'}) $$ 
and 
$$ {\cal M}(O_a,O_{b'})\times_{O_{b'}}( \hat{\cal 
M}(O_{b'},O_{a'}) \times \RR^+). $$ 
Similarly, for 
$$ K\subset \hat{\cal M}(O_a,\theta_0)\times  {\cal 
M}(\theta_0,O_{a'}) $$ 
or 
$$ K\subset {\cal M}(O_a,\theta_1)\times \hat{\cal 
M}(\theta_1,O_{a'}), $$ 
there is a $T_0(K)>0$ and gluing maps 
\beq \label{glue:red:I} \#: K \times [T,\infty) \to {\cal 
M}(O_a,O_{a'}),  \eeq 
which are smooth embeddings for $T\geq T_0(K)$, and with the same 
properties on the  
range and on the orientation. 
\end{thm} 
 
\noindent\underline{Proof:} We know from the results of Lemma 
\ref{noflow} that all the moduli spaces involved are smooth manifolds, 
cut out transversely by the equations. In particular, for any element 
$(x,y)$ in a fibered product as above, the cokernels $Coker({\cal 
L}_x)$ and $Coker({\cal L}_y)$  
are trivial. Thus, we can adapt the same argument used in 
Section 4.2 for the gluing construction, using the result of 
Proposition \ref{glue=fixedpoint}.
 
Namely, we proceed as follows. We show that all broken trajectories of  
the form ${\cal M}(O_a,O_{b'})\times_{O_b}\hat{\cal M}(O_{b'},O_{a'})$ are 
indeed contained in the boundary of ${\cal M}(O_a, O_{a'})$. The  
argument for the components of the form $\hat{\cal 
M}(O_a,O_b)\times_{O_b}{\cal M}(O_b,O_{a'})$ is analogous. 
We need to show that, given $x$ and $\hat y$ in 
${\cal M}(O_a,O_{b'})$ and $\hat{\cal M}(O_{b'},O_{a'})$ respectively, 
we can form a pre-glued approximate solution and perturb it to an 
actual solution as in Theorem \ref{equivgluing}. 
We have temporal gauge representatives 
\[ (A_1(t),\psi_1(t))\in {\cal A}_{k,\delta}(O_a,O_{b'}) \] 
and 
\[ (A_2(t),\psi_2(t))\in {\cal A}_{k,\delta}(O_{b'},O_{a'}). \] 
The solution $y$ is translation invariant, but the solution $x$ 
is not, since the metric is constant only outside the interval 
$[0,1]$. Thus we have to modify slightly the definition of the 
pre-gluing map. Choose slices  
\[ {\cal S}_{\Gamma_{ab'}}\subset {\cal A}_{k,\delta}(O_a,O_{b'}) \] 
and 
\[ {\cal S}_{\Gamma_{b'a'}}\subset {\cal A}_{k,\delta}(O_{b'},O_{a'}), 
\] 
determined by the elements 
\[ \Gamma_{ab'}\in {\cal 
A}_{k,\delta}((A_a,\psi_a),(A_{b'},\psi_{b'}))  \] 
and 
\[ \Gamma_{b'a'}\in {\cal 
A}_{k,\delta}((A_{b'},\psi_{b'}),(A_{a'},\psi_{a'})). \] 
We choose them  
in such a way that we have representatives $x$ and $y$  
in a ball of radius $r$ centered at $\Gamma_{ab'}$ and 
$\Gamma_{b'a'}$ respectively. Then, according to Lemma \ref{slicelemma}, 
there are gauge transformations $\lambda_1^+$ and $\lambda_2^-$ that 
conjugate $x$ and $y$ into $W^+_{\Gamma_{ab'}}$ and 
$W^-_{\Gamma_{b'a'}}$ respectively. For $t\geq T$, there is a gauge 
transformation $\lambda_{b'}\in {\cal G}(Y)$ such that    
we can write $\lambda_1^+ 
x(t)=\lambda_{b'}(A_{b'},\psi_{b'})+(\alpha_1(t),\phi_1(t))$ and  
$\lambda_2^- 
y(t)=\lambda_{b'}(A_{b'},\psi_{b'})+(\alpha_2(t),\phi_2(t))$, where  
$[A_{b'},\psi_{b'}]$ is a point on the orbit $O_{b'}$ and we have 
\[ \lim_{t\to\infty}(A_1(t),\psi_1(t))=\lambda_{b'}(A_{b'},\psi_{b'})= 
\lim_{t\to -\infty} (A_2(t),\psi_2(t)).  \]  
The (right) pre-glued approximate solution is given as 
follows: 
\beq\label{right:pre-glue} 
 x\#^0_T y =\left\{ \begin{array}{lr} \lambda_1^+ x & t\leq T-1\\ 
\lambda_{b'}(A_{b'},\psi_{b'}) +\rho_T^-(t)(\alpha_1(t),\phi_1(t))+& \\ 
\rho^+_T(t)(\alpha_2(t-2T),\phi_2(t-2T))& T-1\leq t\leq T+1 \\  
\lambda_2^- y^{-2T} & t\geq T+1. 
\end{array}\right. \eeq 
Here $\rho^\pm_T(t)$ are smooth cutoff functions with bounded 
derivative, such that $\rho^-_T(t)$ is equal 
to one for $t\leq T-1$ and to zero for $t\geq T$ and $\rho^+_T(t)$ is 
equal to zero for $t\leq T$ and to one for $t\geq T+1$. 
Lemma \ref{inducedgluing} shows that $\#_T^0$ descends to a well 
defined pre-gluing  
\[ \#_T^0: {\cal M}(O_a,O_{b'})\times_{O_{b'}}\hat{\cal 
M}(O_{b'},O_{a'})\to {\cal B}^0(O_a,O_{a'}). \]  
 
We are assuming that $O_a$ and $O_{a'}$ have a free $U(1)$ 
action, thus, according to Lemma \ref{noflow}, Corollary 
\ref{reducibleflow}, and Proposition \ref{relmorseind}, 
we can ensure that, for a generic choice of the metric and 
perturbation we have $Coker({\cal L}_x)=0$ and $Coker({\cal L}_y)=0$. 
Lemma \ref{surj2} and Proposition \ref{approxsurjective2} show that we 
have an isomorphism at the level of the actual kernel of the 
linearization of the pre-glued solution:   
\[ Ker ({\cal L}_{x\#^0_T y})\cong Ker({\cal L}_x)\times 
Ker({\cal L}_y). \] 
This implies that we have $Coker ({\cal L}_{x\#^0_T y})=0$ 
We proceed as in 
Proposition \ref{glue=fixedpoint}  
and we are able to perturb the approximate 
solutions to an actual solutions by the fixed point argument of Theorem  
\ref{equivgluing}. 
 
In the case of pre-gluing broken trajectories in the components 
\[ \hat{\cal M}(O_a,O_b)\times_{O_b}{\cal M}(O_b,O_{a'})\] 
we need to consider the (left) pre-gluing map 
\beq \label{left:pre-glue} 
 x\#^0_{-T} y=\left\{ \begin{array}{lr} \lambda_1^+ x^{2T} &t\leq -T-1\\ 
\lambda_b (A_b,\psi_b)+ \rho^+_{-T}(t)(\alpha_2(t),\phi_2(t)+ & \\ 
\rho_{-T}^-(t)(\alpha_1(t+2T),\phi_1(t+2T))& -T-1\leq t\leq -T+1\\ 
\lambda_2^- y &t\geq -T+1. 
\end{array}\right. \eeq 
Here $\rho^\pm_{-T}(t)$ are smooth cutoff functions with bounded 
derivative, such that $\rho^-_{-T}(t)$ is equal 
to one for $t\leq -T-1$ and to zero for $t\geq -T$ and $\rho^+_{-T}(t)$ is 
equal to zero for $t\leq -T$ and to one for $t\geq -T+1$.  
The rest of the argument is analogous. We have 
$$ Ker ({\cal L}_{x\#^0_T y})\cong Ker({\cal L}_x)\times 
Ker({\cal L}_y) $$ 
and $Coker ({\cal L}_{x\#^0_T y})=0$. We then proceed as in the 
previous case. 
 
\noindent $\diamond$ 
 
We have a similar result for the moduli spaces $\M^P(O_a,O_c)$. 
 
\begin{thm} 
\label{glue:H} 
For $O_a$ and $O_c$ free orbits, consider a compact set  
$$ K\subset \M^P(O_a,O_b)\times_{O_b} \hat\M(O_b,O_c) $$ 
or 
$$ K\subset \hat\M(O_a,O_b)\times_{O_b} \M^P(O_b,O_c). $$ 
There is a bound $T_0(K)$ and gluing maps 
$$ \# : K\times [T_0,\infty) \to \M^P(O_a,O_c) $$ 
that are smooth embeddings, compatible with the orientation induced by 
the product orientation of 
$$ \M^P(O_a,O_b)\times \RR \times \hat\M(O_b,O_c) $$ 
and 
$$ \hat\M(O_a,O_b)\times \RR \times \M^P(O_b,O_c). $$ 
Similarly, for a compact set 
$$ K\subset {\cal M}(O_a,O_{a'})\times_{O_{a'}} {\cal M}(O_{a'},O_c), $$ 
there is a $T_0(K)$ and an orientation preserving gluing map  
\beq \label{glue:gamma} 
 \# : K\subset {\cal M}(O_a,O_{a'})\times_{O_{a'}} {\cal 
M}(O_{a'},O_c) \to {\cal  
M}^P(O_a,O_c) \cap \{ \sigma= 1 \}, \eeq 
that is a smooth embedding. 
Since by construction the moduli space ${\cal M}^P(O_a,O_c)$ satisfies 
$$ {\cal M}^P(O_a,O_c) \cap \{ \sigma \in 
(1-\epsilon, 1] \} \cong  ({\cal M}^P(O_a,O_c) \cap \{ \sigma =1 \}) 
\times (1-\epsilon, 1], $$ 
this provides the necessary collar structure at the boundary 
components 
\beq \label{bound:only} 
{\cal M}(O_a,O_{a'})\times_{O_{a'}} {\cal M}(O_{a'},O_c).  \eeq 
Moreover, any sequence of solutions in 
$\M^P(O_a,O_c)$ converging to the boundary lies eventually in the 
range of the gluing map. 
\end{thm} 
 
\noindent\underline{Proof:} The proof of the first two cases is 
analogous to the case of flow lines analyzed in Section 4, therefore 
we concentrate on the case of the gluing map (\ref{glue:gamma}). 
Consider the moduli space  
$$ {\cal M}_{\sigma=1}^P(O_a,O_c)={\cal M}^P(O_a,O_c) \cap \{ \sigma 
=1 \}, $$ 
with the metric $\gamma_1=g_t \# \tilde g_t$, as in (\ref{gamma1}).  
We define the metric $\gamma_1(R)=g_t \#_R \tilde g_t$, for 
sufficiently large $R\geq R_0$, as 
\beq \label{gamma1R} 
\gamma_1(R)=\left\{ \begin{array}{lr} g_0+dt^2 & t< -2-R\\ 
g_{t+2}+dt^2 & t\in [-2-R,-1-R]\\ 
g_1+dt^2 & t\in [-1-R,1+R] \\ 
\tilde g_{2-t}+dt^2 & t\in [1+R,2+R]\\ 
g_0+dt^2 & t>2+R. 
\end{array}\right. \eeq 
Correspondingly, we have metrics $\gamma_\sigma(R)$, for 
$\sigma\in [0,1]$. We consider the moduli spaces ${\cal M}^P(O_a,O_c)$  
as before, defined with respect to the metrics $\gamma_\sigma(R)$. All  
the properties of these moduli spaces discussed so far, are 
independent of the choice of a fixed $R$. However, to describe in 
detail the appearance of strata of the form (\ref{bound:only}) in the 
compactification, and their gluing properties, we need to work with a 
large $R$. This correspond to stretching a long cylinder $Y\times 
[-1-R,1+R]$ with the constant metric  
$g_1 + dt^2$ inside the manifold $Y\times \RR$. The proof is obtained 
via the following argument. 
 
\begin{prop} 
\label{glue:gamma:proof} 
Consider a solution $x=[\AA,\Psi]$  
in ${\cal M}_{\sigma=1}^P(O_a,O_c)$ 
on the manifold $(Y\times \RR, \gamma_1(R))$. For every $\epsilon >0$, 
there exists $R_0$, such that, if $R \geq R_0$, the restriction of $x$ 
to $Y\times [-2,2]$ is sufficiently close to a solution 
$[A_{a'},\psi_{a'}] \in O_{a'}$, for some critical orbit 
$O_{a'}$. Moreover, given a pair of solutions $(x,y)$ in the fibered 
product (\ref{bound:only}), for $R$ sufficiently large, we can 
consider the pre-glued solution obtained by splicing together the 
solutions $x$ and $y$ cut off by functions $\rho_-$ and $\rho_+$ 
supported in $Y\times (-\infty, -1]$ and $Y\times [1,\infty)$, 
respectively. For $R$ sufficiently large, these approximate solutions 
can be deformed to actual solutions. Thus, we obtain the gluing map 
(\ref{glue:gamma}) with the desired properties.  
\end{prop} 
 
\noindent\underline{Proof of Proposition \ref{glue:gamma:proof}:} the 
proof of the first statement follows closely the argument of 
Proposition 8 of \cite{KM}. The gluing argument again follows very 
closely the argument illustrated in Section 4.2, Theorem 
\ref{equivgluing} in the case of the moduli spaces of flow 
lines. Again, the key feature is the fact that the cokernels 
$Coker({\cal L}_x)$ and $Coker({\cal L}_y)$ 
vanish, hence the gluing argument presents no surprise with respect to  
the case of gluing flow lines. 
 
\noindent $\diamond$ 
 
All the remaining statements of Theorem \ref{glue:H} are proved by the  
same techniques developed in Section 4. 
 
\noindent $\diamond$ 
 
The main problem now is to extend these gluing theorems to the case of  
strata involving a contribution from the singular moduli space 
$\M(\theta_0,\theta_1)$.  
The strategy we develop in the rest of Section 6.1 and in 
Section 6.2 is to identify precisely the contribution of these singular 
components (\ref{singboundary}) to the compactification of the moduli 
spaces ${\cal M}(\theta_0, O_{a'})^*$ and to the analogous case for 
$\M^P(O_a,O_b)^*$. The main technique we employ is   
the obstruction bundle for the gluing construction and a direct 
analysis of the zero set of the canonical obstruction section. 
We aim at showing that all the 
components of the actual boundary in the compactification ${\cal  
M}(O_a, O_{a'})^*$ behave according to the picture described in 
Section 4.3, hence the arguments of the easy case of Theorem 
\ref{metrics} involving the properties of the compactifications 
moduli spaces ${\cal M}(O_a, O_{a'})^*$ and $\M^P(O_a,O_b)^*$ can be 
extended to the case of different chambers as well.  
 
\begin{rem} 
Recall that, in the proof of Theorem 
\ref{metrics}, we need only consider
the case of two sufficiently close metrics and perturbations 
$(g_0,\nu_0)$ and $(g_1,\nu_1)$ in two different chambers.
Thus, we can assume that the path $(g_t,\nu_t)$ 
satisfies $SF(\partial^{g_t}_{\nu_t})=-2$ and the path $(\tilde g_t, 
\tilde \nu_t)$ satisfies $SF(\partial^{\tilde g_t}_{\tilde 
\nu_t})=+2$. Thus, the 
compactified moduli spaces ${\cal M}(O_{a'},O_a)^*$  used in the 
definition of the map $J$ and in the proof of $DJ -JD =0$ all satisfy  
the properties of Section 4.2 and 4.3, by simply rephrasing the same 
arguments. The only problem then remains in the analysis of the moduli 
spaces needed in the definition of $I$, in the 
proof of $ID-DI=0$ and in the construction of the chain homotopy. 
\label{simplify:2} 
\end{rem}

\subsection{Obstruction bundle and gluing theorems} 
 
The obstruction bundle technique we discuss in this Section will both 
complete the proof of the unobstructed gluing theorems, as stated in 
Proposition \ref{glue=fixedpoint}, and identify the precise boundary 
structure in the gluing with obstructions.

Let ${\cal U}(O_a,O_b)$ be the image of ${\cal M}(O_a,O_c)\times_{O_c} {\cal 
M}(O_c,O_b)$ under the pre-gluing map $\#^0$, for a large gluing 
parameter $T>>0$.  
In the case of a translation invariant metric, analyzed in Section 4, 
the statement of the gluing theorem \ref{equivgluing} is equivalent to the 
statement that the subset  
\[ {\cal U}(O_a,O_b)\subset {\cal B}^0_{k,\delta}(O_a,O_b) \] 
can be deformed to ${\cal M}(O_a,O_b)$, whenever ${\cal U}(O_a,O_b)$ 
is obtained by pasting together solutions along $O_b$ with 
$\mu(O_a)>\mu(O_c)>\mu(O_b)$. It is this equivalent statement which we  
are going to verify in the  proof of Proposition \ref{glue=fixedpoint}  
given in this Section, following Lemma \ref{solveqbeta}.  
 
However, if we consider the moduli spaces ${\cal M}(O_{a},O_{a'})$ with 
the metric $\tilde g_t$ that varies along the cylinder, there might be 
topological obstructions that make it impossible to deform the  
set of pre-glued trajectories to actual solutions of (\ref{3SW1P'}) and 
(\ref{3SW2P'}). In other words, given a point 
$[x\#^0_T y]\in {\cal U}(O_{a},O_{a'})$ there might be an obstruction 
to pushing that point onto ${\cal M}(O_{a},O_{a'})$. Typically 
obstructions may arise either because of the presence of non-trivial 
cokernels of the linearizations ${\cal L}_x$ and ${\cal L}_y$, or 
because Sard's theorem is not available. In the problem we are 
considering the obstruction originates in the presence of the 
non-trivial cokernel of Lemma \ref{unique:reducible}. Similar topological 
obstructions to gluing solutions  
arise in other gauge theory problems  \cite{Ta2}, \cite{Ta}. Following 
Taubes' construction of the obstruction bundle, it is possible to 
describe the obstruction completely in terms of the vanishing of a 
canonical section. 
 
We maintain the notation (\ref{f}) as in the proof of the gluing 
theorem: for an element $x\#^0_Ty$, with $x=(\AA_1,\Psi_1)$ and 
$y=(\AA_2,\Psi_2)$, we write 
\[ f(x\#^0_T y)=\left\{\begin{array}{l} 
D_{\AA_1\#_T^0\AA_2}(\Psi_1\#^0_T\Psi_2)+ \rho (\Psi_1\#^0_T\Psi_2) \\ 
F_{\AA_1\#_T^0\AA_2}^+ 
-(\Psi_1\#^0_T\Psi_2)\cdot\overline{(\Psi_1\#^0_T\Psi_2)}  
-i\mu -P_{(\AA_1\#_T^0\AA_2,\Psi_1\#^0_T\Psi_2)} \\ 
G^*_{\Gamma_{aa'}}(\AA_1\#_T^0\AA_2,\Psi_1\#^0_T\Psi_2), \end{array}\right.\] 
for a fixed choice of a slice ${\cal S}_{\Gamma_{aa'}}$ in 
${\cal A}(O_a,O_{a'})$, such that $x\#^0_T y$ is in a ball of radius $r$ 
in ${\cal A}(O_a,O_{a'})$ centered at $\Gamma_{aa'}$. 
If the pre-glued element $x\#^0_T y$ can be deformed to an actual 
solution, then there exists a small $(\alpha,\Phi)\in Ker({\cal 
L}_{x\#^0_T y})^\perp$ such that we have  
\[ [x\#^0_T y+ (\alpha,\Phi)]\in {\cal M}(O_a,O_{a'}), \] 
that is, the following equation is satisfied 
\beq 
f(x\#^0_T y)+{\cal L}_{x\#^0_T y}(\alpha,\Phi) +{\cal N}_{x\#^0_T 
y}(\alpha,\Phi)=0.  
\label{eqalpha} 
\eeq 
Here ${\cal N}$ denotes the nonlinear term in the equation, as in 
Proposition \ref{glue=fixedpoint}. 
\[ {\cal N}_{x\#^0_T y}(\alpha,\Phi)=(\Phi\cdot\bar\Phi+{\cal 
  N}P_{x\#^0_T y}(\alpha,\Phi),\alpha\cdot\Phi), \] 
where the perturbation $P$ of equation (\ref{4SW2P't}) is written as 
sum of a linear and a non-linear term, $P={\cal D}P+{\cal N}P$.  
 
The presence of a non-trival approximate kernel of ${\cal L}^*_{x\#^0_T 
y}$ can generate obstructions to solving equation (\ref{eqalpha}) for 
$(\alpha,\Phi)$. In fact, the hypothesis that ${\cal L}_{x\#^0_T y}$ 
has a trivial approximate cokernel is essential in the proof 
of the gluing theorem (see Lemma \ref{surj2} and Proposition 
\ref{approxsurjective2}):  
the same argument cannot be extended to a case with a non-trivial 
cokernel.  
 
We follow \cite{Ta} and introduce open sets ${\cal U}(\mu)$ of 
elements $[x\#^0_T y]\in {\cal U}(O_{a},O_{a'})$ such that $\mu>0$ is 
not an eigenvalue of the operator ${\cal L}_{x\#^0_T y}{\cal 
L}^*_{x\#^0_T y}$ acting on $L^2_{0,\delta(T)}$ connections and sections.  
There are projection maps $\Pi(\mu, x\#^0_T y)$ onto the span of the 
eigenvectors of ${\cal L}_{x\#^0_T y}{\cal 
L}^*_{x\#^0_T y}$ with eigenvalue smaller than $\mu$. These are smooth 
maps of $[x\#^0_T y]\in {\cal U}(\mu)$. 
 
We have the following result, which gives the necessary eigenvalue 
splitting for the operators ${\cal L}_{x\#^0_T y}{\cal L}_{x\#^0_T 
y}^*$ and ${\cal L}_{x\#^0_T y}^*{\cal L}_{x\#^0_T y}$. 
 
\begin{lem} 
There exists a $\mu_0=\mu(T_0) >0$, such that, for all $T\geq T_0$, 
all the small  
eigenvalues $\mu_T$ of ${\cal L}_{x\#^0_T y}^*{\cal L}_{x\#^0_T y}$ 
(or ${\cal L}_{x\#^0_T y}{\cal L}_{x\#^0_T y}^*$), satisfy $\mu_T < \mu_0$ 
and  $\mu_T \to 0$, and all the other eigenvalues are bounded 
below by $\mu_0$.  
\label{spectral:lemma} 
\end{lem} 
 
\noindent\underline{Proof:} Recall that, from Lemma 
\ref{split:convergence} we know that the number of independent 
eigenvectors of ${\cal L}_{x\#^0_T y}^*{\cal L}_{x\#^0_T y}$ 
or ${\cal L}_{x\#^0_T y}{\cal L}_{x\#^0_T y}^*$ with small 
eigenvalue $\mu_T \to 0$ is at most  
$$\dim Ker({\cal L}_x) + \dim 
Ker({\cal L}_y)$$  
or at most  
$$\dim Coker({\cal L}_x) + \dim 
Coker({\cal L}_y), $$ 
respectively. Moreover, 
by Lemma \ref{apprkercoker}, we know that the spans of the eigenvectors  
of the small eigenvalues have exactly these dimensions, cf. Proposition  
\ref{apprkercoker2}.  
So there exists a $\mu_0$ such that, for $T\geq T_0$, there is 
exactly this number  
of independent eigenvectors with eigenvalues $\mu_T < \mu_0$ and 
$\mu_T\to 0$. In fact, we also know more precise estimates on the rate  
of decay of these small eigenvalues, as specified in Lemma 
\ref{apprkercoker}.  
It remains to be shown that such $\mu_0$ can be chosen 
so that all the remaining eigenvalues are bounded below by $\mu_0$. In  
other words, we have to show that we cannot have other sequences of 
eigenvalues $\mu_T(k)$ satisfying $\mu_T(k)\to \mu_k >0$, as 
$T\to\infty$, but with the limits $\mu_k \to 0$ as $k\to \infty$.  
This follows from Lemma \ref{surj2}. In fact, we have proved that 
the operator ${\cal L}_{x\#^0_T y}^*{\cal L}_{x\#^0_T 
y}$ (or ${\cal L}_{x\#^0_T y}{\cal L}_{x\#^0_T y}^*$) is uniformly 
invertible on the orthogonal complement of the small eigenvalues 
eigenspace.  
Notice that Lemma \ref{apprkercoker} and Lemma \ref{surj2} give 
norms for the operator ${\cal L}_{x\#^0_T y}$ (and by similar 
arguments for the operator ${\cal L}_{x\#^0_T y}^*$) which are 
uniformly bounded in $T$, when the operators are considered as acting 
between Sobolev spaces with the same $\delta$ weight on source and 
target space (or with the same rescaled weight $\delta(T)$ on both 
source and target space). In other words we have 
$$ \| {\cal L}_{x\#^0_T y}^*{\cal L}_{x\#^0_T y} \xi 
\|_{L^2_{1,\delta}} \geq C \| \xi \|_{L^2_{1,\delta}}, $$ 
for 
$$ \xi \in (F_{\#_T}(Ker({\cal L}_x) \times Ker({\cal L}_y)))^\perp 
$$ 
inside  
$$ {\cal T}_{1,\delta}((e_a^+)^{-1}(x_a))\subset {\cal 
T}_{1,\delta}({\cal S}_{\Gamma_{aa'}}), $$   
and 
$$ \| {\cal L}_{x\#^0_T y}{\cal L}_{x\#^0_T y}^* \xi 
\|_{L^2_{0,\delta}} \geq C \| \xi \|_{L^2_{0,\delta}}, $$ 
for 
$$ \xi \in (F_{\#_T}(Coker({\cal L}_x) \times Coker({\cal L}_y)))^\perp 
$$ 
inside ${\cal T}_{0,\delta}({\cal S}_{\Gamma_{aa'}})$. With this 
choices of norms, the constants are all independent of $T$.  
Here $x_a \in O_a$ is the asymptotic value as $t\to -\infty$  of the 
solution $x\in {\cal M}(O_a,O_{a'})$, and $e^+_a: {\cal S}_{\Gamma_{aa'}} 
\to O_a$ is the endpoint map, for the fixed choice of slices. 
 
\noindent $\diamond$ 
 
It is important to notice that   
the difference between using the rescaled $\delta(T)$--norms or 
the original $\delta$--norms in the previous Lemma is reflected in the 
different norm bound on the gluing map $\tilde F_{\#_T}$, as explained 
in the proof of Theorem \ref{equivgluing}.

The kind of analysis of the behavior of small eigenvalues described 
in Lemma \ref{spectral:lemma} has been used 
in similar context in  \cite{Ta2}, \cite{Ta}, with a fuller account in  
\cite{Mrowka}, and in \cite{CLM}. A similar analysis has also been 
used, more recently, in \cite{CMW}, \cite{Fee-Len}.

By Lemma \ref{spectral:lemma}, we now know that there exist $T_0>>0$ and 
$\mu(T_0)>0$, such that, for $T>T_0$ and $0< \mu< \mu(T_0)$, we can 
identify the projection map $\Pi(\mu, x\#^0_T y)$ with projection on 
the approximate cokernel $ApprCoker({\cal L}_{x\#^0_T y})$, as defined 
in Definition \ref{defapprker}. 
 
Consider a compact set $K\subset {\cal U}(O_{a},O_{a'})$. Suppose it can be 
covered with one open set ${\cal U}(\mu)$. If not, we can paste 
together the projections corresponding  
to the finitely many open sets that cover $K$, using a suitable partition 
of unity as in \cite{Ta}, pg. 193 and slightly modifying the argument 
that follows in \cite{Ta}. 
 
Since the element $(\alpha,\Phi)$ is in $Ker({\cal L}_{x\#^0_T 
y})^\perp$, we have 
\[ (\alpha,\Phi)={\cal L}^*_{x\#^0_T y}(\beta,\xi). \] 
Suppose that $[x\#^0_T y+(\alpha,\Phi)]$ is a solution, for a 
non-trivial element $(\alpha,\Phi)$. For $\mu$ small enough, we can 
assume that $(\beta,\xi)$ satisfies  
\[ (\beta,\xi)\in Ker(\Pi(\mu,x\#^0_T y)). \] 
Then $(\beta,\xi)$ solves the equations 
\beq 
\label{eqbeta} 
{\cal L}_{x\#^0_T y}{\cal L}^*_{x\#^0_T y}(\beta,\xi)+ 
(1-\Pi(\mu,x\#^0_T y)) \left( {\cal N}_{x\#^0_T y}{\cal L}^*_{x\#^0_T 
y}(\beta,\xi) + f(x\#^0_T y) \right)=0. 
\eeq 
and 
\beq 
\label{eqbeta2} 
\Pi(\mu,x\#^0_T y) \left( {\cal 
N}_{x\#^0_T y}{\cal L}^*_{x\#^0_T y}(\beta,\xi)+ f(x\#^0_T 
y)\right)=0. 
\eeq 
Lemma \ref{solveqbeta} below ensures that it is always possible to 
find a small solution of equation (\ref{eqbeta}), hence the problem of 
deforming an approximate solution to an actual solution depends on 
whether equation (\ref{eqbeta2}) can also be solved. The latter has a 
geometric interpretation as the section of an obstruction bundle, as 
described in the following (see \cite{D2} Section 4, \cite{Ta2} 
Section 3 and 5, \cite{Ta} Section 6).  
 
\begin{lem} 
There is an $\epsilon>0$ and a constant $C>0$ such that, for any $\mu$ with 
$C\epsilon < \mu < \mu_0$, with $\mu_0$ the least eigenvalue of  
$H_{x\#^0_T y}^0$ on the complement of $ApprCoker({\cal L}_{x\#^0_T y})$,  
the equation (\ref{eqbeta}) has a unique solution  
\[ (\beta,\xi)\in Ker(\Pi(\mu,x\#^0_T y)) \]  
with 
$\| (\beta,\xi)\|_{L^2_{2,\delta}} \leq \epsilon$,  
provided the error term $\| f(x \#_T^0 y) \|$ satisfies an estimate 
$$ \| f(x \#_T^0 y) \| \leq \frac{\mu}{2\epsilon}. $$ 
\label{solveqbeta} 
\end{lem} 
 
\noindent\underline{Proof:} 
For $\mu>0$ and any $(\AA,\Psi)$, the operator 
\[ H^0_{(\AA,\Psi)}={\cal L}_{(\AA,\Psi)}{\cal L}^*_{(\AA,\Psi)} \] 
has a bounded inverse, when restricted to the image of 
$(1-\Pi(\mu,(\AA,\Psi))$, in the space of $L^2_\delta$ connections and  
sections, as discussed in Lemma \ref{surj2}.  
 
Thus, the equation (\ref{eqbeta}) can be rephrased as a fixed point problem 
\beq (\beta,\xi)=-{H^0_{x\#^0_T y}}^{-1} (1-\Pi(\mu,x\#^0_T y)) ({\cal 
N}_{x\#^0_T y}{\cal L}^*_{x\#^0_T y}(\beta,\xi)+f(x\#^0_T 
y)). \label{contraction:map} \eeq 
 
We need to prove that the right hand side is a contraction. For all 
$(\beta,\xi)$ and $(\beta',\xi')$, we have an 
estimate  
\[ \| {H^0_{x\#^0_T y}}^{-1} (1-\Pi(\mu,x\#^0_T y))\left( {\cal N}{\cal 
L}^* (\beta,\xi)-{\cal N}{\cal L}^*(\beta',\xi') \right) 
\|_{L^2_{\delta}} \leq \] 
\[ \frac{1}{\mu} \| (1-\Pi) ({\cal N}{\cal 
L}^* (\beta,\xi)-{\cal N}{\cal L}^*(\beta',\xi') ) 
\|_{L^2_{\delta}}  \leq  \] 
\[ \frac{1}{\mu} \| {\cal N}{\cal 
L}^* (\beta,\xi)-{\cal N}{\cal L}^*(\beta',\xi') \|_{L^2_{\delta}}. \] 
The quadratic form 
$(\phi\cdot\bar\phi,\alpha\cdot\phi)$ satisfies an 
estimate of the form  
\[ \| N(x)-N(y) \|\leq C(\|x\|+\|y\| )\|x-y\| \] 
as required for the contraction method \ref{fixedpoint}. 
The perturbation term also satisfies a similar estimate for  
large enough $T$ because of the assumptions on the perturbation space 
${\cal P}$. Thus we can improve the last estimate to get 
\[ \frac{1}{\mu} \| {\cal N}{\cal 
L}^* (\beta,\xi)-{\cal N}{\cal L}^*(\beta',\xi') \|_{L^2_{\delta}} \leq \] 
\[ \frac{\tilde C}{\mu} (\| {\cal L}^*(\beta,\xi) \| + 
\| {\cal L}^*(\beta',\xi') \|) \| {\cal L}^* 
((\beta,\xi)-(\beta',\xi')) \|_{L^2_{1,\delta}} \leq \] 
\[ \frac{C}{\mu} (\| (\beta,\xi) \|+\| (\beta',\xi')\|) \| 
(\beta,\xi)- (\beta',\xi') \|_{L^2_{2,\delta}}. \] 
 
Thus, we have 
$$  \| {H^0_{x\#^0_T y}}^{-1} (1-\Pi(\mu,x\#^0_T y))\left( {\cal 
N}_{x\#^0_T y} {\cal L}^*_{x\#^0_T y}  (\beta,\xi) + f(x\#^0_T y) 
\right) \| \leq $$ 
$$ \frac{C}{\mu} \| (\beta,\xi) \| + \frac{1}{\mu} \| f(x\#^0_T y) \|. $$ 
 
Suppose given $\epsilon$ and $\mu$ satisfying  
$$ \max\{ 2C\epsilon, \epsilon \} <\mu < \mu_0, $$ 
and assume the the error term also satisfies 
$$ \| f(x \#_T^0 y) \| \leq \frac{\mu}{2\epsilon}. $$ 
Then the map is a contraction on 
the ball $\| (\beta,\xi) \|_{L^2_{\delta}}\leq \epsilon$.  
 
\noindent $\diamond$

\begin{rem} \label{isom:T} 
Lemma \ref{solveqbeta} relies on the identification of the image of  
$1-\Pi(\mu,x\#^0_T y)$, that is, the complement of the small eigenvalues 
eigenspace for the Laplacian $H^0_{x\#^0_T y}$, with the space 
$$  F_{\#_T}(Coker({\cal L}_x) \times Coker({\cal  
L}_y)). $$  
This identification is obtained from Lemma 
\ref{surj2} together with Remark \ref{split:convergence:coker}. In 
particular, recall that, by Lemma \ref{surj2} and Remark 
\ref{split:convergence:coker}, the isomorphism   
$$ F_{\#_T}: Coker{\cal L}_x \times Coker{\cal 
L}_y \stackrel{\cong}{\to} Image(1-\Pi(\mu,x\#^0_T y)) $$ 
has norm bounded uniformly in $T$ if we use weighted norms with 
rescaled weight $\delta(T)$ on the target space and the $\delta$-norms  
on the source spaces. It has norm bounded by $C e^{-\delta T}$ if the 
target space is also endowed with the $\delta$--norm. 
\end{rem} 
 
\begin{rem} \label{free:glue:case} 
In the case of the gluing theorem for flow lines, with approximate 
solutions $x\#^0_T y$ in ${\cal  
U}(O_a,O_c)$, obtained from elements $x\in \hat{\cal M}(O_a,O_b)$ and 
$y\in \hat{\cal M}(O_b,O_c)$, it is 
possible to choose the constant 
$\epsilon$ in Lemma \ref{solveqbeta} independent of $T\geq T_0$. 
In fact, in this case it is possible to choose $\mu \leq \mu_0$, the 
least eigenvalue of $H^0_{x\#^0_T y}$, which in this case is 
independent of $T$, cf. the spectral decomposition of Lemma 
\ref{spectral:lemma}. The constant $C$ is  
also independent of $T\geq T_0$, cf. Lemma \ref{spectral:lemma}. Thus,  
in this case, the projection $\Pi(\mu,x\#^0_T y)$ is trivial, hence 
the gluing result reduces to the following argument that proves 
Proposition \ref{glue=fixedpoint}.  
\end{rem} 
 
A proof of Proposition \ref{glue=fixedpoint} follows from Lemma 
\ref{solveqbeta} in the case the 
cylindrical metric on $Y\times \RR$, and under the assumptions that  
$Coker({\cal L}_x)=0$ and $Coker({\cal L}_y)=0$. 
 
\noindent\underline{Proof of Proposition \ref{glue=fixedpoint}:}  
Let us consider the space ${\cal T}_{1,\delta}({\cal 
S}_{\Gamma_{ac}})$, endowed with the $L^2_{1,\delta}$--norm.  
The initial condition of Remark \ref{fixedpoint} is 
provided by the exponential decay. In fact, we have 
\[ \| f(\AA_1\#_T^0 \AA_2,\Psi_1\#_T^0\Psi_2)\|_{L^2_{1,\delta}} \leq \] 
\[ \| F_{\#_T} \| \cdot ( \|(\alpha_1,\phi_1)|_{Y\times [-1+T,T]} 
\|_{L^2_{1,\delta}} + \|(\alpha_2,\phi_2) |_{Y\times [-T,-T+1]} 
\|_{L^2_{1,\delta}}). \]  
The operator norm $\| F_{\#_T} \|$, for $F_{\#_T}$ acting on the 
$L^2_{1,\delta}$ spaces, is bounded by $C e^{-\delta T}$, as in 
(\ref{estFT}). The terms  
$$\|(\alpha_1,\phi_1)|_{Y\times [-1+T,T]}\|_{L^2_{1,\delta}} \ \ 
\hbox{and} \ \ \|(\alpha_2,\phi_2) |_{Y\times 
[-T,-T+1]}\|_{L^2_{1,\delta}} $$ are bounded by a constant, uniformly 
in $T\geq T_0$, because of the exponential decay of 
$(A_1(t),\psi_1(t))$ and $(A_2(t),\psi_2(t))$ to the endpoints. 
Thus, for all $T\geq T_0$, we have obtained an estimate 
\beq 
\| f(\AA_1\#_T^0 \AA_2,\Psi_1\#_T^0\Psi_2)\|_{1,\delta}\leq C 
e^{-\delta T}.  
\label{gluedecay} 
\eeq 
 
\noindent $\diamond$ 
 
This result can be adapted to prove all the unobstructed gluing 
results, such as the ones stated in Theorems \ref{glue:I} and 
\ref{glue:H}. Now we turn to the more interesting case with 
obstructions.  
 
Using the result of Lemma \ref{solveqbeta}, we can construct the 
obstruction bundle and the canonical section, in the case with 
non-trivial approximate cokernel (i.e. non-trivial small eigenvalue 
eigenspace of $H^0_{x\#^0_T y}$).  
 
\begin{prop} 
There is a local bundle ${\cal E}$ over ${\cal U}(O_a,O_{a'})$ and a 
section $s_{\mu}$ such \index{$s_{\mu}$} 
that the following property holds. A point $x\#^0_T y$ in ${\cal 
U}(O_a,O_{a'})$  
can be deformed to a solution of the flow equations (\ref{3SW1P'}) and 
(\ref{3SW2P'}) iff  
\[ x\#^0_T y\in s_{\mu}^{-1}(0). \] 
\label{obstruction} 
\end{prop} 
 
\noindent\underline{Proof:} 
Consider the section 
\[ s_{\mu} (x\#^0_T y)=\Pi(\mu,x\#^0_T y) \left( {\cal 
N}_{x\#^0_T y}{\cal L}^*_{x\#^0_T y}(\beta,\xi)+ f(x\#^0_T y)\right) \] 
of the vector bundle over ${\cal U}(O_{a},O_{a'})$ with fiber 
$ApprCoker({\cal L}_{x\#^0_T y})$. 
Here $(\beta,\xi)$ is the unique solution of (\ref{eqbeta}).  
 
Since $(\beta,\xi)$ satisfies (\ref{eqbeta}), we have 
\[ [x\#^0_T y +{\cal L}^*_{x\#^0_T y}(\beta,\xi)]\in {\cal M}(O_a,O_{a'}) 
\] 
iff  
\[ s_{\mu}(x\#^0_T y)=0. \] 
 
\noindent $\diamond$

In the case of moduli spaces of solutions of the 
Seiberg--Witten equations on $Y\times \RR$ with the changing metric 
$g_t + dt^2$, as in the construction of the map $I$ and the subsequent  
arguments, 
Proposition \ref{obstruction} characterizes which points of ${\cal 
U}(O_a,O_{a'})$ can be deformed to actual solutions.  
We need a more explicit description of this obstruction bundle in the 
specific problem at hand, namely in the case of approximate solutions 
${\cal U}(\theta_0,O_{a'})$, with the pre-gluing map 
$$ \#^0 : {\cal M}(\theta_0, \theta_1)\times \hat{\cal 
M}(\theta_1,O_{a'})\times [T_0,\infty) \to {\cal U}(\theta_0,O_{a'}) 
$$ 
or ${\cal U}(O_a,\theta_1)$ with the pre-gluing map 
$$ \#^0 : \hat{\cal M}(O_a,\theta_0)\times {\cal M}(\theta_0, 
\theta_1)\times [T_0,\infty) \to {\cal U}(O_a,\theta_1). $$

We first need to derive the equivalent of Remark
\ref{free:glue:case}. Namely, we need a suitable choice of  
$\mu$ and $\epsilon$, so that the map in the fixed point problem 
(\ref{contraction:map}) is a contraction. This choice depends upon 
estimating the error term $\| f(x\#_T^0 y) \|$. We have the following 
result.  
 
\begin{lem} 
\label{singular:glue:case} 
Let $x=[\AA,0]$ be the unique reducible in ${\cal 
M}(\theta_0,\theta_1)$, with 
$\mu(\theta_0)-\mu(\theta_1)=-2$. Consider (right)pre-glued solutions 
$x\#_{T}^0 y$ in ${\cal U}(\theta_0, O_{a'})$, with $y\in \hat{\cal 
M}(\theta_1,O_{a'})$. We have an estimate on the error term $\| 
f(x\#_{T}^0 y) \|$ of the form 
$$ \| f(x\#_{T}^0 y) \|_{L^2_{1,\delta}} \leq \tilde C, $$ 
for some constant independent of $T\geq T_0$, 
satisfying $\tilde C < \mu_0 $. 
Moreover, choose $\tilde C <\mu< \mu_0$, independent of 
$T$, with $\tilde C$ as above, and with $\mu_0$ the  
lower bound on the  
spectrum of the Laplacian $H^0_{x\#_{T}^0 y}$ acting on the complement of 
the space of small eigenvalues, as in Lemma \ref{spectral:lemma}. Then  
for $\max\{ 2C\epsilon, \epsilon \} < \mu$, we obtain that the map of 
(\ref{contraction:map}) is a contraction, and the fixed point problem 
has a unique solution as in Lemma \ref{solveqbeta}. 
The case of (right)pre-gluing maps is completely analogous. 
\end{lem} 
 
\noindent\underline{Proof:}  
Recall that $(\AA,0)$ is the unique solution to $d^+ \AA =\mu$ and $d^* 
\AA=0$. Here we have the 4-dimensional $d^*$, with $*_4$ with respect 
to the metric $g_t + dt^2$, and $\mu$ of the form $\mu=d\nu_t + *_3 
d\nu_t \wedge dt$, with $*_3$ with respect to the metric $g_t$.  
We have  
$$ \| f(x\#_T^0 y) \| \leq C \| F_{\#_T} \| \cdot $$ 
$$ ( \| (\alpha_1,\phi_1) \|_{L^2_{1,\delta}(Y\times [-T-1,-T+1])} + \| 
\AA-\AA_1 \|_{L^2_{1,\delta}(Y\times [T-1,T+1])}). $$ 
Here $\AA_1$ is the solution of $d^+\AA=\mu_1$ and $d^*\AA=0$, with 
$\mu_1=d\nu_1 +*_3 d\nu_1\wedge dt$. Here we have the 4-dimensional 
$d^*$, with $*_4$ with respect to the metric $g_1 + dt^2$ and 
$*_3$ with respect to the metric $g_1$. In the rest of the proof we 
write $d^{*_{g_t}}$ and $d^{*_{g_1}}$ and $*_{g_t}$ and $*_{g_1}$ for the 
4-dimensional $*$-operators with respect to the two different metrics. 
(We hope this will not cause confusion with our previous notation 
$*_4$ and $*_3$.) 
The norm $\| F_{\#_T} \|$ is exponentially small $\| F_{\#_T} \|\leq C 
e^{-\delta T}$. 
The term $\| (\alpha_1,\phi_1) \|_{L^2_{1,\delta}(Y\times [-T-1,-T+1])}$ is 
bounded by a constant because of the exponential decay of 
$(\alpha_1,\phi_1)$, where we write $y=\lambda_1 
(A_1+\nu_1,0)+(\alpha_1,\phi_1)$, as $t\to -\infty$. The estimate for 
the remaining term is obtained as follows. We write 
$$ \AA-\AA_1=d^{*_{g_t}} (\mu_t -\mu_1) + (d^{*_{g_t}}-d^{*_{g_1}}) 
\mu_1. $$   
Thus, we have 
$$ \|  \AA-\AA_1 \|_{L^2_1(Y\times [T-1,T+1])} \leq \| d^{*_{g_t}} 
(\mu_t -\mu_1) \|_{L^2_1(Y\times [T-1,T+1])}  $$  
$$ + \| (d^{*_{g_t}}-d^{*_{g_1}})  \mu_1 \|_{L^2_1(Y\times [T-1,T+1])}. $$ 
The first term in the right hand side goes to zero for $T$ 
sufficiently large, because the path $\nu_t$ ends at $\nu_1$. 
Upon writing  
$$ d^{*_{g_t}}-d^{*_{g_1}}=*_{g_t}\circ d\circ (*_{g_t}-*_{g_1}) + 
(*_{g_t}-*_{g_1})\circ d \circ *_{g_1} $$ 
we estimate the second term by a constant times the norm  
$\| *_{g_t}-*_{g_1} \|$, using the fact that $d$ is invertible on 
$Ker(d^*)$, since $Y$ is a rational homology sphere. 
The difference of the $*$-operators is controlled by the  
difference $| \sqrt{|g_t|}-\sqrt{|g_1|} |$, with  
$|g|$ the absolute value of the determinant of the 
metric. For $(g_0,\nu_0)$ and $(g_1,\nu_1)$  
two sufficiently close points in two different chambers, the path 
$(g_t,\nu_t)$ can be chosen so that this term is small enough, so that  
we obtain an estimate  
$$ \| F_{\#_T} \| \cdot \|  \AA-\AA_1 \|_{L^2_{1,\delta}(Y\times 
[T-1,T+1])} \leq \tilde C, $$  
with $\tilde C < \mu_0$, after combining the resulting estimate 
$$ \|  \AA-\AA_1 \|_{L^2_1(Y\times [T-1,T+1])} \leq C_1, $$ 
which gives 
$$ \| \AA-\AA_1 \|_{L^2_{1,\delta}(Y\times [-T-1,-T+1])} \leq C_1 
e^{\delta T} $$ 
with the estimate $\| F_{\#_T} \|\leq C e^{-\delta T}$. 
 
This gives the required estimate for the error term. 
We then have the estimate 
$$  \| {H^0_{x\#^0_{T} y}}^{-1} (1-\Pi(\mu,x\#^0_{T} y))\left( {\cal 
N}_{x\#^0_{T} y} {\cal L}^*_{x\#^0_{T} y}  (\beta,\xi) + f(x\#^0_{T} y) 
\right) \| \leq $$ 
$$ \frac{C}{\mu} \| (\beta,\xi) \| + \frac{1}{\mu} \| f(x\#^0_T y) \| 
\leq $$ 
$$ \frac{C}{\mu} \epsilon + \frac{\tilde C}{\mu}. $$ 
With the choice of $\mu$ and $\epsilon$ as indicated, the map is a 
contraction. The case of (right)pre-gluing maps is completely analogous. 
 
\noindent $\diamond$ 
 
For $T\geq 
T_0$ sufficiently large, the projection $\Pi(\mu,x\#^0_{T} y)$ is 
identified with the projection on the approximate 
cokernel 
$$ ApprCoker({\cal L}_{x\#^0_T y})\cong Coker({\cal L}_y)\cong \CC. $$ 
 
Thus, the construction of the obstruction bundle 
in \cite{Ta2} can be rephrased as in Theorem 4.53 of 
\cite{D2}. Namely, we obtain the following description of the 
obstruction. 
 
\begin{thm} 
\label{geomobstr} 
Consider the space ${\cal U}(\theta_0,O_{a'})$, obtained by gluing 
$x=(\AA,0)$ in ${\cal M}(\theta_0,\theta_1)$ with $y\in  {\cal 
M}(\theta_1,O_{a'})$.  
There exists a bundle $\tilde{\cal E}$ over  
$ {\cal M}(\theta_1,O_{a'})/U(1)$ with fiber 
\[ Coker({\cal L}_x)=\CC, \] 
with the following property. The set of pre-glued elements $x\#_T^0 y$  
in ${\cal U}(\theta_0,O_{a'})$ that can be deformed to an actual 
solution in $\M(\theta_0,O_{a'})$ form a codimension 2 submanifold of  
$$ \M(\theta_1,a')\cong {\cal M}(\theta_1,O_{a'})/U(1), $$ 
given by the zero set of a transverse section of $\tilde{\cal E}$. 
The case of ${\cal U}(O_a,\theta_1)$ is analogous. 
\end{thm} 
 
\noindent\underline{Proof:} As a result of Lemma 
\ref{singular:glue:case}, we know that, in the case of ${\cal 
U}(\theta_0,O_{a'})$, the obstruction bundle described in Proposition 
\ref{obstruction} is a bundle ${\cal E}$ over ${\cal 
M}(\theta_0,\theta_1)\times  {\cal M}(\theta_1,O_{a'})$ 
with fiber  
$ApprCoker({\cal L}_{x\#^0_T y})$. Since we are 
considering the case with $\mu(\theta_0)-\mu(\theta_1)<0$, according 
to Lemma \ref{noflow}, the moduli space ${\cal M}(\theta_0,\theta_1)$ 
consists of a unique gauge class $[\AA,0]$. Thus we have a 
diffeomorphism  
\[ \pi:  {\cal M}(\theta_1,O_{a'}) \stackrel{\cong}{\to}  
{\cal U}(\theta_0,O_{a'}) \] 
induced by the pre-gluing map. 
According to Corollary \ref{apprkercoker2}, since we can restrict our 
attention to the case $\mu(\theta_0)-\mu(\theta_1)=-2$, we have  
\[ ApprCoker({\cal L}_{x\#^0_T y})\cong Coker({\cal 
L}_x) \cong \CC. \]    
 
Consider the bundle $\tilde{\cal E}$ over $ {\cal 
M}(\theta_1,O_{a'})$ with fiber  
$$ Coker({\cal L}_x)\cong \CC,  $$  
obtained as 
pullback of the obstruction bundle ${\cal E}$ via the diffeomorphism 
$\pi$ of the base spaces. \index{$\tilde s$} 
There is a section $\tilde s=\pi^* s_\mu$ of $\tilde {\cal E}$ that 
corresponds to the canonical section $s_\mu$ of Proposition 
\ref{obstruction}.  
Recall that there is a free $U(1)$-action over the moduli 
space $ {\cal M}(\theta_1,O_{a'})$, hence there is a smooth 
projection to the quotient  
$$ {\cal M}(\theta_1,O_{a'})/U(1)\cong  {\cal 
M}(\theta_1,a'). $$ 
The bundle $\tilde {\cal E}$ can be regarded as the pullback, under 
this quotient map, of a bundle $\bar{\cal E}$ over $ {\cal 
M}(\theta_1,O_{a'})/U(1)$, with fiber  
$$ Coker({\cal L}_x)\cong \CC.  $$ 
The section $\tilde s$ is the pullback \index{$\bar s$} 
of a corresponding section $\bar s$ of $\bar{\cal E}$.  
We are going to proceed as follows. We show 
that, for a generic choice of the 
perturbations $\rho$  
and $P$ of equations (\ref{4SW1P't}) and (\ref{4SW2P't}), this section  
$\bar s$ is a generic section of $\bar{\cal E}$. 
 
Consider the universal pre-glued space 
\[ \tilde{\cal U}(\theta_1,O_{a'})=\{ (\rho,P,x\#^0_T y) \} \] 
with $\rho$ a compactly supported form in $\Lambda^1_c(Y\times \RR)$, 
$P$ a perturbation of ${\cal P}$ as in Definition \ref{calP} (with the 
modified property (1') specified at the beginning of Section 6.1).  
The pre-glued element $x\#^0_T y$ is obtained from 
$x$ and $y$, solutions respectively of the equations (\ref{4SW1P't}), 
(\ref{4SW2P't}) and (\ref{4SW1P'}), (\ref{4SW2P'}), with the 
perturbations $\rho$ and $P$.  
We can extend the bundle $\tilde {\cal E}$ to a local bundle over 
$\tilde{\cal U}(\theta_1,O_{a'})$ with fiber $Coker({\cal 
L}_x)$.  
We can still identify this as the pullback of a local bundle 
$\bar{\cal E}$ on the quotient $\tilde{\cal U}(\theta_1,O_{a'})/U(1)$. 
There is an induced section $\bar s$ of $\bar{\cal E}$ whose pull-back 
agrees with the section $\pi^* s_\mu$ for fixed perturbations $(\rho, P)$. 
We prove that this section $\bar s$ over $\tilde{\cal 
U}(\theta_1,O_{a'})/U(1)$ is transverse to the zero section, hence the 
restriction to a generic $(\rho, P)$ gives a generic section over 
${\cal U}(\theta_1,O_{a'})/U(1)$.   
 
The section $\bar s$ is given by 
\[ \bar s(\rho, P,x\#^0_T y)=\Pi_{Coker({\cal L}_{x})} 
\left( {\cal N}_{(\rho,P,x\#^0_T y)}{\cal L}^*_{(\rho,P,x\#^0_T 
y)}(\beta,\xi)+  
f_{(\rho,P)}(x\#^0_T y)\right). \] 
Suppose given a point $(\rho,P,x\#^0_T y)$ in $\bar 
s^{-1}(0)$. Consider a small variation of the perturbation 
$\rho+\epsilon\eta$. The variation of the term 
$f_{(\rho+\epsilon\eta,P)}(x\#^0_T y)$ is given by 
\[ \eta\cdot \Psi_1\#^0_T\Psi_2= \eta\cdot \rho_{T}^+  
\Psi_2^{-2T}. \]  
Let $\Phi_1$ and $\Phi_2$ in $Ker(D_{\AA-\rho})$ be the generators of 
the 2-dimensional space $Coker({\cal L}_x)$. Consider a small open sets 
$U_i$ where $\Phi_i$ and $\rho_T^+ \Psi_2^{-2T}$ are non-vanishing, 
and almost constant. There exist 1-forms $\eta_i$ supported in  
small neighborhoods of the open sets $U_i$ such that  
\[ \langle \Phi_i, \eta_i\cdot \rho_T^+ \Psi_2^{-2T} \rangle \] 
is non-zero on $U_i$. Thus, we obtain 
\[ \int_{Y\times \RR} \langle \Phi_i, \eta_i\cdot \rho_T^+ 
\Psi_2^{-2T} \rangle dv \neq 0. \] 
Thus, by varying the perturbation $\rho$ alone, it is possible to achieve 
surjectivity of the linearization of the section $\bar s$ onto 
$Coker({\cal L}_x)$.   
 
There is a free $U(1)$-action on the space $ {\cal 
M}(\theta_1,O_{a'})$, whereas the element $[x]=[\AA,0]$ in ${\cal 
M}(\theta_0,\theta_1)$ is fixed by the $U(1)$ action. The section  
\[ \tilde s=\pi^* s_\mu:   {\cal M}(\theta_1,O_{a'}) \to 
\tilde{\cal E} \]  
is invariant under the $U(1)$-action, being the pullback of $\bar s$.  
We have seen that, for a generic choice of the perturbation 
$(\rho,P)$, the section $\bar s$ is a generic 
section.  
Thus, the approximate solutions in ${\cal U}(\theta_0, O_{a'})$ that 
can be glued to actual solutions in $\M(\theta_0, O_{a'})$ are 
identified with the co-dimension 2 $U(1)$-submanifold  
$$ \tilde s^{-1}(0) \subset \M(\theta_1, O_{a'}). $$ 
In other words, we have then proved that the bundle $\bar{\cal E}$ 
over the moduli space $  {\cal M}(\theta_1,O_{a'})/U(1)$ is a model of
the obstruction.  
 
\noindent $\diamond$ 
 
We have the following consequence of Theorem \ref{geomobstr} and Lemma  
\ref{singular:glue:case}. 
 
\begin{prop}\label{gluing:obstr} 
For any compact set 
$$ K\subset \M(\theta_0,\theta_1)\times \tilde s^{-1}(0), $$ 
with 
$\tilde s^{-1}(0) \subset \M(\theta_1, O_{a'})$ 
the zeroes of the obstruction section 
$$ \tilde s^{-1}(0)/U(1)= \bar s^{-1}(0) \subset \M(\theta_1,a'), $$ 
there is an orientation preserving gluing map 
$$ \#: K\subset \M(\theta_0,\theta_1)\times \tilde s^{-1}(0) \to 
\M(\theta_0, O_{a'}) $$ 
that is a smooth embedding. There is a similar gluing map for any 
compact set 
$$ K\subset \tilde s^{-1}(0)\times \M(\theta_0,\theta_1), $$ 
with 
$\tilde s^{-1}(0) \subset \M(O_a,\theta_0)$ 
the zeroes of the obstruction section 
$$ \tilde s^{-1}(0)/U(1)= \bar s^{-1}(0) \subset \M(a,\theta_0). $$ 
\end{prop} 
 
Finally, we have the following result on the singular components in 
the ideal boundary of the moduli spaces $\M(O_a,O_{a'})$. 
 
\begin{thm} 
Assume that $\mu(a) -\mu(\theta_0) \ge 2$ and $\mu(\theta_1)
-\mu(a')\ge 3$, then the contributions of the singular strata   
\beq \label{sing:bound:1} \M(\theta_0,\theta_1)\times
\M(\theta_1,O_{a'}) \eeq 
and 
\beq  \M(O_a,\theta_0)\times \M(\theta_0,\theta_1)
\label{sing:bound:2} \eeq 
to the actual boundary of the compactified moduli spaces 
$$ \M(\theta_0,O_{a'})^* \ \ \hbox{ and } \ \ \M(O_a,\theta_1)^*$$ 
are given by the terms 
\beq\label{sing:bound:I:1}  
\M(\theta_0,\theta_1)\times \partial^{(1)}\left( \tilde 
s^{-1}(0) \cap \M (\theta_1,O_{a'}) \right) 
\eeq 
and 
\beq\label{sing:bound:I:2} 
\partial^{(1)}\left( \tilde s^{-1}(0)\cap 
\M (O_a,\theta_0)\right) \times \M(\theta_0,\theta_1),  
\eeq 
respectively,  
where $\tilde s$ are  the obstruction sections. 
Thus, the compactification of the moduli spaces 
$\M(O_a,O_{a'})$ has the structure of a smooth manifold with corners, 
but in addition to the strata of the form (\ref{boundary:t}), we also 
have the strata (\ref{sing:bound:I:1}) and 
(\ref{sing:bound:I:2}). Thus, we obtain  
\beq 
\begin{array}{lll} 
\partial^{(1)} \M(\theta_0,O_{a'})^* &=&  
\bigcup_{\{a|\mu(O_{b'}) >\mu(O_{a'})\}} ({\cal M}(\theta_0,O_{b'})^* 
\times_{O_{b'}}  
\hat {\cal M}(O_{b'}, O_{a'})^*) \\[1mm] 
&&\bigcup_{\{ b| \mu(O_b) <\mu (\theta_0)  \}}(-\hat  {\cal M}(\theta_0,O_b)^* 
\times_{O_b}   
{\cal M}(O_b,  O_{b'})^*)\\[1mm] 
&&\cup \M(\theta_0,\theta_1)\times \partial^{(1)}\left( \tilde 
s^{-1}(0) \cap \M (\theta_1,O_{a'})^* \right) 
\end{array} 
\label{boundary:t:obstr:1} 
\eeq 
and 
\beq 
\begin{array}{lll} 
\partial^{(1)} \M(O_a,\theta_1)^*&=& 
\bigcup_{\{a'|\mu(O_{a'}) >\mu(\theta_1)\}} ({\cal M}(O_a,O_{a'})^* 
\times_{O_{a'}}  
\hat {\cal M}(O_{a'},\theta_1 )^*) \\[1mm] 
&&\bigcup_{\{ b| \mu(O_b) <\mu (O_a)  \}}(-\hat  {\cal M}(O_a,O_b)^* 
\times_{O_b}   
{\cal M}(O_b, \theta_1 )^*)\\[1mm] 
&& \cup \partial^{(1)}\left( \tilde s^{-1}(0)\cap 
\M (O_a,\theta_0)\right) \times \M(\theta_0,\theta_1). 
\end{array} 
\label{boundary:t:obstr:2} 
\eeq 
\label{noglue} 
\end{thm} 
 
\noindent\underline{Proof:} 
By the results of Proposition \ref{gluing:obstr}, we have gluing maps  
$$ \M(\theta_0,\theta_1)\times \bar s^{-1}(0) \to \M(\theta_0, a') $$ 
and 
$$ \bar s^{-1}(0) \times \M(\theta_0,\theta_1) \to \M(a,\theta_1) $$ 
that are smooth embeddings. The dimension count then implies that 
these gluing maps are diffeomorphisms of $\bar s^{-1}(0)$ to a union 
of connected components of $\M(\theta_0, a') $ or $\M(a,\theta_1) $. 
Similarly, we have diffeomorphisms of $U(1)$-manifolds between 
$\tilde s^{-1}(0)$ and a union of connected components of 
$\M(\theta_0,O_{a'})$ or $\M(O_a,\theta_1) $, induced by the gluing 
maps 
\beq \#: \M(\theta_0,\theta_1)\times \tilde s^{-1}(0) \to 
\M(\theta_0,O_{a'}) \label{glue:tilde:1} \eeq and 
\beq \#: \tilde s^{-1}(0) \times \M(\theta_0,\theta_1)  \to 
\M(O_a,\theta_1). \label{glue:tilde:2}\eeq 
 
Under the gluing maps, the image of  
$$ \M^{bal}(O_a,\theta_0)\cap \tilde s^{-1}(0) \ \ \hbox{ or } \ \ 
\M^{bal}(\theta_1,O_{a'})\cap \tilde s^{-1}(0) $$ 
is a co-dimension 1 submanifold of $\M(\theta_0,O_{a'})$ or 
$\M(O_a,\theta_1) $. This submanifold actually lies in the interior of  
$\M(\theta_0,O_{a'})$ or $\M(O_a,\theta_1) $, and is not part of its 
boundary strata. This can be seen from the fact that the gluing maps 
(\ref{glue:tilde:1}) and (\ref{glue:tilde:2}) provide the collar 
structure around these smooth codimension one embedded submanifolds.   
In other words, this means that any sequence of solutions in 
$\M(\theta_0,O_{a'})$ or $\M(O_a,\theta_1)$ that converges to an 
element in the ideal boundary components (\ref{sing:bound:1}) and
(\ref{sing:bound:2}) is in fact  
already convergent in the interior (top stratum) of $\M(\theta_0,O_{a'})$ or 
$\M(O_a,\theta_1)$. The only contribution of (\ref{sing:bound:1}) and
(\ref{sing:bound:2}) to  
the actual boundary of the compactification then comes from  
the boundary points of 
$$ \M (O_a,\theta_0)\cap \tilde s^{-1}(0) $$ 
and 
$$ \M (\theta_1,O_{a'})\cap \tilde s^{-1}(0). $$ 
This gives the formulae  (\ref{boundary:t:obstr:1}) and 
(\ref{boundary:t:obstr:2}). 
 
\noindent $\diamond$

Similarly, we can analyze the boundary structure of 
$\M^P(\theta_0,O_a)$, in the presence of obstructions. We have the 
following result.  
 
\begin{thm} 
\label{boundary:P:obstr} 
For a generic choice of perturbation, the obstruction section 
$$ \bar s : \M(\theta_1,a) \to {\cal L}_{\theta_1,a}= 
\M(\theta_1,O_a)\times_{U(1)} \CC $$  
defines a codimension 2 submanifold  
 $ \bar s^{-1}(0)\subset \M(\theta_1,a)$, 
which corresponds to a $U(1)$-submanifold of co-dimension 2,  
$ \tilde s^{-1}(0)\subset \M(\theta_1,O_a) $,  
 such that the gluing map 
$$ \#: K\subset \M(\theta_0,\theta_1)\times \tilde s^{-1}(0)\to 
\M^P_{\sigma=1} (\theta_0,O_a) $$ 
is a smooth embedding, for any compact set  
$  K\subset \M(\theta_0,\theta_1)\times \tilde s^{-1}(0). $ 
Thus, the co-dimension one boundary strata of the compactification  
$\M^P(\theta_0,O_a)^*$ are given by 
\beq \label{strataP:2:obstr} 
\begin{array}{rl} 
\partial^{(1)} {\cal M }^P(\theta_0, O_a )^* =& \bigcup_{a'}({\cal 
M}(\theta_0, O_{a'})^* \times   
_{  O_{a'}} {\cal M}(O_{a'}, O_b)^* )\\[1mm] 
& \bigcup  \M(\theta_0,\theta_1)\times ( \tilde 
s^{-1}(0) \cap  \M(\theta_1,O_a) ) \\[1mm]  
& \bigcup_{\{c | \mu(O_c)\ge\mu(\theta_0)\}} ( {\cal M}^P(\theta_0, O_c)^* 
\times_{  O_c}  
\hat {\cal M}(O_c,  O_b)^*)\\[1mm] 
& \bigcup_{\{c| \mu(O_c) \le\mu(\theta_0)\}} (\hat {\cal M}(\theta_0, O_c)^* 
\times_{  O_c}   
{\cal M}^P(O_c, , O_a)^*).\end{array} 
\eeq 
with the orientations given by the gluing theorem. 
\end{thm} 
 
The proof follows from the analysis of the obstruction bundles and 
sections, as in the case of Proposition \ref{gluing:obstr} and Theorem 
\ref{noglue}.

\subsection{Proof of topological invariance} 
 
We gave an argument for the easy case of theorem \ref{metrics}, with 
maps $I$, $J$, and $H$ as in (\ref{I:def}), (\ref{J:def}), and 
(\ref{H:def}). This gives us an isomorphism between  
the equivariant Floer homologies for any two choices $(g_0,\nu_0)$ and  
$(g_1,\nu_1)$ within the same chamber. Now we prove the general case 
 
\noindent\underline{Proof of Theorem \ref{metrics}, Part II: the 
general case.} 
In Section 6.1 and Section 6.2, we analysed the boundary structure  
of the moduli spaces $\M(O_a,O_{a'})^*$ and $\M^P(O_a,O_b)^*$, in the
case of two metrics and  
perturbations $(g_0,\nu_0)$ and  $(g_1,\nu_1)$ in different chambers, 
connected by a path $(g_t,\nu_t)$ satisfying 
$SF(\partial^{g_t}_{\nu_t})=-2$, and a path $(\tilde g_t,\tilde 
\nu_t)$, in the opposite direction, satisfying $SF(\partial^{\tilde 
g_t}_{\tilde \nu_t})= -SF(\partial^{g_t}_{\nu_t})=2$. We assume 
throughout the discussion that the metrics $(g_0,\nu_0)$ and  
$(g_1,\nu_1)$ are ``close enough'', on 
the two sides of the wall. We can assume, similarly, that the paths 
$(g_t,\nu_t)$ and $(\tilde g_t,\tilde \nu_t)$ are also close enough. 
 
Recall that we perform a 
shift of grading in the complex ${C_k}_{U(1)}(Y, g_1,\nu_1)$ by 
setting $\mu(\theta_1)=-SF(\partial^{g_t}_{\nu_t})$, where 
$SF(\partial^{g_t}_{\nu_t})$ is the spectral flow 
of the Dirac operator along the path of reducible solutions 
$[\nu_t,0]$. 
  
The analysis of the obstruction bundles in Section 6.2 implies that 
the boundary structure of the moduli spaces $\M(O_a,O_{a'})$ and 
$\M^P(O_a,O_b)$ is modified by the presence of the zeroes of the 
obstruction sections. This difference determines suitable correction 
terms for the maps $I$ and $H$, so that the  
argument of theorem \ref{metrics} can be adapted to this general case.

\begin{rem} 
The moduli space $\M(\theta_1,\theta_0)$, as a $U(1)$-manifold, consists 
of a disk, containing the fixed point $x_1=[\AA_1, 0]$ (the unique 
solution of $d*\AA_1=0$ and $d^+ \AA_1 =\tilde \mu_t$), and with 
boundary a circle, given either by a component 
$$ \M(\theta_1,O_a)\times_{O_a} \hat\M(O_a,\theta_0) $$ 
with $\mu(O_a)-\mu(\theta_0)=1$, or by a component 
$$ \hat\M(\theta_1,O_a')\times_{O_a'} \M(O_a',\theta_0), $$ 
with $\mu(O_a')-\mu(\theta_0)=0$. 
Thus, we have the following identity that counts the boundary 
components of the framed monopole  
moduli space $\M(\theta_1,\theta_0)$: 
\beq  
\sum_{a_{(1)}} n_{\theta_1 a_{(1)}} n_{a_{(1)} \theta_0}  
+ \sum_{a_{(0)}'} n_{\theta_1 a_{(0)}'} n_{a_{(0)}' \theta_0}   
=1,  
\label{nice:1}  
\eeq 
The sum is over all $O_{a_{(1)}}$ in $\M^0(Y,g_0,\nu_0)$,  
and $O_{a_{(0)}'}$ in $\M^0(Y,g_1,\nu_1)$ of (shifted) index 
$\mu(O_{a_{(1)}})-\mu(\theta_0)=1$ and 
$\mu(O_{a_{(0)}'})-\mu(\theta_0)=0$, respectively.    
\label{disk:component} 
\end{rem}

In particular, notice that we can derive the observation of Remark 
\ref{disk:component} if we make explicit use here of the information on the 
local structure of the moduli space $\M^0 (Y,g_t,\nu_t)$ of critical 
orbits, as the metric $(g_t,\nu_t)$ approaches the wall, This is 
derived in Section 7.3.  
We know that two possibilities arise: either an 
irreducible orbit $O_a$, with $\mu(O_a)-\mu(\theta_0)=1$, disappears  
into the reducible as the metric and perturbation $(g_t,\nu_t)$ hits 
the wall, or an irreducible orbit $O_{a'}$ with 
$\mu(O_{a'})-\mu(\theta_0)=0$ arises from the reducible, as 
$(g_t,\nu_t)$ hits the wall. In the first case, the disk of Remark 
\ref{disk:component} has boundary the circle 
$$ \M(\theta_1,O_a)\times_{O_a} \hat\M(O_a,\theta_0) $$ 
and in the second case it has boundary the circle 
$$ \hat\M(\theta_1,O_a')\times_{O_a'} \M(O_a',\theta_0). $$ 
 
The relation (\ref{nice:1}) then yields the separate identities 
\beq \#\M(\theta_1,O_a)= n_{\theta_1a} =1 \ \ \hbox{ and } \ \ 
\#\hat\M(O_a,\theta_0)= n_{a\theta_0}=1, \label{disk:boundary} \eeq 
for $O_a$ the unique orbit that hits the reducible, and 
$$ \sum_{a_{(1)} \neq a} n_{\theta_1 a_{(1)}} n_{a_{(1)} \theta_0} + \sum_{a_{(0)}'} 
n_{\theta_1 a_{(0)}'}  
n_{a_{(0)}' \theta_0} =0, $$ 
for $\mu(O_{a_{(1)}})-\mu(\theta_0)=1$ and 
$\mu(O_{a_{(0)}'})-\mu(\theta_0)=0$. 
The case of $O_{a'}$ is analogous.

The moduli space $\M(\theta_1,\theta_0)$, containing this separate 
component with a fixed point, is one of the differences with respect 
to the picture for metrics and perturbations within the same chamber. 
Another essential difference is, of course, the presence of the singular 
moduli space  
$\M(\theta_0,\theta_1)$. We have already seen, in the analysis of the 
obtructions how this moduli space plays an essential role. However, 
what we wish to point out here is that the basic asymmetry between the  
moduli spaces $\M(\theta_0,\theta_1)$ and $\M(\theta_1,\theta_0)$ is 
what calls for correction terms in the maps $I$ and $H$, but not in 
the map $J$. 
 
In fact, consider first the action of the map $J$, defined as in 
(\ref{J:def}) on the generator $\Omega^n \otimes 1_{\theta_1}$ in the 
Floer complex for $(g_1,\nu_1)$. We have 
$$ J(\Omega^n \otimes 1_{\theta_1})=  
\sum_{a_{(1)}} n_{\theta_1 a_{(1)}} \Omega^n \otimes 1_{a_{(1)}}, $$ 
with the sum over all $O_{a_{(1)}}$ with 
$\mu(O_{a_{(1)}})-\mu(\theta_0)=1$.

Notice that the existence of the extra disk component in the boundary 
of $\M(\theta_1,\theta_0)$, as in remark \ref{disk:component}, does 
not affect the identity $JD-DJ=0$, in fact, even though now the count 
of boundary terms in $\M(\theta_0,\theta_1)$ satisfies (\ref{nice:1}), 
we still do get 
$$ \langle \Omega^n\otimes 1_{\theta_0}, JD-DJ (\Omega^n\otimes 
1_{\theta_1}) \rangle =0. $$ 
In fact, the fibrations 
$$ \M(O_{a_{(0)}'},\theta_0) \to \theta_0, $$ 
$$ \hat\M(O_{a_{(1)}},\theta_0) \to \theta_0, $$ 
have 1-dimensional fibers, and the pushforward of a zero-form 
is trivial.  
The identity $JD-DJ=0$ at all the other components follows the 
argument given in the proof of the easy case of Theorem \ref{metrics}, 
without any modification.  
 
Now, instead, consider the case of the map $I$ acting on 
$\Omega^n\otimes 1_{\theta_0}$. The first difference we notice, with
respect to the model case of metrics and perturbations in the same
chamber, is that the moduli space $\M(\theta_0,\theta_1)$ consisting
of the unique point  $x_0= [\AA_0,0]$ defines a non-trivial
``pull-back push-forward'' acting on the zero form $1_{\theta_0}$.
To account for this moduli space, we 
have to assume the existence here of an extra component of the map $I$  
connecting $\Omega^n\otimes 1_{\theta_0}$ to $\Omega^{n-1}\otimes 
1_{\theta_1}$, where the drop of degree in $\Omega$ accounts for the 
change of grading of the reducible point, so that the map $I$ can be 
of degree zero. Thus, we have a new component 
\beq \langle \Omega^{n-1}\otimes 1_{\theta_1},  I (\Omega^n\otimes 
1_{\theta_0}) \rangle =1 \label{tildeI:1}. \eeq 
Notice that the necessity of the additional term (\ref{tildeI:1})  
in the map $I$ can be made clear by looking at the following 
example.  
 
\noindent{\bf Example.} 
Consider a model case where the moduli space $\M^0(g_0,\nu_0)$ 
consists solely of the fixed point $\theta_0$ and of an irreducible 
orbit $O_a$ in degree one. The Floer complex for $(g_0,\nu_0)$ then 
has generators $\Omega^n\otimes 1_{\theta_0}$ in even degrees $p=2n 
\geq 0$, $\Omega^n\otimes \eta_a$ in odd degrees $p=2n+1 \geq 1$, and  
$\Omega^{n-1} \otimes 1_a$ in even degrees $p=2n\geq 2$. The boundary 
operator $D$ has a component 
$$ \Omega^n\otimes \eta_a \mapsto - \Omega^{n-1} \otimes 1_a + n_{a0}
\Omega^n\otimes 1_{\theta_0}. $$    
In the wall crossing with 
$\mu(\theta_0)-\mu(\theta_1)=-2$, the orbit $O_a$ disappears into the 
reducible, and the Floer complex for $(g_1,\nu_1)$ only has the fixed 
point as generator. In this case, we have the equivariant Floer homology 
$$ SWH_{p,U(1)}(Y, (g_0,\nu_0))= \left\{ \begin{array}{lr} \RR [ 
\Omega^{n-1} \otimes 1_a ] & p=2n\geq 2 \\ 
0 & \hbox{otherwise,} \end{array} \right. $$ 
with $[\Omega^{n-1} \otimes 1_a ]=[\Omega^n\otimes 1_{\theta_0}]$ for 
all $n\geq 1$, and $[1\otimes 1_{\theta_0}]=0$. 
The Floer complex for $(g_1,\nu_1)$ has $\Omega^n\otimes 1_{\theta_1}$  
as unique generators, and no boundary components, hence the Floer 
homology is 
$$ SWH_{p,U(1)}(Y, (g_1,\nu_1))=\left\{ \begin{array}{lr} \RR 
[\Omega^n\otimes 1_{\theta_1}] & p=2n+2 \geq 2 \\ 
0 & \hbox{otherwise,} \end{array} \right. $$ 
after degree shift. 
The map $I$ that maps $\Omega^n \otimes 1_{\theta_0} \to \Omega^{n-1} 
\otimes 1_{\theta_1}$ gives the desired isomorphism. 
Notice that such examples can in fact be realized geometrically, for 
instance when considering a metric of positive scalar curvature on 
$S^3$ as $g_1$. 
 
\noindent $\diamond$ 
 
It becomes then clear that, if we change the map $I$ with an extra 
component as  
in (\ref{tildeI:1}), to account for the moduli space 
$\M(\theta_0,\theta_1)$, we need to add further correction terms to 
the original map $I$ on other generators, so that the identity $  
ID-DI=0$ continues to be satisfied. In assigning the necessary 
correction terms, we need to take into account the different structure  
of the compactification of the moduli spaces $\M(O_a,O_{a'})$, with 
the boundary strata (\ref{boundary:t:obstr:1}) and 
(\ref{boundary:t:obstr:2}). Studying the boundary information of $\M
(\theta_0, a_{(-2)}')$, we know that there is a contribution of the 
singular gluing from the monopoles in
\[ 
s^{-1}(0)   
\cap \M (\theta_1, a_{(-2)}'), 
\]
which contributes to the expression
\[ 
(ID -DI ) (\Omega^n\otimes 1_{\theta_0}), 
\] 
in addition to the ordinary components 
$n_{\theta_0 a_{(-1)}'}\Omega^n \otimes 1_{a_{(-1)}'}$, where 
$n_{\theta_0 a_{(-1)}'}$  effectively counts the monopoles 
from the zeros of the obstruction section over 
$\M ( \theta_1, a_{(-1)}')$.  This is precisely the correction term
which is needed in order to obtain the identity $ID-DI=0$ once we take 
into account the presence of the component (\ref{tildeI:1})
originates from the presence of  the moduli space 
$\M(\theta_0,\theta_1)$ consisting of the unique point 
$x_0= [\AA_0,0]$.

Similarly, we need to 
introduce correction terms to the map $H$ so that it continues to be a  
chain homotopy, satisfying $id-JI =DH+HD$, with respect to the 
modified map $I$ and with the modified structure of the 
compactification of $\M^P(O_a,O_b)$. 
 
In order to introduce the correction terms for the maps $I$ and $H$, 
we need a preliminary discussion on some identities derived from the 
counting of zeroes of the obstruction sections. 
 
Consider the case of the 2-dimensional moduli space $\M(a,\theta_0)$, 
with $ a$ an irreducible critical point 
 of index 2 and $\theta_0$ the unique reducible point. 
 
In Section 5.3, in Lemma \ref{rel:eul:class}, we described the 
invariant $m_{ac}$ as the relative Euler class of the complex line
bundle ${\cal L}_{ac}$ over 
$\M(a,c)$. Consider now the case of the 2-dimensional moduli space 
$\M(a,\theta_0)$, with $O_a$ a free orbit of index 2 and $\theta_0$ 
the fixed point in $\M^0(Y,g_0,\nu_0)$. We also introduced the
invariant $m_{a\theta_0}$ in Section 5.3, as 
as the relative Euler class of the associated complex line 
bundle over $\M (a, \theta_0)$, according to Lemma 
\ref{tens:line} and Lemma \ref{zeroes:s:ac}. As we point out in the
following Remark, we can use the obstruction section to obtain the
necessary trivialization and compute the invariant $m_{a\theta_0}$ as
relative Euler class of ${\cal L}_{a\theta_0}$.
 
\begin{rem}\label{sec:obstr:induction}
We use the same choice of trivialization $s_{ab}$ over $\M(a,b)$, for 
all the free orbits $O_a$ and $O_b$ of relative index 1, as determined  
by the trivialization $\varphi$ of Lemma \ref{rel:eul:class}. For 
$\M(b,\theta_0)$ with $\mu(O_b)-\mu(\theta_0)=1$,  
we set $s_{b,\theta_0}\neq 0$ to be the obstruction section  
$$ s_{b,\theta_0}=\bar s : \M(b,\theta_0) \to \M(O_b,\theta_0)\times_{U(1)} 
Coker( {\cal L}_x)\cong {\cal L}_{b\theta_0}. $$ 
 
We then set $s_{a,\theta_0}$ over $\M(a,\theta_0)$, 
with $\mu(O_a)-\mu(\theta_0)\geq 2$ to be the obstruction section 
\beq s_{a,\theta_0}=\bar s : \M(a,\theta_0) \to \M(O_a,\theta_0)\times_{U(1)} 
Coker( {\cal L}_x)\cong {\cal L}_{a\theta_0}. \label{sec:ac:obst} \eeq 
This choice satisfies the requirement of the class of sections 
specified above.  Over a moduli spaces $\M(a,\theta_0)$, with 
$\mu(O_a)-\mu(\theta_0)\geq 2$, then 
the obstruction section 
$$ s_{a,\theta_0}=\bar s : \M(a,\theta_0) \to \M(O_a,\theta_0)\times_{U(1)} 
Coker( {\cal L}_x)\cong {\cal L}_{a\theta_0} $$ 
is homotopic to  
$$  \pi_1^*s_{ab}\otimes \pi_2^*s_{b,\theta_0} $$ 
over all the submanifolds 
$$ \M(a,b)\times \M(b,\theta_0) $$ 
of $\M(a,\theta_0)$, where we identify 
$$ {\cal L}_{a,\theta_0}=\pi^*_1 {\cal L}_{ab}\otimes \pi_2^* {\cal 
L}_{b,\theta_0}, $$ 
with $\pi_1$ and $\pi_2$ the projections on the two factors in 
$$ \M(a,b)\times \M(b,\theta_0), $$ 
with $\mu(O_a)>\mu(O_b)>\mu(\theta_0)$,  
and  section $s_{b,\theta_0}$ is the obstruction section over 
$\M(b,\theta_0)$, and the section  
$s_{ab}$ is a transverse section of the line bundle 
$$ {\cal L}_{ab}=\M(O_a,O_b)\times_{U(1)} \CC $$ 
over $\M(a,b)$ as discussed in Section 5.3. 
\end{rem}
  
In summary, we have obtained the following lemma which claims that 
the relative Euler class $m_{a_{(2)} \theta_0}$ and  
$m_{\theta_1 a_{(-1)}'}$ can be calculated by counting 
the zeros of the obstruction section  
over $\M (a_{(2)}, \theta_0)$  
and $\M (\theta_1, a_{(-1)}')$ respectively. 
 
\begin{lem} 
For $\mu( a_{(2)})-\mu(\theta_0)=2$, the relative Euler 
class $m_{a_{(2)}\theta_0} $ is given by 
$$ \# s_{a_{(2)}\theta_0}^{-1}(0)$$ 
where $s_{a_{(2)}\theta_0}$ is the obstruction section  
over $\M (a_{(2)}, \theta_0)$, similarly, for $\mu(\theta_1) 
-\mu (a_{(-1)}')=3$, the relative Euler 
class $m_{\theta_1 a_{(-1)}'}$  
is given by 
$$ \# s_{ \theta_0a_{(-1)}'}^{-1}(0)$$ 
where $s_{ \theta_0a_{(-1)}'}$  is the obstruction section  
over $\M (\theta_1, a_{(-1)}')$. 
\label{sec:obstr:euler} 
\end{lem}

As a consequence  of 
Remark \ref{sec:obstr:induction} and Lemma \ref{sec:obstr:euler}, 
we have the following identities.  
 
\begin{lem} 
Let $O_{a_{(p)}}$ denote a free orbit in the moduli space 
$\M^0(Y,g_0,\nu_0)$, of index 
$\mu(O_{a_{(p)}})-\mu(\theta_0)=p$. Let $O_{a_{(q)}'}$ denote a free 
orbit in the moduli space  
$\M^0(Y,g_1,\nu_1)$, of (shifted) index 
$\mu(O_{a_{(q)}'})-\mu(\theta_0)=q$. We have the following 
properties.  
\begin{enumerate}  
\item Consider the obstruction bundle over $\M (a_{(2)}, \theta_0) $.
The counting of the zeroes of the obstruction, that is, of the
flowlines in $\M (a_{(2)}, \theta_0) $  glued
to the singular reducible $x=[\AA^0,0]$ in $\M(\theta_0,\theta_1)$,
is given by   
\[  
m_{a_{(2)}\theta_0}=\# \bigl(\bar s^{-1} (0) \cap \M (a_{(2)}, 
\theta_0) \bigr).   
\]  
This satisfies 
\beq  
\sum_{a_{(1)}} m_{a_{(3)}a_{(1)}} n_{a_{(1)} \theta_0}  
- \sum_{a_{(2)}} n_{a_{(3)}a_{(2)}}m_{a_{(2)} \theta_0}=0.  
\label{nice0:1}  
\eeq  
\item The counting of the zeroes of the obstruction   
section $\bar s$ of the obstruction bundle   
over $\M (\theta_1, a_{(-1)}')$, which counts the gluing of the 
singular reducible $x=[\AA^0,0]$ in $\M(\theta_0,\theta_1)$ to flowlines  
in $\M(\theta_1, a_{(-1)}')$, is given by   
\[  
m_{\theta_1 a_{(-1)}'}= \# \bigl(\bar s^{-1} (0) \cap \M (\theta_1, 
a_{(-1)}')  \bigr).   
\]  
This satisfies 
\beq  
\sum_{a_{(-1)}'} m_{\theta_1 a_{(-1)}'} n_{a_{(-1)}'  a_{(-2)}'} - 
\sum_{a_{(0)}'} n_{\theta_1 a_{(0)}'} m_{a_{(0)}'  a_{(-2)}'} = 0.  
 \label{nice0:2}  
\eeq  
\end{enumerate}  
\label{obstr:sec:euler} 
\end{lem}  
 
\noindent\underline{Proof:} 
In this case, we know that, for dimensional reasons, we have  
$$ \bar s^{-1} (0) \cap \M (a_{(1)}, \theta_0) = \emptyset. $$ 
Thus, by the previous discussion, and the results of Lemma \ref{tens:line}, 
Lemma \ref{zeroes:s:ac}, and Remark \ref{sec:obstr:induction}, we know 
that the boundary of the 1-dimensional manifold 
$\bigl(\bar s^{-1} (0) \cap \M (a_{(3)}, \theta_0) \bigr)$    
consists of the set 
$$ \bigcup_{a_{(2)}} \hat\M (a_{(3)}, a_{(2)}) \times  
\bigl(\bar s^{-1} (0) \cap \M (a_{(2)}, \theta_0) \bigr) $$  
$$ \cup \bigcup_{a_{(1)}}  
\bigl( - s^{-1}_{a_{(2)}a_{(1)}} (0) \cap (\M (a_{(2)}, a_{(1)}) \times   
\hat \M (a_{(1)}, \theta_0) ) \bigr). $$ 
This implies that we have the identity 
\[  
\sum_{a_{(1)}} m_{a_{(3)}a_{(1)}} n_{a_{(1)} \theta_0}  
-  
\sum_{a_{(2)}} n_{a_{(3)}a_{(2)}}m_{a_{(2)} \theta_0}=0,  
\]  
with the sign denoting the different orientation.
The remaining case is analogous. 
 
\noindent $\diamond$

Now we define the necessary modifications to the maps $I$ and $H$. 
 
\begin{defin} 
We modify the maps $I, H$ on   
$ \Omega^n \otimes 1_{\theta_0}$ and   
$\Omega^{n-1} \otimes 1_{a_{(1)}} $ as follows:  
\begin{enumerate}  
\item  
\beq \label{I:change:1} \begin{array}{rl} I (\Omega^n \otimes 
1_{\theta_0})=&  \Omega^{n-1} \otimes 1_{\theta_1} \\[2mm] 
& + \sum_{a_{(-1)}'} 
 n_{\theta_0a_{(-1)}'} \Omega^n\otimes 1_{a_{(-1)}'} \\[2mm] 
& + \sum_{a_{(0)}'} n_{\theta_1 a_{(0)}'}  \Omega^n\otimes 
\eta_{a'_{(0)}};  \end{array} \eeq  
\item  
\beq \label{I:change:2} \begin{array}{rl}  
I (\Omega^{n-1} \otimes 1_{a_{(1)}} )=&  n_{a_{(1)}\theta_0} 
\Omega^{n-1} \otimes 1_{\theta_1}   \\[2mm]  
& + \sum_{a_{(1)}'} n_{a_{(1)}a_{(1)}'} 
 \Omega^{n-1} \otimes 1_{a'_{(1)}}; \end{array} \eeq 
\item 
\beq \label{H:change}\begin{array}{rl} 
H(\Omega^n \otimes 1_{\theta_0})=& \sum_{a_{(0)}} n^P_{\theta_0 a_{(0)}} 
\Omega^n \otimes 1_{a(0)} \\[2mm] 
& + \sum_{ a_{(1)}} n_{\theta_1a_{(1)}}   \Omega^n \otimes 
\eta_{a_{(1)}}. \end{array} \eeq 
\end{enumerate}  
\label{def:I:H:change} 
\end{defin}

In order to show that the relations $ID=DI$ and $id-JI=DH+HD$ are 
still satisfied, we only need to check explicitly all the terms that 
are directly affected by the presence of the correction terms. 
 
\begin{lem} 
We have 
\begin{enumerate}  
\item  $ID =DI$  
on $\Omega^n \otimes 1_{\theta_0}$ and  
$\Omega^{n-1} \otimes 1_{a_{(1)}}$;  
\item $ID =DI $ on $\Omega^{n } \otimes \eta_{a_{(3)}}$;  
\item $ID =DI $ on $\Omega^{n } \otimes \eta_{a_{(1)}}$;  
\end{enumerate}  
\label{check:ID=DI} 
\end{lem} 
 
\noindent\underline{Proof:}  We prove the claim by direct analysis of 
the boundary strata of the various moduli spaces involved. 
We have 
\[  
ID (\Omega^n \otimes 1_{\theta_0}) =  
\sum_{a_{(-2)}}  n_{\theta_0 a_{(-2)}} n_{a_{(-2)}a_{(-2)}'}   
\Omega^n \otimes 1_{a_{(-2)}'}\]  
and 
\[  
DI (\Omega^n \otimes 1_{\theta_0}) =  
- \sum_{a_{(-1)}'}   n_{\theta_0 a_{(-1)}'} n_{a_{(-1)}'a_{(-2)}'}  
\Omega^n \otimes 1_{a_{(-2)}'} \] \[ +  
\sum_{a_{(0)}'}  n_{\theta_1 a_{(0)}'}m_{a_{(0)}'a_{(-2)}'}  
\Omega^n \otimes 1_{a_{(-2)}'}.  
\]  
 
The first claim then follows, since we know that the co-dimension one 
boundary of  $\M (\theta_0, O_{a_{(-2)}'})$ consists of  
\[ \bigcup_{a_{(-2)}} \hat \M (\theta_0, O_{a_{(-2)}}) 
\times_{O_{a_{(-2)}}} \M (O_{a_{(-2)}},O_{a_{(-2)}'})   \] 
\[ \cup \bigcup_{a_{(-1)}'} \M (\theta_0, O_{a_{(-1)}'}) \times  \hat\M 
(O_{a_{(-1)}}, O_{a_{(-2)}'})  \] 
\[ \cup \M(\theta_0,\theta_1) \times \bigl( \tilde s^{-1}(0)\cap (\hat 
\M (\theta_1, O_{a_{(0)}'})  \times_{O_{a_{(0)}'}} \M (O_{a_{(0)}'}, 
O_{a_{(-2)}'}) \bigr). \]    
 
In fact, according to Theorem \ref{noglue}, we have boundary strata as  
in (\ref{boundary:t:obstr:1}), 
$$ \M(\theta_0,\theta_1)\times \partial^{(1)}\left( \tilde 
s^{-1}(0) \cap \M (\theta_1,O_{a_{(-2)}'})^* \right). $$ 
Here $\tilde s$ is the obstruction section.  
By Remark \ref{sec:obstr:induction}, we have 
\[ \begin{array}{lll} 
&&\partial^{(1)}\left( \tilde 
s^{-1}(0) \cap \M (\theta_1,O_{a_{(-2)}'})^* \right) \\[1mm] 
&=& \tilde 
s^{-1}(0)\cap  
\bigl( \bigcup_{a_{(0)}'}\hat 
\M (\theta_1, O_{a_{(0)}'})^*  \times_{O_{a_{(0)}'}} \M (O_{a_{(0)}'}, 
O_{a_{(-2)}'})^*\bigr)\\[1mm] 
&&  \cup \tilde 
s^{-1}(0)\cap \bigl( \bigcup_{a_{(-1)}'} \hat\M (\theta_1, 
O_{a_{(-1)}'})^*\times_{O_{a_{(-1)}'}} \M 
(O_{a_{(-1)}'},O_{a_{(-2)}'})^*\bigr).  
\end{array} 
\] 
 Note that the coefficient 
$n_{ \theta_0 a_{(-1)}'}$ effectively counts the monopoles 
from the zeros of the obstruction section over 
$\M ( \theta_1, a_{(-1)}')$,  
 so we are left with the remaining counting, 
\[ 
\sum_{a_{(0)}'}  n_{\theta_1 a_{(0)}'}m_{a_{(0)}'a_{(-2)}'}  
\Omega^n \otimes 1_{a_{(-2)}'}, 
\] 
 which proves the claim.

In order to prove Claim (2), we compute 
\[  
ID (\Omega^n \otimes 1_{a_{(1)}}) = - \sum_{a_{(0)}}  
n_{a_{(1)}a_{(0)}} n_{a_{(0)} a_{(0)}'}  
\Omega^n \otimes 1_{a_{(0)}'}  
\]  
and   
\[  
DI (\Omega^n \otimes 1_{a_{(1)}}) =   
 - \sum_{a_{(1)}'}  
n_{a_{(1)}a_{(1)}'} n_{a_{(1)}' a_{(0)}'}\Omega^n \otimes 1_{a_{(0)}'}  
+ n_{a_{(1)}\theta_0} n_{\theta_1a_{(0)}'}  
\Omega^n \otimes 1_{a_{(0)}'}  
\]  
Then the claim follows from the boundary structure of $\M (a_{(1)}, 
a_{(0)}')$. The zeroes of the obstruction section contribute the term 
$$ n_{a_{(1)}\theta_0} n_{\theta_1a_{(0)}'} $$  
in the counting of the boundary points. 
 
The argument for Claims (3) and (4) is analogous. 
  
\noindent{$\diamond$}  
 
Now we need to check the effect of the correction terms on the 
identity $1-JI =HD +HD$.

\begin{lem} \label{check:1} 
We have the identity 
$1-JI =DH+HD$ on $\Omega^n \otimes 1_{\theta_0}$.  
\end{lem} 
 
\noindent\underline{Proof:}  By direct calculation, we obtain  
\[  
\begin{array}{lll}  
&& (1-JI) (\Omega^n \otimes 1_{\theta_0}) \\[2mm]  
&=& \Omega^n \otimes 1_{\theta_0} -  
 J (\Omega^{n-1} \otimes 1_{\theta_1}  + 
\sum_{a_{(-1)}'}  n_{\theta_0a_{(-1)}'}  
\Omega^n\otimes 1_{a_{(-1)}'} \\[2mm]  && 
+ \sum_{a_{(0)}'}    n_{\theta_1 a_{(0)}'}  \Omega^n\otimes 
\eta_{a'_{(0)}})\\[2mm]   
&=& \Omega^n \otimes 1_{\theta_0} - \sum_{a_{(1)}}  n_{\theta_1 a_{(1)}}  
 \Omega^{n-1}\otimes 1_{a_{(1)}} \\[2mm]  
&&-   
\bigl( \sum_{a'_{(-1)}} n_{\theta_0a'_{(-1)}}n_{a_{(-1)}'a_{(-1)}}  
+ \sum_{a'_{(0)}} n_{\theta_1a'_{(0)}}m_{a_{(0)}'a_{(-1)}}\bigr) 
\Omega^{n}\otimes  
1_{a_{(-1)}} \\[2mm]   
&&-  
\sum_{a'_{(0)}} 
n_{\theta_1a'_{(0)}}n_{a_{(0)}'a_{(0)}}\Omega^{n}\otimes 
\eta_{a_{(0)}}   
-\sum_{a'_{(0)}} n_{\theta_1a'_{(0)}}n_{a_{(0)}'\theta_0} \Omega^n 
\otimes 1_{\theta_0}   
\end{array}  
\] 
 
We also get 
\[  
\begin{array}{lll}  
&& (DH+HD) (  (\Omega^n \otimes 1_{\theta_0}) \\[2mm]  
&=& D(  \sum_{a_{(0)}} n^P_{\theta_0 a_{(0)}}  
\Omega^n \otimes 1_{a(0)} + 
 \sum_{a_{(1)}}    n_{\theta_1a_{(1)}}\Omega^n \otimes \eta_{a_{(1)}}
) \\[2mm]  && +   
H (  \sum_{a_{(-2)}}  n_{\theta_0a_{(-2)}} \Omega^n \otimes 
1_{a_{(-2)}})  \\[2mm]   
&=& \sum_{a_{(-1)}} \bigl ( \sum_{a_{(1)}} n_{\theta_1 a_{(1)}} 
m_{a_{(1)} a_{(-1)}}-   
\sum_{a_{(0)}} n^P_{\theta_0 a_{(0)}}n_{a_{(0)}a_{(-1)}} \\[2mm]      
&&- \sum_{a_{(-2)}} n_{\theta_0a_{(-2)}}n^P_{a_{(-2)}a_{(-1)}} 
\bigr)\Omega^n \otimes 1_{a(-1)}    
 \\[2mm]  
&& + \sum_{a_{(1)},a_{(0)}} n_{\theta_1a_{(1)}}n_{a_{(1)}a_{(0)}}    
 \Omega^{n}\otimes \eta_{a_{(0)}}- \sum_{a_{(1)}}  
n_{\theta_1a_{(1)}}\Omega^{n-1} \otimes 1_{a_{(1)}} \\[2mm]  
&&  
+ \sum_{a_{(1)}} n_{\theta_1a_{(1)}}n_{a_{(1)}\theta_0}  
\Omega^n \otimes 1_{\theta_0}.  
\end{array}  
\]   
 
Let us first check the coefficient of $\Omega^n\otimes 
1_{\theta_0}$. Equating the coefficients of $\Omega^n\otimes 
1_{\theta_0}$ in the two expressions above 
yields the identity 
\[  
\sum_{a_{(1)}} n_{\theta_1 a_{(1)}} n_{a_{(1)} \theta_0}  
+ \sum_{a_{(1)}'} n_{\theta_1 a_{(0)}'} n_{a_{(0)}' \theta_0}   
=1.  
\]  
This identity is satisfied, since it is exactly the 
counting of (\ref{nice:1}). 
 
We then check the coefficient of $ \Omega^{n}\otimes 
\eta_{a_{(0)}}$. We need to prove the identity 
\[ \sum_{a_{(1)}} 
n_{\theta_1 a_{(1)}} n_{a_{(1)}a_{(0)}}  + 
\sum_{a_{(0)}'} n_{\theta_1 a'_{(0)}} n_{a_{(0)}'a_{(0)}} =0. \] 
This follows from the boundary structure of $\M(\theta_1, O_{a_{(0)}})^*$,  
given by  
$$ \bigcup_{a_{(1)}} \M(\theta_1, O_{a_{(1)}})^* \times_{O_{a_{(1)}}}  
\hat \M (O_{a_{(1)}}, O_{a_{(0)}})^*  $$  
$$ \cup \bigcup_{a'_{(0)}} \hat \M \theta_1, 
O_{a_{(0)}'})^*\times_{O_{a_{(0)}'}}  
\M (O_{a_{(0)}'}, O_{a_{(0)}} )^*. $$ 
 
In order to compare the coefficients of $\Omega^n \otimes 1_{a(-1)} $,  
we need to prove the identity 
\beq  
\begin{array}{lll}  
&& \sum_{a_{(-1)}'} 
n_{\theta_0 a_{(-1)}'}n_{a_{(-1)}'a_{(-1)}} -  
\sum_{a_{(0)}} n^P_{\theta_0 a_{(0)}} n_{a_{(0)}a_{(-1)}} \\[2mm] 
&-& \sum_{a_{(-2)}}  n_{\theta_0a_{(-2)}} n^P_{a_{(-2)}a_{(-1)}} \\[2mm]  
&+& \sum_{a'_{(0)}} n_{\theta_1a'_{(0)}}m_{a_{(0)}'a_{(-1)}}  
+\sum_{a_{(1)}} n_{\theta_1a_{(1)}}  m_{a_{(1)} a_{(-1)}} =0  
\end{array}  
\label{vanish:a_1}  
\eeq  
 
Among these terms we can isolate a contribution 
\[ \sum_{a_{(-1)}'} 
n_{\theta_0 a_{(-1)}'} n_{a_{(-1)}'a_{(-1)}} - \sum_{a_{(0)}} 
n^P_{\theta_0 a_{(0)}}n_{a_{(0)}a_{(-1)}}       - \sum_{a_{(-2)}}  
 n_{\theta_0a_{(-2)}}n^P_{a_{(-2)}a_{(-1)}},  
\] 
which is the contribution of the boundary terms in 
$\M^P (\theta_0, O_{a_{(-1)}})^*$ obtained by gluing  
co-dimension one boundary strata along irreducible critical orbits.   
The term 
$$ \sum_{a_{(-1)}'}  n_{\theta_0 a_{(-1)}'} n_{a_{(-1)}' a_{(-1)}} $$  
counts the contribution  of  the special gluing 
\[ \M (\theta_0, \theta_1) \times \M (\theta_1, a_{(-1)}') \times  
\M (a_{(-1)}', a_{(-1)}), \] 
where the coefficient 
$$ n_{\theta_0 a_{(-1)}'} =\# \bar s^{-1}(0) $$ 
counts the zeroes of the obstruction section in $\M (\theta_1, 
a_{(-1)}')$.   
 
According to Theorem \ref{boundary:P:obstr}, the counting of the 
remaining boundary  
components for $\M^P (\theta_0, a_{(-1)})$ is given by the gluing 
$$ \M (\theta_0, \theta_1) \times s^{-1}(0), $$ 
with 
\[  s^{-1}(0)\subset \M (\theta_1, a_{(-1)}).\] 
Recall that, in this case, $s^{-1}(0)$ is zero-dimensional. 
 
We claim that the zeroes of this obstruction section are counted precisely 
by the expression  
\[ \sum_{a_{(0)}'} 
 n_{\theta_1 a_{(0)}'}m_{a_{(0)}'a_{(-1)}}  
+\sum_{a_{(1)}}  n_{\theta_1 a_{(1)}}  m_{a_{(1)} a_{(-1)}}.   
\]  
In fact, we can describe the  
obstruction bundle as the line bundle 
\[  
 \M (\theta_1, O_{a_{(-1)}} )^* \times _{U(1)} Coker ({\cal L}_{x}),  
\]  
over $\M (\theta_1,a_{(-1)})^*$, with $x=[\AA_0,0]$ the singular 
reducible in $\M(\theta_0,\theta_1)$, and $Coker ({\cal 
L}_{x})=\CC$. With the choice of trivilizations, 
the zeros of the obstruction section  are 
localized in some compact set in  
$$ \bigcup_{a_{(0)}'} \M (\theta_1, a_{(0)}') \times \M ( 
a'_{(0)},  a_{(-1)}) $$ 
$$ \cup \bigcup_{a_{(1)}} \M(\theta_1,a_{(1)})\times  \M 
(a_{(1)},a_{(-1)}). $$ 
 
Thus, counting the zeros of this obstruction section, which gives the
relative Euler class of the associated
complex line bundle over $\M (\theta_1, a_{(1)})$,
contributes 
$$ \begin{array}{lr} \#  s^{-1}(0)=&  \sum_{a_{(0)}'} n_{\theta_1, 
a_{(0)}'} m_{a_{(0)}'a_{(-1)}} \\[2mm] 
 & +\sum_{a_{(1)}} n_{\theta_1 a_{(1)}}  m_{a_{(1)} a_{(-1)}}, 
\end{array}  
$$ 
 
Then the vanishing condition of (\ref{vanish:a_1}) follows  
from $\# (\partial \M^P (\theta_0, a_{(-1)}) ) =0$.  
Thus, we have completed the proof of the identity $1-JI =DH+HD$  
on the generator $\Omega^n \otimes 1_{\theta_0}$.  
 
\noindent{$\diamond$}  
 
We proceed to check the remaining identities. 
 
\begin{lem} \label{check:2} 
We have the identity $1-JI =DH+HD$ on $\Omega^n \otimes 1_{a_{(1)}}$. 
\end{lem} 
 
\noindent\underline{Proof:} We compute 
\[ \begin{array}{lll}  
&&(1-JI) (\Omega^n \otimes 1_{a_{(1)}}) \\[2mm]  
&=& \Omega^n \otimes 1_{a_{(1)}} -  J\bigl( n_{a_{(1)}\theta_0} 
\Omega^{n-1}  
\otimes 1_{\theta_1}  + \sum_{a_{(1)}'} n_{a_{(1)}a_{(1)}'} 
\Omega^{n-1} \otimes 1_{a'_{(1)}}  
 \bigr)\\[2mm]  
&=&  \Omega^n \otimes 1_{a_{(1)}} -  \sum_{\tilde a_{(1)}} 
( \sum_{a_{(1)}'} n_{a_{(1)}a_{(1)}'} n_{a_{(1)}'\tilde a_{(1)}}  
+  n_{a_{(1)}\theta_0} n_{\theta_1 \tilde a_{(1)}}   )  
\Omega^n \otimes 1_{\tilde a_{(1)}}   
\end{array}  
\]  
\[ \begin{array}{lll} && 
(DH+HD) (\Omega^n \otimes 1_{a_{(1)}}) \\[2mm]   
&=&  (\sum_{a_{(2)}} 
 n^P_{a_{(1)}a_{(2)}} n_{a_{(2)}\tilde a_{(1)}} + \sum_{a_{(0)}} 
n_{a_{(1)}a_{(0)}} n^P_{a_{(0)}\tilde a_{(1)}})  
\Omega^n \otimes 1_{\tilde a_{(1)}}. \end{array} 
\]  
 
If we compare the coefficient of each term in $DH+HD-1+JI$, we see that  
the coefficient of $\Omega^n \otimes 1_{a_{(1)}}$ is given by  
\[  
\sum_{a_{(2)}} n^P_{a_{(1)}a_{(2)}} n_{a_{(2)}  a_{(1)}}   
+ \sum_{a_{(0)}}   
 n_{a_{(1)}a_{(0)}} n^P_{a_{(0)}  a_{(1)}}  \]
\[  + \sum_{a_{(1)}'}   
n_{a_{(1)}a_{(1)}'} n_{a_{(1)}'  a_{(1)}} -1   
+  n_{a_{(1)}\theta_0} n_{\theta_1   a_{(1)}}.   
\]  
Notice that the sum 
\[ \sum_{a_{(2)}} n^P_{a_{(1)}a_{(2)}} n_{a_{(2)}  a_{(1)}}   
+ \sum_{a_{(0)}}   
 n_{a_{(1)}a_{(0)}} n^P_{a_{(0)} a_{(1)}}  
 + \sum_{a_{(1)}'}   
n_{a_{(1)}a_{(1)}'} n_{a_{(1)}'  a_{(1)}} -1  
 \] 
is the algebraic counting  of the  boundary points in  
$\M^P(a_{(1)},  a_{(1)})^*$ that correspond to boundary components of 
the form 
\[  
\bigcup_{a_{(2)}}   
\M^P ( a_{(1)},  a_{(2)}  ) \times   
\hat \M (a_{(2)},  a_{(1)}) \cup \bigcup_{ a_{(0)}}  
\hat \M (a_{(1)}, a_{(0)}) \times \M^P (a_{(0)}, a_{(1)})  
\]  
\[  
\cup \bigcup_{a_{(1)}'} \M (a_{(1)}, a_{(1)}') \times 
   \M  ( a_{(1)}' ,  a_{(1)}) \cup \{ -a\}.  
\]  
 
Similarly, for $\tilde a_{(1)} \neq   a_{(1)}$,   
the coefficient of $\Omega^n \otimes 1_{ \tilde a_{(1)}}$ is given by  
\[  
\sum_{a_{(2)}} n^P_{a_{(1)}a_{(2)}} n_{a_{(2)} \tilde   a_{(1)}}   
+ \sum_{a_{(0)}}   
 n_{a_{(1)}a_{(0)}} n^P_{a_{(0)}\tilde a_{(1)}}  
 + \sum_{a_{(1)}'}   
n_{a_{(1)}a_{(1)}'} n_{a_{(1)}'\tilde a_{(1)}} +  n_{a_{(1)}\theta_0} 
n_{\theta_1 \tilde a_{(1)}}.    
\]  
 
Again, notice that the sum 
\[  
\sum_{a_{(2)}} n^P_{a_{(1)}a_{(2)}} n_{a_{(2)} \tilde  a_{(1)}}   
+ \sum_{a_{(0)}}   
 n_{a_{(1)}a_{(0)}} n^P_{a_{(0)} \tilde a_{(1)}}  
 + \sum_{a_{(1)}'}   
n_{a_{(1)}a_{(1)}'} n_{a_{(1)}' \tilde  a_{(1)}}   
 \]  
is the algebraic counting of the  boundary points of  
$\M^P ( a_{(1)}, \tilde a_{(1)})^*$ which correspond to boundary 
components  of the form 
\[  
\bigcup_{a_{(2)}}   
\M^P ( a_{(1)},  a_{(2)}  ) \times   
\hat \M (a_{(2)}, \tilde  a_{(1)}) \cup \bigcup_{ a_{(0)}}  
\hat \M (a_{(1)}, a_{(0)}) \times \M^P ( a_{(0)},\tilde a_{(1)})\]  
\[  
\cup \bigcup_{a_{(1)}'} \M ( a_{(1)}, a_{(1)}') \times 
   \M ( a_{(1)}' , \tilde  a_{(1)}) .  
\] 
 
We only need to prove that the counting of the remaining boundary 
components of $\M^P(a_{(1)}, \tilde a_{(1)})$ is given by  
\[  
n_{a_{(1)}\theta_0} n_{\theta_1 \tilde a_{(1)}},   
\]  
which counts the contribution of the zeros of the obstruction section,  
that is, of those monopoles of $\M^P(a_{(1)}, \tilde a_{(1)})$ which 
are obtained by gluing the singular reducible $x\in 
\M(\theta_0,\theta_1)$ with  
$$ \M (a_{(1)}, \theta_0) \times U(1)\times \M (\theta_1, \tilde  
a_{(1)}). $$  
 
Recall that the non-equivariant gluing 
$$ \# : \bar s^{-1}(0)\subset \M (a_{(1)}, \theta_0) \times U(1)\times \M 
(\theta_1, \tilde  a_{(1)}) \to \M^P (a_{(1)},\tilde  a_{(1)}) $$ 
corresponds to the equivariant gluing in the framed moduli spaces 
$$ \# : \tilde s^{-1}(0)\subset \M (O_{a_{(1)}}, \theta_0) \times \M 
(\theta_1, O_{\tilde a_{(1)}}) \to \M^P( O_{a_{(1)}},O_{\tilde 
a_{(1)}}). $$ 
 
The pull-back and push-forward map   
$$(e_{\tilde a_{(1)}}^-)_* (e^+_{ a_{(1)}})^*(\eta_{a_{(1)}})$$  
defines the relative Euler characteristic number  
on the associated line bundle  
of   
\[  
{\cal L}_{a_{(1)}, \tilde a_{(1)}} =  
(\M (O_{a_{(1)}}, \theta_0) \times \M (\theta_1, O_{\tilde a_{(1)}}) ) 
\times _{U(1)} \CC.   
\]  
This gives  
\[  
(e_{\tilde a_{(1)}}^-)_* (e^+_{a_{(1)}})^*(\eta_{a_{(1)}})  
= n_{a_{(1)}\theta_0} n_{\theta_1 \tilde a_{(1)}}. \]  
 
Consider  the obstruction bundle    
\[  
\bigl( \M (O_{a_{(1)}}, \theta_0) \times \M (\theta_0, \theta_1)  
\times \M (\theta_1, O_{\tilde a_{(1)}})  
\bigr) \times _{U(1)} Coker({\cal L}_x) \] 
over the space   
\[  
\M (a_{(1)}, \theta_0) \times \M (\theta_1, \tilde  
a_{(1)})\times U(1),  
\]  
with $Coker({\cal L}_x)=\CC$, for $x\in \M(\theta_0,\theta_1)$ the 
singular point $x=[\AA_0,0]$. 
Counting the zeros of the obstruction section 
\[  
\# \bigl ( \bar s^{-1}(0) \cap   
(\M (a_{(1)}, \theta_0) \times \M (\theta_1, \tilde  
a_{(1)})\times U(1) )\bigr)  
\]   
gives another computation of the same relative Euler class above, 
that is, 
\[ \begin{array}{lll}  
&&\# \bigl ( \bar s^{-1}(0) \cap   
( \M (a_{(1)}, \theta_0) \times \M (\theta_1, \tilde  
a_{(1)})\times U(1) )\bigr)  \\[2mm]  
&=&   
n_{a_{(1)}\theta_0} n_{\theta_1 \tilde a_{(1)}}.  
\end{array}  
\]  
 
Thus,  the components of the boundary of $\M ^P ( a_{(1)}, \tilde 
a_{(1)})^*$ which come from the gluing with the singular $x\in 
\M(\theta_0,\theta_1)$ contribute a term 
\[  
n_{a_{(1)}\theta_0} n_{\theta_1 \tilde a_{(1)}}  
\]  
to the counting of the boundary points. This completes  
the proof of the Lemma.  
 
\noindent{$\diamond$}  
 
Finally, we have to check the following. 
 
\begin{lem} 
We have the identity $1-JI =DH+HD$ on $\Omega^n \otimes 
\eta_{a_{(1)}}$.  \label{check:3} 
\end{lem} 
 
\noindent\underline{Proof:} Direct calculation of the terms 
$(1-JI)(\Omega^n \otimes \eta_{a_{(1)}})$ and $(DH+HD)(\Omega^n 
\otimes \eta_{a_{(1)}})$ shows that the conditions required in order 
to have the same coefficients on all the generators are precisely the 
conditions already verified in the case of Lemma \ref{check:2}. 
 
\noindent{$\diamond$}

Notice  how clearly this argument of topological invariance
breaks down for the non-equivariant Floer homology. The
invariance within the same chamber is still verified: in fact, no
substantial changes are necessary in that first part of the proof, 
in order to adapt it to the case of the non-equivariant Floer
homology. However, as wee see clearly from the structure of this
second part of the proof, the general argument for the proof of
Theorem \ref{metrics}, for metrics and perturbations in two different
chambers, relies essentially on the contribution of the reducible
points, in order to construct the chain
map $I$ and chain homotopy $H$, as discussed in Definition
\ref{def:I:H:change}, Lemma \ref{check:ID=DI}, Lemma 
\ref{check:1},  Lemma \ref{check:2}, and  Lemma \ref{check:3}.
The example presented at the beginning of Section 6.3 also clarifies
why the argument cannot be adapted to the non-equivariant Floer
homology. 

\section{Wall crossing formula for the Casson-type invariant}

We want to compare the equivariant Floer homology with the ordinary
Floer homology in the cases where the latter is defined, i.e. when
$b^1(Y)$ is non-trivial \cite{Ma} or when $Y$ is a homology sphere,
\cite{Wa}. 

In the case when $b^1(Y)$ is non-trivial, we expect to find that the
equivariant Floer homology, which is computed by considering framed
moduli spaces, is isomorphic to the ordinary Floer homology computed
in the unframed space. In fact this is the analogue of the well
known result for equivariant cohomology of a finite dimensional
manifold, where, if the action of the group is free, then
the equivariant cohomology is just the ordinary cohomology of the
quotient, $H^*_G(M;\RR)\cong H^*(M/G;\RR)$ as $H^*(BG;\RR)$-modules.

In the case of a homology sphere, instead, we expect to find an
exact sequence that connects the equivariant Floer homology with the
ordinary Floer homology and an extra copy of $\RR [\Omega]$ that
corresponds to the unique reducible solution that has been removed in
the computation of the non-equivariant Floer homology.

Recall that we have the explicit description of the boundary operator
in the equivariant Floer complex, as analyzed in Section 5,
Proposition \ref{boundaryterms}, which gives the formula
(\ref{equivD}), 
\[ D: \begin{array}{ccc}
\Omega^n \otimes 1_a &\mapsto & -n_{ab}\Omega^n \otimes 1_b \\
                        &        & \\
\Omega^n \otimes \eta_a &\mapsto & (n_{ab}\Omega^n \otimes \eta_b)\oplus
                                   (m_{ac}\Omega^n \otimes 1_c) \\
                        &        & \oplus(-\Omega^{n-1}\otimes 1_a).
\end{array} \]

\subsection{Comparison with the non-equivariant Floer Homology}

Now we can define a chain homomorphism that maps the equivariant
to the non-equivariant complex.

Let us first work in the case with
no reducible solution (i.e. with $b^1(Y)>0$). In this case for each
$O_a$ that appears in the equivariant complex we have a generator $\RR
a$ that appears in the non-equivariant complex (coefficients in
$\RR$). 

Now we define the chain map
$$ i_{k} : C_{k,U(1)}(Y)\to C_{k}(Y),$$
so that it satisfies $\partial_{k}i_{k}=i_{k-1}D_k$. 
Let $i_{k}$ act on the generators as follows 
\beq
i_{k} :\bigoplus_{\mu(a)+j=k} \Omega_{j,U(1)}(O_a)
\to \sum_{\mu(O_{a})=k} \RR a, 
\label{imap}
\eeq
$$ i_{k} (\Omega^n\otimes 1_a)=0, $$
for all values of $n$ and $\mu(O_a)$,
$$ i_{k} (1\otimes \eta_a)=a, $$
if $\mu(O_a)=k$, and in all other cases
$$ i_{k} (\Omega^n\otimes \eta_a)=0.$$

This means that the map $i_k$ kills all the generators in degree $k$
that are not the generator of the equivariant homology of some orbit $O_a$
of degree $k$.

With this definition it is clear that $i_*$ is a chain map.
Thus it defines a sub-complex of
$C_{*,U(1)}(Y)$ given by $Q_*=Ker(i_*)$ with the restriction of the
boundary operator $D$.

\begin{thm}
If there are no reducible solutions (i.e. $b^1(Y)>0$) the map $i_*$
defined in (\ref{imap}) induces an isomorphism in cohomology,
\[ SWH_{*,U(1)}(Y)\cong SWH_*(Y). \]
\label{isomorphism}
\end{thm}

\noindent\underline{Proof:}
The complexes  $C_{k,U(1)}$, $Q_k$, and $C_k$ all have a filtration by
index. For $C_{k, U(1)}$ the filtration is given by
\beq
C_{k,U(1)}(n)=\bigoplus_{\mu(O_a)+j=k, \mu(O_a)\leq n} \Omega_{j,
  U(1)}(O_a). 
\label{equivfiltr}
\eeq

The complex $Q_*$ is written as
\beq
\label{Qk}
Q_k = \bigoplus_{\mu(a)+j=k, j\geq 1} \Omega_{j,U(1)}(O_a).
\eeq 

It has a filtration by index of the form
\beq
Q_k(n)=\bigoplus_{\mu(a)+j=k, j\geq 1, \mu(a)\leq n}\Omega_{j,U(1)}(O_a).
\label{Qfiltr}
\eeq

On the other hand also the non-equivariant complex has a filtration by
index of the form
\beq
\label{Cfiltr}
C_k(n)=\bigoplus_{\mu(a)=k\leq n} \RR a.
\eeq

Thus we can look at the spectral
sequences associated to the filtrations and prove that $i$ induces a
map of spectral sequences and an isomorphism of the $E^1$-terms of the
spectral sequences associated to the filtration of $C_{*,U(1)}$ and of $C_*$.
Thus we get the resulting isomorphism of the $E^\infty$ terms, i.e. of 
the homology of $C_{*,U(1)}$ and of $C_*$.

\begin{lem} \label{spectr:seq1}
Let $E_{U(1)}$, $E_Q$, and $E$ be the spectral sequences associated to
the filtration of the complexes $C_{U(1)}$, $Q$ and $C$ respectively.
The chain map $i_*$ induces a map of spectral sequences. 
Moreover, in the case when $b^1(Y)>0$ the map $i_*$ induces an
isomorphism of the $E^1$-terms 
\[ E^1_{kl, U(1)}\cong E^1_{kl}. \]
\end{lem}

\noindent\underline{Proof:}

Consider the filtration of $C_{*,U(1)}$. The $E^0_{kl}$ terms
of the spectral sequence are given by
\[ E^0_{kl, U(1)}= C_{k+l,U(1)}(k)/C_{k+l, U(1)}(k-1). \]
{}From (\ref{equivfiltr}) we get
\[ E^0_{kl, U(1)}= \bigoplus_{\mu(O_a)=k}\Omega_{l,U(1)}(O_a), \]
and the differentials $E^0_{kl, U(1)}\to E^0_{kl-1, U(1)}$ are just
given by the differential of the equivariant complex on each fixed
orbit $O_a$. Thus the $E^1$-term of this spectral sequence is given by
\[ E^1_{kl, U(1)}=H_{k+l}(E^0_{k*, U(1)}), \]
\[ E^1_{kl, U(1)}=\bigoplus_{\mu(O_a)=k}H_{k+l, U(1)}(O_a). \]
But since the $O_a$ are irreducible orbits, with a free $U(1)$ action,
the equivariant homology is concentrated in degree zero,
\[ E^1_{k0, U(1)}=\bigoplus_{\mu(O_a)=k}\RR a, \]
\[ E^1_{kl, U(1)}=0 \]
for $l\neq 0$.

Now let us consider the filtration of the non-equivariant complex $C_*$.
{}From (\ref{Cfiltr}) we get
\[ E^0_{k0}= \bigoplus_{\mu(O_a)=k} \RR a \]
and $E^0_{kl}=0$ for $l\neq 0$. Thus in this case the only terms that
survive in the $E^1$ is in degree $l=0$ and is
\[ E^1_{k0}=\bigoplus_{\mu(O_a)=k} \RR a. \]

\noindent $\diamond$

Thus the map $i_*$ induces an isomorphism on the $E^1$-terms of the
two spectral sequences, hence on the $E^\infty$-terms, namely on the
homology 
\[ H_*(C_{k,U(1)}(Y),D_k)\stackrel{i}{\cong}
H_*(C_k(Y),\partial_k). \]
This completes the proof of theorem \ref{isomorphism}.

\noindent $\diamond$

Now let us consider the case when the manifold $Y$ is an integral
homology sphere. In this case there are also terms in the equivariant
complex that come from the reducible solution $\theta=[\nu,0]$. We
assume that $\theta$ has index zero, $\mu(\theta)=0$.

\begin{thm}
Let $Y$ be an integral homology sphere. Then there is an exact
sequence
\[ \cdots\to H_{k,U(1)}(\theta)\to SWH_{k,U(1)}(Y)\stackrel{i}{\to}
SWH_k(Y)\stackrel{\Delta}{\to}H_{k-1,U(1)}(\theta)\cdots \]
\label{exact}
\end{thm}

In the equivariant complex in degree $k$ we have an extra generator
$\Omega^k\otimes \theta$. The boundary maps that come from the
equivariant complex associated to the degenerate orbit $\nu$ with
the trivial action of $U(1)$ are trivial: in fact the
equivariant homology of a point is
\[ H_{*,U(1)}(\theta)=\RR[\Omega]=H_*(BU(1), \RR). \]
However, there are non-trivial boundary maps that hit the generators
$\Omega^n\otimes \theta$. These can be described as follows.
Suppose $O_a$ is the orbit of an irreducible solution with index
$\mu(O_a)=1$. Then we have a moduli space $\hat{\cal M}(O_a,\theta)$
that is 1-dimensional and that fibers over $\theta$ with a
1-dimensional fiber. Thus the pullback-pushforward map acts as
\[ \eta_a\mapsto (e_\theta^-)_*(e_a^+)^*\eta_a=m_{a\theta}, \]
where $m_{a\theta}$ is the integration along the 1-dimensional fiber
of the 1-form $(e_a^+)^*\eta_a$. This gives rise to a component of the
boundary map of the form
\[ 1\otimes\eta_a\mapsto m_{a\theta}1\otimes \theta. \]
Moreover, there is a non-trivial boundary map that comes from the
moduli spaces that connect the reducible to generators with the orbits
with  $\mu(O_a)=-2$.

Now the map $i_k$ is defined as before with the additional condition
that it kills the extra generator $\Omega^k\otimes \theta$.

\begin{lem} \label{spectr:seq2}
Let $Y$ be a homology sphere. Then the homology of the complex $Q_*$
is just the equivariant homology of the point $\theta$,
\[ H_*(Q)\cong H_{*,U(1)}(\theta)=\RR [\Omega]. \]
\end{lem}

\noindent\underline{Proof:}
The complex $Q_*$ contains the extra generator $\Omega^k\otimes \theta$
in degree $k$ and this generator appears in all levels of the
filtration $Q_k(n)$ for any $n\geq 0$, since $\theta$ is of degree
zero.

Thus if we look at the spectral sequence associated to the filtration
of the complex $Q_*$ we find
\[ E^0_{kl,Q}= Q_{k+l}(k)/Q_{k+l}(k-1), \]
that is, for $k>0$ and $l\geq 0$,
\[ E^0_{kl,Q}=\bigoplus_{\mu(O_b)=k-1} \Omega_{l+1,U(1)}(O_b). \]
This complex is clearly acyclic because the differentials are just the
equivariant differentials for each orbit and no generator survives in
homology because we are counting only $l+1\geq 1$, hence the terms
$1\otimes \eta_b$ are suppressed (in fact they are not in $Ker(i_*)$).
On the other hand for $k=0$ we get
\[ E^0_{0l,Q}=\RR <\Omega^l\otimes \theta>, \]
with trivial differentials, so that the $E^1$-terms are
\[ E^1_{kl,Q}=H_{k+l, U(1)}(\theta). \]

This means that the homology of the complex $Q_*$ is actually
$H_{*,U(1)}(\theta)=\RR[\Omega]$. 

\noindent $\diamond$

Thus if we consider
the long exact sequence induced by the short sequence
$$ 0\to Q_k\to C_{k,U(1)}\stackrel{i}{\to} C_k\to 0 $$ 
we have

\[ \cdots\to H_{k,U(1)}(\theta)\to SWH_{k,U(1)}(Y)\stackrel{i}{\to}
SWH_k(Y)\to\cdots \]

This proves theorem \ref{exact}. 

\noindent $\diamond$

The connecting homomorphism $\Delta$ in the above long exact sequences  is
particularly interesting. 

\begin{prop}\label{connect:homom:prop}
Suppose given a representative
$\sum_a x_a a$ in
$SWH_{2k+1}(Y)$. We assume that $Y$ is a homology sphere. Then the connecting
homomorphism in the long exact sequence is 
\[
\bigtriangleup_{k} : \qquad SWH_{2k+1} \to H_{2k, U(1)}(\theta )= \RR \Omega^k,
\label{diagonalmap}
\]
where
\beq
\bigtriangleup_{k}(\sum_a x_a a) =\sum x_a m_{ac} m_{ce}\cdots
m_{\alpha' \alpha} n_{\alpha\theta} \Omega^k\otimes \theta. 
\label{representative}
\eeq
Here the sum is understood over all the repeated indices, that is over
all critical points with indices $\mu(O_a)=2k+1$, $\mu(O_c)=2k-1$,
$\mu(O_{\alpha'})=3$, $\mu(O_\alpha)=1$. 
\end{prop}

\noindent\underline{Proof:} 
The map is defined by the standard diagram chase and by adding 
boundary terms in order to find a representative of the form
(\ref{representative}), as illustrated in the following diagram. Sums
over repeated indices are understood.

$$
\spreaddiagramrows{-1pc}
\spreaddiagramcolumns{-1.75pc}
\diagram
& & \vdots \dto & & & \vdots \dto & & & \vdots \dto \\
0 \rrto & & Q_{2k+1} \xto[0,3] \xto[3,0] & & & C_{2k+1,U(1)}
\xto[0,3]^{\qquad i}
\xto[3,0] & & & C_{2k+1} \rrto \xto[3,0] & & 0 \\
& & & & & & x_a~1\otimes\eta_a \xdashed[3,0]|<\stop|>\tip & & & x_a~a
\xdashed[0,-3]|<\stop|>\tip \\
\\
0 \rrto & & Q_{2k} \xto[0,3] \ddto & & & C_{2k,U(1)} \xto[0,3]^{\qquad i} \ddto
& & & C_{2k} \rrto \ddto & & 0 \\
& & & \Delta_k(x_a~a) & & & x_am_{ac}~1\otimes 1_c
\xdashed[0,-3]|<\stop|>\tip \\
& & \vdots & & & \vdots & & & \vdots
\enddiagram
$$

For a cycle
$\sum _{a: \mu(O_a)= 2k+1} x_a a $ we have $\sum _{a: \mu (a)= 2k+1} x_a
n_{ab} = 0 $ for any $b$ with $\mu(O_b) = 2k$. The element $\sum_a x_a
a$ has a preimage $\sum x_a 1\otimes \eta_a$ under $i_{2k+1}$. The
image of this element under the equivariant boundary is given by 
$$\sum x_a D(1\otimes \eta_a)=\sum x_a n_{ab} 1\otimes \eta_b +\sum
x_a m_{ac} 1\otimes 1_c, $$
where the first term is zero due to our assumption on $\sum x_a a$.
The element $\Delta_k(x_a a)$ is unchanged if we add to $\sum x_a
1\otimes \eta_a$ an element in the kernel of $i_{2k+1}$. Adding
the element $\sum_{a,c} x_a m_{ac} \Omega\otimes \eta_c$, which is in the
kernel of $i$, we get the following diagram

$$
\spreaddiagramrows{-1pc}
\spreaddiagramcolumns{-1.75pc}
\diagram
& & \vdots \dto & & & \vdots \dto & & & \vdots \dto \\
0 \rrto & & Q_{2k+1} \xto[0,3] \xto[3,0] & & & C_{2k+1,U(1)}
\xto[0,3]^{\qquad i}
\xto[3,0] & & & C_{2k+1} \rrto \xto[3,0] & & 0 \\
& & & & & & 
\Text{$x_a~1\otimes\eta_a$ \\$+ x_a m_{ac}~\Omega\otimes \eta_c$}
\xdashed[3,0]|<\stop|>\tip & & & x_a~a
\xdashed[0,-3]|<\stop|>\tip \\
\\
0 \rrto & & Q_{2k} \xto[0,3] \ddto & & & C_{2k,U(1)} \xto[0,3]^{\qquad i} \ddto
& & & C_{2k} \rrto \ddto & & 0 \\
& & & \Delta_k(x_a~a) & & & x_a m_{ac} m_{ce}~\Omega\otimes 1_e
\xdashed[0,-3]|<\stop|>\tip \\
& & \vdots & & & \vdots & & & \vdots
\enddiagram
$$
where we have
\[ D(x_a 1\otimes \eta_a + x_a m_{ac} \Omega\otimes \eta_c)= x_a
m_{ac} 1\otimes 1_c + \]
\[ x_a m_{ac} (-1\otimes 1_c +n_{cd} \Omega \otimes
\eta_d + m_{ce} \Omega\otimes 1_e)= x_a m_{ac} m_{ce} \Omega\otimes 1_e \]
with sums over repeated indices. The last equality follows from the
identities $\sum m_{ac}n_{cd}+n_{ab}m_{bd}=0$ and $\sum x_a
n_{ab}m_{bd}=0$. 
We can iterate the procedure. In the following step we add a term
$\sum x_a m_{ac}m_{ce} \Omega^2 \otimes \eta_e$ to the preimage of
$\sum x_a a$. The corresponding image under the boundary of the
equivariant complex is $\sum x_a m_{ac}m_{ce}m_{eg} \Omega^2\otimes
1_g$, where $\mu(O_c)-\mu(O_e)=\mu(O_e)-\mu(O_g)=2$. The procedure can
be iterated until the reducible point $\theta$ is hit. Contributions
from other irreducible critical orbits $O_p$ with
$\mu(O_p)-\mu(\theta)=0$ are killed in finitely many steps, iterating
the same procedure, since the complex is finitely generated and they
eventually hit the lowest index critical points. Thus,  
the resulting image under the connecting homomorphism 
$\bigtriangleup_{k} (\sum _{a: \mu{a}= 2k+1} x_a a)$ is given by
\[
\sum x_a m_{ac}m_{ce}\cdots m_{\alpha' \alpha} m_{\alpha\theta} \Omega ^k.
\]

\noindent $\diamond$
 
As in the case of instanton Floer theory \cite{Taubes}, one expects
$Ker(\Delta_k)\cong Im(i_{2k+1})$ to be the part of the Floer homology
where the relative invariants of four-manifolds with boundary live.

\subsection{Wall-crossing formula: the algebraic picture}

In this section, we will apply the equivariant Seiberg-Witten-Floer
homology theory to study the dependence of the metric for the Casson-type
invariant \cite{CMWZ2} \cite{Wa} of a homology sphere. In
order to define the Casson-type invariant, 
we choose a metric whose ordinary Dirac operator has trivial kernel.
The metrics whose ordinary Dirac operator has 
non-trivial kernel form a chamber structure as proven in Theorem
\ref{codim:1} and Theorem \ref{wallstrata}.
The usual cobordism argument can be adopted to prove that
the Casson-type invariant is constant in each chamber. The aim of this
section is to get a wall-crossing formula for a path of metrics and
perturbations that crosses the wall.

Denote by $\lambda_{SW}(Y, g, \nu)$ the Casson-type invariant for the
metric and perturbation $(g,\nu)$:
$\lambda_{SW}(Y, g,\nu)$ is the Euler characteristic of the non-equivariant
Seiberg-Witten-Floer homology. Recall that this Floer homology \cite{Wa} is
defined by removing the reducible critical point, the trivial
solution $\theta=[\nu,0]$. We have
\beq
\lambda_{SW}(Y, g,\nu) = \sum_{k} (-1)^k \dim SWH_k(Y, g,\nu).
\label{casson}
\eeq

In this Section 7.2 we derive the wall crossing formula under the
following assumption. 
Fix metrics and perturbations $(g_0,\nu_0)$, $(g_1,\nu_1)$ in two
different chambers, with the property that there exists an open set of 
paths $(g_t,\nu_t)$ connecting $(g_0,\nu_0)$ to $(g_1,\nu_1)$ such
that $(g_t,\nu_t)$ hits a co-dimension one wall only once,
transversely. 
This situation is certainly verified if the two points $(g_0,\nu_0)$
and $(g_1,\nu_1)$ are close enough points on either side of a wall
${\cal W}_1$, in the notation of Theorem \ref{wallstrata}, that is, a
wall of metrics and perturbations satisfying $Ker(\partial^{g}_\nu)=\CC$.
We also assume that, along the path $(g_t,\nu_t)$ we have
$SF(\partial^{g_t}_{\nu_t})= 2$, that is, that we have
$\mu(\theta_0)-\mu(\theta_1)=2$ in the notation used in Section 6,
which determines a global grade shift between the equivariant Floer
complexes for $(g_0,\nu_0)$ and for $(g_1,\nu_1)$.

Using the topological invariance 
of the equivariant Seiberg-Witten-Floer homology, we obtain the following
isomorphism:
\[
SWH_{k, U(1)}(Y, g_0,\nu_0) \cong SWH _{k+2, U(1)}(Y, g_1,\nu_1).
\]
In fact, the necessary degree shift is computed as follows: if we set
$\mu(\theta_0)=0$, hence $\mu(\theta_1)=-2$, and 
the map $I$ has degree zero, we obtain that the generators $\Omega^n
\otimes 1_{\theta_0}$ and $\Omega^{n+1}\otimes 1_{\theta_1}$ have the
same degree.

First of all, we express the Casson-type invariant (\ref{casson}) in
terms of some 
alternating sum of the equivariant  Seiberg-Witten-Floer homology groups.
We consider the Casson-type invariant for metric an perturbation $(g_0,\nu_0)$.
The exact sequences given in Theorem \ref{exact}, relating the
equivariant to the non-equivariant Seiberg-Witten-Floer 
homology, and the fact that we have $H_{*, U(1)}(\theta ) = \RR
[\Omega]$ give us the following result.

\begin{prop} 
\label{wallcross1}
(1) For $k< 0 $,
\[ SWH_k (Y, g_0,\nu_0)  \cong SWH_{k, U(1)}(Y, g_0,\nu_0). \]
For $k\ge 0$, we have the following exact sequences
\[
0\to SWH_{2k+1, U(1)}(Y, g_0,\nu_0) \to SWH_{2k+1}(Y, g_0,\nu_0)\to
\RR \Omega^k  
\to
\]\[
\to SWH_{2k, U(1)}(Y, g_0,\nu_0) \to SWH_{2k}(Y, g_0,\nu_0)\to 0. 
\]

Thus the dimensions are related by 
\[
\dim SWH_{2k} - \dim  SWH_{2k+1} =\dim SWH_{2k, U(1)}- \dim  SWH_{2k+1, U(1)}
-1. \]
\end{prop}

This gives the wall-crossing formula for the Casson-type invariant.

\begin{thm}
\label{wallcross2}
Suppose given two metrics and perturbations $(g_0,\nu_0)$ and
$(g_1,\nu_1)$ in two different chambers and a 
generic path $(g_t,\nu_t)$ that connects them and that crosses the
wall ${\cal W}$ once. Assume that the crossing happens at a generic
point of ${\cal W}$, so that the relative Morse index with respect to the
reducible solution decreases by 2 across the wall, i.e. that we have
$\mu(\theta_0)-\mu(\theta_1)=2$. Then the
Casson-type invariant changes by
\[ \lambda_{SW}(Y, g_1,\nu_1)= \lambda_{SW}(Y, g_0,\nu_0) -1. \]
\end{thm}  

\noindent\underline{Proof:}
  We can assume that we are in one of the following two cases for the
non-equivariant Floer homology group.

\noindent\underline{Case 1:} There exists an integer $N$ 
such that $SWH_{p}(Y, g_0,\nu_0) =0$
for all $p \ge 2N$ but $SWH_{2N-1}(Y, g_0,\nu_0) \ne 0$.

\noindent\underline{Case 2:} here exists an integer $N$ such that
$SWH_{p}(Y, g_0,\nu_0) =0$ for all $p \ge 2N +1$ but $SWH_{2N}(Y,
g_0,\nu_0) \ne 0$. 

For case 1, we can assume $N \ge 0$,
 otherwise the non-equivariant Floer homology groups
are isomorphic to the equivariant Floer homology groups by Proposition
\ref{wallcross1}. Therefore, from the exact sequences in Proposition
\ref{wallcross1}, we get 

\begin{itemize}

\item  $SWH_{p, U(1)} (Y, g_0,\nu_0) = \left\{
  \begin{array}{ll}
0 \qquad & p > 2N, \  \hbox{p is odd}\\[2mm]
\RR \Omega^{m} \qquad & p= 2m \ge 2N.
\end{array}\right.$ 

\item
\[ \begin{array} {ll}
& \dim SWH_{2N-2}(Y, g_0,\nu_0)-\dim  SWH_{2N-1}(Y, g_0,\nu_0) \\[2mm]
 =& \dim SWH_{2N-2, U(1)}(Y, g_0,\nu_0) - \dim  SWH_{2N-1, U(1)}(Y,
 g_0,\nu_0) -1 \\[2mm] 
& \dim SWH_{2N-4}(Y, g_0) -\dim  SWH_{2N-3}(Y, g_0,\nu_0)\\[2mm] 
=& \dim SWH_{2N-4, U(1)}(Y, g_0,\nu_0)- \dim  SWH_{2N-3, U(1)}(Y,
g_0,\nu_0) -1 \\[2mm] 
 & \cdots \\[2mm]
& \dim SWH_{0}(Y, g_0,\nu_0) -\dim  SWH_{1}(Y, g_0,\nu_0) \\[2mm]
= & \dim SWH_{0, U(1)}(Y, g_0,\nu_0)- \dim  SWH_{1, U(1)}(Y, g_0,\nu_0) -1. 
\end{array} \]

\item For $k< 0 $, $SWH_{k} (Y, g_0,\nu_0) = SWH_{k, U(1)}(Y, g_0,\nu_0)$.  
 \end{itemize}

Given the above information, we can calculate the Casson-type invariant for the
metric and perturbation $(g_0,\nu_0)$ as,
\[
\begin{array}{lll}
\lambda_{SW}(Y, g_0,\nu_0) &=& \sum_k (-1)^k dim SWH_k(Y, g_0,\nu_0)\\[2mm]
&=&  \sum_{k <N} (dim SWH_{2k, U(1)} (Y, g_0,\nu_0) \\[2mm]
& &  - \dim SWH_{2k+1, U(1)}(Y, g_0,\nu_0))-N.
\end{array}
\]

{}From the isomorphism $SWH_{k, U(1)}(Y, g_1,\nu_1) = SWH_{k-2,
  U(1)}(Y, g_0,\nu_0)$, 
we have
$$ SWH_{p, U(1)} (Y, g_1,\nu_1) =  \left\{
  \begin{array}{ll}
0 \qquad & p > 2N+2, \  \hbox{p is odd}\\[2mm]
\RR \Omega^{m} \qquad & p= 2m \ge 2N+2.
\end{array}\right.$$
If we apply the above isomorphisms to the exact sequences in
Proposition \ref{wallcross1}, we can see that the Casson-type
invariant for metric and perturbation $(g_1,\nu_1)$ can be rewritten as 
\[
\begin{array}{lll}
\lambda_{SW}(Y, g_1,\nu_1) &=& \sum_k (-1)^k \dim SWH_k(Y, g_1,\nu_1)\\[2mm]
&=&  \sum_{k <N+1}  (\dim SWH_{2k, U(1)} (Y, g_1,\nu_1) \\[2mm]
& &    - \dim SWH_{2k+1, U(1)}(Y, g_1,\nu_1) )- (N+1) \\[2mm]
&=& \lambda_{SW}(Y, g_0,\nu_0) -1 
\end{array}
\]

For case 2, similarly, we have 
\begin{itemize}

\item  $SWH_{p, U(1)} (Y, g_0,\nu_0) = \left\{
  \begin{array}{ll}
0 \qquad & p > 2N, \  \hbox{p is odd}\\[2mm]
\RR \Omega^{m} \qquad & p= 2m \ge 2N +2 .
\end{array}\right.$

\item
\[ \begin{array} {ll}
& \dim SWH_{2N}(Y, g_0,\nu_0) = \dim SWH_{2N, U(1)}(Y, g_0,\nu_0) - 1 \\[2mm] 
& \dim SWH_{2N-2}(Y, g_0,\nu_0)-\dim  SWH_{2N-1}(Y, g_0,\nu_0) \\[2mm]
 =& \dim SWH_{2N-2, U(1)}(Y, g_0,\nu_0) - \dim  SWH_{2N-1, U(1)}(Y,
 g_0,\nu_0) -1 \\[2mm] 
& \dim SWH_{2N-4}(Y, g_0,\nu_0) -\dim  SWH_{2N-3}(Y, g_0,\nu_0)\\[2mm]
=& \dim SWH_{2N-4, U(1)}(Y, g_0,\nu_0)- \dim  SWH_{2N-3, U(1)}(Y,
g_0,\nu_0) -1 \\[2mm] 
 & \cdots \\[2mm]
& \dim SWH_{0}(Y, g_0,\nu_0) -\dim  SWH_{1}(Y, g_0,\nu_0) \\[2mm]
= & \dim SWH_{0, U(1)}(Y, g_0,\nu_0)- \dim  SWH_{1, U(1)}(Y, g_0,\nu_0) -1
\end{array} \]

\item For $k< 0 $, $SWH_{k} (Y, g_0,\nu_0) = SWH_{k, U(1)}(Y, g_0,\nu_0)$.
 \end{itemize}

Then the Casson-type invariant for $(g_0,\nu_0)$ is
\[
\begin{array}{lll}
\lambda_{SW}(Y, g_0,\nu_0) &=& \sum_k (-1)^k \dim SWH_k(Y, g_0,\nu_0)\\[2mm]
&=&  \sum_{k <N} (\dim SWH_{2k, U(1)} (Y, g_0,\nu_0) \\[2mm]
& & - \dim SWH_{2k+1, U(1)}(Y, g_0,\nu_0))-N
\\[2mm]
& & + \dim SWH_{2N, U(1)}(Y, g_0,\nu_0) -1. \end{array}
\]

{}From the isomorphism $SWH_{k, U(1)}(Y, g_1,\nu_1) = SWH_{k-2,
  U(1)}(Y, g_0,\nu_0)$, 
we know that
$$ SWH_{p, U(1)} (Y, g_1,\nu_1) =  \left\{
  \begin{array}{ll}
0 \qquad & p > 2N+2, \  \hbox{p is odd}\\[2mm]
\RR \Omega^{m} \qquad & p= 2m \ge 2N+2. 
\end{array}\right.$$

If we apply the above isomorphisms to the exact sequences in
Proposition \ref{wallcross1} 
again, we can see that the Casson-type invariant for metric and
perturbation $(g_1,\nu_1)$ is
 \[
\begin{array}{lll}
&&\lambda_{SW}(Y, g_1,\nu_1)\\[2mm]
 &=& \sum_k (-1)^k \dim SWH_k(Y, g_1,\nu_1)\\[2mm]
&=&  \sum_{k <N+1}  (\dim SWH_{2k, U(1)} (Y, g_1,\nu_1) - \dim
SWH_{2k+1, U(1)}(Y, g_1,\nu_1) 
)\\[2mm]
&& - (N+1) + \dim SWH_{2N+2. U(1)}(Y, g_1,\nu_1)- 1 \\[2mm]
& =& \lambda_{SW}(Y, g_0,\nu_0) -1 
\end{array}
\]

Thus, we have proven the wall-crossing formula:

\[
\lambda_{SW}(Y, g_1,\nu_1)= \lambda_{SW}(Y, g_0,\nu_0) -1. \]

\noindent $\diamond$

Now we are interested in generalizing the argument of Theorem
\ref{wallcross2} to the case of a path $(g_t,\nu_t)$ that crosses the
wall structure ${\cal W}$ at a point which lies in a stratum ${\cal
W}_n$ of higher codimension.

If we know that a stratum ${\cal W}_n$ of metrics and perturbations
satisfying $Ker(\partial^g_\nu)=\CC^n$ is obtained as the transverse 
intersection of $n$ strata ${\cal W}_1^{i_1} \cap \cdots \cap {\cal
W}_1^{i_n}$, where every ${\cal W}_1^{i_k}$ consists of metrics with
$Ker(\partial^g_\nu)=\CC$, then a path $(g_t,\nu_t)$ that crosses
${\cal W}$ at a point in ${\cal W}_n$ can be deformed to a path that
crosses each ${\cal W}_1^{i_k}$ once transversely. In this case, the
wall crossing formula simply follows by applying repeatedly Theorem
\ref{wallcross2}. However, we do not really need the assumption on the
structure of ${\cal W}$ near a stratum of higher codimension. In fact, 
it is enough to know that the complex spectral flow
\index{$SF(\partial^{g_t}_{\nu_t})$} 
$SF_{\CC}(\partial^{g_t}_{\nu_t})=\frac{1}{2}
SF(\partial^{g_t}_{\nu_t})$ is equal to $\pm n$ along the path
$(g_t,\nu_t)$ that crosses a point on ${\cal W}_n$. In that case, we
can follow the same argument in the proof of Theorem \ref{wallcross2}, 
but starting with a grade shift of $2n$ between the equivariant Floer
complexes for $(g_0,\nu_0)$ and $(g_1,\nu_1)$. We obtain the following 
result. 

\begin{prop} \label{wallcross:n}
Let $(g_0,\nu_0)$ and $(g_1,\nu_1)$ be two metrics and perturbations 
in two different chambers. Suppose given a 
path $(g_t,\nu_t)$ joining them that crosses the
wall ${\cal W}$ once transversely at a point of a
stratum ${\cal W}_n$ of codimension $2n-1$. The 
relative Morse index with respect to the
reducible solution decreases by 2n across the wall, that is, that we
have $\mu(\theta_0)-\mu(\theta_1)=2n$. Then the
Casson-type invariant changes by
\[ \lambda_{SW}(Y, g_1,\nu_1)= \lambda_{SW}(Y, g_0,\nu_0) -n. \]
\end{prop}

Since the argument of Theorem \ref{wallcross2} depends only on the
counting of the grade shift between the Floer complexes for
$(g_0,\nu_0)$ and $(g_1,\nu_1)$ given by the spectral flow
$SF(\partial^{g_t}_{\nu_t})$ , together with the proof of
topological invariance of the equivariant Floer homology (up to this
grade shift), we can formulate the result under these more general
hypothesis. 

\begin{prop} \label{wallcross:general}
Let $(g_0,\nu_0)$ and $(g_1,\nu_1)$ be two metrics and perturbations 
in two different chambers. Suppose given a path $(g_t,\nu_t)$
joining them that crosses the wall ${\cal W}$ transversely in finitely
many points. Then the Casson invariant changes by
\beq \lambda_{SW}(Y, g_1,\nu_1)= \lambda_{SW}(Y, g_0,\nu_0)
-SF_{\CC}(\partial^{g_t}_{\nu_t}), \label{wall:gen} \eeq
where $SF_{\CC}(\partial^{g_t}_{\nu_t})=\frac{1}{2}
SF(\partial^{g_t}_{\nu_t})$ is the complex spectral flow of the Dirac
operator along the path of reducible solutions $[\nu_t,0]$.
\end{prop} 

\noindent\underline{Proof:} If the spectral flow along the path
$(g_t,\nu_t)$ is given by $SF(\partial^{g_t}_{\nu_t})$, the
topological invariance of the equivariant Floer homology gives
\[
SWH_{k, U(1)}(Y, g_0,\nu_0) \cong SWH _{k+SF(\partial^{g_t}_{\nu_t}),
U(1)}(Y, g_1,\nu_1). 
\]
We can then follow the steps of the proof of Theorem \ref{wallcross2}
and compare the ranks of the Floer groups and the counting of the
Euler characteristic. This can be done by induction on
$|SF(\partial^{g_t}_{\nu_t})|$. The result is the formula
(\ref{wall:gen}). 

\noindent $\diamond$

The wall crossing formula in the  case of $J$-invariant perturbations
constructed by W.Chen 
\cite{Chen} can also be derived with the same method. This gives rise to
the following wall crossing result. 

\begin{corol} \label{Jwallcross}
Consider the invariant $\lambda_{SW}(Y, g, f)$ where $f$ is the
$J$-invariant perturbation of \cite{Chen} Prop.2.6. 
Given two metrics and perturbations
$(g_0,f_0)$ and $(g_1,f_1)$ in two different chambers and a generic
path $(g_t,f_t)$ that crosses the wall once with  
$Ker(\partial^g_f)=\HH$. The relative Morse index with respect to the
reducible solution decreases by 4 across the wall, namely, we have
$\mu(\theta_0)-\mu(\theta_1)=4$. Then the 
Casson-type invariant changes by
\[ \lambda_{SW}(Y,g_1,f_1)= \lambda_{SW}(Y,g_0,f_0) -2. \]
\end{corol}

\noindent\underline{Proof:} Again the main issue is the change of
grading of the equivariant Floer homology induced by the spectral flow 
of the Dirac operator along the path $(g_t,f_t)$. This time, since we
are using $J$--invariant perturbations, the Dirac operator
$\partial^g_f$ is quaternion linear, hence, for the chosen path the
spectral flow satisfies $SF(\partial_{f_t}^{g_t})=4$. This implies
that there is a degree shift
$$ SWH_{k, U(1)}(Y, g_0,f_0) \cong SWH _{k+4,
U(1)}(Y, g_1, f_1). $$
By applying the previous Proposition \ref{wallcross:general}, we
obtain the result.

\noindent $\diamond$

\subsection{Wall-crossing formula: the geometric picture}

In this section we re-derive, in a more geometric way, the wall crossing
formula for a homology three-sphere $Y$ that we proved
algebraically in the previous section. We analyze the local
structure of the parameterized moduli space. A geometric
proof of the wall-crossing formula has been also worked out by \cite{Lim}. 
 
Let ${\cal M}^*(g,\nu)$ denote the irreducible part of the 
moduli space ${\cal M}$ for the metric and perturbation $(g,\nu)$.
Given a family of metrics and perturbation $(g_t, \nu_t)$ with $(t
\in [-1, 1])$, the moduli spaces ${\cal M}^*(g_{-1},\nu_{-1})$ and
${\cal M}^*(g_1,\nu_1)$ are cobordant as long as the path
$(g_t,\nu_t)$ does not cross the wall, that is the co-dimension one
subspace $W$ in the space of metrics and perturbations
\[ W=\{ (g,\nu ) | Ker (\partial_{\nu}^g)\neq 0 \}. \]

Suppose the path $(g_t,\nu_t)$ crosses the wall $W$ just once at
$t=0$. Generically, $Ker (\partial^{g_0}_{\nu_0}) \cong \CC$. We want
to analyze the local structure of the parameterized moduli space
\[ \bar {\cal M} = \{ {\cal M}(g_t,\nu_t) \times \{ t\} | t \in [-1, 1]\}
\]
at the reducible point $\vartheta_0=(\theta_0,0)$, where
$\theta_0=[\nu_0,0]$ is 
the class of the reducible solution of (\ref{extremP}) with the metric and
perturbation $(g_0,\nu_0)$. There is a family
of reducibles $\vartheta_t$ in $\bar {\cal M}$.  Let
$\bar {\cal M}^*$ be the irreducible set in $\bar {\cal M}$,
${\cal U}$ be a sufficiently small neighborhood of $\theta_0$ in 
$\bar {\cal M}$, and ${\cal U} ^*$ be the irreducible part of ${\cal
  U}$. 

We construct a bundle over neighborhood of $\vartheta_0$ in ${\cal
  A}\times [-1,1]$, together with a section $\varsigma$ such that 
\[ {\cal U} ^*=(\varsigma^{-1}(0)-\{(\vartheta_t,t)\})/{\cal G}. \]

\begin{lem}
\label{slicemodel}
1) The slice of the ${\cal G}/ U(1)$-action at a point $(A_0, 0)$
is $V_{(A_0, 0)} =  Ker (d^*) \times \Gamma _{L^2}(S)$,

(2) The slice of the ${\cal G}/ U(1)$-action at a point $(A, \psi)$ is
\[
V_{(A, \psi)} = \{ (\alpha, \phi) | d^*(\alpha) - 2i Im \langle \phi ,
\psi \rangle 
\hbox{is a constant function on}\  Y.\} 
\]

(3) For $(A, \psi )$ close to $(A_0, 0)$ there is an isomorphism 
\[
\lambda_{(A, \psi )} :\qquad V_{(A, \psi )}  \to V_{(A_0, 0)}
\]
\end{lem}

\noindent\underline{Proof:} 
Properties (1) and (2) follow by direct computation. For (3), choose
$(\alpha,\phi)$ in $V_{(A, \psi)}$ and define $\lambda_{(A, \psi
  )}(\alpha, \phi)$ to be 
\[
(\alpha - 2d\xi_{(\alpha, \phi)}, \xi_{(\alpha, \phi)} \psi + \phi)
\]
where $\xi_{(\alpha, \phi)}$ is the unique solution of the following equations:
\[
\left\{ \begin{array}{l}
2d^*d\xi_{(\alpha, \phi)} = d^* \alpha\\[2mm]
\displaystyle{\int_Y} \xi_{(\alpha, \phi)} dv = 0
\end{array}
\right.
\]
Direct computation shows that $\lambda_{(A, \psi)}$ is an isomorphism. 

\noindent $\diamond$

The above Lemma shows that we obtain a locally trivial vector bundle $V$ 
over the space of connections and sections ${\cal A}$ endowed with a
$U(1)$-action.  

Define the section $\varsigma$ 
\[ \varsigma: {\cal A}\times [-1,1] \to V \]
to be
\[ \varsigma(A, \psi, t) = \lambda_{(A, \psi)} (*_{g_t} (F_A - d \nu_t) - 
\sigma (\psi, \psi), \partial_A\psi).
\]

Near $\theta_0$, we know that ${\cal U} = \varsigma^{-1}(0) / {\cal
  G}$. 
Therefore, the local structure of ${\cal U}^*$ at $\theta_0$ is given
by the Kuranishi model of $\varsigma^{-1}(0) / {\cal G}$ at
$\theta_0$.  

Suppose $(A_t,\psi_t)$ is an element in ${\cal U}^*$.
Consider a formal expansion at $\vartheta_0$ of the form
\[ A_t=\nu_t +t\alpha_1 +t^2\alpha_2 +\cdots, \]
\[ \psi_t= t\psi_1 +t^2\psi_2+\cdots . \]

The section $\varsigma$ is approximated by 
\[ *_{g_0}d(\nu_t +t\alpha_1 +t^2\alpha_2 +\cdots)-*_{g_0}d\nu_t
-\sigma(t\psi_1 +t^2\psi_2+\cdots,t\psi_1 +t^2\psi_2+\cdots), \]
\[ \partial^{g_0}_{\nu_t}(t\psi_1 +t^2\psi_2+\cdots)+(t\alpha_1
+t^2\alpha_2 +\cdots)\cdot (t\psi_1 +t^2\psi_2+\cdots), \]
where we are perturbing in a neighborhood of the wall $W$ just by
changing the perturbation and fixing the metric $g_0$. 
    
The zero set of the section therefore determines the conditions
$*d\alpha_1=0$ and $d^*\alpha_1=0$, which imply $\alpha_1=0$ on a
homology sphere. Moreover, we have $d^*\alpha_2=0$ and 
\[ *d\alpha_2=\sigma(\psi_1,\psi_1). \]
On the kernel of $d^*$ the operator $*d$ is invertible, hence we have
\[ \alpha_2=(*d)^{-1} \sigma(\psi_1,\psi_1). \]

The Kuranishi model near $\vartheta_0$ is given by a
$U(1)$-equivariant map
\[ {\cal S}: \RR\times Ker(\partial^{g_0}_{\nu_0})\to
CoKer(\partial^{g_0}_{\nu_0}), \]
where $U(1)$ acts on $Ker(\partial^{g_0}_{\nu_0})\cong
CoKer(\partial^{g_0}_{\nu_0})\cong \CC$ by the natural multiplication
on $\CC$.

There exists a sufficiently small 
$\delta > 0$ such that, for $t \in [-\delta , \delta]$, we have that
$\partial^{g_0}_{\nu_t}$ has exactly one small eigenvalue $\lambda(t)$ with
eigenvector $\phi_t$ and with $\lambda(0)=0 $, that is
\[
\partial^{g_0}_{\nu_t} \phi_t = \lambda(t) \phi_t.
\] 

This implies that, if $\lambda'(0) >0$, then the spectral flow of 
$\partial^{g_0}_{\nu_t}$ for $(t \in [-1, 1])$ is $1$ and, if
$\lambda'(0) < 0$, the spectral flow of $\partial^{g_0}_{\nu_t}$ for
$(t \in [-1, 1])$ is $-1$. 

The map ${\cal S}$ is given by
\[ {\cal S}:\RR\times\CC\to\CC, \]
\[ {\cal S}(t,w\phi)=\Pi_{Ker(\partial^{g_0}_{\nu_0})}(\partial_{A_t}
w\phi). \]
Here we assume that $\phi$ is a spinor in $Ker(\partial^{g_0}_{\nu_0})$
with $\| \phi\|=1$, so that $Ker(\partial^{g_0}_{\nu_0})\cong \CC\phi$.
Consider the expression
\[ <\partial^{g_0}_{\nu_t}\phi,\phi>=z(t). \]
Notice that we have $z'(0)=\lambda'(0)$, in fact, we write formally
$\lambda(t)\sim t\lambda'(0)$, $\phi_t\sim\phi+t\phi_1$ and the Dirac
operator $\partial^{g_0}_{\nu_t}\sim\partial^{g_0}_{\nu_0}+tC$, where
$\nu_t\sim \nu_0+t\nu_1$ and $C$ acts as Clifford multiplication by
$\nu_1$. We can write the first order term in the relation
$\partial^{g_0}_{\nu_t}\phi_t=\lambda(t)\phi_t$ as  
\[ t<\partial^{g_0}_{\nu_0}\phi_1,\phi> +t<C\phi,\phi>
=t\lambda'(0) +\lambda(0). \]
Here the term $<\partial^{g_0}_{\nu_0}\phi_1,\phi>=
<\phi_1,\partial^{g_0}_{\nu_0}\phi>$ vanishes, and also
$\lambda(0)=0$. Thus, we have the relation
\[ <C\phi,\phi>=\lambda'(0). \]
On the other hand, we have 
\[ <C\phi,\phi>=z'(0) \]
from the expansion of $\partial^{g_0}_{\nu_t}\phi=z(t)\phi$.

Thus the map ${\cal S}$ can be rewritten as
\[ {\cal S}(t,w\phi)=z(t)w\phi +t^2 <\alpha_2\phi,\phi> w\phi +
O(t^3) \]
\[ = w\phi \left( z(t)+t^2 r^2
<(*d)^{-1}\sigma(\phi,\phi),\sigma(\phi,\phi)> \right)+ O(t^3). \]
Here we use the fact that the first order term of the Dirac equation
gives $\partial^{g_0}_{\nu_0}\psi_1=0$, therefore
$\psi_1=re^{i\theta}\phi$ and
$\sigma(\psi_1,\psi_1)=r^2\sigma(\phi,\phi)$. 

The term
\[ \gamma(Y,g_0,\nu_0)=<(*d)^{-1}\sigma(\phi,\phi),\sigma(\phi,\phi)>
\]
is a constant that only depends on the manifold and on
$\theta_0$. An inductive argument shows that, if $\gamma(Y,g_0,\nu_0)$
vanishes, then all the forms $\alpha_i$ in the formal expansion of
$A_t$ must also vanish identically. Thus, we can assume that
$\gamma(Y,g_0,\nu_0)\neq 0$. 

Notice that we have
 $$\quad \RR \times (Ker(\partial^{g_0}_{\nu_0})- \{0\}) /U(1) =  \RR
\times \RR^+.$$
The irreducible part of $\varsigma^{-1}(0) / {\cal G}$ is
tangent to $\{ 0\} \times \RR^+$ as $t$ approaches $0$, as we see
in the following.

The difference $\lambda$ between the Casson-type invariant at $t=\pm
\delta$ can be evaluated by counting the number (with sign) of
oriented lines in $\varsigma^{-1}(0) / {\cal G}$, with
$t\in[-\delta,\delta]$, that are tangent to
$\{ 0\} \times \RR^+ \times \{ 0\}$.
Here we identify ${\cal U}^* $ with the set $({\cal S}^{-1}(0)- \{
w=0\})/U(1)$.  The sign of the wall crossing term is determined by the 
section ${\cal S}$, as follows.
The zero set $({\cal S}^{-1}(0)- \{
w=0\})/U(1)$ is given by the condition
\[ t=-\frac{\lambda'(0)}{r^2 \gamma(Y,g_0,\nu_0)}. \]
Thus, we have one line in ${\cal U}^* $ which is counted with
the orientation determined by the sign of $-\gamma(Y,g_0,\nu_0)$ and the
spectral flow. 
Suppose that we have $\lambda'(0)>0$,  then the
spectral flow is $SF_{{\CC}}(\partial^{g_t}_{\nu_t})=1$  on the path $t\in
[-1 ,1]$.
If we have $\gamma(Y,g_0,\nu_0)>0$, then there is a unique irreducible
solution, which contributes a $+1$ to the invariant, that flows into
the reducible as $t\to 0$, with $t<0$. If we have
$\gamma(Y,g_0,\nu_0)<0$, then a unique irreducible that contributes a
$-1$ to the invariant approaches the reducible as $t\to 0$, $t>0$.
This gives the wall crossing term
\[ \lambda(Y,g_1)=\lambda(Y,g_{-1})-  SF_{{\CC}}(\partial^{g_t}_{\nu_t}). \]

This provides a geometric interpretation of the wall-crossing
formula that we derived algebraically in the
previous section from the exact sequences.

\vspace{.5cm}

\noindent{\bf Acknowledgements}. We are deeply grateful to T. Mrowka
for the many invaluable comments and suggestions. 
We thank  A. Carey and M. Rothenberg for the many useful
conversations. 
We are grateful to L. Nicolaescu for having suggested the argument of
Proposition \ref{orientation}
and to R.G. Wang for having suggested the use of condition
(\ref{int=0}).  
In the early stages of this work useful were comments of R. Brussee,
G. Mati\'c, R. Mazzeo, and J.W. Morgan.
We thank the referee for the very detailed and useful comments and for
many valuable suggestions on how to improve various sections of the
manuscript. Parts of this work were carried out during visits of the
two authors to the Max Planck Institut f\"ur Mathematik in Bonn, and
visits of the first author to the University of Adelaide. We thank
these institutions for the kind hospitality and for support.

\input index.ind

%\vspace{3cm}
\pagebreak

\noindent Matilde Marcolli

\noindent {\em Department of Mathematics 2-275, MIT, Cambridge MA 02139,
  USA}

\noindent matilde@math.mit.edu 

\vspace{1cm}

\noindent Bai-Ling Wang

\noindent {\em Department of Pure Mathematics, The University of
  Adelaide, Adelaide SA 5005, Australia}

\noindent bwang@maths.adelaide.edu.au

\end{document}